\documentclass[mathpazo]{aamm}
\renewcommand\thisnumber{x}
\renewcommand\thisyear {200x}
\renewcommand\thismonth{xxx}
\renewcommand\thisvolume{xx}
\usepackage{enumerate}
\usepackage{bm}
\newtheorem{rmk}{Remark}
\newtheorem{thm}{Theorem}

\usepackage{wrapfig}
\usepackage{multicol}
\usepackage{verbatim}
\usepackage{amssymb}
\usepackage{amsmath}
\usepackage{enumitem}
\usepackage{bm}
\usepackage{color}
\usepackage{exscale}
\usepackage{relsize}
\usepackage{graphicx}
\allowdisplaybreaks[4]
\usepackage{epstopdf}
\begin{document}
	
	\markboth{Hao Dong and Yifei Ding and Jiale Linghu}{Higher-order multi-scale computational method}
	\title{Higher-order multi-scale computational method and its convergence analysis for hygro-thermo-mechanical coupling problems of quasi-periodic composite structures}
	
	\author[Hao Dong and Yifei Ding and Jiale Linghu and Yufeng Nie and Yaochuang Han]{Hao Dong\affil{1,2}\comma\corrauth and Yifei Ding\affil{1} and Jiale Linghu\affil{1} and Yufeng Nie\affil{3} and Yaochuang Han\affil{4}}
	\address{\affilnum{1}\ School of Mathematics and Statistics, Xidian University, Xi'an 710071, PR China\\
		\affilnum{2}\ Xi'an Key Laboratory of Information Network Optimization and Mathematical Methods, Xi'an 710071, PR China\\
		\affilnum{3}\ School of Mathematics and Statistics, Northwestern Polytechnical University, Xi'an 710129, PR China\\
		\affilnum{4}\ School of Mathematical Sciences, Luoyang Normal University, Luoyang 471934, PR China}
	\emails{{\tt donghao@mail.nwpu.edu.cn} (Hao Dong)}
	
	%
	%
	%
	
	\begin{abstract}
		This paper proposes a novel higher-order multi-scale (HOMS) computational method, which is highly targeted for efficient, high-accuracy and low-computational-cost simulation of hygro-thermo-mechanical (H-T-M) coupling problems in quasi-periodic composite structures. The first innovation of this work is that the establishment of the high-accuracy multi-scale model incorporating the higher-order correction terms for H-T-M coupling problems of quasi-periodic composite structures. The second innovation of this work is that the error analyses in the point-wise and integral senses are rigorously derived for multi-scale asymptotic solutions. Especially from the point-wise error analysis, the primary impetus for current study to develop the HOMS approach for quasi-periodic composite structures is illustrated. Furthermore, an high-accuracy multi-scale numerical algorithm is developed based on finite element method, while corresponding convergent analysis is also obtained. Finally, extensive numerical experiments are conducted to validate the computational performance of the proposed HOMS computational approach, demonstrating not only exceptional numerical accuracy, but also reduced computational cost.
	\end{abstract}
	
	
	
	\maketitle
	
	\section{Introduction}
	The outstanding physical properties of composite materials, including high specific stiffness, superior specific strength, lightweight characteristics, excellent corrosion resistance, high-temperature tolerance, thermal stability, among others, have led to their extensive utilization in cutting-edge engineering fields such as aeronautics, space systems, architectural engineering, mechanical manufacturing, optoelectronics, and electronic packaging \cite{R1}. However, when subjected to manufacturing process variations or fatigue damage accumulation, the composites with periodic geometric configurations exhibit a degradation of periodicity in their physical and mechanical parameters, resulting in macroscopic position-dependent material properties. Such special materials are classified as quasi-periodic composites \cite{R2,R3,R4}. Functionally graded materials (FGMs), a class of inhomogeneous composites characterized by continuous gradient variations, are recognized as a representative type of quasi-periodic composites \cite{R5}.
	
	With the extensive expansion of engineering applications, composite materials and structures routinely endure simultaneous high-heat and high-moisture conditions, while also withstanding substantial external loads \cite{R6,R7,R8,R9,R10,R45,R11,R12}. From a microscopic perspective, the absorption of large amounts of moisture or other liquids by composite structures can weaken or even break the chemical bond connection inside component materials. At the same time, as the ambient temperature continues to increase, the chemical bonds inside component structures become abnormally active due to obtaining more energy. These factors can change the basic mechanical properties of composite structures. From a macroscopic perspective, due to the thermal expansion and contraction, and hygroscopic swelling and shrinkage of composite structures, the overall sizes of composite structures will change, inducing severe hygrothermal residual stress. Moreover, excessive thermal and hygroscopic deformation may also lead to the overall performance degradation and load-bearing capacity loss, or even structural failure of composite structures. Hence, it is of great theoretical and engineering values for structural design, damage assessment and life prediction of composite structures \cite{R13}. In the field of electronic packaging technology, significant breakthroughs have been made in the multiphysics coupling, namely H-T-M coupling problems within packaging composite structures. Liu et al. \cite{R6} developed a nonlinear finite element (FE) model to simulate these coupled phenomena, successfully forecasting delamination, deformation, and fracture behavior in plastic-encapsulated devices subjected to combined mechanical, moisture, and thermal loading. Kessentini et al. \cite{R7} developed a H-T-M coupling model, revealing the influence of mechanical tension on both the coefficient of thermal expansion and moisture-driven stresses within multi-layer
	bonded assemblies. Moleiro et al. \cite{R8} formulated a mixed layerwise framework for H-T-M static analysis of multilayered plates, covering three advanced composite types: hybrid laminates, fibre-metal systems, and sandwich configurations. Zheng et al. \cite{R9} derived a second-order two-scale (SOTS) asymptotic model for periodic composites under H-T-M coupling loading and validated its effectiveness through numerical simulations. Meski et al. \cite{R10} disclosed the critical impact of nonlinear hygro-thermo-mechanical coupled loading on flexural responses of functionally graded sandwich beams, employing a novel quasi-3D  higher-order shear deformation framework. This work demonstrated the inherently nonlinear nature of mechanical responses under complex multiphysics coupling conditions. Zhang et al. \cite{R45} developed a H-T-M coupled model for RC bridge piers in plateau environments, successfully predicting temperature-induced stress evolution and cracking risks under extreme climate conditions. Hirwani et al. \cite{R11,R12} developed a micromechanical FE approach to numerically investigate the time-evolving nonlinear deflection behavior of delaminated composite shells subjected to simultaneous hygroscopic, thermal, and mechanical loading. These research advances demonstrate that the influence of hygro-thermal environments on the mechanical behavior of composites has remained a key focus for scientists and engineers. Therefore, investigating the H-T-M coupling response of quasi-periodic composite structures holds significant theoretical and engineering values.
	
	In recent years, the advancement of materials science and technology has driven the continuous emergence of advanced composite structures, which typically exhibit complex multi-scale characteristics. The characterization of physical and mechanical properties in composite structures commonly involves solving initial-boundary value problems governed by partial differential equations (PDE) exhibiting highly oscillating coefficients \cite{R2,R13,R14}. For such multi-scale problems, analytical solutions are generally unavailable \cite{R15}, necessitating the use of computational methods for numerical solutions. In order to achieve effective simulation of composite structures, computational methods should resolve the challenging issues inherent in the multi-scale nature of composite structures. Capturing highly oscillatory information and generating high-fidelity solutions at the smallest scale demands prohibitive computational resources. This leads to significantly reduced efficiency and even failure to numerically converge for conventional numerical approaches, including the finite difference method (FDM), finite element method (FEM), and finite volume method (FVM), etc \cite{R13,R14}. In the past five decades, theoretical value and engineering demand have spurred the development of diverse multi-scale frameworks, including asymptotic homogenization method (AHM) \cite{R2}, heterogeneous multi-scale method (HMM) \cite{R16}, variational multi-scale method (VMS) \cite{R17}, multi-scale finite element method (MsFEM) \cite{R18}, generalized multiscale finite element method (GMsFEM, which extends MsFEM by systematically constructing multiscale basis functions to capture sub-grid heterogeneity)\cite{R46,R47,R48,R49}, localized orthogonal decomposition method (LOD) \cite{R19}, multi-scale eigenfunction method (MEM) \cite{R20}, and variational asymptotic homogenization (VAM) \cite{R21}, etc. Current literature indicates that multi-scale computational methodologies for quasi-periodic composite structures remain relatively underdeveloped. The pioneering work in \cite{R2} by Lions et al. developed a homogenization framework governing second-order elliptic equations with quasi-periodic coefficients, thereby establishing the theoretical foundation for multi-scale computation of quasi-periodic composite structures. Cui and Cao \cite{R21} systematically developed a first-order two-scale model for elastic boundary value problems in quasi-periodic composite structures. Their derivation of several key estimates for displacement, stress and strain energy provided the mathematical foundation for subsequent finite element numerical implementations. However, practical engineering simulation find that classical homogenized and lower-order multi-scale approaches are inadequate in numerical accuracy when dealing with the physical fields exhibiting significantly local oscillations, thereby driving the development of innovative HOMS methods. In the past three decades, Cui and his research team systematically established a class of HOMS approaches for precisely and efficiently simulating the thermal, mechanical and multiphysics behaviors of composite structures, as shown in references \cite{R3,R4,R22,R23,R24,R25,R26,R27,R28,R29,R30} for more details. Especially for quasi-periodic composite structures, Su and Cui et al. \cite{R4,R28} established a SOTS method for elliptic and elastic boundary value problems, while proposing high-precision approximate solutions accompanied by rigorous error analysis. Dong, Ma and Cui et al. \cite{R3,R29,R30} developed HOMS approaches to damped wave propagation problems and dynamic thermo-mechanical coupling problems, obtaining both pointwise and integral error estimates for HOMS asymptotic solutions while developing corresponding multi-scale numerical algorithms. These researches demonstrate that HOMS methods exhibit both high-accuracy computational performance and strong applicability in predicting the performance of quasi-periodic composite structures, and also showing significant potential for engineering applications. To summarize, systematic research on multi-scale simulation and analysis for H-T-M coupling problems of quasi-periodic composite structures remains notably lacking. However, widespread engineering demands strongly prompt continued research about this challenging issue.
	
	To effectively address this H-T-M coupling problems of quasi-periodic composite structures, this study proposes a HOMS computational model and corresponding multi-scale numerical algorithm, which can preserve the local balance of concerned physical quantities enabling high-accuracy multi-scale simulation. The paper is structured as follows: Section \ref{sec:2} constructs a HOMS computational model for H-T-M coupling problems in quasi-periodic composite structures via multi-scale asymptotic analysis. The novel higher-order correction terms are rigorously derived to enrich the numerical accuracy for H-T-M coupling simulation in quasi-periodic solid structures. Section \ref{sec:3} performs point-wise convergence analysis to quantitatively compare the accuracy between lower-order multi-scale (LOMS) and HOMS solutions. The point-wise analysis theoretically demonstrates the superior capability of the higher-order computational model in capturing highly oscillatory features at micro-scale, while rigorously deriving explicit convergence rates in integral sense. Section \ref{sec:4} develops a two-stage multi-scale algorithm based on FEM and provides corresponding convergence analysis. The proposed multi-scale algorithm comprise with off-line microscale computation, and on-line macroscale and multi-scale computation. Section \ref{sec:5} validates both the superiority of the proposed HOMS method and the necessity of incorporating higher-order correction terms through comprehensive 2D and 3D numerical experiments. Section \ref{sec:6} concisely summarizes the principal contributions of this study, while delineating critical pathways for subsequent research advancement.
	
	To enhance notational conciseness, the Einstein summation convention is employed throughout this work.
	
	\section{Novel higher-order multi-scale computational model}
	\label{sec:2}
	\subsection{Problem setting and governing equations}
	\label{sec:21}
\begin{figure}[!htb]
		\centering
		\begin{minipage}[c]{0.4\textwidth}
			\centering
			\includegraphics[width=50mm]{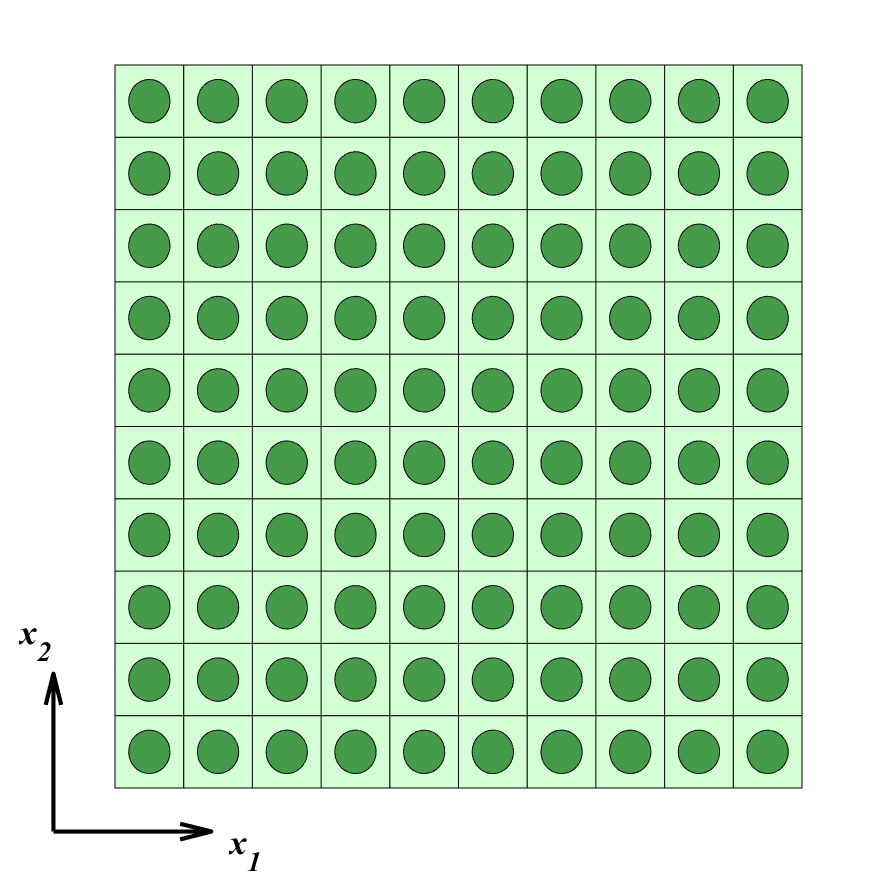}\\
			(a)
		\end{minipage}
		\begin{minipage}[c]{0.4\textwidth}
			\centering
			\includegraphics[width=50mm]{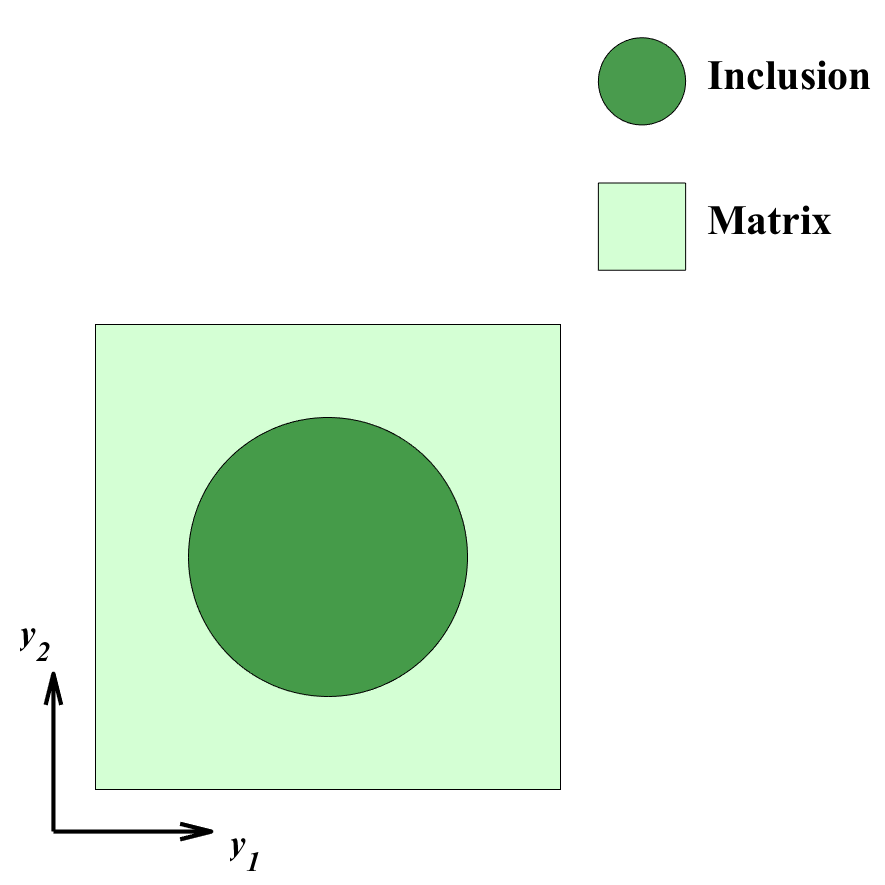}\\
			(b)
		\end{minipage}
		\caption{The schematic of composite structure: (a) composite structure $\Omega$; (b) PUC $Y$.}\label{f1:2D}
	\end{figure}

	Inspired by the H-T-M coupling model of periodic composites in \cite{R7,R8,R9}, the following governing equations for H-T-M coupling problems of quasi-periodic composite structures can be formulated over domain $\Omega$, where $\Omega\in\mathbb{R}^n(n=2,3)$ is a bounded convex domain with the Lipschitz continuous boundary $\partial\Omega=\Gamma_{T}\cup\Gamma_{q}\cup\Gamma_{c}\cup\Gamma_{d}\cup\Gamma_{u}\cup\Gamma_{\sigma}$, where these boundary parts are pairwise disjoint. The domain is formed by the repetition of periodic unit cell (PUC) $Y$, as shown in Fig.\hspace{1mm}\ref{f1:2D}.
	\begin{equation}
		\label{eq:2.1}
		\begin{cases}
			\begin{aligned}
				& -\frac{\partial}{\partial x_{i}}\Bigl(k_{ij}\bigl(\mathbf{x}, \frac{\mathbf{x}}{\epsilon}\bigr) \frac{\partial T^{\epsilon}(\mathbf{x})}{\partial x_{j}}\Bigr) = h(\mathbf{x}),\;\text{in }\;\Omega, \\
				& -\frac{\partial}{\partial x_{i}}\Bigl(g_{ij}\bigl(\mathbf{x}, \frac{\mathbf{x}}{\epsilon}\bigr) \frac{\partial c^{\epsilon}(\mathbf{x})}{\partial x_{j}}\Bigr) = m(\mathbf{x}),\;\text{in }\;\Omega, \\
				& -\frac{\partial}{\partial x_{j}}\Bigl[D_{ijkl}(\mathbf{x}, \frac{\mathbf{x}}{\epsilon}) \big( \frac{\partial u_{k}^{\epsilon}(\mathbf{x})}{\partial x_{l}}\!-\!\alpha_{kl}(\mathbf{x}, \frac{\mathbf{x}}{\epsilon}) T^{\epsilon}(\mathbf{x})\!-\!\beta_{kl}(\mathbf{x}, \frac{\mathbf{x}}{\epsilon}) c^{\epsilon}(\mathbf{x}) \big)\Bigr]\!=\!f_{i}(\mathbf{x}),\;\text{in }\;\Omega,\\
				& T^{\epsilon}(\mathbf{x}) = \overline{T}(\mathbf{x}),\;\text{on }\;\Gamma_{T},\\
				& k_{ij}\bigl(\mathbf{x}, \frac{\mathbf{x}}{\epsilon}\bigr) \frac{\partial T^{\epsilon}(\mathbf{x})}{\partial x_{j}}n_i=\overline{q}(\mathbf{x}),\;\text{on }\;\Gamma_{q},\\
				& c^{\epsilon}(\mathbf{x}) = \overline{c}(\mathbf{x}),\;\text{on }\;\Gamma_{c},\\
				& g_{ij}\bigl(\mathbf{x}, \frac{\mathbf{x}}{\epsilon}\bigr) \frac{\partial c^{\epsilon}(\mathbf{x})}{\partial x_{j}}n_i=\overline{d}(\mathbf{x}),\;\text{on }\;\Gamma_{d},\\
				& \bm{u}^{\epsilon}(\mathbf{x}) = \overline{\bm{u}}(\mathbf{x}),\;\text{on }\;\Gamma_{{u}},\\
				& \Bigl[D_{ijkl}(\mathbf{x}, \frac{\mathbf{x}}{\epsilon}) \big( \frac{\partial u_{k}^{\epsilon}(\mathbf{x})}{\partial x_{l}}\! -\!\alpha_{kl}(\mathbf{x}, \frac{\mathbf{x}}{\epsilon}) T^{\epsilon}(\mathbf{x})\!-\!\beta_{kl}(\mathbf{x}, \frac{\mathbf{x}}{\epsilon}) c^{\epsilon}(\mathbf{x}) \big)\Bigr] n_{j}\!=\!\overline{\sigma}_{i}(\mathbf{x}),\;\text{on }\;\Gamma_{\sigma},
			\end{aligned}
		\end{cases}
	\end{equation}
	where $T^{\epsilon}(\mathbf{x})$, $c^{\epsilon}(\mathbf{x})$ and $\bm{u}^{\epsilon}(\mathbf{x})$ are the undetermined temperature increment, moisture and displacement fields; $\overline{T}(\mathbf{x})$, $\overline{c}(\mathbf{x})$ and $\overline{\bm{u}}(\mathbf{x})$ are the prescribed temperature increment, moisture and displacement on the domain boundaries $\Gamma_{T}\cup\Gamma_{c}\cup\Gamma_{u}$ with $meas(\Gamma_{T})>0$, $meas(\Gamma_{c})>0$ and $meas(\Gamma_{u})>0$; $\overline{q}(\mathbf{x})$, $\overline{d}(\mathbf{x})$ and $\overline{\sigma}_{i}(\mathbf{x})$ are the prescribed heat flux, moisture flux and traction on the domain boundaries $\Gamma_{q}\cup\Gamma_{d}\cup\Gamma_{\sigma}$, where $n_j$ denotes the $j$-th component of the unit normal vector at a given point on the domain boundaries. Furthermore, $h(\mathbf{x})$, $m(\mathbf{x})$, and $f_{i}(\mathbf{x})$ are the internal heat source, internal moisture source and body force, respectively. Moreover, $\displaystyle\{k_{ij}(\mathbf{x},\frac{\mathbf{x}}{\epsilon})\}$ and $\displaystyle\{g_{ij}(\mathbf{x},\frac{\mathbf{x}}{\epsilon})\}$ are the second order thermal conductivity tensor and moisture diffusion tensor; $\displaystyle\{D_{i j k l}(\mathbf{x},\frac{\mathbf{x}}{\epsilon})\}$ is the fourth order elastic tensor; $\displaystyle\{\alpha_{k l}(\mathbf{x}, \frac{\mathbf{x}}{\epsilon})\}$ and  $\displaystyle\{\beta_{k l}(\mathbf{x},\frac{\mathbf{x}}{\epsilon})\}$ are the second order thermal expansion tensor and moisture expansion tensor. Clearly, these material parameters depend on both the macroscopic slow variable $\mathbf{x}$ and the microscopic fast variable $\displaystyle\frac{\mathbf{x}}{\epsilon}$, where $\epsilon$ represents the characteristic length of the PUC $Y=[0,1]^n$, manifesting the quasi-periodic characteristic of the heterogeneous structures under investigation.
	
	To implement subsequent multi-scale modeling and theoretical analysis, we define $\displaystyle\mathbf{y} = \frac{\mathbf{x}}{\epsilon}$ as microscopic coordinates of PUC $Y$. Based on this definition, $\displaystyle k_{ij}(\mathbf{x},\frac{\mathbf{x}}{\epsilon })$, $\displaystyle g_{ij}(\mathbf{x},\frac{\mathbf{x}}{\epsilon })$, $\displaystyle D_{ijkl}(\mathbf{x},\frac{\mathbf{x}}{\epsilon })$, $\displaystyle \alpha _{kl}(\mathbf{x},\frac{\mathbf{x}}{\epsilon })$ and $\displaystyle\beta _{kl}(\mathbf{x},\frac{\mathbf{x}}{\epsilon })$ can be rewritten as $k_{ij}(\mathbf{x},\mathbf{y})$, $g_{ij}(\mathbf{x},\mathbf{y})$, $D_{ijkl}(\mathbf{x},\mathbf{y})$, $\alpha_{kl}(\mathbf{x},\mathbf{y})$ and $\beta_{kl}(\mathbf{x},\mathbf{y})$, respectively. Next, we make the following fundamental assumptions:
	\begin{enumerate}
		\item[(A)] The material parameters $k_{ij}$, $g_{ij}$, $D_{ijkl}$, $\alpha_{kl}$ and $ \beta_{kl}$ belong to $L^\infty(\Omega)$. Moreover, functions $k_{ij}(\mathbf{x},\mathbf{y})$, $g_{ij}(\mathbf{x},\mathbf{y})$, $D_{ijkl}(\mathbf{x},\mathbf{y})$, $\alpha_{kl}(\mathbf{x},\mathbf{y})$ and $\beta_{kl}(\mathbf{x},\mathbf{y})$ are $1$-periodic in microscopic variable $\mathbf{y}$.
		\item[(B)] The material coefficients $k_{ij}$, $g_{ij}$, $D_{ijkl}$, $\alpha_{kl}$, and $\beta_{kl}$ are symmetric and uniformly elliptic for all vectors $\bm{\xi}=( \xi_i)\in \mathbb{R}^n $ and  symmetric matrices $\{\eta_{ij}\}\in\mathbb{R}^{n \times n}$, which means material coefficients satisfy
		\begin{displaymath}
			\begin{aligned}
				&k_{ij}=k_{ji},\; \underline{\gamma} | \bm{\xi} |^2 \leq k_{ij} (\mathbf{x},\mathbf{y} ) \xi_i \xi_j \leq \overline{\gamma} | \bm{\xi} |^2,\\
				&g_{ij}=g_{ji},\; \underline{\gamma} | \bm{\xi} |^2 \leq g_{ij} (\mathbf{x},\mathbf{y} ) \xi_i \xi_j \leq \overline{\gamma} | \bm{\xi} |^2,\\		
				&D_{ijkl} = D_{ijlk} = D_{klij}, \;\underline{\gamma} \eta_{ij} \eta_{ij} \leq D_{ijkl} (\mathbf{x},\mathbf{y} ) \eta_{ij} \eta_{kl} \leq \overline{\gamma} \eta_{ij} \eta_{ij},\\	
				&\alpha_{kl}=\alpha_{lk},\;\underline{\gamma}| \bm{\xi} |^2 \leq \alpha_{kl} (\mathbf{x},\mathbf{y} ) \xi_i \xi_j \leq \overline{\gamma} | \bm{\xi} |^2,\\
				&\beta_{kl}=\beta_{lk},\;\underline{\gamma}| \bm{\xi} |^2 \leq \beta_{kl} (\mathbf{x},\mathbf{y} ) \xi_i \xi_j \leq \overline{\gamma} | \bm{\xi} |^2,\\
			\end{aligned}
		\end{displaymath}
		where $\underline{\gamma}$ and $\overline{\gamma}$ are two positive constants independent of $\epsilon$.
		\item[(C)] $h(\mathbf{x})$, $m(\mathbf{x})$ and $f_i(\mathbf{x}) \in L^2(\Omega)$. $\overline{T}(\mathbf{x})$ and $\overline{c}(\mathbf{x}) \in H^1(\Omega)$, and $\overline{\bm{u}}(\mathbf{x}) \in (H^1(\Omega))^n$. $\overline{q}(\mathbf{x})$, $\overline{d}(\mathbf{x})$ and $\overline{\sigma}_{i}(\mathbf{x})\in L^2(\Omega)$.
	\end{enumerate}
	
	\subsection{HOMS modeling for quasi-periodic composite structures}
	\label{sec:22}
	Considering the relation of macro- and micro-coordinates in composite structures, the multi-scale chain rule achieves as below
	\begin{equation}
		\label{eq:2.2}
		\frac{\partial \Phi^{\epsilon} (\mathbf{x})}{\partial x_i}=\frac{\partial \Phi(\mathbf{x},\mathbf{y})}{\partial x_i}+\frac{1}{\epsilon }\frac{\partial \Phi^{\epsilon} (\mathbf{x},\mathbf{y})}{\partial y_i},
	\end{equation}
	which forms the cornerstone for subsequent multi-scale modeling.
	
	Next, performing a standard way, the succeeding multi-scale asymptotic expansions are derived for the exact solutions $T^{\epsilon}(\mathbf{x})$, $c^{\epsilon}(\mathbf{x})$ and $u_{i}^{\epsilon}(\mathbf{x})$ in terms of $\epsilon$
	\begin{equation}
		\label{eq:2.3}
		\begin{cases}
			\begin{aligned}
				&T^{\epsilon}(\mathbf{x}) = T^{(0)}(\mathbf{x},\mathbf{y}) + \epsilon T^{(1)}(\mathbf{x},\mathbf{y}) + \epsilon^{2} T^{(2)}(\mathbf{x},\mathbf{y}) + O(\epsilon^{3}), \\
				&c^{\epsilon}(\mathbf{x}) = c^{(0)}(\mathbf{x},\mathbf{y}) + \epsilon c^{(1)}(\mathbf{x},\mathbf{y}) + \epsilon^{2} c^{(2)}(\mathbf{x},\mathbf{y}) + O(\epsilon^{3}), \\
				&u_{i}^{\epsilon}(\mathbf{x}) = u_{i}^{(0)}(\mathbf{x},\mathbf{y}) + \epsilon u_{i}^{(1)}(\mathbf{x},\mathbf{y}) + \epsilon^{2} u_{i}^{(2)}(\mathbf{x},\mathbf{y}) + O(\epsilon^{3}),
			\end{aligned}
		\end{cases}
	\end{equation}
	where $T^{(0)}$, $c^{(0)}$ and $u_{i}^{(0)}$ are defined as zeroth-order expansion terms, $T^{(1)}$, $c^{(1)}$ and $u_{i}^{(1)}$ are defined as first-order (lower-order) asymptotic terms, and $T^{(2)}$, $c^{(2)}$ and $u_{i}^{(2)}$ are defined as second-order (higher-order) asymptotic terms.

	Then, substituting \eqref{eq:2.3} into multi-scale equations \eqref{eq:2.1} and expanding all spatial derivatives by chain rule \eqref{eq:2.2}, we hence derive a sequence of equations by balancing the two sides in terms of $\epsilon$.
	\begin{equation}
		\label{eq:2.4}
		O(\epsilon^{-2}):
		\begin{cases}
			\begin{aligned}
				& \frac{\partial}{\partial y_i}\Bigl(k_{ij} \frac{\partial T^{(0)}}{\partial y_j}\Bigr) = 0, \\
				& \frac{\partial}{\partial y_i}\Bigl(g_{ij} \frac{\partial c^{(0)}}{\partial y_j}\Bigr) = 0, \\
				& \frac{\partial}{\partial y_j}\Bigl(D_{ijkl} \frac{\partial u_k^{(0)}}{\partial y_l}\Bigr) = 0.
			\end{aligned}
		\end{cases}
	\end{equation}
	\begin{equation}
		\label{eq:2.5}
		O(\epsilon^{-1}):
		\begin{cases}
			\begin{aligned}
				& \frac{\partial}{\partial x_i}\Bigl(k_{ij} \frac{\partial T^{(0)}}{\partial y_j}\Bigr)+ \frac{\partial}{\partial y_i}\Bigl[k_{ij} \bigl(\frac{\partial T^{(0)}}{\partial x_j} + \frac{\partial T^{(1)}}{\partial y_j}\bigr)\Bigr] = 0, \\
				& \frac{\partial}{\partial x_i}\Bigl(g_{ij} \frac{\partial c^{(0)}}{\partial y_j}\Bigr) + \frac{\partial}{\partial y_i}\Bigl[g_{ij} \bigl(\frac{\partial c^{(0)}}{\partial x_j} + \frac{\partial c^{(1)}}{\partial y_j}\bigr)\Bigr] = 0, \\
				& \frac{\partial}{\partial x_j}\Bigl(D_{ijkl} \frac{\partial u_k^{(0)}}{\partial y_l}\Bigr) + \frac{\partial}{\partial y_j}\Bigl[D_{ijkl} \bigl(\frac{\partial u_k^{(0)}}{\partial x_l} + \frac{\partial u_k^{(1)}}{\partial y_l}\bigr)\Bigr] \\
				& - \frac{\partial}{\partial y_j}\Bigl[D_{ijkl} \bigl(\alpha_{kl} T^{(0)} + \beta_{kl} c^{(0)}\bigr)\Bigr] = 0.
			\end{aligned}
		\end{cases}
	\end{equation}
	\begin{equation}
		\label{eq:2.6}
		O(\epsilon^0):
		\begin{cases}
			\begin{aligned}
				&\!\frac{\partial}{\partial x_i}\Bigl[k_{ij} \bigl(\frac{\partial T^{(0)}}{\partial x_j} + \frac{\partial T^{(1)}}{\partial y_j}\bigr)\Bigr] + \frac{\partial}{\partial y_i}\Bigl[k_{ij} \bigl(\frac{\partial T^{(1)}}{\partial x_j} + \frac{\partial T^{(2)}}{\partial y_j}\bigr)\Bigr] + h = 0, \\
				&\!\frac{\partial}{\partial x_i}\Bigl[g_{ij} \bigl(\frac{\partial c^{(0)}}{\partial x_j} + \frac{\partial c^{(1)}}{\partial y_j}\bigr)\Bigr] + \frac{\partial}{\partial y_i}\Bigl[g_{ij} \bigl(\frac{\partial c^{(1)}}{\partial x_j} + \frac{\partial c^{(2)}}{\partial y_j}\bigr)\Bigr] + m = 0, \\
				&\!\frac{\partial}{\partial x_j}\Bigl[D_{ijkl} \bigl(\frac{\partial u_k^{(0)}}{\partial x_l} + \frac{\partial u_k^{(1)}}{\partial y_l}\bigr)\Bigr] + \frac{\partial}{\partial y_j}\Bigl[D_{ijkl} \bigl(\frac{\partial u_k^{(1)}}{\partial x_l} + \frac{\partial u_k^{(2)}}{\partial y_l}\bigr)\Bigr] \\
				&\!-\!\frac{\partial}{\partial x_j}\Bigl[D_{ijkl} \bigl( \alpha_{kl} T^{(0)} + \beta_{kl} c^{(0)}\bigr)\Bigr] -\frac{\partial}{\partial y_j}\Bigl[D_{ijkl} \bigl(\alpha_{kl} T^{(1)}\!+\!\beta_{kl} c^{(1)}\bigr)\Bigr]\!+\!f_i\!=\!\!0.
			\end{aligned}
		\end{cases}
	\end{equation}
	From $O({\epsilon^{-2}})$-order equations \eqref{eq:2.4}, we deduce that $T^{(0)}$, $c^{(0)}$ and $u_{i}^{(0)}$ are independent of the microscopic variable $\mathbf{y}$, namely
	\begin{equation}
		\label{eq:2.7}
		T^{(0)}(\mathbf{x}, \mathbf{y}) = T^{(0)}(\mathbf{x}),\;c^{(0)}(\mathbf{x}, \mathbf{y}) = c^{(0)}(\mathbf{x}), \;u_i^{(0)}(\mathbf{x},\mathbf{y}) = u_i^{(0)}(\mathbf{x}).
	\end{equation}
	After that, by substituting \eqref{eq:2.7} into $O({\epsilon^{-1}})$-order equations \eqref{eq:2.5} and employing the linearity of \eqref{eq:2.5}, the first-order correction terms $T^{(1)}$, $c^{(1)}$ and $u_{i}^{(1)}$ were found to have a linear dependence on the gradients of macroscopic homogenized fields $T^{(0)}$, $c^{(0)}$ and $u_{i}^{(0)}$ as given by
	\begin{equation}
		\label{eq:2.8}
		\begin{cases}
			\begin{aligned}
				T^{(1)}(\mathbf{x},\mathbf{y}) &= \mathcal{H}_{\alpha_1}(\mathbf{x},\mathbf{y}) \frac{\partial T^{(0)}(\mathbf{x})}{\partial x_{\alpha_1}}, \\
				c^{(1)}(\mathbf{x},\mathbf{y}) &= \mathcal{L}_{\alpha_1}(\mathbf{x},\mathbf{y}) \frac{\partial c^{(0)}(\mathbf{x})}{\partial x_{\alpha_1}}, \\
				u_i^{(1)}(\mathbf{x},\mathbf{y}) &= \mathcal{X}_{ih}^{\alpha_1}(\mathbf{x},\mathbf{y}) \frac{\partial u_h^{(0)}(\mathbf{x})}{\partial x_{\alpha_1}} - \mathcal{M}_i(\mathbf{x},\mathbf{y}) T^{(0)}(\mathbf{x}) - \mathcal{N}_i(\mathbf{x},\mathbf{y}) c^{(0)}(\mathbf{x}),
			\end{aligned}
		\end{cases}
	\end{equation}
	where $\mathcal{H}_{\alpha_1}$, $\mathcal{L}_{\alpha_1}$, $\mathcal{X}_{ih}^{\alpha_1}$, $\mathcal{M}_i$ and $\mathcal{N}_i$ are referred as the first-order cell functions defined in PUC $Y$.
	\begin{rmk}
		Crucially, the first-order cell functions exhibit quasi-periodicity with explicit dependence on the macroscopic coordinate $\mathbf{x}$, which acts as a varying parameter. This represents a key distinction from classical periodic composites.	
	\end{rmk}
	After combining \eqref{eq:2.7} and \eqref{eq:2.8} with $O({\epsilon^{-1}})$-order equations \eqref{eq:2.5}, simplification and calculation yield the following equations  subject to homogeneous Dirichlet boundary condition.
	\begin{equation}
		\label{eq:2.9}
		\begin{cases}
			\begin{aligned}
				& \frac{\partial}{\partial y_i}\Bigl(k_{ij}(\mathbf{x},\mathbf{y}) \frac{\partial \mathcal{H}_{\alpha_1}}{\partial y_j}\Bigr) = -\frac{\partial k_{i\alpha_1}(\mathbf{x},\mathbf{y})}{\partial y_i},\quad \mathbf{y} \in Y, \\
				& \mathcal{H}_{\alpha_1}(\mathbf{x},\mathbf{y}) = 0, \quad \mathbf{y} \in \partial Y.
			\end{aligned}
		\end{cases}
	\end{equation}
	\begin{equation}
		\label{eq:2.10}
		\begin{cases}
			\begin{aligned}
				& \frac{\partial}{\partial y_i}\Bigl(g_{ij}(\mathbf{x},\mathbf{y}) \frac{\partial \mathcal{L}_{\alpha_1}}{\partial y_j}\Bigr) = -\frac{\partial g_{i\alpha_1}(\mathbf{x},\mathbf{y})}{\partial y_i}, \quad \mathbf{y} \in Y, \\
				& \mathcal{L}_{\alpha_1}(\mathbf{x},\mathbf{y}) = 0, \quad \mathbf{y} \in \partial Y.
			\end{aligned}
		\end{cases}
	\end{equation}
	\begin{equation}
		\label{eq:2.11}
		\begin{cases}
			\begin{aligned}
				& \frac{\partial}{\partial y_j}\Bigl(D_{ijkl}(\mathbf{x},\mathbf{y}) \frac{\partial \mathcal{X}_{kh}^{\alpha_1}}{\partial y_l}\Bigr) = -\frac{\partial D_{ijh\alpha_1}(\mathbf{x},\mathbf{y})}{\partial y_j}, \quad \mathbf{y} \in Y, \\
				& \mathcal{X}_{kh}^{\alpha_1}(\mathbf{x},\mathbf{y}) = 0, \quad \mathbf{y} \in \partial Y.
			\end{aligned}
		\end{cases}
	\end{equation}
	\begin{equation}
		\label{eq:2.12}
		\begin{cases}
			\begin{aligned}
				& \frac{\partial}{\partial y_j}\Bigl(D_{ijkl}(\mathbf{x},\mathbf{y}) \frac{\partial \mathcal{M}_k}{\partial y_l}\Bigr) = -\frac{\partial \bigl(D_{ijkl}(\mathbf{x},\mathbf{y}) \alpha_{kl}(\mathbf{x},\mathbf{y})\bigr)}{\partial y_j}, \quad \mathbf{y} \in Y, \\
				& \mathcal{M}_k(\mathbf{x},\mathbf{y}) = 0, \quad \mathbf{y} \in \partial Y.
			\end{aligned}
		\end{cases}
	\end{equation}	
	\begin{equation}
		\label{eq:2.13}
		\begin{cases}
			\begin{aligned}
				& \frac{\partial}{\partial y_j}\Bigl(D_{ijkl}(\mathbf{x},\mathbf{y}) \frac{\partial \mathcal{N}_k}{\partial y_l}\Bigr) = -\frac{\partial \bigl(D_{ijkl}(\mathbf{x},\mathbf{y}) \beta_{kl}(\mathbf{x},\mathbf{y})\bigr)}{\partial y_j}, \quad \mathbf{y} \in Y, \\
				& \mathcal{N}_k(\mathbf{x},\mathbf{y}) = 0, \quad \mathbf{y} \in \partial Y.
			\end{aligned}
		\end{cases}
	\end{equation}
	Subsequently, we integrate both sides of $O({\epsilon^{0}})$-order equations \eqref{eq:2.6} over the unit cell $Y$ with respect to microscopic variable $\mathbf{y}$ and apply Gauss divergence theorem, yielding the macroscopic homogenized equations related with \eqref{eq:2.1} as below
	\begin{equation}
		\label{eq:2.14}
		\begin{cases}
			\begin{aligned}
				&-\frac{\partial}{\partial x_i}\Bigl(\hat{k}_{ij}(\mathbf{x}) \frac{\partial T^{(0)}}{\partial x_j}\Bigr) = h,\; \text{in }\; \Omega,\\
				&-\frac{\partial}{\partial x_i}\Bigl(\hat{g}_{ij}(\mathbf{x}) \frac{\partial c^{(0)}}{\partial x_j}\Bigr) = m,\; \text{in }\; \Omega,\\
				&-\frac{\partial}{\partial x_j}\Bigl(\hat{D}_{ijkl}(\mathbf{x}) \frac{\partial u_{k}^{(0)}}{\partial x_l}
				- \hat{A}_{ij}(\mathbf{x}) T^{(0)}
				- \hat{B}_{ij}(\mathbf{x}) c^{(0)}\Bigr) = f_i,\; \text{in }\;\Omega,\\
				& T^{(0)}(\mathbf{x}) = \overline{T}(\mathbf{x}),\;\text{on }\;\Gamma_{T},\\
				& \hat{k}_{ij}\bigl(\mathbf{x}\bigr) \frac{\partial T^{(0)}(\mathbf{x})}{\partial x_{j}}n_i=\overline{q}(\mathbf{x}),\;\text{on }\;\Gamma_{q},\\
				& c^{(0)}(\mathbf{x}) = \overline{c}(\mathbf{x}),\;\text{on }\;\Gamma_{c},\\
				& \hat{g}_{ij}\bigl(\mathbf{x}\bigr) \frac{\partial c^{(0)}(\mathbf{x})}{\partial x_{j}}n_i=\overline{d}(\mathbf{x}),\;\text{on }\;\Gamma_{d},\\
				& \bm{u}^{(0)}(\mathbf{x}) = \overline{\bm{u}}(\mathbf{x}),\;\text{on }\;\Gamma_{{u}},\\
				& \Bigl[\hat{D}_{ijkl}(\mathbf{x}) \frac{\partial u_{k}^{(0)}(\mathbf{x})}{\partial x_l}
				- \hat{A}_{ij}(\mathbf{x}) T^{(0)}(\mathbf{x})
				- \hat{B}_{ij}(\mathbf{x}) c^{(0)}(\mathbf{x})\Bigr] n_{j}=\overline{\sigma}_{i}(\mathbf{x}),\;\text{on }\;\Gamma_{\sigma},
			\end{aligned}
		\end{cases}
	\end{equation}
	where the homogenized material parameters at macro-scale are defined as below:
	\begin{equation}
		\label{eq:2.15}
		\begin{aligned}
			& \hat{k}_{ij}(\mathbf{x})
			= \frac{1}{|Y|} \int_{Y} \Bigl( k_{ij}(\mathbf{x}, \mathbf{y}) + k_{ik}(\mathbf{x}, \mathbf{y}) \frac{\partial \mathcal{H}_{j}}{\partial y_{k}} \Bigr) dY, \\
			& \hat{g}_{ij}(\mathbf{x})
			= \frac{1}{|Y|} \int_{Y} \Bigl( g_{ij}(\mathbf{x}, \mathbf{y}) + g_{ik}(\mathbf{x}, \mathbf{y}) \frac{\partial \mathcal{L}_{j}}{\partial y_{k}} \Bigr) dY, \\
			& \hat{D}_{ijkl}(\mathbf{x})
			= \frac{1}{|Y|} \int_{Y} \Bigl( D_{ijkl}(\mathbf{x}, \mathbf{y}) + D_{ijmn}(\mathbf{x}, \mathbf{y}) \frac{\partial \mathcal{X}_{mk}^{l}}{\partial y_{n}} \Bigr) dY, \\
			& \hat{A}_{ij}(\mathbf{x})
			= \frac{1}{|Y|} \int_{Y} \Bigl( D_{ijkl}(\mathbf{x}, \mathbf{y}) \alpha_{kl}(\mathbf{x}, \mathbf{y}) + D_{ijkl}(\mathbf{x}, \mathbf{y}) \frac{\partial \mathcal{M}_{k}}{\partial y_{l}} \Bigr) dY, \\
			& \hat{B}_{ij}(\mathbf{x})
			= \frac{1}{|Y|} \int_{\Gamma} \Bigl( D_{ijkl}(\mathbf{x}, \mathbf{y}) \beta_{kl}(\mathbf{x}, \mathbf{y}) + D_{ijkl}(\mathbf{x}, \mathbf{y}) \frac{\partial \mathcal{N}_{k}}{\partial y_{l}} \Bigr) dY.
		\end{aligned}
	\end{equation}
	\begin{rmk}
		Stemming from the approach as outlined in \cite{R23,R31,R32}, it can be proved that $\underline{\varsigma} | \bm{\xi} |^2\leq \hat{k}_{ij}(\mathbf{x}) \xi _i\xi_j \leq\overline{\varsigma} | \bm{\xi} |^2$, $\underline{\varsigma} | \bm{\xi} |^2\leq \hat{g}_{ij}(\mathbf{x}) \xi _i\xi_j \leq\overline{\varsigma} | \bm{\xi} |^2$, $\underline{\varsigma} \eta_{ij} \eta_{ij} \leq \hat{D}_{ijkl} (\mathbf{x}) \eta_{ij} \eta_{kl} \leq \overline{\varsigma} \eta_{ij} \eta_{ij}$, $\underline{\varsigma} | \bm{\xi} |^2\leq \hat{A}_{ij}(\mathbf{x}) \xi _i\xi_j \leq\overline{\varsigma} | \bm{\xi} |^2$ and $\underline{\varsigma} | \bm{\xi} |^2\leq \hat{B}_{ij}(\mathbf{x}) \xi _i\xi_j \leq\overline{\varsigma} | \bm{\xi} |^2$, where $\underline{\varsigma}$ and $\overline{\varsigma}$ are two positive constants independent of $\epsilon$.
	\end{rmk}
	
	Furthermore, substituting the terms $h(\mathbf{x})$, $m(\mathbf{x})$ and $f_i(\mathbf{x})$ in $O({\epsilon^{0}})$-order equations \eqref{eq:2.6} with their equivalent definitions in macroscopic homogenized equations \eqref{eq:2.14}, we can formulate the following equalities.
	\begin{equation}
		\label{eq:2.16}
		\begin{aligned}
			\frac{\partial}{\partial y_i} \Bigl( k_{ij} \frac{\partial T^{(2)}}{\partial y_j} \Bigr)
			& = \frac{\partial T^{(0)}}{\partial x_{\alpha_1}} \Bigl[ \frac{\partial \hat{k}_{i\alpha_1}}{\partial x_i}
			- \frac{\partial k_{i\alpha_1}}{\partial x_i} - \frac{\partial}{\partial x_i} \bigl( k_{ij} \frac{\partial \mathcal{H}_{\alpha_1}}{\partial y_j} \bigr) - \frac{\partial}{\partial y_i} \bigl( k_{ij} \frac{\partial \mathcal{H}_{\alpha_1}}{\partial x_j} \bigr) \Bigr] \\
			& + \frac{\partial^2 T^{(0)}}{\partial x_{\alpha_1} \partial x_{\alpha_2}} \Bigl[ \hat{k}_{\alpha_1\alpha_2} - k_{\alpha_1\alpha_2} - k_{\alpha_1j} \frac{\partial \mathcal{H}_{\alpha_2}}{\partial y_j}
			- \frac{\partial}{\partial y_i} \bigl( k_{i\alpha_2} \mathcal{H}_{\alpha_1} \bigr) \Bigr].
		\end{aligned}
	\end{equation}
	\begin{equation}
		\label{eq:2.17}
		\begin{aligned}
			\frac{\partial}{\partial y_i} \Bigl( g_{ij} \frac{\partial c^{(2)}}{\partial y_j} \Bigr)
			& = \frac{\partial c^{(0)}}{\partial x_{\alpha_1}} \Bigl[ \frac{\partial \hat{g}_{i\alpha_1}}{\partial x_i}
			- \frac{\partial g_{i\alpha_1}}{\partial x_i} - \frac{\partial}{\partial x_i} \bigl( g_{ij} \frac{\partial \mathcal{L}_{\alpha_1}}{\partial y_j} \bigr) - \frac{\partial}{\partial y_i} \bigl( g_{ij} \frac{\partial \mathcal{L}_{\alpha_1}}{\partial x_j} \bigr) \Bigr] \\
			& + \frac{\partial^2 c^{(0)}}{\partial x_{\alpha_1} \partial x_{\alpha_2}} \Bigl[ \hat{g}_{\alpha_1\alpha_2} - g_{\alpha_1\alpha_2} - g_{\alpha_1j} \frac{\partial \mathcal{L}_{\alpha_2}}{\partial y_j}
			- \frac{\partial}{\partial y_i} \bigl( g_{i\alpha_2} \mathcal{L}_{\alpha_1} \bigr) \Bigr].
		\end{aligned}
	\end{equation}
	\begin{equation}
		\label{eq:2.18}
		\begin{aligned}	
			\frac{\partial}{\partial y_{j}}\Bigl(D_{ijkl}\frac{\partial u_{k}^{(2)}}{\partial y_{l}}\Bigr)
			&\!=\!-\frac{\partial u_{h}^{(0)}}{\partial x_{\alpha_{1}}}\!\Bigl[\!\frac{\partial D_{ijh\alpha_{1}}}{\partial x_{j}}\!-\!\frac{\partial\hat{D}_{ijh\alpha_{1}}}{\partial x_{j}} \!+\!\frac{\partial}{\partial x_{j}}\bigl(D_{ijkl}\frac{\partial\mathcal{X}_{kh}^{\alpha_{1}}}{\partial y_{l}}\bigr) \!+\!\frac{\partial}{\partial y_{j}}\bigl(D_{ijkl}\frac{\partial\mathcal{X}_{kh}^{\alpha_{1}}}{\partial x_{l}}\bigr)\!\Bigr]\! \\
			& \!-\!\frac{\partial^2u_h^{(0)}}{\partial x_{\alpha_1}\partial x_{\alpha_2}}\Bigl[D_{i\alpha_1h\alpha_2} -\hat{D}_{i\alpha_1h\alpha_2} +D_{i\alpha_1kl}\frac{\partial\mathcal{X}_{kh}^{\alpha_2}}{\partial y_l} +\frac{\partial}{\partial y_j}\bigl(D_{ijk\alpha_2}\mathcal{X}_{kh}^{\alpha_1}\bigr)\Bigr] \\
			& \!+\!T^{(0)}\Bigl[\frac{\partial}{\partial x_{j}}\bigl(D_{ijkl}\frac{\partial \mathcal{M}_{k}}{\partial y_{l}} +D_{ijkl}\alpha_{kl} -\hat{A}_{ij}\bigr) +\frac{\partial}{\partial y_{j}}\bigl(D_{ijkl}\frac{\partial \mathcal{M}_{k}}{\partial x_{l}}\bigr)\Bigr]\\
			& \!+\!\frac{\partial T^{(0)}}{\partial x_{\alpha_1}}\!\Bigl[\!D_{i\alpha_1kl}\bigl(\frac{\partial \mathcal{M}_{k}}{\partial y_l}\!+\!\alpha_{kl}\bigr)\!-\!\hat{A}_{i\alpha_1}\!+\!\frac{\partial}{\partial y_{j}}\bigl(D_{ijk\alpha_{1}}\mathcal{M}_{k}\!+\!D_{ijkl}\alpha_{kl}\mathcal{H}_{\alpha_{1}}\bigr)\!\Bigr]\! \\
			& \!+\!c^{(0)}\Bigl[\frac{\partial}{\partial x_{j}}\bigl(D_{ijkl}\frac{\partial \mathcal{N}_{k}}{\partial y_{l}}+D_{ijkl}\beta_{kl}-\hat{B}_{ij}\bigr) +\frac{\partial}{\partial y_{j}}\bigl(D_{ijkl}\frac{\partial \mathcal{N}_{k}}{\partial x_{l}}\bigr)\Bigr] \\
			& \!+\!\frac{\partial c^{(0)}}{\partial x_{\alpha_{1}}}\!\Bigl[\!D_{i\alpha_{1}kl}\bigl(\frac{\partial \mathcal{N}_{k}}{\partial y_{l}}\!+\!\beta_{kl}\bigr) \!-\!\hat{B}_{i\alpha_{1}}\!+\!\frac{\partial}{\partial y_{j}}\bigl(D_{ijk\alpha_{1}}\mathcal{N}_{k}\!+\!D_{ijkl}\beta_{kl}\mathcal{L}_{\alpha_{1}}\bigr)\!\Bigr]\!.
		\end{aligned}
	\end{equation}
	Based on equalities \eqref{eq:2.16}-\eqref{eq:2.18}, we can establish the detailed expressions for second-order correction terms as below
	\begin{equation}
		\label{eq:2.19}
		\begin{cases}
			\begin{aligned}
				T^{(2)}(\mathbf{x},\mathbf{y}) &=\mathcal{H}_{\alpha_1\alpha_2}(\mathbf{x},\mathbf{y})\frac{\partial^2T^{(0)}(\mathbf{x})}{\partial x_{\alpha_1}\partial x_{\alpha_2}}+\mathcal{R}_{\alpha_1}(\mathbf{x},\mathbf{y})\frac{\partial T^{(0)}(\mathbf{x})}{\partial x_{\alpha_1}}, \\			
				c^{(2)}(\mathbf{x},\mathbf{y}) &=\mathcal{L}_{\alpha_1\alpha_2}(\mathbf{x},\mathbf{y})\frac{\partial^2c^{(0)}(\mathbf{x})}{\partial x_{\alpha_1}\partial x_{\alpha_2}}+\mathcal{S}_{\alpha_1}(\mathbf{x},\mathbf{y})\frac{\partial c^{(0)}(\mathbf{x})}{\partial x_{\alpha_1}}, \\
				u_i^{(2)}(\mathbf{x},\mathbf{y}) &=\mathcal{P}_{ih}^{\alpha_1\alpha_2}(\mathbf{x},\mathbf{y})\frac{\partial^2u_h^{(0)}(\mathbf{x})}{\partial x_{\alpha_1}\partial x_{\alpha_2}}+\mathcal{Q}_{ih}^{\alpha_1}(\mathbf{x},\mathbf{y})\frac{\partial u_h^{(0)}(\mathbf{x})}{\partial x_{\alpha_1}}+\mathcal{W}_i(\mathbf{x},\mathbf{y})T^{(0)}(\mathbf{x}) \\
				& +\mathcal{Z}_i^{\alpha_1}(\mathbf{x},\mathbf{y})\frac{\partial T^{(0)}(\mathbf{x})}{\partial x_{\alpha_1}}+\mathcal{F}_i(\mathbf{x},\mathbf{y})c^{(0)}(\mathbf{x})+\mathcal{G}_i^{\alpha_1}(\mathbf{x},\mathbf{y})\frac{\partial c^{(0)}(\mathbf{x})}{\partial x_{\alpha_1}},
			\end{aligned}
		\end{cases}
	\end{equation}
	where $\mathcal{H}_{\alpha_1\alpha_2}$, $\mathcal{R}_{\alpha_1}$, $\mathcal{L}_{\alpha_1\alpha_2}$, $\mathcal{S}_{\alpha_1}$, $\mathcal{P}_{ih}^{\alpha_1\alpha_2}$, $\mathcal{Q}_{ih}^{\alpha_1}$, $\mathcal{W}_i$, $\mathcal{Z}_i^{\alpha_1}$, $\mathcal{F}_i$ and $\mathcal{G}_i^{\alpha_1}$ are the second-order cell functions defined in PUC $Y$.
	
	Afterwards, the substitution of \eqref{eq:2.19} into \eqref{eq:2.16}-\eqref{eq:2.18} yields a series of equations attaching homogeneous Dirichlet boundary conditions for solving second-order cell functions respectively.
	\begin{equation}
		\label{eq:2.20}
		\begin{cases}
			\begin{aligned}
				& \frac{\partial}{\partial y_i}\Bigl(k_{ij}(\mathbf{x},\mathbf{y})\frac{\partial \mathcal{H}_{\alpha_1\alpha_2}}{\partial y_j}\Bigr)=\hat{k}_{\alpha_1\alpha_2}-k_{\alpha_1\alpha_2}-k_{\alpha_1j}\frac{\partial \mathcal{H}_{\alpha_2}}{\partial y_j} -\frac{\partial}{\partial y_i}\bigl(k_{i\alpha_2}\mathcal{H}_{\alpha_1}\bigr),\quad \mathbf{y}\in Y, \\
				&\mathcal{H}_{\alpha_1\alpha_2}(\mathbf{x},\mathbf{y})=0,\quad \mathbf{y}\in\partial Y.
			\end{aligned}
		\end{cases}
	\end{equation}
	\begin{equation}
		\label{eq:2.21}
		\begin{cases}
			\begin{aligned}
				& \frac{\partial}{\partial y_i}\Bigl(k_{ij}(\mathbf{x},\mathbf{y})\frac{\partial \mathcal{R}_{\alpha_1}}{\partial y_j}\Bigr)\!=\!\frac{\partial\hat{k}_{i\alpha_1}}{\partial x_i}\!-\!\frac{\partial k_{i\alpha_1}}{\partial x_i}\!-\!\frac{\partial}{\partial x_i}\bigl(k_{ij}\frac{\partial \mathcal{H}_{\alpha_1}}{\partial y_j}\bigr) \!-\!\frac{\partial}{\partial y_i}\bigr(k_{ij}\frac{\partial \mathcal{H}_{\alpha_1}}{\partial x_j}\bigr),\;\; \mathbf{y}\in Y, \\
				& \mathcal{R}_{\alpha_1}(\mathbf{x},\mathbf{y})=0,\quad \mathbf{y}\in\partial Y.
			\end{aligned}
		\end{cases}
	\end{equation}
	\begin{equation}
		\label{eq:2.22}
		\begin{cases}
			\begin{aligned}
				& \frac{\partial}{\partial y_i}\Bigl(g_{ij}(\mathbf{x},\mathbf{y})\frac{\partial \mathcal{L}_{\alpha_1\alpha_2}}{\partial y_j}\Bigr)=\hat{g}_{\alpha_1\alpha_2}-g_{\alpha_1\alpha_2}-g_{\alpha_1j}\frac{\partial \mathcal{L}_{\alpha_2}}{\partial y_j} -\frac{\partial}{\partial y_i}\bigl(g_{i\alpha_2}\mathcal{L}_{\alpha_1}\bigl),\quad \mathbf{y}\in Y,\\
				& \mathcal{L}_{\alpha_1\alpha_2}(\mathbf{x},\mathbf{y}){=}0,\quad\mathbf{y}\in\partial Y.
			\end{aligned}
		\end{cases}
	\end{equation}
	\begin{equation}
		\label{eq:2.23}
		\begin{cases}
			\begin{aligned}
				& \frac{\partial}{\partial y_i}\Bigl(g_{ij}(\mathbf{x},\mathbf{y})\frac{\partial \mathcal{S}_{\alpha_1}}{\partial y_j}\Bigr)\!=\!\frac{\partial\hat{g}_{i\alpha_1}}{\partial x_i}\!-\!\frac{\partial g_{i\alpha_1}}{\partial x_i}\!-\!\frac{\partial}{\partial x_i}\bigl(g_{ij}\frac{\partial \mathcal{L}_{\alpha_1}}{\partial y_j}\bigr) \!-\!\frac{\partial}{\partial y_i}\bigl(g_{ij}\frac{\partial \mathcal{L}_{\alpha_i}}{\partial x_j}\bigr),\;\; \mathbf{y}\in Y, \\
				& \mathcal{S}_{\alpha_1}(\mathbf{x},\mathbf{y})=0,\quad \mathbf{y}\in\partial Y.
			\end{aligned}
		\end{cases}
	\end{equation}
	\begin{equation}
		\label{eq:2.24}
		\begin{cases}
			\begin{aligned}
				& \!\frac{\partial}{\partial y_j}\!\Bigl(\!D_{ijkl}(\mathbf{x},\mathbf{y})\frac{\partial \mathcal{P}_{kh}^{\alpha_1\alpha_2}}{\partial y_l}\!\Bigr)\!\!=\!\hat{D}_{i\alpha_1h\alpha_2}\!\!-\!D_{i\alpha_1h\alpha_2}\!\!-\!D_{i\alpha_1kl}\frac{\partial\mathcal{X}_{kh}^{\alpha_2}}{\partial y_l}\!-\!\frac{\partial}{\partial y_j}\!\bigl(\!D_{ijk\alpha_2}\mathcal{X}_{kh}^{\alpha_1}\!\bigr),\; \mathbf{y}\!\in\!Y\!,\\
				& \!\mathcal{P}_{kh}^{\alpha_{k}\alpha_{2}}(\mathbf{x},\mathbf{y})=0,\quad\mathbf{y}\in\partial Y.
			\end{aligned}
		\end{cases}
	\end{equation}
	\begin{equation}
		\label{eq:2.25}
		\begin{small}
			\begin{cases}
				\begin{aligned}
					&\frac{\partial}{\partial y_j}\!\Bigl(\!D_{ijkl}(\mathbf{x},\mathbf{y})\frac{\partial \mathcal{Q}_{kh}^{\alpha_1}}{\partial y_l}\!\Bigr)\!\!=\!\frac{\partial\hat{D}_{ijha_1}}{\partial x_j}\!-\!\frac{\partial D_{ijha_1}}{\partial x_j} \!-\!\frac{\partial}{\partial x_j}\!\Bigl(\!D_{ijkl}\frac{\partial\mathcal{X}_{kh}^{\alpha_1}}{\partial y_l}\!\Bigr)\!\!-\!\frac{\partial}{\partial y_j}\!\Bigl(\!D_{ijkl}\frac{\partial\mathcal{X}_{kh}^{\alpha_1}}{\partial x_l}\!\Bigr)\!,\; \mathbf{y}\!\in\!Y\!,\! \\
					& \mathcal{Q}_{kh}^{\alpha_1}(\mathbf{x},\mathbf{y})=0,\quad\mathbf{y}\in\partial Y.
				\end{aligned}
			\end{cases} 		
		\end{small}
	\end{equation}
	\begin{equation}
		\label{eq:2.26}
		\begin{cases}
			\begin{aligned}
				&\frac{\partial}{\partial y_j}\!\Bigl(\!D_{ijkl}(\mathbf{x},\mathbf{y})\frac{\partial \mathcal{W}_k}{\partial y_l}\Bigr)\!=\!
				\frac{\partial}{\partial y_j}\!\Bigl(\!D_{ijkl}\frac{\partial \mathcal{M}_k}{\partial x_l}\!\Bigr)\!\!+\!\frac{\partial}{\partial x_j}\!\Bigl(\!D_{ijkl}\bigl(\frac{\partial \mathcal{M}_k}{\partial y_l}\!+\!\alpha_{kl} \bigr)\!-\!\hat{A}_{ij}\!\Bigr)\!,\; \mathbf{y}\in Y, \\
				&\mathcal{W}_k(\mathbf{x},\mathbf{y})=0,\quad \mathbf{y}\in\partial Y.
			\end{aligned}
		\end{cases}
	\end{equation}
	\begin{equation}
		\label{eq:2.27}
		\begin{small}
			\begin{cases}
				\begin{aligned}
					&\!\frac{\partial}{\partial y_j}\!\Bigl(\!D_{ijkl}(\mathbf{x},\mathbf{y})\frac{\partial \mathcal{Z}_k^{\alpha_1}}{\partial y_l}\!\Bigr)\!\!=\!D_{i\alpha_1kl}\!\Bigl(\!\frac{\partial \mathcal{M}_k}{\partial y_l}\!+\!\alpha_{kl}\!\Bigr)\!\!-\!\hat{A}_{i\alpha_1} \!+\!\frac{\partial}{\partial y_j}\!\bigl(\!D_{ijk\alpha_1}\mathcal{M}_k \!+\!D_{ijkl}\alpha_{kl}\mathcal{H}_{\alpha_1}\!\bigr)\!,\;\!\mathbf{y}\!\in\!Y\!,\!\\
					&\!\mathcal{Z}_k^{\alpha_1}(\mathbf{x},\mathbf{y})=0,\quad \mathbf{y}\in\partial Y.
				\end{aligned}
			\end{cases}
		\end{small}
	\end{equation}
	\begin{equation}
		\label{eq:2.28}
		\begin{cases}
			\begin{aligned}
				& \frac{\partial}{\partial y_j}\!\Bigl(\!D_{ijkl}(\mathbf{x},\mathbf{y})\frac{\partial \mathcal{F}_k}{\partial y_l}\!\Bigr)\!\!=\!
				\frac{\partial}{\partial y_j}\Bigl(D_{ijkl}\frac{\partial \mathcal{N}_k}{\partial x_l}\Bigr) \!+\!\frac{\partial}{\partial x_j}\Bigl(D_{ijkl}\bigl( \frac{\partial \mathcal{N}_k}{\partial y_l}\!+\!\beta_{kl}\bigr)\!-\!\hat{B}_{ij}\!\Bigr)\!,\;\; \mathbf{y}\in Y, \\
				& \mathcal{F}_k(\mathbf{x},\mathbf{y})=0,\quad \mathbf{y}\in\partial Y.
			\end{aligned}
		\end{cases}
	\end{equation}
	\begin{equation}
		\label{eq:2.29}
		\begin{small}
			\begin{cases}
				\begin{aligned}
					&\frac{\partial}{\partial y_j}\!\Bigl(\!D_{ijkl}(\mathbf{x},\mathbf{y})\frac{\partial \mathcal{G}_k^{\alpha_1}}{\partial y_l}\!\Bigr)\!\!=\!D_{i\alpha_1kl}\!\Bigl(\!\frac{\partial \mathcal{N}_k}{\partial y_l}\!+\!\beta_{kl}\!\Bigr)\!\!-\!\hat{B}_{i\alpha_1} \!+\!\frac{\partial}{\partial y_j}\!\bigl(\!D_{ijk\alpha_1}\mathcal{N}_k\!+\!D_{ijkl}\beta_{kl}\mathcal{L}_{\alpha_1}\!\bigr)\!,\;\!\mathbf{y}\!\in\!Y\!,\!\\
					&\mathcal{G}_k^{\alpha_1}(\mathbf{x},\mathbf{y})=0, \quad \mathbf{y}\in\partial Y.
				\end{aligned}
			\end{cases}
		\end{small}
	\end{equation}
	\begin{rmk}	
		According to Refs. \cite{R23,R31,R32}, the homogeneous Dirichlet boundary condition  may substitute for the classical periodic boundary condition for the first-order cell problems \eqref{eq:2.9}-\eqref{eq:2.13} and second-order cell problems \eqref{eq:2.20}-\eqref{eq:2.29}, when the material property parameters satisfy geometric symmetry and regularity assumptions. 
	\end{rmk}
	\begin{rmk}
		By the Lax-Milgram theorem and assumption (B), the auxiliary cell problems \eqref{eq:2.9}-\eqref{eq:2.13} and \eqref{eq:2.20}-\eqref{eq:2.29} admit a unique solution for arbitrary macroscopic coordinates $\mathbf{x}$.
	\end{rmk}
	In conclusion, we establish the LOMS solutions for multi-scale hygro-thermo-mechanical problems \eqref{eq:2.1} as below
	\begin{equation}
		\label{eq:2.30}
		\begin{aligned}
			T^{(1,\epsilon)}(\mathbf{x})&=T^{(0)}(\mathbf{x},\mathbf{y}) + \epsilon T^{(1)}(\mathbf{x},\mathbf{y}) \\
			&=T^{(0)}(\mathbf{x})+\epsilon \mathcal{H}_{\alpha_1}(\mathbf{x},\mathbf{y})\frac{\partial T^{(0)}(\mathbf{x})}{\partial x_{\alpha_1}}.
		\end{aligned}
	\end{equation}
	\begin{equation}
		\label{eq:2.31}
		\begin{aligned}
			c^{(1,\epsilon)}(\mathbf{x})& = c^{(0)}(\mathbf{x},\mathbf{y}) + \epsilon c^{(1)}(\mathbf{x},\mathbf{y}) \\
			&=c^{(0)}(\mathbf{x})+\epsilon \mathcal{L}_{\alpha_1}(\mathbf{x},\mathbf{y})\frac{\partial c^{(0)}(\mathbf{x})}{\partial x_{\alpha_1}}.
		\end{aligned}
	\end{equation}
	\begin{equation}
		\label{eq:2.32}
		\begin{aligned}
			u_{i}^{(1,\epsilon)}(\mathbf{x})&= u_{i}^{(0)}(\mathbf{x},\mathbf{y}) + \epsilon u_{i}^{(1)}(\mathbf{x},\mathbf{y})\\
			&=u_{i}^{(0)}(\mathbf{x})\!+\!\epsilon \Bigl(\mathcal{X}_{ih}^{\alpha_{1}}(\mathbf{x},\mathbf{y})\frac{\partial u_{h}^{(0)}(\mathbf{x})}{\partial x_{\alpha_{1}}}\!-\!\mathcal{M}_{i}(\mathbf{x},\mathbf{y})T^{(0)}(\mathbf{x})\!-\!\mathcal{N}_{i}(\mathbf{x},\mathbf{y})c^{(0)}(\mathbf{x})\Bigr). \end{aligned}
	\end{equation}
	Furthermore, the HOMS solutions for multi-scale problem \eqref{eq:2.1} as follows
	\begin{equation}
		\label{eq:2.33}
		\begin{aligned}
			T^{(2,\epsilon)}(\mathbf{x})&=T^{(0)}(\mathbf{x},\mathbf{y}) + \epsilon T^{(1)}(\mathbf{x},\mathbf{y}) + \epsilon^{2} T^{(2)}(\mathbf{x},\mathbf{y})\\
			&=T^{(0)}(\mathbf{x})+\epsilon \mathcal{H}_{\alpha_1}(\mathbf{x},\mathbf{y})\frac{\partial T^{(0)}(\mathbf{x})}{\partial x_{\alpha_1}}\\
			&+\epsilon^2\Bigl(\mathcal{H}_{\alpha_1\alpha_2}(\mathbf{x},\mathbf{y})\frac{\partial^2T^{(0)}(\mathbf{x})}{\partial x_{\alpha_1}\partial x_{\alpha_2}}+\mathcal{R}_{\alpha_1}(\mathbf{x},\mathbf{y})\frac{\partial T^{(0)}(\mathbf{x})}{\partial x_{\alpha_1}}\Bigr).
		\end{aligned}
	\end{equation}
	\begin{equation}
		\label{eq:2.34}
		\begin{aligned}
			c^{(2,\epsilon)}(\mathbf{x})& = c^{(0)}(\mathbf{x},\mathbf{y}) + \epsilon c^{(1)}(\mathbf{x},\mathbf{y}) + \epsilon^{2} c^{(2)}(\mathbf{x},\mathbf{y})\\
			&=c^{(0)}(\mathbf{x})+\epsilon \mathcal{L}_{\alpha_1}(\mathbf{x},\mathbf{y})\frac{\partial c^{(0)}(\mathbf{x})}{\partial x_{\alpha_1}}\\
			&+\epsilon^2\Bigl(\mathcal{L}_{\alpha_1\alpha_2}(\mathbf{x},\mathbf{y})\frac{\partial^2c^{(0)}(\mathbf{x})}{\partial x_{\alpha_1}\partial x_{\alpha_2}}+\mathcal{S}_{\alpha_1}(\mathbf{x},\mathbf{y})\frac{\partial c^{(0)}(\mathbf{x})}{\partial x_{\alpha_1}}\Bigr).
		\end{aligned}
	\end{equation}
	\begin{equation}
		\label{eq:2.35}
		\begin{aligned}
			u_{i}^{(2,\epsilon)}(\mathbf{x})&= u_{i}^{(0)}(\mathbf{x},\mathbf{y}) + \epsilon u_{i}^{(1)}(\mathbf{x},\mathbf{y}) + \epsilon^{2} u_{i}^{(2)}(\mathbf{x},\mathbf{y})\\
			&=u_{i}^{(0)}(\mathbf{x})+\epsilon \Bigl(\mathcal{X}_{ih}^{\alpha_{1}}(\mathbf{x},\mathbf{y})\frac{\partial u_{h}^{(0)}(\mathbf{x})}{\partial x_{\alpha_{1}}}-\!\mathcal{M}_{i}(\mathbf{x},\mathbf{y})T^{(0)}(\mathbf{x})-\!\mathcal{N}_{i}(\mathbf{x},\mathbf{y})c^{(0)}(\mathbf{x})\Bigr) \\
			&+\epsilon^{2}\Bigl(\mathcal{P}_{ih}^{\alpha_{1}\alpha_{2}}(\mathbf{x},\mathbf{y})\frac{\partial^{2}u_{h}^{(0)}(\mathbf{x})}{\partial x_{\alpha_{1}}\partial x_{\alpha_{2}}}+\mathcal{Q}_{ih}^{\alpha_{1}}(\mathbf{x},\mathbf{y})\frac{\partial u_{h}^{(0)}(\mathbf{x})}{\partial x_{\alpha_{1}}}+\mathcal{W}_{i}(\mathbf{x},\mathbf{y})T^{(0)}(\mathbf{x}) \\
			& +\mathcal{Z}_{i}^{\alpha_{1}}(\mathbf{x},\mathbf{y})\frac{\partial T^{(0)}(\mathbf{x})}{\partial x_{\alpha_{1}}}+\mathcal{F}_{i}(\mathbf{x},\mathbf{y})c^{(0)}(\mathbf{x})+\mathcal{G}_{i}^{\alpha_{1}}(\mathbf{x},\mathbf{y})\frac{\partial c^{(0)}(\mathbf{x})}{\partial x_{\alpha_{1}}}\Bigr).
		\end{aligned}
	\end{equation}
	
	In the engineering damage assessment of composite materials, engineers and scientists developed a simplified model to quantify the damage degree of composites. This simplified model introduces a scalar damage parameter $\omega(\mathbf{x})$ to evaluate the damage degree and assumes that the material parameters possess the macro-micro separation properties \cite{R33,R34,R35}, namely
	\begin{equation}
		\label{eq:2.36}
		\begin{aligned}
			& k_{ij}(\mathbf{x},\mathbf{y})=\omega(\mathbf{x})k_{ij}^*(\mathbf{y}),\quad g_{ij}(\mathbf{x},\mathbf{y})=\omega(\mathbf{x})g_{ij}^*(\mathbf{y}),\quad D_{ijkl}(\mathbf{x},\mathbf{y})=\omega(\mathbf{x})D_{ijkl}^*(\mathbf{y}), \\
			& \alpha_{kl}(\mathbf{x},\mathbf{y})=\omega(\mathbf{x})\alpha_{kl}^*(\mathbf{y}),\quad \beta_{kl}(\mathbf{x},\mathbf{y})=\omega(\mathbf{x})\beta_{kl}^*(\mathbf{y}).
		\end{aligned}
	\end{equation}
	Based on this new model, engineers can employ a single unit cell as basic unit for damage evaluation and reduce the computational cost of multi-scale simulation of quasi-periodic composite structures. Furthermore, the HOMS asymptotic solutions for the multi-scale problem \eqref{eq:2.1} under this novel model are detailed in Appendix A.
	
	\section{The error analyses of multi-scale asymptotic solutions}
	\label{sec:3}
	This section conducts rigorous point-wise and integral error analyses pertaining to multi-scale asymptotic solutions. Before giving the detailed error analyses, we define the residual functions $T_\Delta^{(1,\epsilon)}$, $c_\Delta^{(1,\epsilon)}$ and $u_{\Delta i}^{(1,\epsilon)}$ for LOMS solutions as below
	\begin{equation}
		\label{eq:3.1}
		\begin{aligned}
			&T_\Delta^{(1,\epsilon)}(\mathbf{x})=T^\epsilon-T^{(1,\epsilon)},\;\;c_\Delta^{(1,\epsilon)}(\mathbf{x})=c^\epsilon-c^{(1,\epsilon)},\;\;u_{\Delta i}^{(1,\epsilon)}(\mathbf{x})=u_i^\epsilon-u_i^{(1,\epsilon)}.
		\end{aligned}
	\end{equation}
	Additionally, the residual functions $T_\Delta^{(2,\epsilon)}$, $c_\Delta^{(2,\epsilon)}$ and $u_{\Delta i}^{(2,\epsilon)}$ are denoted for HOMS solutions as below
	\begin{equation}
		\label{eq:3.2}
		\begin{aligned}
			&T_\Delta^{(2,\epsilon)}(\mathbf{x})=T^\epsilon-T^{(2,\epsilon)},\;\;c_\Delta^{(2,\epsilon)}(\mathbf{x})=c^\epsilon-c^{(2,\epsilon)},\;\;u_{\Delta i}^{(2,\epsilon)}(\mathbf{x})=u_i^\epsilon-u_i^{(2,\epsilon)}.
		\end{aligned}
	\end{equation}
	
	\subsection{The error analysis in the point-wise sense}
	\label{sec:31}
	Firstly, substituting the residual functions \eqref{eq:3.1} into the multi-scale equations \eqref{eq:2.1} yields the residual equations for LOMS solutions presented below
	\begin{equation}
		\label{eq:3.3}
		\begin{cases}
			\begin{aligned}
				& -\frac{\partial}{\partial x_i}\Bigl(k_{ij}(\mathbf{x}, \mathbf{y}) \frac{\partial T_{\Delta}^{(1,\epsilon)}}{\partial x_j}\Bigr)
				= \mathcal{A}_0(\mathbf{x}, \mathbf{y}) + \epsilon \mathcal{A}_1(\mathbf{x}, \mathbf{y}), & \text{in } \Omega, \\
				& -\frac{\partial}{\partial x_i}\Bigl(g_{ij}(\mathbf{x}, \mathbf{y}) \frac{\partial c_{\Delta}^{(1,\epsilon)}}{\partial x_j}\Bigr)
				= \mathcal{B}_0(\mathbf{x}, \mathbf{y}) + \epsilon \mathcal{B}_1(\mathbf{x}, \mathbf{y}), & \text{in } \Omega, \\
				& -\frac{\partial}{\partial x_j}\Bigl[D_{ijkl}(\mathbf{x}, \mathbf{y})
				\bigl(\frac{\partial u_{\Delta k}^{(1,\epsilon)}}{\partial x_l}
				-\alpha_{kl}(\mathbf{x}, \mathbf{y}) T_{\Delta}^{(1,\epsilon)}
				-\beta_{kl}(\mathbf{x}, \mathbf{y}) c_{\Delta}^{(1,\epsilon)}\bigr)\Bigr] \\
				& = \mathcal{C}_{0i}(\mathbf{x}, \mathbf{y}) + \epsilon \mathcal{C}_{1i}(\mathbf{x}, \mathbf{y}), & \text{in } \Omega.
			\end{aligned}
		\end{cases}
	\end{equation}
	Secondly, by substituting the residual functions \eqref{eq:3.2} into multi-scale equations \eqref{eq:2.1}, we derive the residual equations for HOMS solutions as below
	\begin{equation}
		\label{eq:3.4}
		\begin{cases}
			\begin{aligned}
				&-\frac{\partial}{\partial x_i} \Bigl( k_{ij}(\mathbf{x}, \mathbf{y}) \frac{\partial T^{(2,\epsilon)}_{\Delta}}{\partial x_j} \Bigr)
				= \epsilon \mathcal{U}(\mathbf{x}, \mathbf{y}), & \text{in } \Omega, \\
				&-\frac{\partial}{\partial x_i} \Bigl( g_{ij}(\mathbf{x}, \mathbf{y}) \frac{\partial c^{(2,\epsilon)}_{\Delta}}{\partial x_j} \Bigr)
				= \epsilon \mathcal{V}(\mathbf{x}, \mathbf{y}), & \text{in } \Omega, \\
				&-\frac{\partial}{\partial x_j} \!\Bigl[\! D_{ijkl}(\mathbf{x}, \mathbf{y})
				\bigl( \frac{\partial u^{(2,\epsilon)}_{\Delta k}}{\partial x_l}
				\!-\!\alpha_{kl}(\mathbf{x}, \mathbf{y}) T^{(2,\epsilon)}_{\Delta}
				\!-\!\beta_{kl}(\mathbf{x}, \mathbf{y}) c^{(2,\epsilon)}_{\Delta}\bigr) \!\Bigr]\!\!=\!\epsilon \mathcal{J}_i(\mathbf{x}, \mathbf{y}),& \text{in } \Omega.
			\end{aligned}
		\end{cases}
	\end{equation}
	In residual equations \eqref{eq:3.3} and \eqref{eq:3.4}, the specific expressions of functions $ \mathcal{A}_0$, $\mathcal{A}_1$, $\mathcal{B}_0$, $\mathcal{B}_1$, $\mathcal{C}_{0i}$, $\mathcal{C}_{1i}$, $\mathcal{U}$, $\mathcal{V}$ and $\mathcal{J}_i$ are exhibited in Appendix B of the present study because of their lengthy forms.
	
	We now present the principal conclusions of the point-wise error analysis: the residual equations \eqref{eq:3.3} demonstrate that the LOMS solutions exhibit $O(1)$-order residuals, whereas the residual equations \eqref{eq:3.4} clearly illustrate that the HOMS solutions achieve superior $O(\epsilon)$-order residuals. Since the $\epsilon$-independent terms $\mathcal{A}_0$, $\mathcal{B}_0$ and $\mathcal{C}_{0i}$ in \eqref{eq:3.3} do not vanish as $\epsilon \to 0$, the LOMS solutions fail to preserve the locally physical balance of multi-scale governing equations. In contrast, through introducing the higher-order correction terms, the HOMS solutions rigorously maintain the locally physical balance in the original governing equations \eqref{eq:2.1} while attaining $O(\epsilon)$-order point-wise convergence. Consequently, even for small but finite $\epsilon$, the HOMS solutions maintain sufficient accuracy for engineering simulations while capturing microscopic oscillatory behaviors in quasi-periodic composites. This forms the primary impetus for developing the proposed HOMS methodology.
	
	\subsection{The error analysis in the integral sense}
	\label{sec:32}
	In order to obtain the optimal error estimations in the integral sense, we postulate three fundamental assumptions regarding the multi-scale problem \eqref{eq:2.1} as follows:
	\begin{enumerate}
		\item[(i)] The domain $\Omega$ is the union of the entire periodic cells, expressed as $\bar{\Omega}=\cup_{\mathbf{z}\in I_\epsilon}\epsilon(\mathbf{z}+\bar{Y})$ with the index set $I_{\epsilon}=\{\mathbf{z}=(z_{1},\ldots,z_{n})\in Z^{n},\epsilon(\mathbf{z}+\bar{Y})\subset\bar{\Omega}\}$. 
		\item[(ii)] Let $\Delta_1,\cdots,\Delta_n(n=2,3)$ denote the middle hyperplanes of the PUC $Y$. We assume reference cell $Y$ are symmetric with respect to $\Delta_1,\cdots,\Delta_n$.
		\item[(iii)] The error estimations for the HOMS solutions are presented for multi-scale H-T-M problems with pure Dirichlet boundary conditions.
	\end{enumerate}
	\begin{lemma}
		\label{lem:3.1}
		Defining three differential operators $\displaystyle\sigma_{TY}(\chi)=n_i k_{ij}(\mathbf{x},\mathbf{y})\frac{\partial \chi}{\partial y_j}$,  $\displaystyle\sigma_{cY}(\chi)=n_i g_{ij}(\mathbf{x},\mathbf{y})\frac{\partial \chi}{\partial y_j}$ and $\displaystyle\sigma_{iY}(\bm{\phi})=n_jD_{ijkl}(\mathbf{x},\mathbf{y})\frac{\partial \phi_{k}}{\partial x_{l}}$, and employing the identical approach in Refs. \cite{R23,R31,R32}, then the normal derivatives of all auxiliary cell functions are continuous on the boundary of PUC $Y$ on the basis of foregoing assumptions (A)-(B) and (ii).
	\end{lemma}
	\begin{thm}
		\label{thm:1}
		Let $T^{\epsilon}$, $c^{\epsilon}$ and $\bm{u}^{\epsilon}$ be the weak solutions of the H-T-M coupling problems \eqref{eq:2.1}, and let $T^{(0)}$, $c^{(0)}$ and $\bm{u}^{(0)}$ be the solutions of the corresponding homogenized problem \eqref{eq:2.14}. Under these assumptions (A)-(C) and (i)-(iii), we derive the following global error estimates for the HOMS solutions if $T^{(0)}$, $c^{(0)}\in H^{4}(\Omega)$ and $\bm{u}^{(0)}\in (H^{4}(\Omega))^n$.
		\begin{equation}
			\label{eq:3.5}
			\begin{aligned}
				\|T^\epsilon-T^{(2,\epsilon)}\|_{H^1(\Omega)}=\|T_{\Delta}^{(2,\epsilon)}\|_{H^1(\Omega)}\leq C\epsilon,
			\end{aligned}
		\end{equation}
		\begin{equation}
			\label{eq:3.6}
			\begin{aligned}
				\|c^\epsilon-c^{(2,\epsilon)}\|_{H^1(\Omega)}=\|c_{\Delta}^{(2,\epsilon)}\|_{H^1(\Omega)}\leq C\epsilon,
			\end{aligned}
		\end{equation}
		\begin{equation}
			\label{eq:3.7}
			\begin{aligned}
				\|\bm{u}^{\epsilon}-\bm{u}^{(2,\epsilon)}\|_{(H^1(\Omega))^n}=\|\bm{u}_{\Delta}^{(2,\epsilon)}\|_{(H^1(\Omega))^n}\leq C\epsilon,
			\end{aligned}
		\end{equation}
		where $C$ represents a positive constant that depends on $\Omega$ but remains independent of the small periodic parameter $\epsilon$. Throughout this work, we systematically represent all generic constants uniformly by $C$ without distinction.\\
	\end{thm}
	
	\textbf{Proof.}\hspace{1mm} Firstly, we define three differential operators $\displaystyle\mathcal{K}_{T}(T)=-\frac{\partial}{\partial x_{i}}\Bigl(k_{ij}(\mathbf{x},\mathbf{y})\frac{\partial T}{\partial x_{j}}\Bigr)$, $\displaystyle \mathcal{K}_{c}(c)=-\frac{\partial}{\partial x_{i}}\Bigl(g_{ij}(\mathbf{x},\mathbf{y})\frac{\partial c}{\partial x_{j}}\Bigr)$ and $\displaystyle \mathcal{K}_{u}^{i}(\bm{u})=-\frac{\partial}{\partial x_{j}}\Bigl(D_{ijkl}(\mathbf{x},\mathbf{y})\frac{\partial u_{k}}{\partial x_{l}}\Bigr)$ for temperature increment, moisture and displacement fields respectively for simplifying the subsequent proof process.
	
	After that, we obtain that  $T_\Delta^{(2,\epsilon)}$, $c_\Delta^{(2,\epsilon)}$ and $\bm{u}_{\Delta}^{(2,\epsilon)}$ are the weak solutions of the following boundary value problems, respectively
	\begin{equation}
		\label{eq:3.8}
		\begin{cases}
			\begin{aligned}
				& \mathcal{K}_T(T_{\Delta}^{(2,\epsilon)}) = \epsilon F_0^T + \epsilon \frac{\partial}{\partial x_i}F_i^T, \quad \text{in } \Omega, \\
				& T_{\Delta}^{(2,\epsilon)} = 0, \quad \text{on } \partial \Omega,
			\end{aligned}
		\end{cases}
	\end{equation}
	\begin{equation}
		\label{eq:3.9}
		\begin{cases}
			\begin{aligned}
				& \mathcal{K}_c(c_{\Delta}^{(2,\epsilon)}) = \epsilon F_0^c + \epsilon \frac{\partial}{\partial x_i}F_i^c, \quad \text{in } \Omega, \\
				& c_{\Delta}^{(2,\epsilon)} = 0, \quad \text{on } \partial \Omega,
			\end{aligned}
		\end{cases}
	\end{equation}
	\begin{equation}
		\label{eq:3.10}
		\begin{cases}
			\begin{aligned}
				& \mathcal{K}_{u}^{i}(\bm{u}_{\Delta}^{(2,\epsilon)})=-\frac{\partial}{\partial x_{j}}\bigl(D_{ijkl}(\mathbf{x},\mathbf{y})\alpha_{kl}(\mathbf{x},\mathbf{y})T_{\Delta}^{(2,\epsilon)}+D_{ijkl}(\mathbf{x},\mathbf{y})\beta_{kl}(\mathbf{x},\mathbf{y})c_{\Delta}^{(2,\epsilon)}\bigr) \\
				& +\epsilon F_{i0}^{u}+\epsilon\frac{\partial}{\partial x_{j}}F_{ij}^{u}+\epsilon F_{i0}^{uT}+\epsilon\frac{\partial}{\partial x_{j}}F_{ij}^{uT}+\epsilon F_{i0}^{uc}+\epsilon\frac{\partial}{\partial x_{j}}F_{ij}^{uc},\;\;\text{in }\Omega, \\
				& \bm{u}_{\Delta}^{(2,\epsilon)}=0,\;\;\text{on }\partial\Omega,
			\end{aligned}
		\end{cases}
	\end{equation}
	where the specific expressions for $F_{0}^{T}$, $F_{i}^{T}$, $F_{0}^{c}$, $F_{i}^{c}$, $F_{i0}^{u}$, $F_{ij}^{u}$, $F_{i0}^{uT}$, $F_{ij}^{uT}$, $F_{i0}^{uc}$, and $F_{ij}^{uc}$ can be readily derived and are exhibited in Appendix C of the present study due to their lengthy forms.
	
	Since all first-order and second-order auxiliary cell functions, along with their partial derivatives with respect to macroscopic variable $\mathbf{x}$ and microscopic variable $\mathbf{y}$, belong to $L^{2}(Y)$, the following inequalities hold \cite{R22,R36,R37,R42,R44}		
	\begin{equation}
		\label{eq:3.11}
		\sum_{k=0}^n\|F_k^T\|_{L^2(\Omega)}\leq C\|T^{(0)}\|_{H^{4}(\Omega)},
	\end{equation}
	\begin{equation}
		\label{eq:3.12}
		\sum_{k=0}^n\|F_k^c\|_{L^2(\Omega)}\leq C\|c^{(0)}\|_{H^{4}(\Omega)},
	\end{equation}
	\begin{equation}
		\label{eq:3.13}
		\begin{cases}
			\begin{aligned}
				& \sum_{j=0}^n\lVert F_{ij}^u\rVert_{L^2(\Omega)}\leq C\|\bm{u}^{(0)}\|_{(H^{4}(\Omega))^n}, \\
				& \sum_{j=0}^n\lVert F_{ij}^{uT}\rVert_{L^2(\Omega)}\leq C\|T^{(0)}\|_{H^{3}(\Omega)}, \\
				& \sum_{j=0}^n\lVert F_{ij}^{uc}\rVert_{L^2(\Omega)}\leq C\|c^{(0)}\|_{H^{3}(\Omega)}.
			\end{aligned}
		\end{cases}
	\end{equation}
	Subsequently, utilizing the uniform ellipticity of equations \eqref{eq:3.8}-\eqref{eq:3.9}, along with inequalities \eqref{eq:3.11}-\eqref{eq:3.12}, we can easily conclude that the error estimates \eqref{eq:3.5} and \eqref{eq:3.6} for the temperature increment and moisture fields in Theorem \ref{thm:1} hold.
	
	Finally, employing the uniform ellipticity of equation \eqref{eq:3.10}, applying lemma \ref{lem:3.1}, and substituting the estimates \eqref{eq:3.5}-\eqref{eq:3.6} and \eqref{eq:3.13} into equation \eqref{eq:3.10}, we establish the displacement error estimate \eqref{eq:3.7} with the same $O(\epsilon)$ convergence.
	
	\section{Multi-scale finite element algorithm}
	\label{sec:4}
	This section details the multi-scale algorithm for solving the hygro-thermo-mechanical problems \eqref{eq:2.1}. It is important to note that the auxiliary cell functions of quasi-periodic composite structures exhibit explicit dependence on the macroscopic coordinate $\mathbf{x}$. To obtain the final HOMS solutions, the same auxiliary cell functions must be computed at numerous distinct $\mathbf{x}$, which consumes significant computational  time. Fortunately, it can be demonstrated that all auxiliary cell functions exhibit continuity in $\mathbf{x}$ \cite{R3}. Therefore, computation of the auxiliary cell functions is required solely for discrete representative points $\mathbf{x}_I$ in the domain $\Omega$, enabling their determination at other required macro-coordinates through interpolation. In the following, we present the numerical algorithm, consisting of off-line and on-line stages, for efficiently simulating the multi-scale problem \eqref{eq:2.1}.
	
	\subsection{Off-line microscale computation}
	\label{sec:41}
	\begin{enumerate}
		\item[(1)] Determine the geometric configuration of PUC $Y=[0,1]^n$ in $\mathbb{R}^n (n=2,3)$, and generate a family of triangular $(n=2)$ or tetrahedral $(n=3)$ finite element meshes $J^{h_0}(Y) = \{K\}$ of $Y$, where $h_0=\max_K \{ h_K\}$. Then denote the linear conforming finite element space $V_{h_0}(Y)=\{\nu\in C^{0}(\bar{Y}):\nu|_{\partial Y}=0,\nu|_{K}\in P_{1}(K)\}\subset H_{0}^{1}(Y)$ for auxiliary cell problems.
		\item[(2)] Use FEM to solve the first-order cell problems \eqref{eq:2.9}-\eqref{eq:2.13} on $V_{h_0}(Y)$ corresponding to different representative macro-scale points $\mathbf{x}_I\in\Omega$, and obtain $\mathcal{H}_{\alpha_1}(\mathbf{x}_I,\mathbf{y})$, $\mathcal{L}_{\alpha_1}(\mathbf{x}_I,\mathbf{y})$, $\mathcal{X}_{ih}^{\alpha_1}(\mathbf{x}_I,\mathbf{y})$, $\mathcal{M}_i(\mathbf{x}_I,\mathbf{y})$ and $\mathcal{N}_i(\mathbf{x}_I,\mathbf{y})$. Then evaluate the macroscopic homogenized material parameters $\hat{k}_{ij}(\mathbf{x}_I)$,  $\hat{g}_{ij}(\mathbf{x}_I)$, $\hat{D}_{ijkl}(\mathbf{x}_I)$,  $\hat{A}_{ij}(\mathbf{x}_I)$ and $\hat{B}_{ij}(\mathbf{x}_I)$ via formula \eqref{eq:2.15}.
		\item[(3)] Employing the same mesh as first-order cell problems, we solve the second-order cell problems \eqref{eq:2.20}-\eqref{eq:2.29} on $V_{h_0}(Y)$ corresponding to different macro-scale points $\mathbf{x}_I$ via FEM, and obtain $\mathcal{H}_{\alpha_1\alpha_2}(\mathbf{x}_I,\mathbf{y})$, $\mathcal{R}_{\alpha_1}(\mathbf{x}_I,\mathbf{y})$,  $\mathcal{L}_{\alpha_1\alpha_2}(\mathbf{x}_I,\mathbf{y})$,  $\mathcal{S}_{\alpha_1}(\mathbf{x}_I,\mathbf{y})$, $\mathcal{P}_{ih}^{\alpha_1\alpha_2}(\mathbf{x}_I,\mathbf{y})$, $\mathcal{Q}_{ih}^{\alpha_1}(\mathbf{x}_I,\mathbf{y})$, $\mathcal{W}_i(\mathbf{x}_I,\mathbf{y})$,  $\mathcal{Z}_i^{\alpha_1}(\mathbf{x}_I,\mathbf{y})$, $\mathcal{F}_i(\mathbf{x}_I,\mathbf{y})$ and $\mathcal{G}_i^{\alpha_1}(\mathbf{x}_I,\mathbf{y})$. When solving  $\mathcal{R}_{\alpha_1}(\mathbf{x}_I,\mathbf{y})$, $\mathcal{S}_{\alpha_1}(\mathbf{x}_I,\mathbf{y})$, $\mathcal{Q}_{ih}^{\alpha_1}(\mathbf{x}_I,\mathbf{y})$, $\mathcal{W}_i(\mathbf{x}_I,\mathbf{y})$ and $\mathcal{F}_i(\mathbf{x}_I,\mathbf{y})$, it should be mentioned that the partial derivatives with respect to the macroscopic variable $\mathbf{x}$ in the variational formulations are approximated by using the central difference scheme according to their explicit expressions.
	\end{enumerate}
	\subsection{On-line macroscale and multi-scale computation}
	\label{sec:42}
	\begin{enumerate}
		\item[(1)] Let $J^{h_1}(\Omega)= \{e\}$ be a triangular or tetrahedral  finite element mesh of the macroscopic homogenized region $\Omega$, where $h_1=\max_e \{ h_e \}$. Then define the linear conforming finite element spaces $V^T_{h_1}(\Omega)=\{\nu\in C^0(\bar{\Omega}):\nu|_{\Gamma_{T}}=0,\nu|_e\in P_1(e)\}\subset H^1(\Omega)$, $V^c_{h_1}(\Omega)=\{\nu\in C^0(\bar{\Omega}):\nu|_{\Gamma_{c}}=0,\nu|_e\in P_1(e)\}\subset H^1(\Omega)$ and $V^u_{h_1}(\Omega)=\{\nu\in C^0(\bar{\Omega}):\nu|_{\Gamma_{u}}=0,\nu|_e\in P_1(e)\}\subset H^1(\Omega)$ for temperature increment, moisture and displacement fields, respectively.
		\item[(2)] The macroscopic homogenized material parameters at each node of homogenized domain $\Omega$ can be calculated using an interpolation approach. Then, by solving the macroscopic homogenized problem \eqref{eq:2.14} in the finite element spaces $V^T_{h_1}(\Omega)$, $V^c_{h_1}(\Omega)$ and $V^u_{h_1}(\Omega)$, the macroscopic homogenized solutions $T^{(0)}$, $c^{(0)}$ and $u_{i}^{(0)}$ for multi-scale problem \eqref{eq:2.1} are obtained.
		\item[(3)] For any point $\mathbf{x}\in\Omega$, the linear interpolation technique is applied to compute the values of first-order and second-order cell functions alongside macroscopic homogenized solutions. Their partial derivatives $\displaystyle \frac{\partial T^{(0)}}{\partial x_{\alpha_1}}$, $\displaystyle \frac{\partial^2T^{(0)}}{\partial x_{\alpha_1}\partial x_{\alpha_2}}$, $\displaystyle \frac{\partial c^{(0)}}{\partial x_{\alpha_1}}$, $\displaystyle \frac{\partial^2c^{(0)}}{\partial x_{\alpha_1}\partial x_{\alpha_2}}$, $\displaystyle \frac{\partial u_h^{(0)}}{\partial x_{\alpha_1}}$ and $\displaystyle \frac{\partial^2u_h^{(0)}}{\partial x_{\alpha_1}\partial x_{\alpha_2}}$ are approximated through the average technique on relative elements \cite{R38,R32,R31,R39}. Subsequently, the temperature increment field $T^{(2,\epsilon)}$, moisture field $c^{(2,\epsilon)}$, and displacement field $\bm{u}^{(2,\epsilon)}$ are obtained from equations \eqref{eq:2.33}-\eqref{eq:2.35}. Furthermore, the higher-order interpolation and post-processing techniques can still be employed to yield the high-precision HOMS solutions \cite{R32,R40,R41}.
	\end{enumerate}
	
	\subsection{Error estimation of multi-scale finite element algorithm}
	\label{sec:43}
	The total error of the proposed multi-scale algorithm comprises not only the error arising from multi-scale modeling, but also two additional components: the numerical error from solving auxiliary cell problems \eqref{eq:2.9}-\eqref{eq:2.13} and \eqref{eq:2.20}-\eqref{eq:2.29}, and that from solving the macroscopic homogenized problem \eqref{eq:2.14}, using FEM. Before giving the detailed error estimation, we prepare some lemmas in advance.
	\begin{lemma}
		\label{lem:4.1}
		Let $\mathcal{H}_{\alpha_1}^{h_0}$, $\mathcal{L}_{\alpha_1}^{h_0}$, $\mathcal{X}_{ih}^{\alpha_1, h_0}$, $\mathcal{M}_{i}^{h_0}$, $\mathcal{N}_{i}^{h_0}$, $\mathcal{H}_{\alpha_1\alpha_2}^{h_0}$, $\mathcal{R}_{\alpha_1}^{h_0}$, $\mathcal{L}_{\alpha_1\alpha_2}^{h_0}$, $\mathcal{S}_{\alpha_1}^{h_0}$, $\mathcal{P}_{ih}^{\alpha_1\alpha_2, h_0}$, $\mathcal{Q}_{ih}^{\alpha_1, h_0}$, $\mathcal{W}_{i}^{h_0}$,  $\mathcal{Z}_i^{\alpha_1, h_0}$, $\mathcal{F}_{i}^{h_0}$ and $\mathcal{G}_i^{\alpha_1,h_0}$ denote the corresponding FE solutions of the first-order and second-order cell functions, respectively. If all microscopic cell functions belong to $H^2({Y})$ for any fixed $\mathbf{x}$, then there holds the following inequality
		\begin{equation}
			\label{eq:4.1}
			\begin{aligned}
				&\left\|\mathcal{H}_{\alpha_1}^{h_0}(\mathbf{x},\mathbf{y})-\mathcal{H}_{\alpha_1}(\mathbf{x},\mathbf{y})\right\|_{H^m(Y)}\leq Ch_0^{2-m}\left\|\mathcal{H}_{\alpha_1}\right\|_{H^2(Y)},
			\end{aligned}
		\end{equation}
		where $m=0,1$ and $C$ denotes the finite element estimate constant independent of $h_0$ and dependent on $Y$. Moreover, other microscopic cell functions have the similar error estimates to the $\mathcal{H}_{\alpha_1}^{h_0}$.
	\end{lemma}
	$\mathbf{Proof:}$ By employing the classical finite element theory, the above inequalities are easily obtained.
	\begin{lemma}
		\label{lem:4.2}
		Denote $\hat k_{ij}^{h_0}(\mathbf{x})$, $\hat g_{ij}^{h_0}(\mathbf{x})$, $\hat D_{ijkl}^{h_0}(\mathbf{x})$, $\hat A_{ij}^{h_0}(\mathbf{x})$ and $\hat B_{ij}^{h_0}(\mathbf{x})$ be the FE approximation of the corresponding macroscopic homogenized parameters, the following results hold
		\begin{equation}
			\label{eq:4.2}
			\begin{aligned}
				\left|\hat k_{ij}^{h_0}(\mathbf{x})-\hat k_{ij}(\mathbf{x})\right|\leq Ch_0^2\left\| \mathcal{H}_{i}\right\|_{H^2(Y)}\left\| \mathcal{H}_{j}\right\|_{H^2(Y)},\;
				\underline{\kappa} | \bm{\xi} |^2 \leq \hat k_{ij}^{h_0} \xi_i \xi_j \leq \overline{\kappa} | \bm{\xi} |^2,
			\end{aligned}
		\end{equation}
		\begin{equation}
			\label{eq:4.3}
			\begin{aligned}
				\left|\hat g_{ij}^{h_0}(\mathbf{x})-\hat g_{ij}(\mathbf{x})\right|\leq Ch_0^2\left\| \mathcal{L}_{i}\right\|_{H^2(Y)}\left\| \mathcal{L}_{j}\right\|_{H^2(Y)},\;
				\underline{\kappa} | \bm{\xi} |^2 \leq \hat g_{ij}^{h_0} \xi_i \xi_j \leq \overline{\kappa} | \bm{\xi} |^2,
			\end{aligned}
		\end{equation}
		\begin{equation}
			\label{eq:4.4}
			\begin{aligned}
				\left|\hat D_{ijkl}^{h_0}(\mathbf{x})-\hat D_{ijkl}(\mathbf{x})\right|\leq Ch_0^2\left\| \mathcal{X}_{ih}^{j}\right\|_{H^2(Y)}\left\| \mathcal{X}_{kh}^{l}\right\|_{H^2(Y)},\;
				\underline{\kappa} \eta_{ij} \eta_{ij}  \leq \hat D_{ijkl}^{h_0} \eta_{ij} \eta_{kl} \leq \overline{\kappa} \eta_{ij} \eta_{ij},
			\end{aligned}
		\end{equation}
		\begin{equation}
			\label{eq:4.5}
			\begin{aligned}
				\left|\hat A_{ij}^{h_0}(\mathbf{x})-\hat A_{ij}(\mathbf{x})\right|\leq Ch_0^2\left\| \mathcal{M}_{i}\right\|_{H^2(Y)}\left\| \mathcal{M}_{j}\right\|_{H^2(Y)},\;
				\underline{\kappa} | \bm{\xi} |^2 \leq \hat A_{ij}^{h_0} \xi_i \xi_j \leq \overline{\kappa} | \bm{\xi} |^2,
			\end{aligned}
		\end{equation}
		\begin{equation}
			\label{eq:4.6}
			\begin{aligned}
				\left|\hat B_{ij}^{h_0}(\mathbf{x})-\hat B_{ij}(\mathbf{x})\right|\leq Ch_0^2\left\| \mathcal{N}_{i}\right\|_{H^2(Y)}\left\| \mathcal{N}_{j}\right\|_{H^2(Y)},\;
				\underline{\kappa} | \bm{\xi} |^2 \leq \hat B_{ij}^{h_0} \xi_i \xi_j \leq \overline{\kappa} | \bm{\xi} |^2,
			\end{aligned}
		\end{equation}
		where $C$ is a constant independent of $h_0$.
	\end{lemma}
	$\mathbf{Proof:}$ By employing the definitions of macroscopic homogenized material parameters in \eqref{eq:2.15}, assumption (B) and lemma \ref{lem:4.1}, it follows that
	\begin{equation}
		\label{eq:4.7}
		\begin{aligned}
			&\left|\hat k_{ij}^{h_0}(\mathbf{x})-\hat k_{ij}(\mathbf{x})\right|\\
			&=\left|\frac{1}{|{Y}|}{\int_{{Y}}}\big(k_{ij} + k_{ik} \frac{\partial \mathcal{H}_{j}^{h_0}}{\partial y_{k}}\big)d{Y}-\frac{1}{|{Y}|}{\int_{{Y}}}\big(k_{ij} + k_{ik} \frac{\partial \mathcal{H}_{j}}{\partial y_{k}}\big)d{Y}\right|\\
			&=\left|\frac{1}{|{Y}|}{\int_{{Y}}}{k_{ik}{\frac{\partial \big(\mathcal{H}_j^{h_0}-\mathcal{H}_j\big)}{\partial y_{k}}}}d{Y}\right|\\
			&=\frac{1}{|{Y}|}\left|-{\int_{{Y}}}\frac{\partial \mathcal{H}_i}{\partial y_{\alpha_1}}k_{\alpha_1\alpha_2}{\frac{\partial }{\partial y_{\alpha_2}}\big(\mathcal{H}_j^{h_0}-\mathcal{H}_j\big)}d{Y}\right|\\
			&=\frac{1}{|{Y}|}\left|{\int_{{Y}}}\frac{\partial \mathcal{H}_i^{h_0}}{\partial y_{\alpha_1}}k_{\alpha_1\alpha_2}{\frac{\partial }{\partial y_{\alpha_2}}\big(\mathcal{H}_j^{h_0}-H_j\big)}d{Y}-{\int_{{Y}}}\frac{\partial \mathcal{H}_i}{\partial y_{\alpha_1}}k_{\alpha_1\alpha_2}{\frac{\partial }{\partial y_{\alpha_2}}\big(\mathcal{H}_j^{h_0}-\mathcal{H}_j\big)}d{Y}\right|\\
			&=\frac{1}{|{Y}|}\left|{\int_{{Y}}}\frac{\partial}{\partial y_{\alpha_1}}\big(\mathcal{H}_i^{h_0}-\mathcal{H}_i\big)k_{\alpha_1\alpha_2}{\frac{\partial }{\partial y_{\alpha_2}}\big(\mathcal{H}_j^{h_0}-\mathcal{H}_j\big)}d{Y}\right|\\
			&\leq C\left\|\mathcal{H}_i^{h_0}-\mathcal{H}_i\right\|_{H^1({Y})}\left\|\mathcal{H}_j^{h_0}-\mathcal{H}_j\right\|_{H^1({Y})}\leq Ch_0^2\left\|\mathcal{H}_i\right\|_{H^2({Y})}\left\|\mathcal{H}_j\right\|_{H^2({Y})}.
		\end{aligned}
	\end{equation}
	Furthermore, choosing a sufficiently small $h_0>0$ satisfies
	\begin{equation}
		\label{eq:4.8}
		Ch_0^2\left\|\mathcal{H}_i(\mathbf{x}, \mathbf{y})\right\|_{H^2({Y})}\left\| \mathcal{H}_j(\mathbf{x}, \mathbf{y})\right\|_{H^2({Y})}\leq\underline{\varsigma}/2.
	\end{equation}
	Hence, we can verify that the lower bound in \eqref{eq:4.2} holds
	\begin{equation}
		\label{eq:4.9}
		\hat k_{ij}^{h_0}(\mathbf{x})\xi_i\xi_j=\hat k_{ij}(\mathbf{x})\xi_i\xi_j+\big[\hat k_{ij}^{h_0}(\mathbf{x})-\hat k_{ij}(\mathbf{x})\big]\xi_i\xi_j\geq(\underline{\varsigma}-\underline{\varsigma}/2)\xi_i\xi_i=\underline{\kappa}|\bm{\xi}|^2,
	\end{equation}
	where $\underline{\kappa}=\underline{\varsigma}/2$ is a constant independent of $h_0$. Moreover, the upper bound in \eqref{eq:4.2} is easily derived when setting $\overline{\kappa}=\overline{\varsigma}+\underline{\varsigma}/2$. Finally, following the similar way, we can obtain the results \eqref{eq:4.3}-\eqref{eq:4.6}.
	
	As shown in lemmas \ref{lem:4.1} and \ref{lem:4.2}, the values of macroscopic homogenized material parameters $\hat k_{ij}^{h_0}$, $\hat g_{ij}^{h_0}$, $\hat D_{ijkl}^{h_0}$, $\hat A_{ij}^{h_0}$ and $\hat B_{ij}^{h_0}$ depend on the finite element computations of
	first-order cell functions $\mathcal{H}_{\alpha_1}^{h_0}$, $\mathcal{L}_{\alpha_1}^{h_0}$, $\mathcal{X}_{ih}^{\alpha_1, h_0}$, $\mathcal{M}_{i}^{h_0}$ and $\mathcal{N}_{i}^{h_0}$. Therefore, in practice, we need to numerically solve the modified homogenized equations as below
	\begin{equation}
		\label{eq:4.10}
		\begin{cases}
			\begin{aligned}
				&-\frac{\partial}{\partial x_i}\Bigl(\hat{k}_{ij}^{h_0}(\mathbf{x}) \frac{\partial T^{(0, h_0)}}{\partial x_j}\Bigr) = h,\; \text{in }\; \Omega,\\
				&-\frac{\partial}{\partial x_i}\Bigl(\hat{g}_{ij}^{h_0}(\mathbf{x}) \frac{\partial c^{(0, h_0)}}{\partial x_j}\Bigr) = m,\; \text{in }\; \Omega,\\
				&-\frac{\partial}{\partial x_j}\Bigl(\hat{D}_{ijkl}^{h_0}(\mathbf{x}) \frac{\partial u_{k}^{(0, h_0)}}{\partial x_l}
				- \hat{A}_{ij}^{h_0}(\mathbf{x}) T^{(0, h_0)}
				- \hat{B}_{ij}^{h_0}(\mathbf{x}) c^{(0, h_0)}\Bigr) = f_i,\; \text{in }\;\Omega,\\
				& T^{(0, h_0)}(\mathbf{x}) = \overline{T}(\mathbf{x}),\;\text{on }\;\Gamma_{T},\\
				& \hat{k}_{ij}^{h_0}\bigl(\mathbf{x}\bigr) \frac{\partial T^{(0, h_0)}(\mathbf{x})}{\partial x_{j}}n_i=\overline{q}(\mathbf{x}),\;\text{on }\;\Gamma_{q},\\
				& c^{(0, h_0)}(\mathbf{x}) = \overline{c}(\mathbf{x}),\;\text{on }\;\Gamma_{c},\\
				& \hat{g}_{ij}^{h_0}\bigl(\mathbf{x}\bigr) \frac{\partial c^{(0, h_0)}(\mathbf{x})}{\partial x_{j}}n_i=\overline{d}(\mathbf{x}),\;\text{on }\;\Gamma_{d},\\
				& \bm{u}^{(0, h_0)}(\mathbf{x}) = \overline{\bm{u}}(\mathbf{x}),\;\text{on }\;\Gamma_{{u}},\\
				& \!\Bigl[\!\hat{D}_{ijkl}^{h_0}(\mathbf{x}) \frac{\partial u_{k}^{(0, h_0)}(\mathbf{x})}{\partial x_l}
				- \hat{A}_{ij}^{h_0}(\mathbf{x}) T^{(0, h_0)}(\mathbf{x})
				- \hat{B}_{ij}^{h_0}(\mathbf{x}) c^{(0, h_0)}(\mathbf{x})\!\Bigr]\! n_{j}=\overline{\sigma}_{i}(\mathbf{x}),\;\text{on }\;\Gamma_{\sigma}.
			\end{aligned}
		\end{cases}
	\end{equation}
	\begin{lemma}
		\label{lem:4.3}
		Let $T^{(0, h_0)}$, $c^{(0, h_0)}$ and $u_{i}^{(0, h_0)}$ represent the exact solutions of the revised macroscopic homogenized equations \eqref{eq:4.10}, and $T^{(0, h_0, h_1)}$, $c^{(0, h_0, h_1)}$ and $u_{i}^{(0, h_0, h_1)}$ denote the corresponding finite element solutions of the revised macroscopic homogenized equations \eqref{eq:4.10}, the following estimates hold
		\begin{equation}
			\label{eq:4.11}
			\begin{aligned}
				\| T^{(0,h_0,h)}-T^{(0)} \|_{H^1(\Omega)}\leq C(h_0^2 + h_1),
			\end{aligned}
		\end{equation}
		\begin{equation}
			\label{eq:4.12}
			\begin{aligned}
				\| c^{(0,h_0,h)}-c^{(0)} \|_{H^1(\Omega)}\leq C(h_0^2 + h_1),
			\end{aligned}
		\end{equation}
		\begin{equation}
			\label{eq:4.13}
			\|\bm{u}^{(0,h_0,h)}-\bm{u}^{(0)}\|_{(H^1(\Omega))^n}\leq C(h_0^2 + h_1),
		\end{equation}
		where $C$ is a constant independent of $h_0$ and $h_1$.
	\end{lemma}
	\textbf{Proof:}\hspace{1mm}Firstly, subtracting the homogenized thermal equation in \eqref{eq:2.14} from corresponding thermal equation in \eqref{eq:4.10}, we can obtain
	\begin{equation}
		\label{eq:4.14}
		\begin{aligned}
			- \frac{\partial }{{\partial {x_i}}}\Big[ {{\hat k_{ij}^{h_0}}(\mathbf{x})\frac{{\partial \big(T^{(0,h_0)}-T^{(0)}\big)}}{{\partial {x_j}}}}\Big]=-\frac{\partial }{{\partial {x_i}}}\Big[ {{\big(\hat k_{ij}(\mathbf{x})-\hat k_{ij}^{h_0}(\mathbf{x})\big)}\frac{{\partial T^{(0)}}}{{\partial {x_j}}}}\Big].
		\end{aligned}
	\end{equation}
	Furthermore, multiplying on both sides of equality \eqref{eq:4.14} by $T^{(0,h_0)}-T^{(0)}$ and integrating on $\Omega$, it follows that
	\begin{equation}
		\label{eq:4.15}
		\begin{aligned}
			&\int_{\Omega}{{\hat k_{ij}^{h_0}}(\mathbf{x})\frac{{\partial \big(T^{(0,h_0)}-T^{(0)}\big)}}{{\partial {x_j}}}}\frac{{\partial \big(T^{(0,h_0)}-T^{(0)}\big)}}{{\partial {x_i}}}d\Omega\\
			&=\int_{\Omega}{{\big(\hat k_{ij}(\mathbf{x})-\hat k_{ij}^{h_0}(\mathbf{x})\big)}\frac{{\partial T^{(0)}}}{{\partial {x_j}}}}\frac{{\partial \big(T^{(0,h_0)}-T^{(0)}\big)}}{{\partial {x_i}}}d\Omega.
		\end{aligned}
	\end{equation}
	Recalling the inequality in \eqref{eq:4.2} and employing Poincar$\rm{\acute{e}}$-Friedrichs inequality and Cauchy-Schwarz inequality, we shall naturally derive the following inequality from equality \eqref{eq:4.15}.
	\begin{equation}
		\label{eq:4.16}
		\begin{aligned}
			\| T^{(0,h_0)}-T^{(0)} \|_{H^1(\Omega)}\leq Ch_0^2.
		\end{aligned}
	\end{equation}
	Subsequently, according to classical finite element theory, we derive the following estimate
	\begin{equation}
		\label{eq:4.17}
		\|T^{(0,h_0,h_1)}- T^{(0,h_0)} \|_{H^1(\Omega)} \leq Ch_1.
	\end{equation}
	Combining the obtained inequalities \eqref{eq:4.16} and \eqref{eq:4.17}, and employing the triangle inequality, we can prove that the following inequality holds, namely the error estimate \eqref{eq:4.11} holds
	\begin{equation}
		\label{eq:4.18}
		\begin{aligned}
			&\|T^{(0,h_0,h_1)} - T^{(0)}\|_{H^1(\Omega)} \\
			&\leq \|T^{(0,h_0,h_1)} -T^{(0,h_0)} \|_{H^1(\Omega)}+ \|T^{(0,h_0)}-T^{(0)}  \|_{H^1(\Omega)}  \\
			&\leq C(h_0^2 + h_1).
		\end{aligned}
	\end{equation}
	Finally, following the similar way, we can obtain the results \eqref{eq:4.12} and \eqref{eq:4.13}.
	
	In what follows, rigorous convergence analysis will be conducted for the proposed multi-scale FE algorithm.
	\begin{thm}
		\label{thm:2}
		Let $T^{(\epsilon,h_0,h_1)}, c^{(\epsilon,h_0,h_1)}$  and $u_{i}^{(\epsilon,h_0,h_1)}$ be the corresponding multi-scale finite element solutions of the HOMS solutions as below
		\begin{equation}
			\label{eq:4.19}
			\begin{aligned}
				T^{(\epsilon,h_0,h_1)}(\mathbf{x}) & =T^{(0, h_0, h_1)}(\mathbf{x}) +\epsilon \mathcal{H}_{\alpha_1}^{h_0}(\mathbf{x},\mathbf{y})\frac{\partial T^{(0, h_0, h_1)}(\mathbf{x})}{\partial x_{\alpha_1}}\\
				& +\epsilon^2\Bigl(\mathcal{H}_{\alpha_1\alpha_2}^{h_0}(\mathbf{x},\mathbf{y})\frac{\partial^2T^{(0, h_0, h_1)}(\mathbf{x})}{\partial x_{\alpha_1}\partial x_{\alpha_2}}+\mathcal{R}_{\alpha_1}^{h_0}(\mathbf{x},\mathbf{y})\frac{\partial T^{(0, h_0, h_1)}(\mathbf{x})}{\partial x_{\alpha_1}}\Bigr),
			\end{aligned}
		\end{equation}
		\begin{equation}
			\label{eq:4.20}
			\begin{aligned}
				c^{(\epsilon,h_0,h_1)}(\mathbf{x}) & =c^{(0, h_0, h_1)}(\mathbf{x})+\epsilon \mathcal{L}_{\alpha_1}^{h_0}(\mathbf{x},\mathbf{y})\frac{\partial c^{(0, h_0, h_1)}(\mathbf{x})}{\partial x_{\alpha_1}}\\
				& +\epsilon^2\Bigl(\mathcal{L}_{\alpha_1\alpha_2}^{h_0}(\mathbf{x},\mathbf{y})\frac{\partial^2c^{(0, h_0, h_1)}(\mathbf{x})}{\partial x_{\alpha_1}\partial x_{\alpha_2}}+\mathcal{S}_{\alpha_1}^{h_0}(\mathbf{x},\mathbf{y})\frac{\partial c^{(0, h_0, h_1)}(\mathbf{x})}{\partial x_{\alpha_1}}\Bigr),
			\end{aligned}
		\end{equation}
		\begin{equation}
			\label{eq:4.21}
			\begin{aligned}
				u_{i}^{(\epsilon,h_0,h_1)}(\mathbf{x}) & =u_{i}^{(0, h_0, h_1)}(\mathbf{x})+\epsilon\Bigl(\mathcal{X}_{ih}^{\alpha_{1}, h_0}(\mathbf{x},\mathbf{y})\frac{\partial u_{h}^{(0, h_0, h_1)}(\mathbf{x})}{\partial x_{\alpha_{1}}}-\mathcal{M}_{i}^{h_0}(\mathbf{x},\mathbf{y})T^{(0, h_0, h_1)}(\mathbf{x}) \\
				& -\mathcal{N}_{i}^{h_0}(\mathbf{x},\mathbf{y})c^{(0, h_0, h_1)}(\mathbf{x})\Bigr) +\epsilon^{2}\Bigl(\mathcal{P}_{ih}^{\alpha_{1}\alpha_{2}, h_0}(\mathbf{x},\mathbf{y})\frac{\partial^{2}u_{h}^{(0, h_0, h_1)}(\mathbf{x})}{\partial x_{\alpha_{1}}\partial x_{\alpha_{2}}} \\
				& +\mathcal{Q}_{ih}^{\alpha_{1}, h_0}(\mathbf{x},\mathbf{y})\frac{\partial u_{h}^{(0, h_0, h_1)}(\mathbf{x})}{\partial x_{\alpha_{1}}} +\mathcal{W}_{i}^{h_0}(\mathbf{x},\mathbf{y})T^{(0, h_0, h_1)}(\mathbf{x}) \\
				& +\mathcal{Z}_{i}^{\alpha_{1}, h_0}(\mathbf{x},\mathbf{y})\frac{\partial T^{(0, h_0, h_1)}(\mathbf{x})}{\partial x_{\alpha_{1}}} +\mathcal{F}_{i}^{h_0}(\mathbf{x},\mathbf{y})c^{(0, h_0, h_1)}(\mathbf{x})+\mathcal{G}_{i}^{\alpha_{1}, h_0}(\mathbf{x},\mathbf{y})\frac{\partial c^{(0, h_0, h_1)}(\mathbf{x})}{\partial x_{\alpha_{1}}}\Bigr),
			\end{aligned}
		\end{equation}
		then holds when assuming $T^{(0)}$, $T^{(0,h_0)}$ and $T^{(0,h_0,h_1)}$ $\in H^{4}(\Omega)$, $c^{(0)}$, $c^{(0,h_0)}$ and $c^{(0,h_0,h_1)}$ $\in H^{4}(\Omega)$, $\bm{u}^{(0)}$, $\bm{u}^{(0,h_0)}$ and $\bm{u}^{(0,h_0,h_1)}$ $\in (H^{4}(\Omega))^n$.
		\begin{equation}
			\label{eq:4.22}
			\|T^\epsilon(\mathbf{x})-T^{(\epsilon,h_0,h_1)}(\mathbf{x})\|_{H^1(\Omega)}\leq C\bigl(\epsilon+h_0^2+h_1 \bigr),
		\end{equation}
		\begin{equation}
			\label{eq:4.23}
			\|c^\epsilon(\mathbf{x})-c^{(\epsilon,h_0,h_1)}(\mathbf{x})\|_{H^1(\Omega)}\leq C\bigl(\epsilon+h_0^2+h_1 \bigr),
		\end{equation}
		\begin{equation}
			\label{eq:4.24}
			\|\bm{u}^\epsilon(\mathbf{x})-\bm{u}^{(\epsilon,h_0,h_1)}(\mathbf{x})\|_{(H^1(\Omega))^n}\leq C\bigl(\epsilon+h_0^2+h_1 \bigr).
		\end{equation}
	\end{thm}
	
	\textbf{Proof:}\hspace{1mm}Firstly, by employing the triangle inequality, we establish the following inequalities
	\begin{equation}
		\label{eq:4.25}
		\|T^\epsilon-T^{(\epsilon,h_0,h_1)}\|_{H^1(\Omega)}\leq \|T^\epsilon-T^{(2,\epsilon)}\|_{H^1(\Omega)}+\|T^{(2,\epsilon)}-T^{(\epsilon,h_0,h_1)}\|_{H^1(\Omega)},
	\end{equation}
	\begin{equation}
		\label{eq:4.26}
		\|c^\epsilon-c^{(\epsilon,h_0,h_1)}\|_{H^1(\Omega)}\leq \|c^\epsilon-c^{(2,\epsilon)}\|_{H^1(\Omega)}+\|c^{(2,\epsilon)}-c^{(\epsilon,h_0,h_1)}\|_{H^1(\Omega)},
	\end{equation}
	\begin{equation}
		\label{eq:4.27}
		\|\bm{u}^\epsilon-\bm{u}^{(\epsilon,h_0,h_1)}\|_{(H^1(\Omega))^n}\leq \|\bm{u}^\epsilon-\bm{u}^{(2,\epsilon)}\|_{(H^1(\Omega))^n}+\|\bm{u}^{(2,\epsilon)}-\bm{u}^{(\epsilon,h_0,h_1)}\|_{(H^1(\Omega))^n}.
	\end{equation}
	Next, by combining the definition of the multi-scale FE solutions given in \eqref{eq:4.19}-\eqref{eq:4.21} with the triangle inequality, the following inequalities	sequentially hold	
	\begin{equation}
		\label{eq:4.28}
		\begin{aligned}
			&\|T^{(2,\epsilon)}-T^{(\epsilon,h_0,h_1)}\|_{H^1(\Omega)}\\
			&\leq \| T^{(0)} - T^{(0,h_0,h_1)} \|_{H^1(\Omega)} + \Big\| \epsilon \mathcal{H}_{\alpha_1}\frac{\partial T^{(0)}}{\partial x_{\alpha_1}} - \epsilon \mathcal{H}_{\alpha_1}^{h_0}\frac{\partial T^{(0,h_0,h_1)}}{\partial x_{\alpha_1}} \Big\|_{H^1(\Omega)} \\
			&+ \Big\| \epsilon^2 \mathcal{H}_{\alpha_1\alpha_2}\frac{\partial^2T^{(0)}}{\partial x_{\alpha_1}\partial x_{\alpha_2}} - \epsilon^2 \mathcal{H}_{\alpha_1\alpha_2}^{h_0}\frac{\partial^2T^{(0,h_0,h_1)}}{\partial x_{\alpha_1}\partial x_{\alpha_2}} \Big\|_{H^1(\Omega)} \\
			&+ \Big\| \epsilon^2 \mathcal{R}_{\alpha_1}\frac{\partial T^{(0)}}{\partial x_{\alpha_1}} - \epsilon^2 \mathcal{R}_{\alpha_1}^{h_0}\frac{\partial T^{(0,h_0,h_1)}}{\partial x_{\alpha_1}} \Big\|_{H^1(\Omega)} \\
			&\triangleq M_1^T + M_2^T + M_3^T + M_4^T.
		\end{aligned}
	\end{equation}
	\begin{equation}
		\label{eq:4.29}
		\begin{aligned}
			&\|c^{(2,\epsilon)}-c^{(\epsilon,h_0,h_1)}\|_{H^1(\Omega)}\\
			&\leq \| c^{(0)} - c^{(0,h_0,h_1)} \|_{H^1(\Omega)} + \Big\| \epsilon \mathcal{L}_{\alpha_1}\frac{\partial c^{(0)}}{\partial x_{\alpha_1}} - \epsilon \mathcal{L}_{\alpha_1}^{h_0}\frac{\partial c^{(0,h_0,h_1)}}{\partial x_{\alpha_1}} \Big\|_{H^1(\Omega)} \\
			&+ \Big\| \epsilon^2 \mathcal{L}_{\alpha_1\alpha_2}\frac{\partial^2c^{(0)}}{\partial x_{\alpha_1}\partial x_{\alpha_2}} - \epsilon^2 \mathcal{L}_{\alpha_1\alpha_2}^{h_0}\frac{\partial^2c^{(0,h_0,h_1)}}{\partial x_{\alpha_1}\partial x_{\alpha_2}} \Big\|_{H^1(\Omega)} \\			
			& + \Big\| \epsilon^2 \mathcal{S}_{\alpha_1}\frac{\partial c^{(0)}}{\partial x_{\alpha_1}} - \epsilon^2 \mathcal{S}_{\alpha_1}^{h_0}\frac{\partial c^{(0,h_0,h_1)}}{\partial x_{\alpha_1}} \Big\|_{H^1(\Omega)} \\
			&\triangleq M_1^c + M_2^c + M_3^c + M_4^c.
		\end{aligned}
	\end{equation}
	\begin{equation}
		\label{eq:4.30}
		\begin{aligned}
			&\|u_{i}^{(2,\epsilon)}-u_{i}^{(\epsilon,h_0,h_1)}\|_{H^1(\Omega)}\\
			&\leq \| u_{i}^{(0)} - u_{i}^{(0,h_0,h_1)} \|_{H^1(\Omega)} + \Big\| \epsilon \mathcal{X}_{ih}^{\alpha_{1}}\frac{\partial u_{h}^{(0)}}{\partial x_{\alpha_{1}}} - \epsilon \mathcal{X}_{ih}^{\alpha_{1},h_0}\frac{\partial u_{h}^{(0,h_0,h_1)}}{\partial x_{\alpha_{1}}} \Big\|_{H^1(\Omega)} \\			
			&+ \| \epsilon \mathcal{M}_{i}T^{(0)} - \epsilon \mathcal{M}_{i}^{h_0} T^{(0,h_0,h_1)} \|_{H^1(\Omega)}  + \| \epsilon \mathcal{N}_{i}c^{(0)} - \epsilon \mathcal{N}_{i}^{h_0} c^{(0,h_0,h_1)} \|_{H^1(\Omega)} \\			
			&+ \Big\| \epsilon^2 \mathcal{P}_{ih}^{\alpha_{1}\alpha_{2}}\frac{\partial^{2}u_{h}^{(0)}}{\partial x_{\alpha_{1}}\partial x_{\alpha_{2}}} - \epsilon^2 \mathcal{P}_{ih}^{\alpha_{1}\alpha_{2},h_0}\frac{\partial^{2}u_{h}^{(0,h_0,h_1)}}{\partial x_{\alpha_{1}}\partial x_{\alpha_{2}}} \Big\|_{H^1(\Omega)} \\		
			& + \Big\| \epsilon^2 \mathcal{Q}_{ih}^{\alpha_{1}}\frac{\partial u_{h}^{(0)}}{\partial x_{\alpha_{1}}} - \epsilon^2 \mathcal{Q}_{ih}^{\alpha_{1},h_0}\frac{\partial u_{h}^{(0,h_0,h_1)}}{\partial x_{\alpha_{1}}} \Big\|_{H^1(\Omega)} \!+\! \| \epsilon^2 \mathcal{W}_{i}T^{(0)} - \epsilon^2 \mathcal{W}_{i}^{h_0} T^{(0,h_0,h_1)} \|_{H^1(\Omega)} \\			
			&+ \Big\| \epsilon^2 \mathcal{Z}_{i}^{\alpha_{1}}\frac{\partial T^{(0)}}{\partial x_{\alpha_{1}}} - \epsilon^2 \mathcal{Z}_{i}^{\alpha_{1},h_0}\frac{\partial T^{(0,h_0,h_1)}}{\partial x_{\alpha_{1}}} \Big\|_{H^1(\Omega)}  + \| \epsilon^2 \mathcal{F}_{i}c^{(0)} - \epsilon^2 \mathcal{F}_{i}^{h_0} c^{(0,h_0,h_1)} \|_{H^1(\Omega)} \\
			& + \Big\| \epsilon^2 \mathcal{G}_{i}^{\alpha_{1}}\frac{\partial c^{(0)}}{\partial x_{\alpha_{1}}} - \epsilon^2 \mathcal{G}_{i}^{\alpha_{1},h_0}\frac{\partial c^{(0,h_0,h_1)}}{\partial x_{\alpha_{1}}} \Big\|_{H^1(\Omega)} \\
			&\triangleq M_1^u + M_2^u + M_3^u + M_4^u + M_5^u + M_6^u + M_7^u + M_8^u + M_9^u + M_{10}^u.
		\end{aligned}
	\end{equation}
	
	Furthermore, we continue to evaluate each term in the inequalities \eqref{eq:4.28}-\eqref{eq:4.30}. According to reference \cite{R43}, these error terms reveal three distinct categories:
	\begin{enumerate}
		\item[(I)] $M_1^T$, $M_1^c$, $M_1^u$, which are all $\epsilon^0$-order terms;
		\item[(II)] $M_2^T$, $M_2^c$, $M_2^u$, $M_3^u$, $M_4^u$, which are all $\epsilon^1$-order terms;
		\item[(III)] $M_3^T$, $M_3^c$, $M_4^T$, $M_4^c$, $M_5^u$, $M_6^u$, $M_7^u$, $M_8^u$, $M_9^u$, $M_{10}^u$, which are all $\epsilon^2$-order terms.
	\end{enumerate}
	Thus, according to the lemma \ref{lem:4.3}, the following estimate is derived for type-I term $M_1^T$
	\begin{equation}
		\label{eq:4.31}
		\begin{aligned}
			M_1^T= \|T^{(0)} - T^{(0,h_0,h_1)}\|_{H^1(\Omega)}\leq C(h_0^2 + h_1).
		\end{aligned}
	\end{equation}
	Using the above approach, we obtain the same error estimates for other type-I terms.
	
	For type-II terms, with $M_2^T$ as an example, we apply the triangle inequality yielding
	\begin{equation}
		\label{eq:4.32}
		\begin{aligned}
			M_2^T & =\Big\| \epsilon \mathcal{H}_{\alpha_1}\frac{\partial T^{(0)}}{\partial x_{\alpha_1}} - \epsilon \mathcal{H}_{\alpha_1}^{h_0}\frac{\partial T^{(0,h_0,h_1)}}{\partial x_{\alpha_1}} \Big\|_{H^1(\Omega)}\\
			&\leq \Big\| \epsilon \bigl( \mathcal{H}_{\alpha_1} - \mathcal{H}_{\alpha_1}^{h_0}\bigr)\frac{\partial T^{(0)}}{\partial x_{\alpha_1}} \Big\|_{H^1(\Omega)} + \Big\| \epsilon \mathcal{H}_{\alpha_1}^{h_0} \Bigl( \frac{\partial T^{(0)}}{\partial x_{\alpha_1}} - \frac{\partial T^{(0, h_0)}}{\partial x_{\alpha_1}} \Bigr) \Big\|_{H^1(\Omega)} \\
			& + \Big\| \epsilon \mathcal{H}_{\alpha_1}^{h_0} \Bigl( \frac{\partial T^{(0, h_0)}}{\partial x_{\alpha_1}} - \frac{\partial T^{(0, h_0, h_1)}}{\partial x_{\alpha_1}} \Bigr) \Big\|_{H^1(\Omega)} \\
			&\triangleq e_1 + e_2 + e_3.
		\end{aligned}
	\end{equation}
	After that, by employing lemma \ref{lem:4.1} and inequalities \eqref{eq:4.16} and \eqref{eq:4.17}, and applying the Cauchy-Schwarz inequality, the following inequalities are obtained
	\begin{equation}
		\label{eq:4.33}
		e_1^2 \leq C \epsilon^2 h_0^2,\; e_2^2 \leq C \epsilon^2 h_0^4,\; e_3^2 \leq C \epsilon^2 h_1^2.
	\end{equation}
	Hence, on the basis of \eqref{eq:4.32} and \eqref{eq:4.33}, we derive the following error estimate for $M_2^T$
	\begin{equation}
		\label{eq:4.34}
		M_2^T \leq C\epsilon(h_0+h_0^2+h_1).
	\end{equation}
	Utilizing the above approach, the same error estimates hold for other type-II terms.
	
	Moreover, for type-III terms, taking $M_3^T$ as an example, the following error estimate yields by employing lemma \ref{lem:4.1} and the inequalities \eqref{eq:4.16} and \eqref{eq:4.17}, and applying the Cauchy-Schwarz inequality
	\begin{equation}
		\label{eq:4.35}
		M_3^T \leq C\epsilon^2 (h_0+h_0^2+h_1).
	\end{equation}
	Employing the above method, the same error estimates are obtained for other type-III terms.
	
	Finally, the inequality \eqref{eq:4.22} is easily derived by combining the proposed inequalities \eqref{eq:3.5}, \eqref{eq:4.25}, \eqref{eq:4.31}, \eqref{eq:4.34} and \eqref{eq:4.35} while neglecting higher-order error terms. Applying the same proof procedures, we can verify the inequalities \eqref{eq:4.23} and \eqref{eq:4.24}. In a summary, the proof of theorem \ref{thm:2} is complete.
	
	\section{Numerical examples and results}
	\label{sec:5}
	This section presents extensive numerical examples to validate the proposed HOMS computational model along with its numerical algorithm. All numerical experiments are conducted on the same computer equipped with an Intel Core i7-13650HX processor (2.60 GHz) and 24.0 GB RAM, and all numerical simulations are performed based on Freefem++ software.
	
	Given the difficulty in obtaining exact solutions for the multi-scale problem \eqref{eq:2.1}, we substitute $T^{\epsilon}, c^{\epsilon}$ and $\bm{u}^{\epsilon}$ with corresponding high-resolution FEM solutions $T_e, c_e$ and $\bm{u}_e$, which serve as reference solutions. Moreover, defining $\|\cdot\|_{L^{2}}$ and $\mid\cdot\mid _{H^{1}}$ as the $L^2$ norm and $H^1$ semi-norm respectively, the $L^2$ norm measures the global, integrated error in the solution, which is crucial for validating the macroscopic response, while the $H^1$ semi-norm measures the error in the gradient, which is the key metric for assessing the method's ability to capture microscopic oscillatory behavior. Then, the relative errors in the $L^2$ norm for $T^{(0)}$, $T^{(1,\epsilon)}$, $T^{(2,\epsilon)}$, $c^{(0)}$, $c^{(1,\epsilon)}$, $c^{(2,\epsilon)}$, $\bm{u}^{(0)}$, $\bm{u}^{(1,\epsilon)}$ and $\bm{u}^{(2,\epsilon)}$ are defined as $TerrorL^20$, $TerrorL^21$, $TerrorL^22$, $cerrorL^20$, $cerrorL^21$, $cerrorL^22$, $\bm{u}errorL^20$, $\bm{u}errorL^21$ and $\bm{u}errorL^22$, respectively. The relative errors in the $H^1$ semi-norm are defined as $TerrorH^10$, $TerrorH^11$, $TerrorH^12$, $cerrorH^10$, $cerrorH^11$, $cerrorH^12$, $\bm{u}errorH^10$, $\bm{u}errorH^11$ and $\bm{u}errorH^12$, respectively.
	
	\subsection{Example 1: 2D quasi-periodic composite structure}
	\label{sec:51}
	In this example, a 2D quasi-periodic composite structure is investigated. The detailed macroscopic structure $\Omega$ and microscopic unit cell $Y$ are shown in Fig.\hspace{1mm}\ref{f1:2D}, where $\Omega= (x_1,x_2)= [0,1]^2 \mathrm{cm}^2$ and small periodic parameter $\epsilon=1/10$.
	
	Numerical experiments will be conducted on two cases of quasi-periodic composite structures: with and without scale-separated material parameters. First, for quasi-periodic material parameters exhibiting scale-separated properties, they can be expressed as $a(\mathbf{x},\mathbf{y}) = \tilde{a}(\mathbf{x}) \cdot \hat{a}(\mathbf{y})$ in case 1, with specific parameter values given in Table\hspace{1mm}\ref{t1}. The associated weight function is defined as $\psi({\mathbf{x}}) = 5 + \sin(4\pi x_{1}) + \sin(4\pi x_{2})$. Second, when the material parameters possess scale-coupling properties without scale-separation, their material parameters can be expressed as $a(\mathbf{x},\mathbf{y}) = \tilde{a}(\mathbf{x}) + \hat{a}(\mathbf{y})$ in case 2, with detailed parameters provided in Table\hspace{1mm}\ref{t2}. The weight function for this case is $\psi(\mathbf{x}) = (x_{1} - 0.5)^{2} \cdot (x_{2} - 0.5)^{2}$.
	\begin{table}[!htb]{\caption{Material property parameters with scale-separation.}\label{t1}}
		\centering
		\begin{tabular}{cccc}
			\hline
			Property & Matrix & Inclusion & $\tilde{a}(x)$\\
			\hline
			Young's modulus $E\bigl(\mathrm{GPa}\bigr)$ & 10.0     & 1.0 & $\psi(\mathbf{x})$\\
			Poisson's ratio $\nu$ & 0.30   & 0.25 & 1.0 \\
			Thermal conductivity $k_{ij}\bigl(\mathrm{W / (m K})\bigr)$ & 100.0   & 1.0 & $\psi(\mathbf{x})$\\
			Moisture diffusion $g_{ij}\bigl(\mathrm{10^{-12}m^2/s}\bigr)$ & 1.0 & 0.02 &$\psi(\mathbf{x})$\\
			Thermal expansion $\alpha_{kl}\bigl(10^{-6}/\mathrm{K}\bigr)$ & 10.0 & 0.1 & $\psi(\mathbf{x})$\\
			Moisture expansion $\beta_{kl}\bigl(10^{-3}/(\mathrm{kg/m^3})\bigr)$ & 1.0 & 0.02 & $\psi(\mathbf{x})$\\
			\hline
		\end{tabular}
	\end{table}
	\begin{table}[!htb]{\caption{Material property parameters without scale-separation.}\label{t2}}
		\centering
		\begin{tabular}{cccc}
			\hline
			Property & Matrix & Inclusion & $\tilde{a}(x)$\\
			\hline
			Young's modulus $E\bigl(\mathrm{GPa}\bigr)$ & 10.0     & 1.0 & $0.5\psi(\mathbf{x})$\\
			Poisson's ratio $\nu$ & 0.30   & 0.25 & 0.0 \\
			Thermal conductivity $k_{ij}\bigl(\mathrm{W / (m K})\bigr)$ & 100.0   & 1.0 & $0.005\psi(\mathbf{x})$\\
			Moisture diffusion $g_{ij}\bigl(\mathrm{10^{-12}m^2/s}\bigr)$ & 1.0 & 0.02 &$0.01\psi(\mathbf{x})$\\
			Thermal expansion $\alpha_{kl}\bigl(10^{-6}/\mathrm{K}\bigr)$ & 10.0 & 0.1 & $0.005\psi(\mathbf{x})$\\
			Moisture expansion $\beta_{kl}\bigl(10^{-3}/(\mathrm{kg/m^3})\bigr)$ & 1.0 & 0.02 & $0.01\psi(\mathbf{x})$\\
			\hline
		\end{tabular}
	\end{table}
	
	Moreover, the heat source, moisture source, and body forces are given by $h=500 \mathrm{J/(cm^{3}\cdot s)}$, $m=500 \mu\mathrm{ g/(cm^{3}\cdot s)}$, and $(f_1,f_2)=(1000,1000) \mathrm{N/cm^3}$, respectively. The boundary conditions on $\partial\Omega$ are prescribed as $\overline{T}=273.15 \mathrm{K}$, $\overline{c}=0 \mathrm{g/cm^3}$ and $\overline{\bm{u}}=0 \mathrm{cm}$.
	
	Conducting the HOMS method and the precise FEM for the investigated 2D quasi-periodic composite structures respectively, the computational cost for this example is presented in Table\hspace{1mm}\ref{t3}, including the numbers of FEM nodes and elements, as well as the computational time required to numerical experiments for two cases with scale-separated parameters and those without scale-separated parameters.
	\begin{table}[!htb]{\caption{Comparison of computational cost.}\label{t3}}
		\centering
		\begin{tabular}{cccc}
			\hline
			& Cell equations & Homogenized equations & Multi-scale equations\\
			\hline
			FEM nodes & 461 & 2601 & 35761\\
			FEM elements & 840 & 5000 & 70800 \\
			\hline
			& & HOMS method & precise FEM\\
			\hline
			\multicolumn{2}{c}{Computing time for case 1} & 5.520s & 13.424s \\
			\multicolumn{2}{c}{Computing time for case 2} & 315.386s & 16.198s \\
			\hline
		\end{tabular}
	\end{table}
	
	As shown in Table\hspace{1mm}\ref{t3}, the HOMS method significantly reduces computational resource requirements compared to the precise FEM. For material parameters exhibiting scale-separation, the HOMS method decreases computational time and enhances efficiency. However, without scale-separation, the HOMS method exhibits marginally longer computational times than precise FEM in this steady-state problem. For time-dependent problems, the proposed HOMS method achieves progressively higher computational efficiency, when auxiliary cell functions are precomputed and stored. Therefore, the HOMS method shall outperform precise FEM in time-dependent dynamic simulations. Furthermore, although the HOMS approach exhibits apparently lower efficiency than precise FEM in this example, its computational time grows at a slower rate as the scale of microscopic unit cells increases. This comparative advantage will be validated in subsequent Example\hspace{1mm}3. In a summary, the HOMS method substantially reduces computational cost compared to the precise FEM, which is of great significance in engineering computations.
	
\begin{figure}[!htb]
		\centering
		\begin{minipage}[c]{0.24\textwidth}
			\centering
			\includegraphics[width=\linewidth]{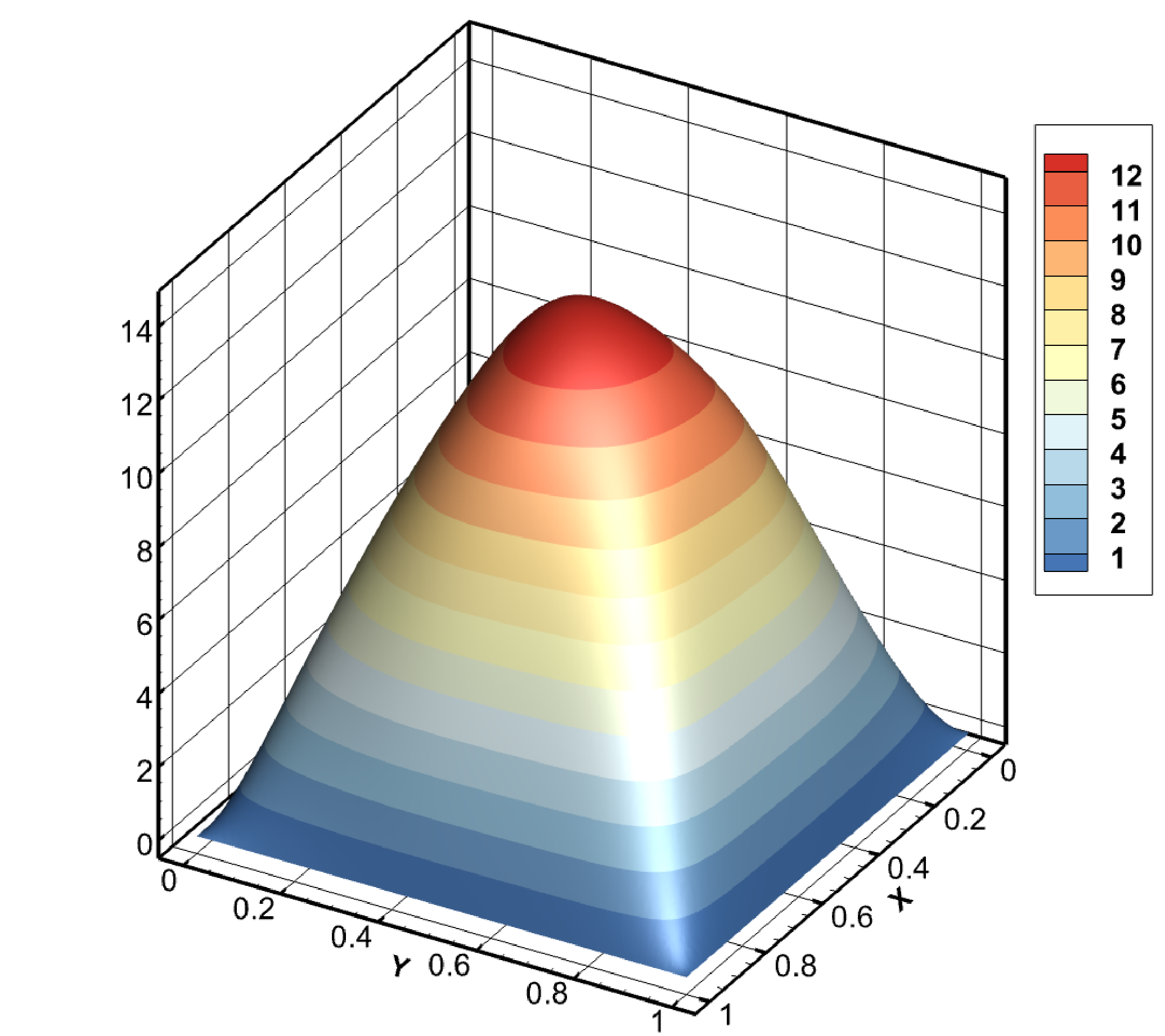}
			(a)
		\end{minipage}
		\hfill
		\begin{minipage}[c]{0.24\textwidth}
			\centering
			\includegraphics[width=\linewidth]{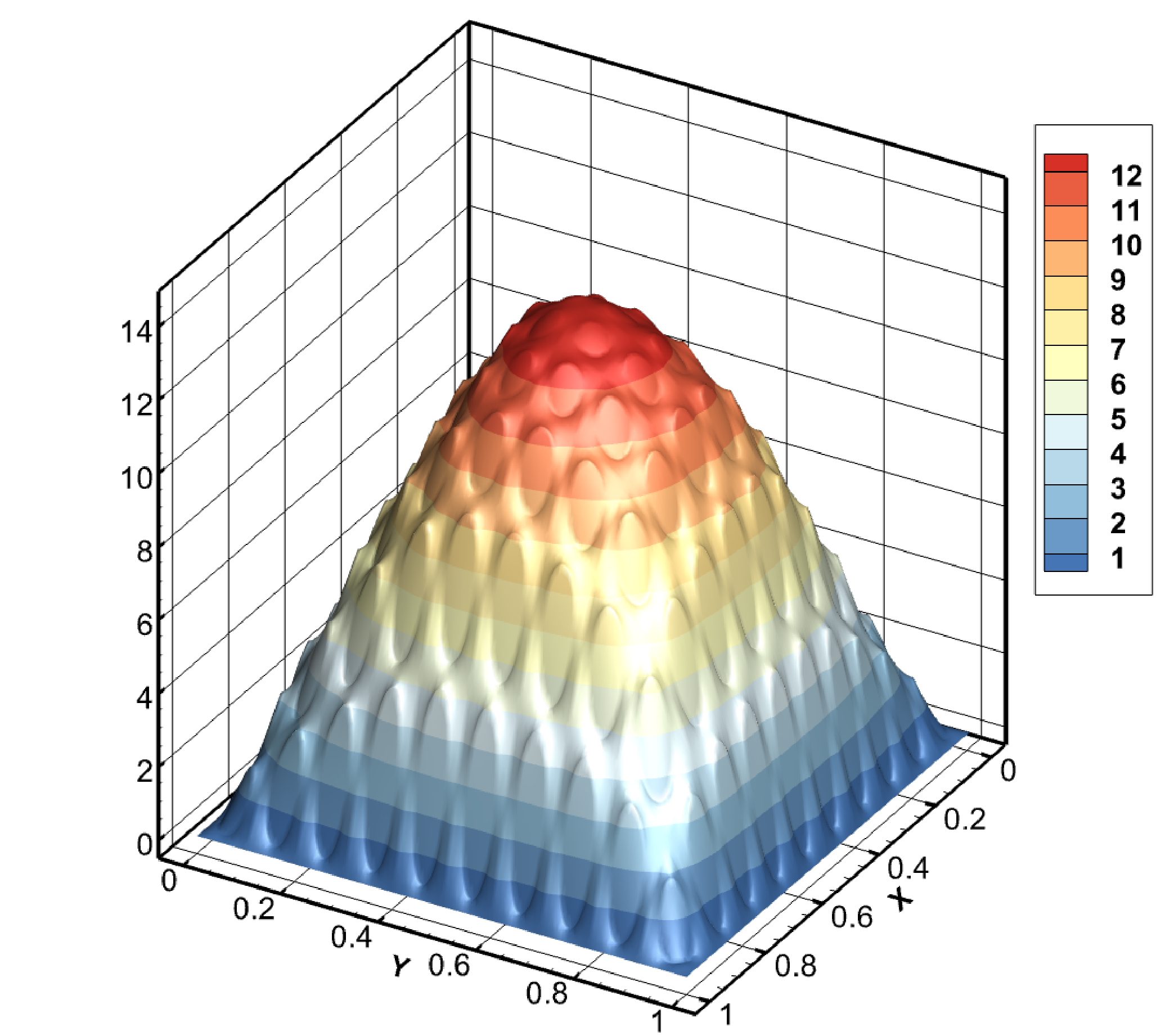}
			(b)
		\end{minipage}
		\hfill
		\begin{minipage}[c]{0.24\textwidth}
			\centering
			\includegraphics[width=\linewidth]{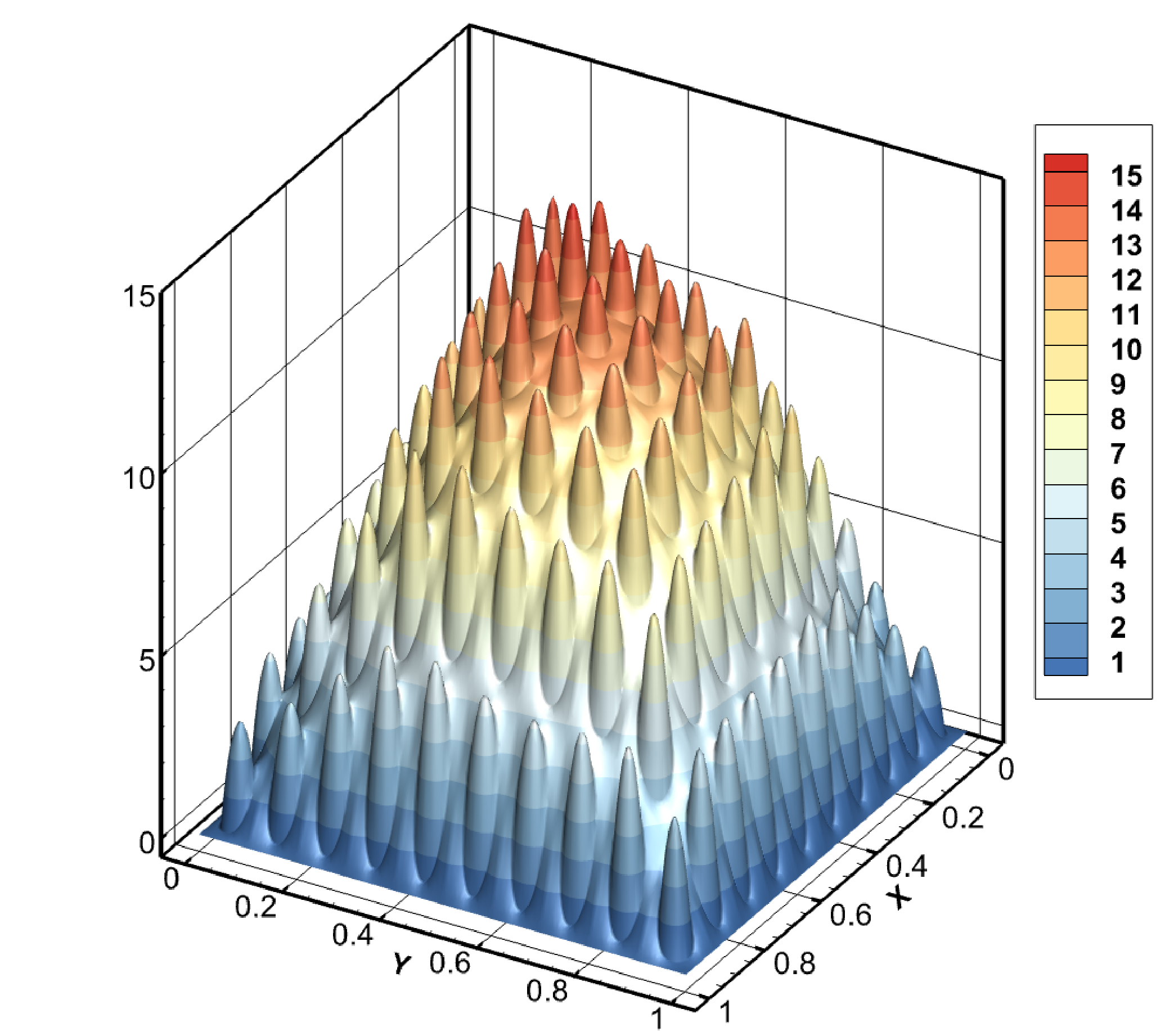}
			(c)
		\end{minipage}
		\hfill
		\begin{minipage}[c]{0.24\textwidth}
			\centering
			\includegraphics[width=\linewidth]{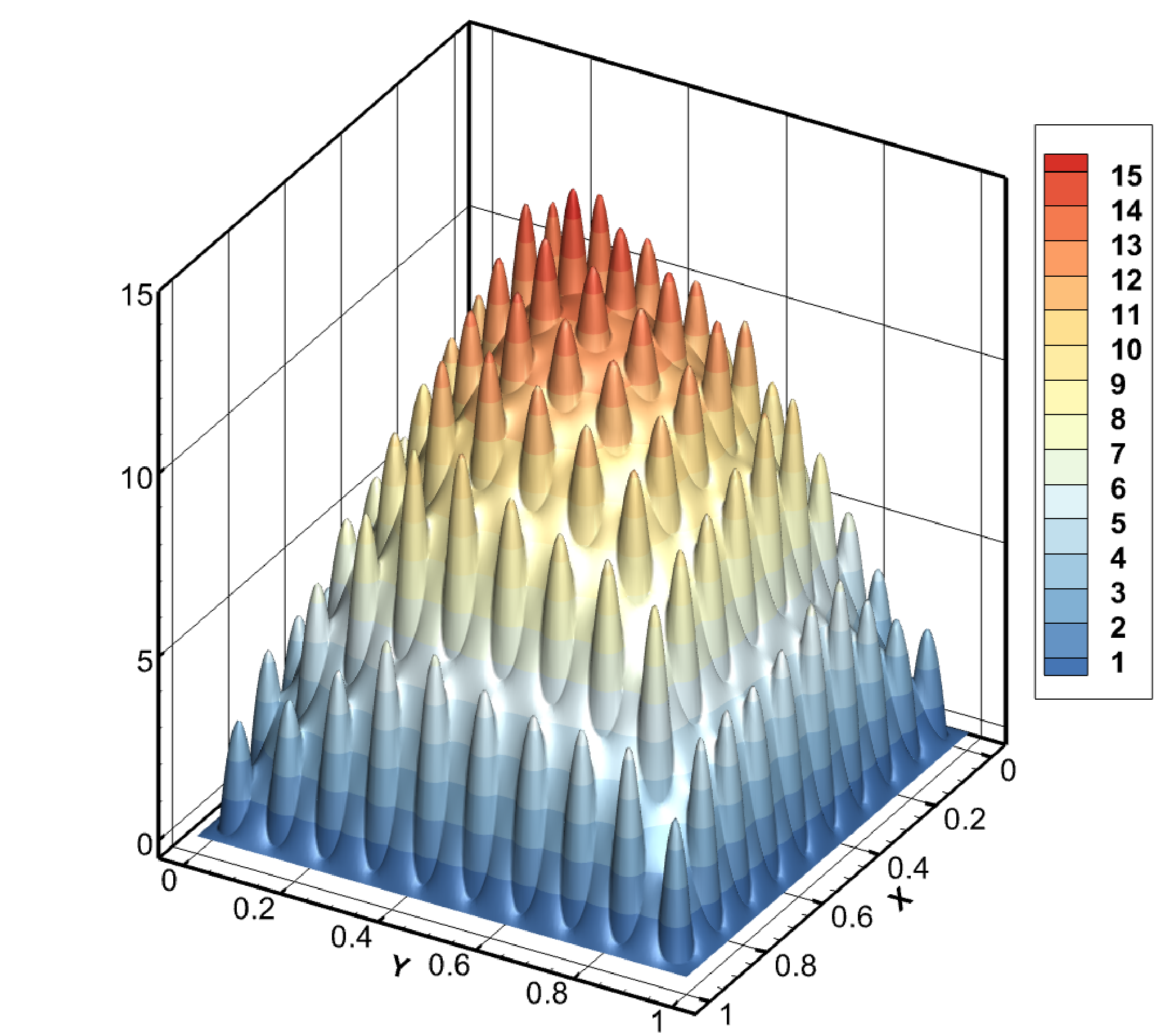}
			(d)
		\end{minipage}
		\caption{Temperature increment field with scale separation: (a) $T^{(0)}$; (b) $T^{(1,\epsilon)}$; (c) $T^{(2,\epsilon)}$; (d) $T^{\epsilon}$.}\label{f2}
	\end{figure}

\begin{figure}[!htb]
		\centering
		\begin{minipage}[c]{0.24\textwidth}
			\centering
			\includegraphics[width=\linewidth]{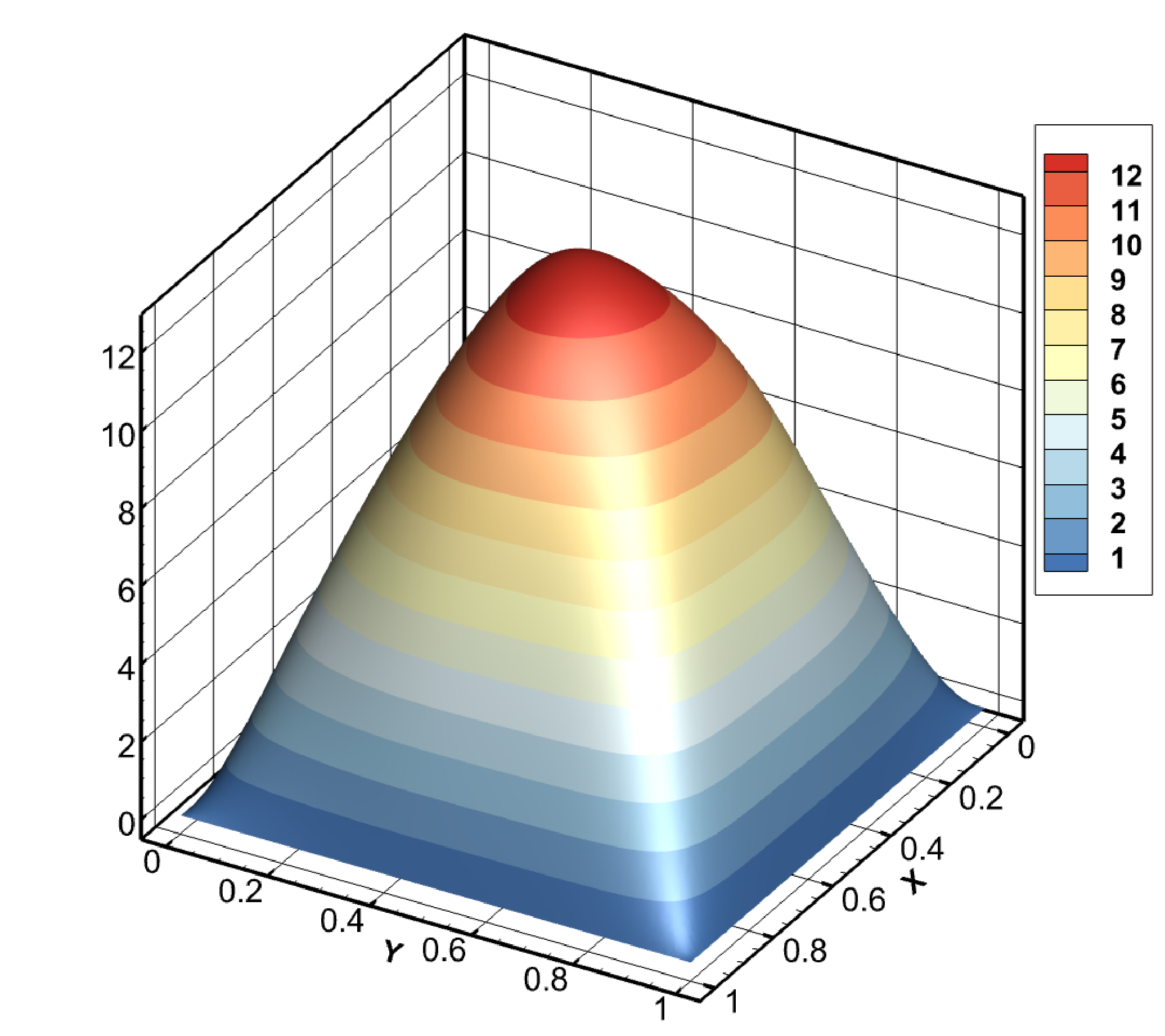}
			(a)
		\end{minipage}
		\hfill
		\begin{minipage}[c]{0.24\textwidth}
			\centering
			\includegraphics[width=\linewidth]{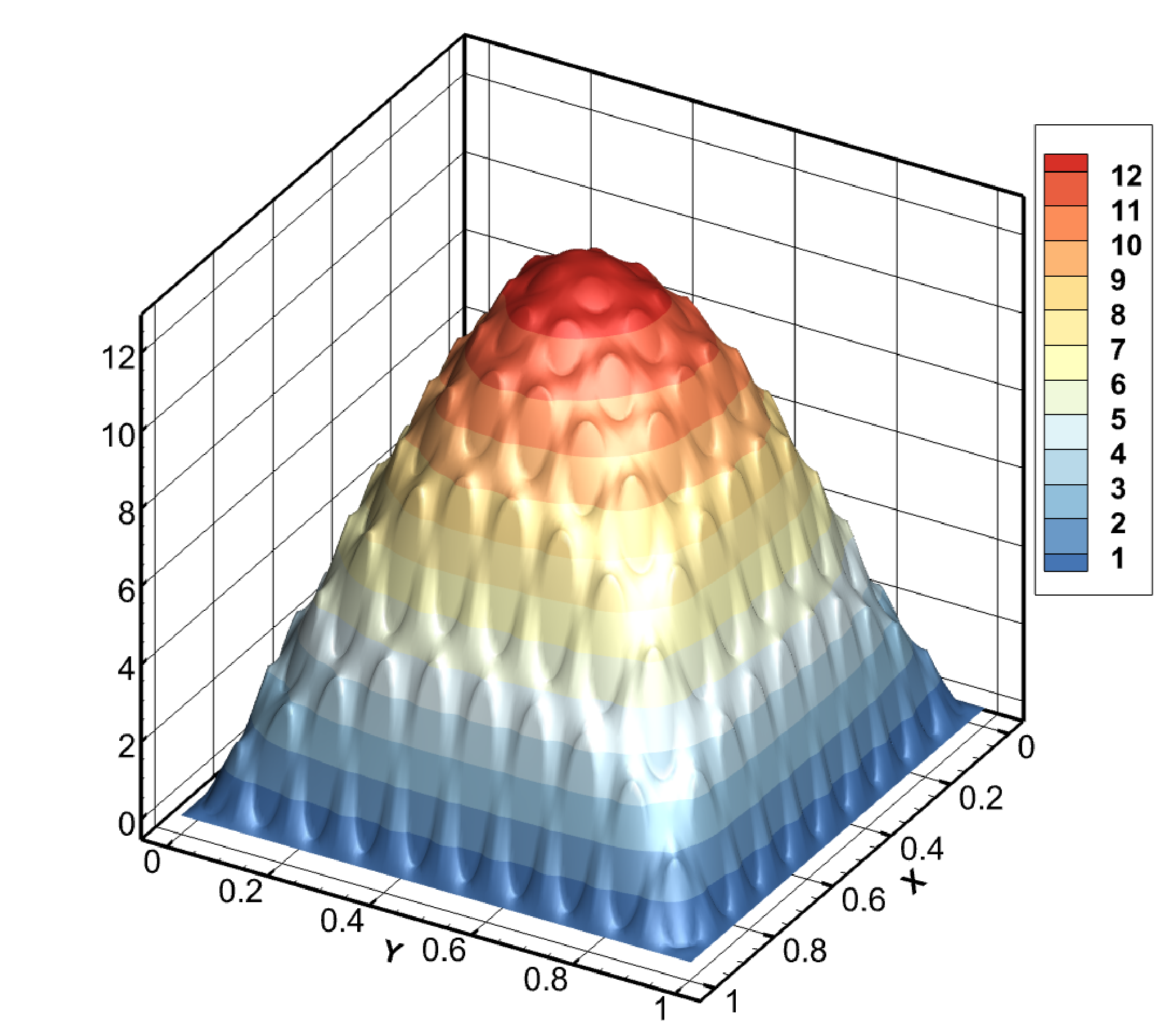}
			(b)
		\end{minipage}
		\hfill
		\begin{minipage}[c]{0.24\textwidth}
			\centering
			\includegraphics[width=\linewidth]{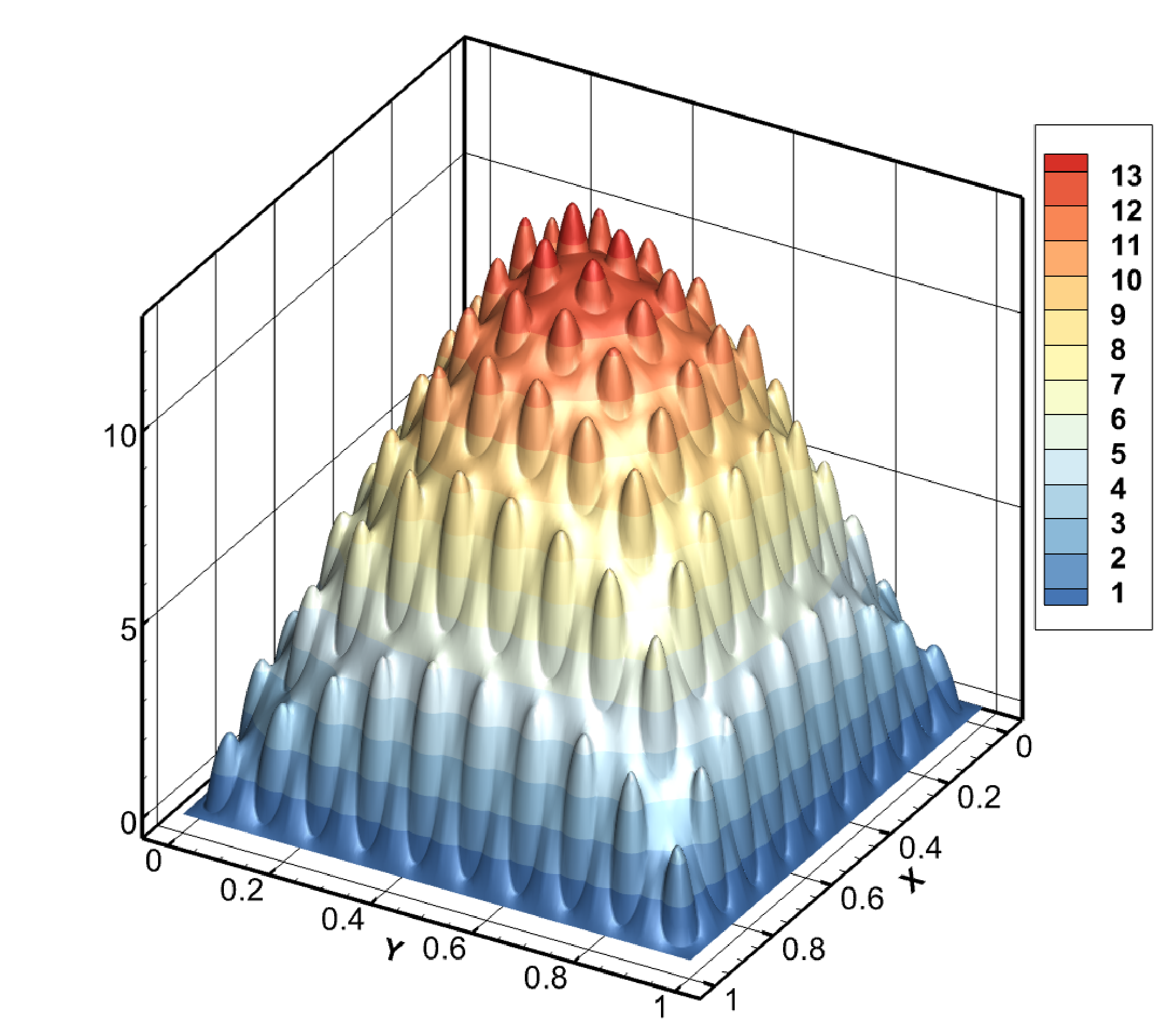}
			(c)
		\end{minipage}
		\hfill
		\begin{minipage}[c]{0.24\textwidth}
			\centering
			\includegraphics[width=\linewidth]{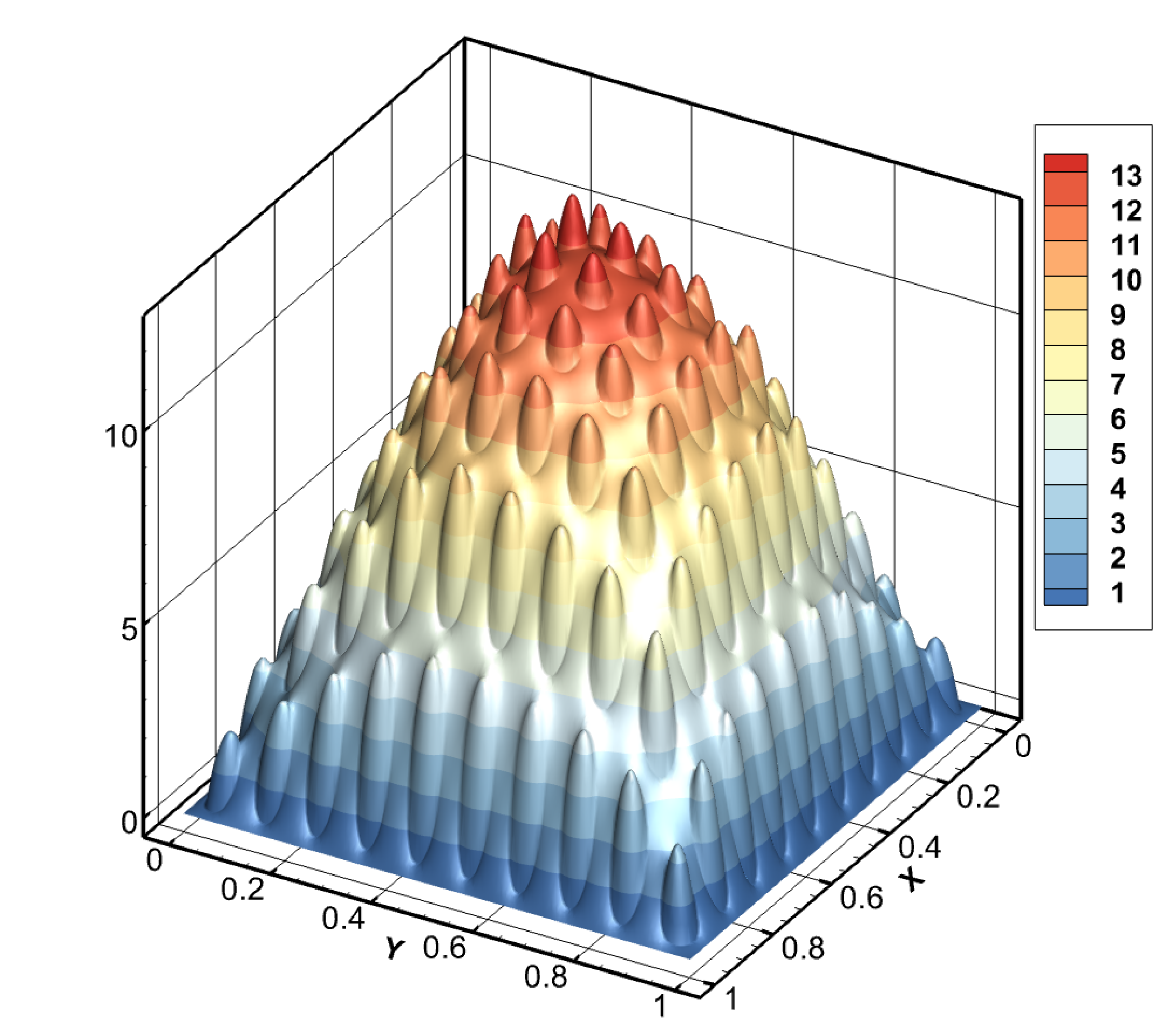}
			(d)
		\end{minipage}
		\caption{Moisture field with scale separation: (a) $c^{(0)}$; (b) $c^{(1,\epsilon)}$; (c) $c^{(2,\epsilon)}$; (d) $c^{\epsilon}$.}\label{f3}
	\end{figure}

\begin{figure}[!htb]
		\centering
		\begin{minipage}[c]{0.24\textwidth}
			\centering
			\includegraphics[width=\linewidth]{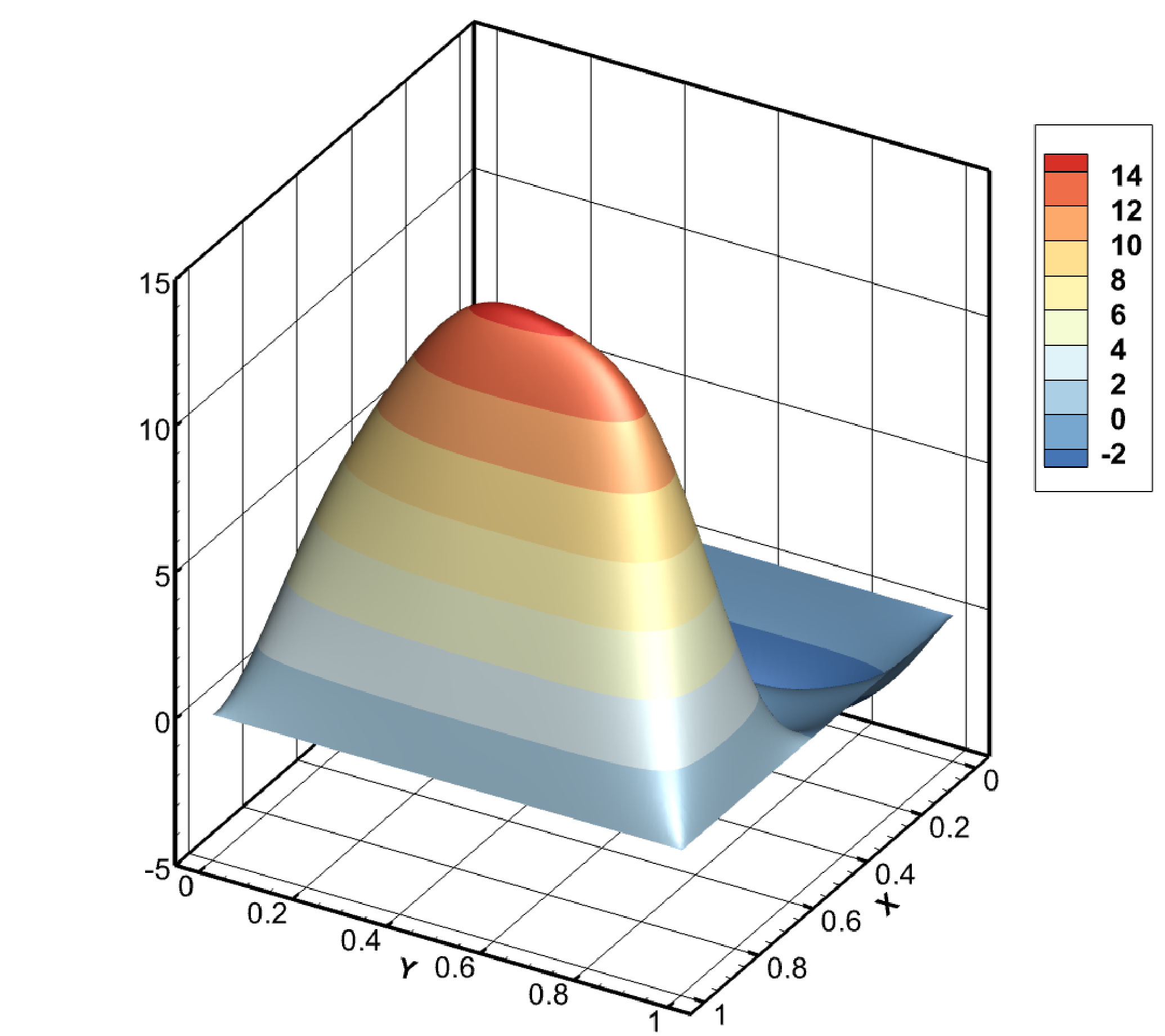}
			(a)
		\end{minipage}
		\hfill
		\begin{minipage}[c]{0.24\textwidth}
			\centering
			\includegraphics[width=\linewidth]{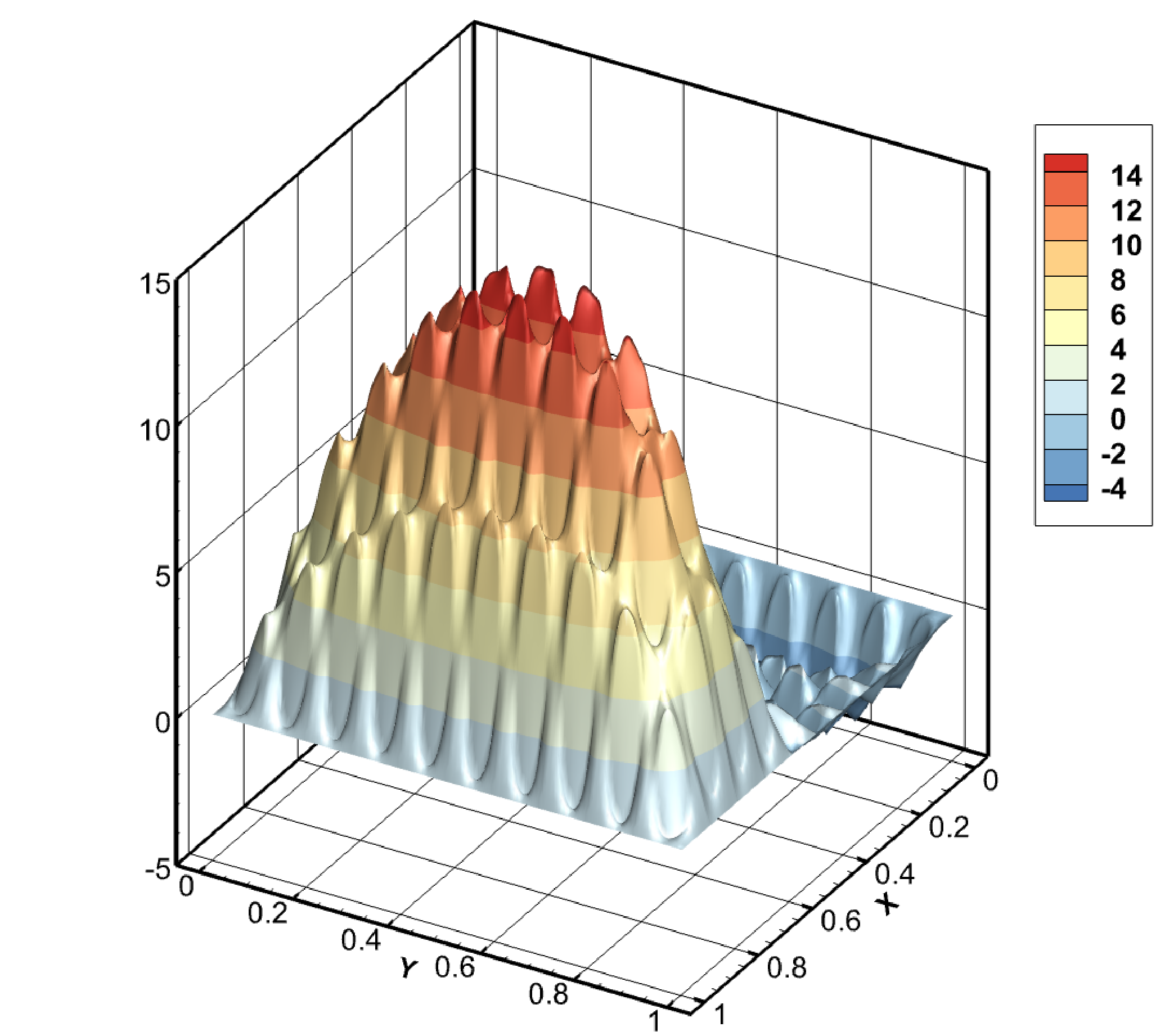}
			(b)
		\end{minipage}
		\hfill
		\begin{minipage}[c]{0.24\textwidth}
			\centering
			\includegraphics[width=\linewidth]{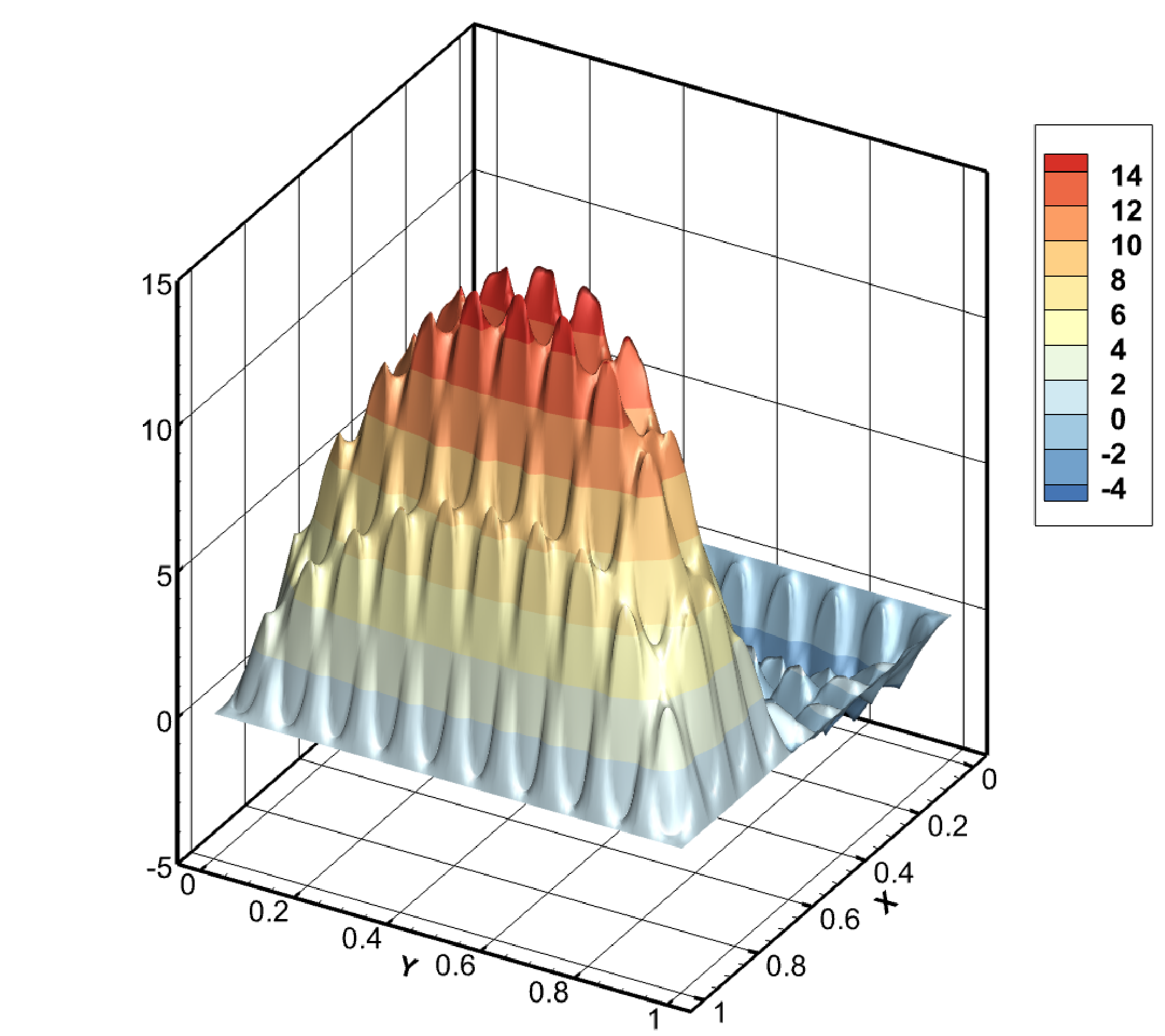}
			(c)
		\end{minipage}
		\hfill
		\begin{minipage}[c]{0.24\textwidth}
			\centering
			\includegraphics[width=\linewidth]{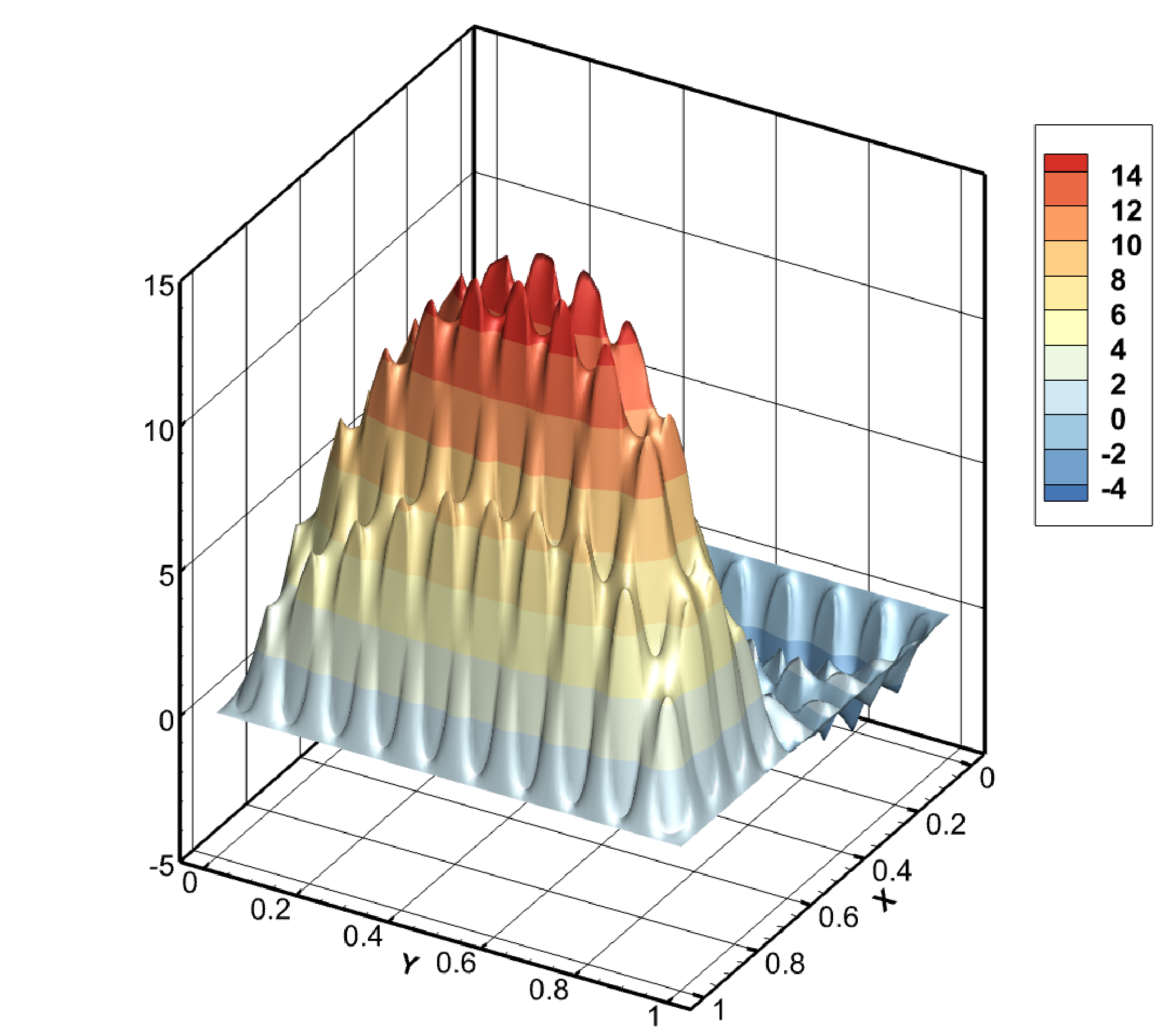}
			(d)
		\end{minipage}
		\caption{First displacement field component with scale separation: (a) $u_1^{(0)}$; (b) $u_1^{(1,\epsilon)}$; (c) $u_1^{(2,\epsilon)}$; (d) $u_1^{\epsilon}$.}\label{f4}
	\end{figure}

\begin{figure}[!htb]
		\centering
		\begin{minipage}[c]{0.24\textwidth}
			\centering
			\includegraphics[width=\linewidth]{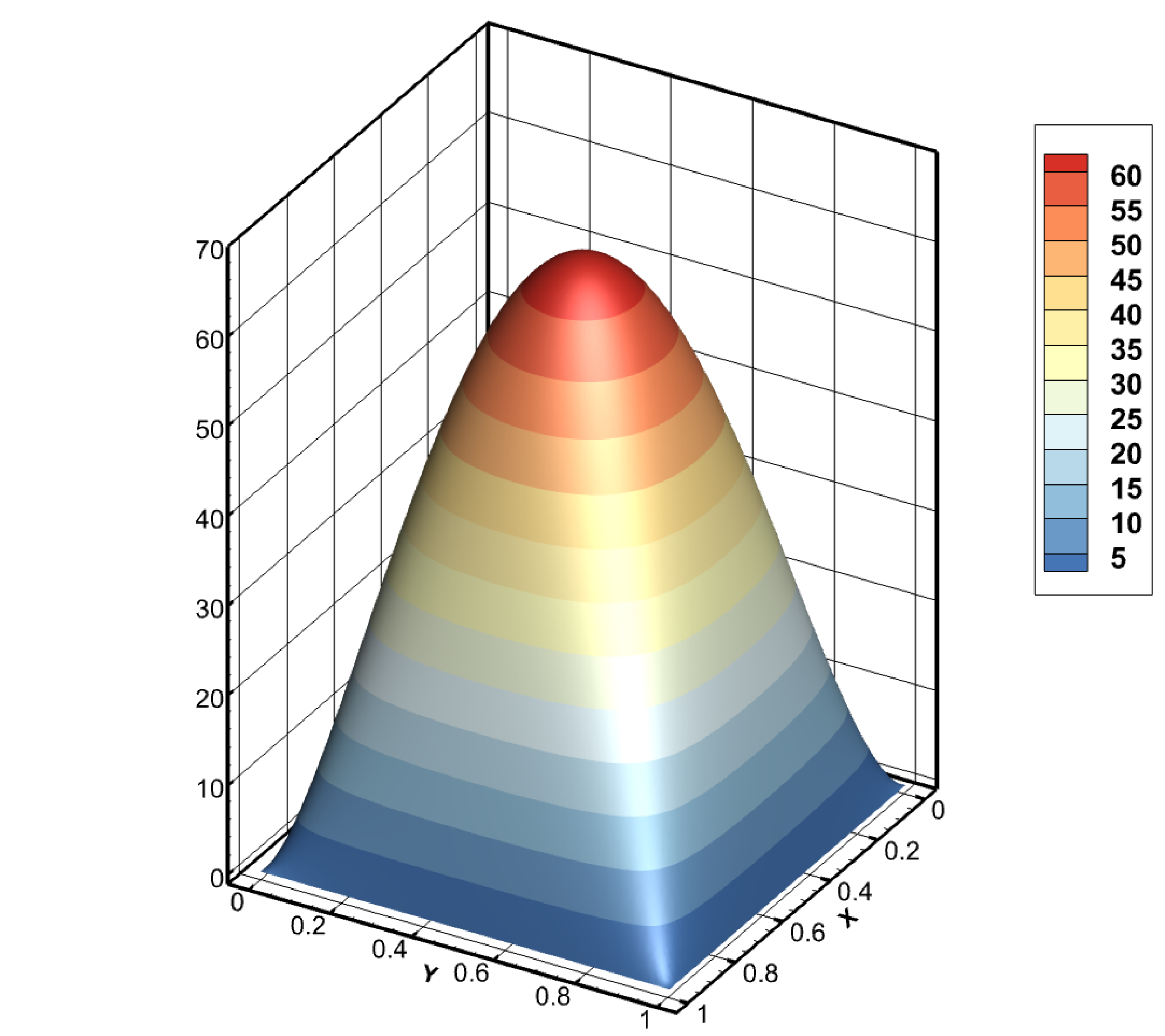}
			(a)
		\end{minipage}
		\hfill
		\begin{minipage}[c]{0.24\textwidth}
			\centering
			\includegraphics[width=\linewidth]{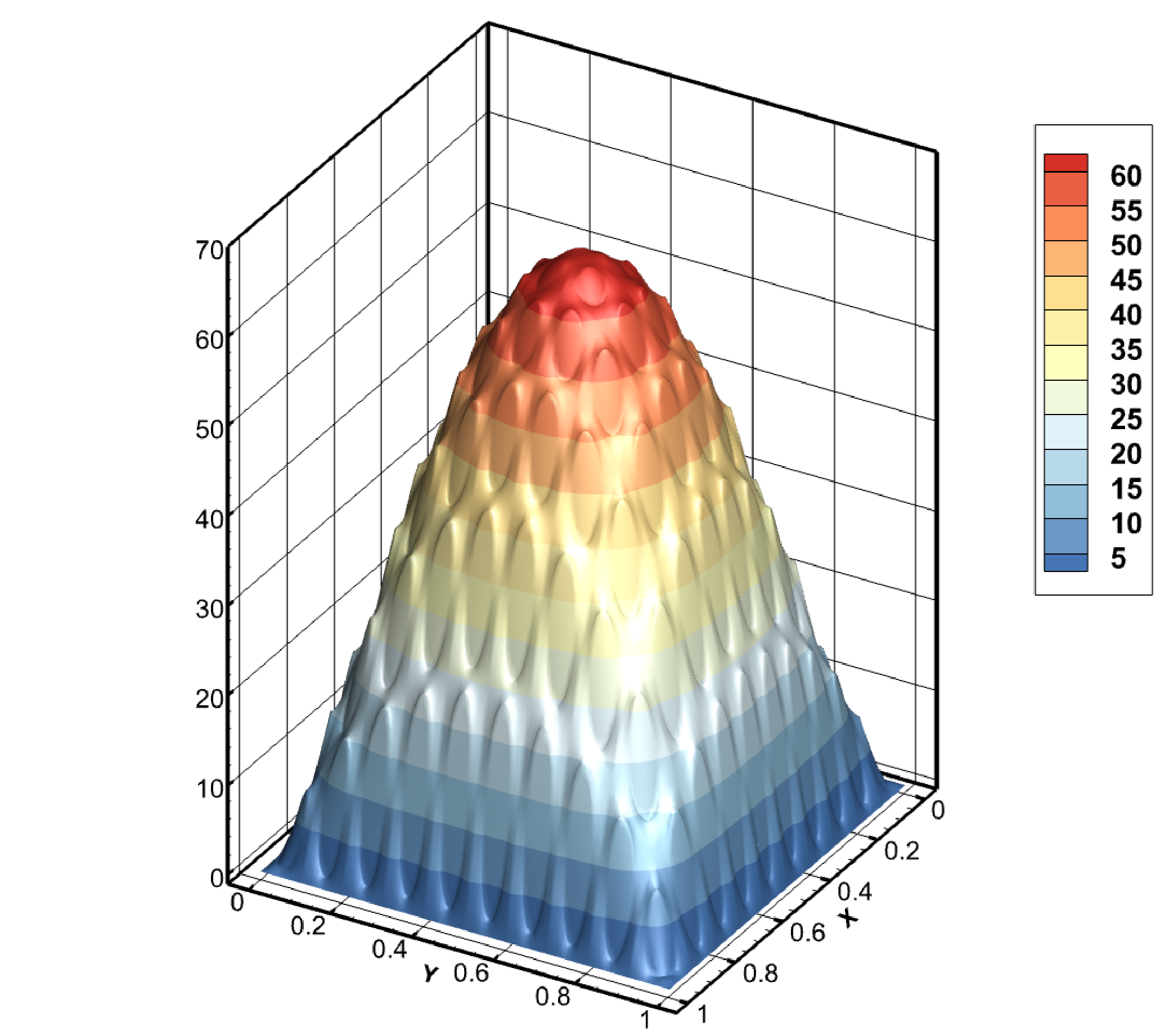}
			(b)
		\end{minipage}
		\hfill
		\begin{minipage}[c]{0.24\textwidth}
			\centering
			\includegraphics[width=\linewidth]{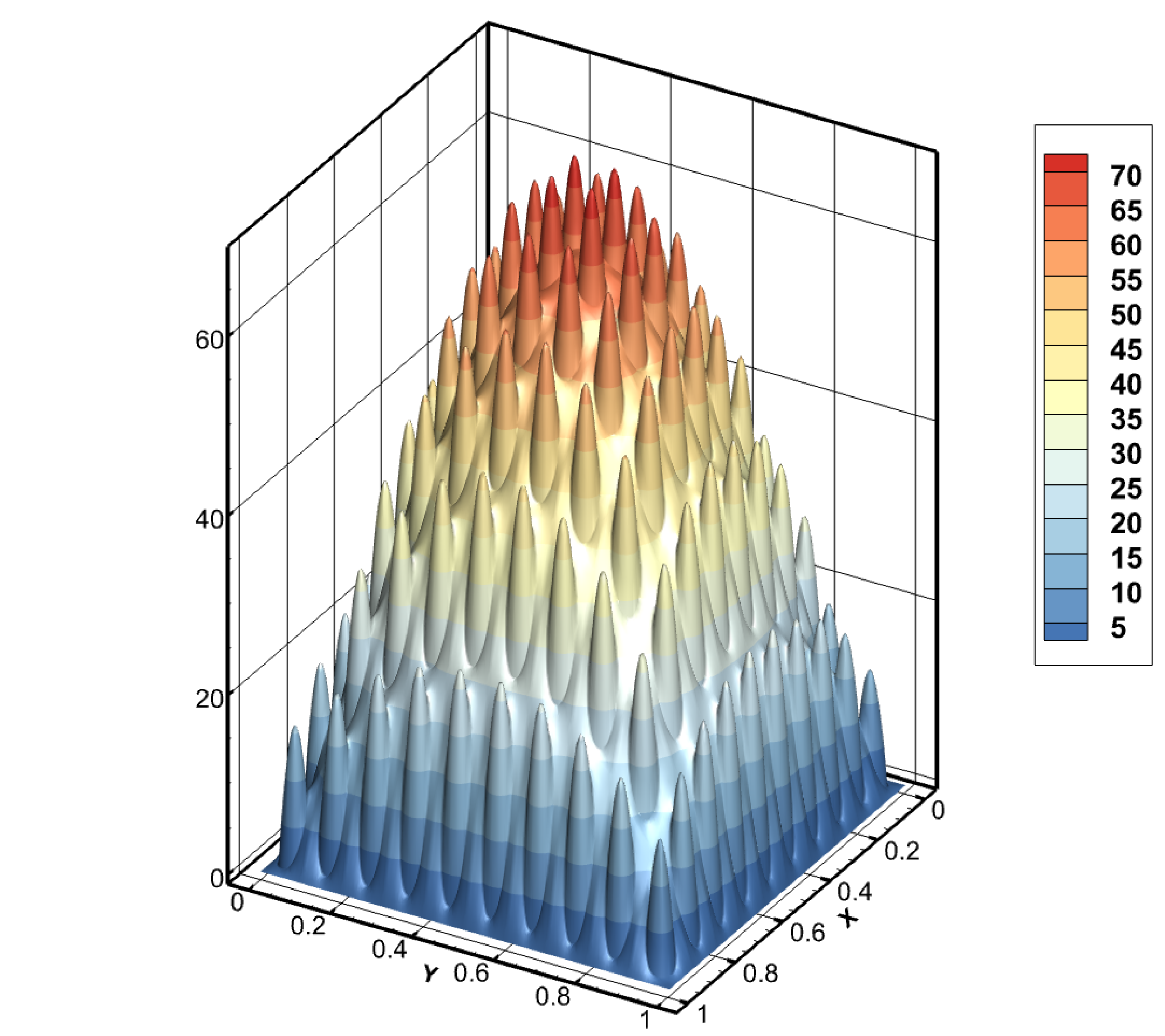}
			(c)
		\end{minipage}
		\hfill
		\begin{minipage}[c]{0.24\textwidth}
			\centering
			\includegraphics[width=\linewidth]{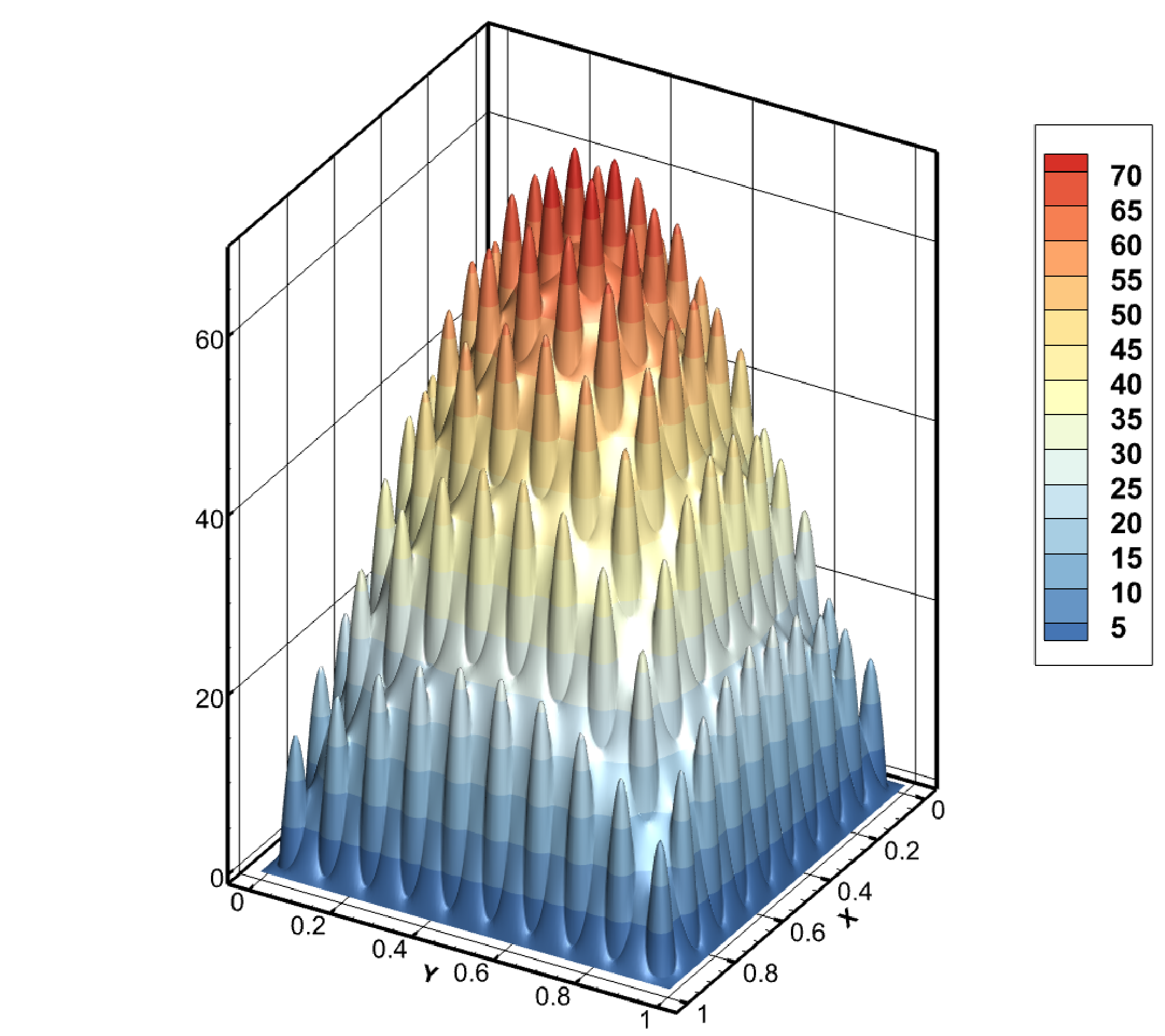}
			(d)
		\end{minipage}
		\caption{Temperature increment field under scale coupling: (a) $T^{(0)}$; (b) $T^{(1,\epsilon)}$; (c) $T^{(2,\epsilon)}$; (d) $T^{\epsilon}$.}\label{f6}
	\end{figure}

\begin{figure}[!htb]
		\centering
		\begin{minipage}[c]{0.24\textwidth}
			\centering
			\includegraphics[width=\linewidth]{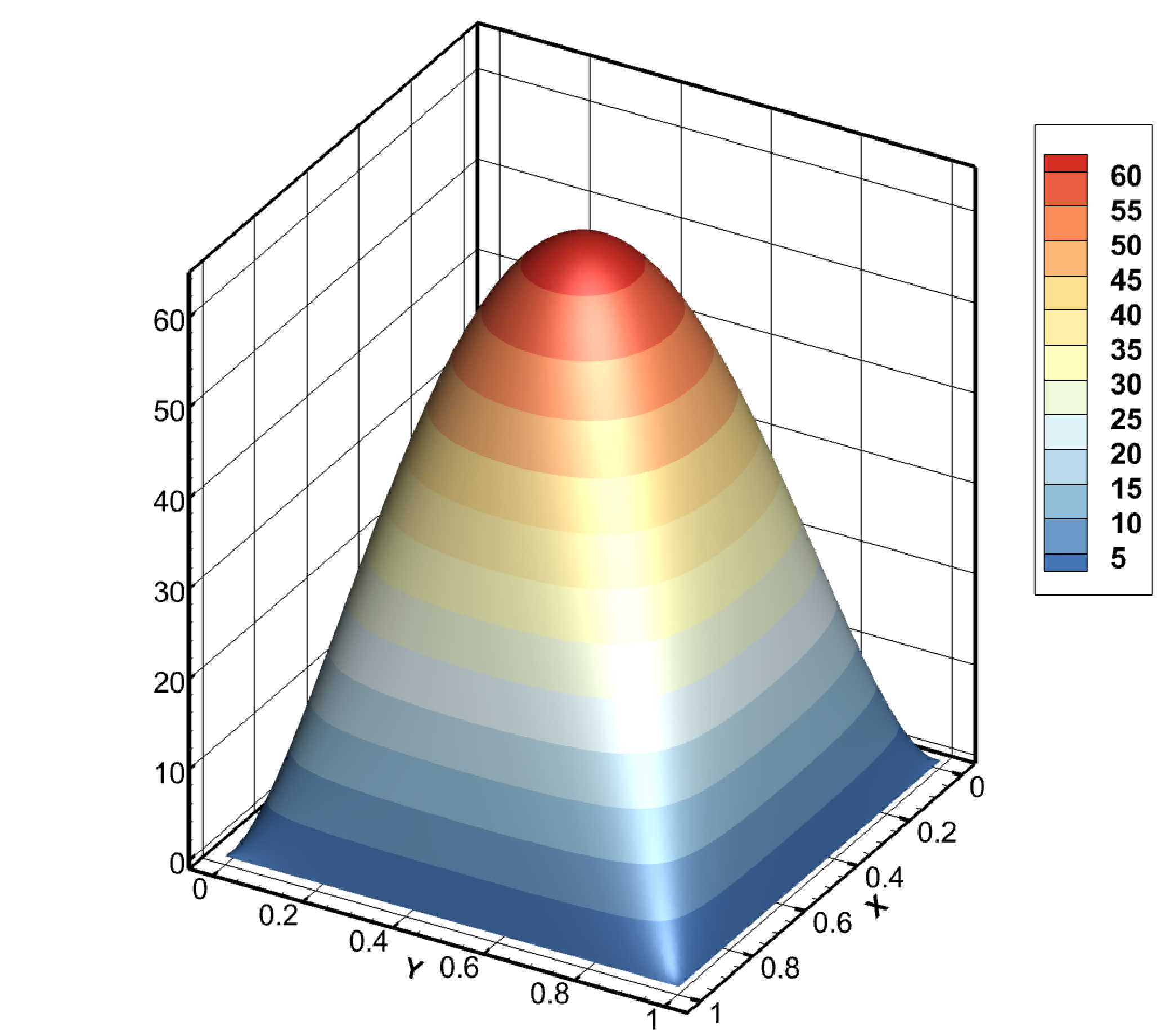}
			(a)
		\end{minipage}
		\hfill
		\begin{minipage}[c]{0.24\textwidth}
			\centering
			\includegraphics[width=\linewidth]{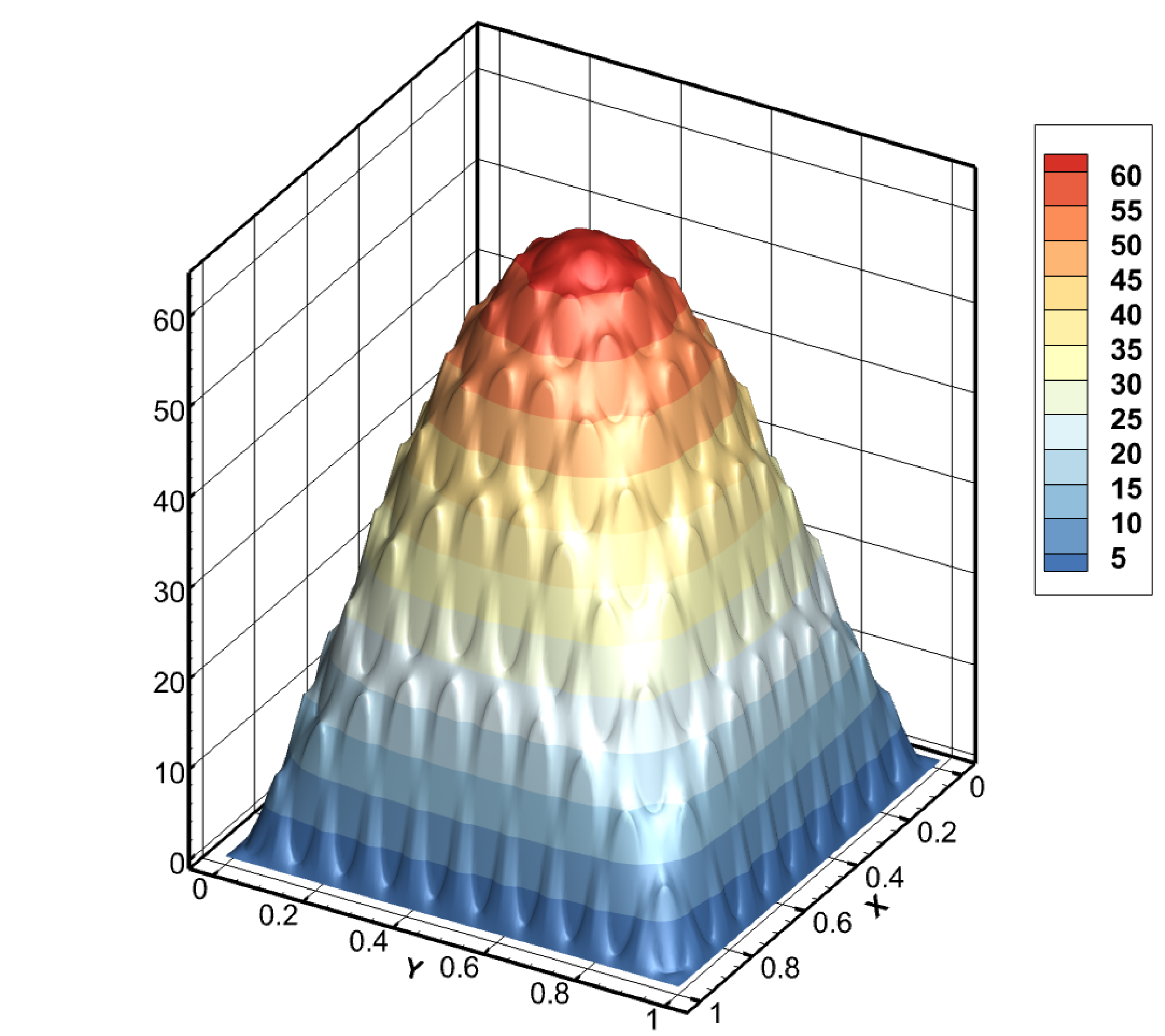}
			(b)
		\end{minipage}
		\hfill
		\begin{minipage}[c]{0.24\textwidth}
			\centering
			\includegraphics[width=\linewidth]{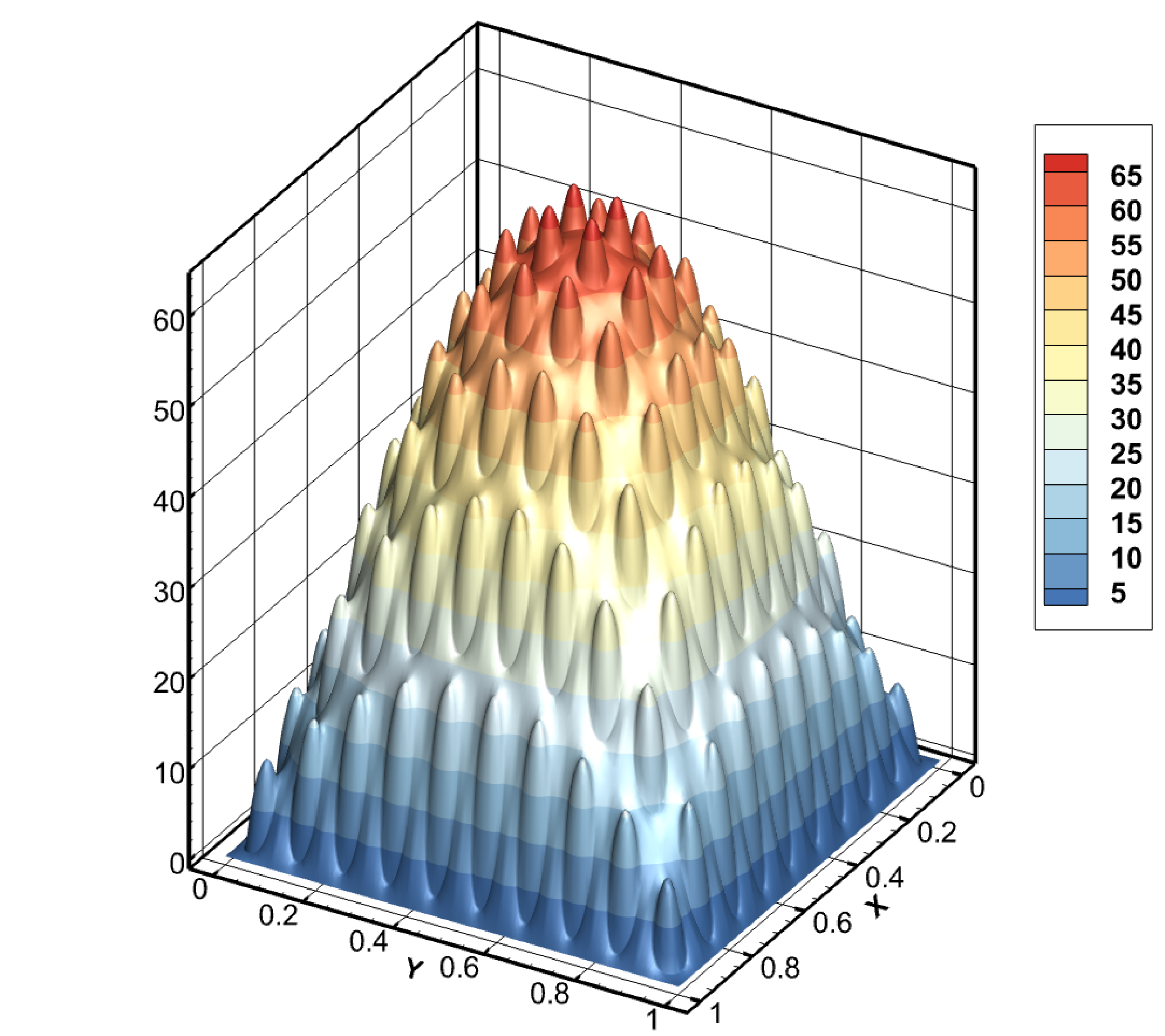}
			(c)
		\end{minipage}
		\hfill
		\begin{minipage}[c]{0.24\textwidth}
			\centering
			\includegraphics[width=\linewidth]{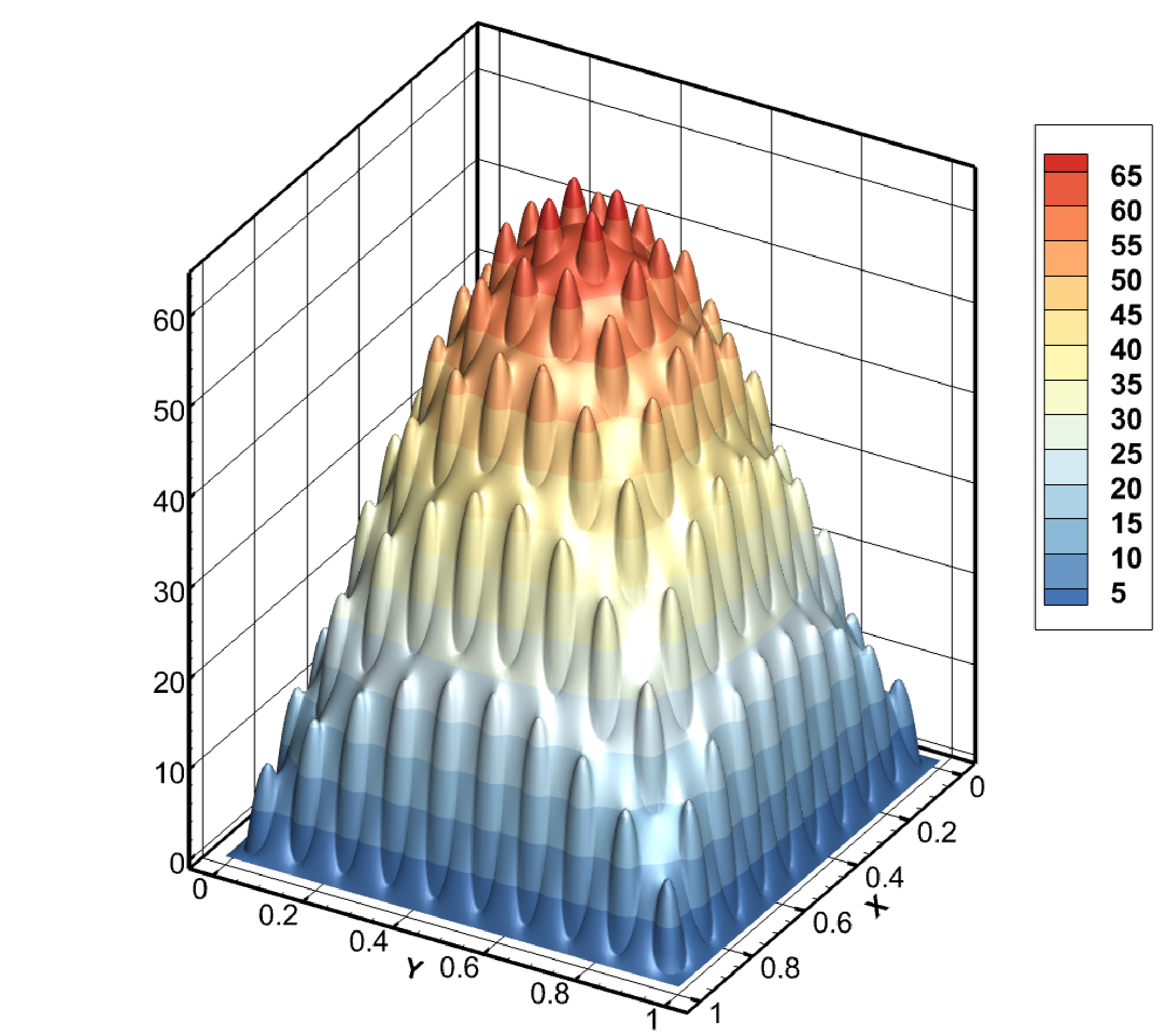}
			(d)
		\end{minipage}
		\caption{Moisture field under scale coupling: (a) $c^{(0)}$; (b) $c^{(1,\epsilon)}$; (c) $c^{(2,\epsilon)}$; (d) $c^{\epsilon}$.}\label{f7}
	\end{figure}

\begin{figure}[!htb]
		\centering
		\begin{minipage}[c]{0.24\textwidth}
			\centering
			\includegraphics[width=\linewidth]{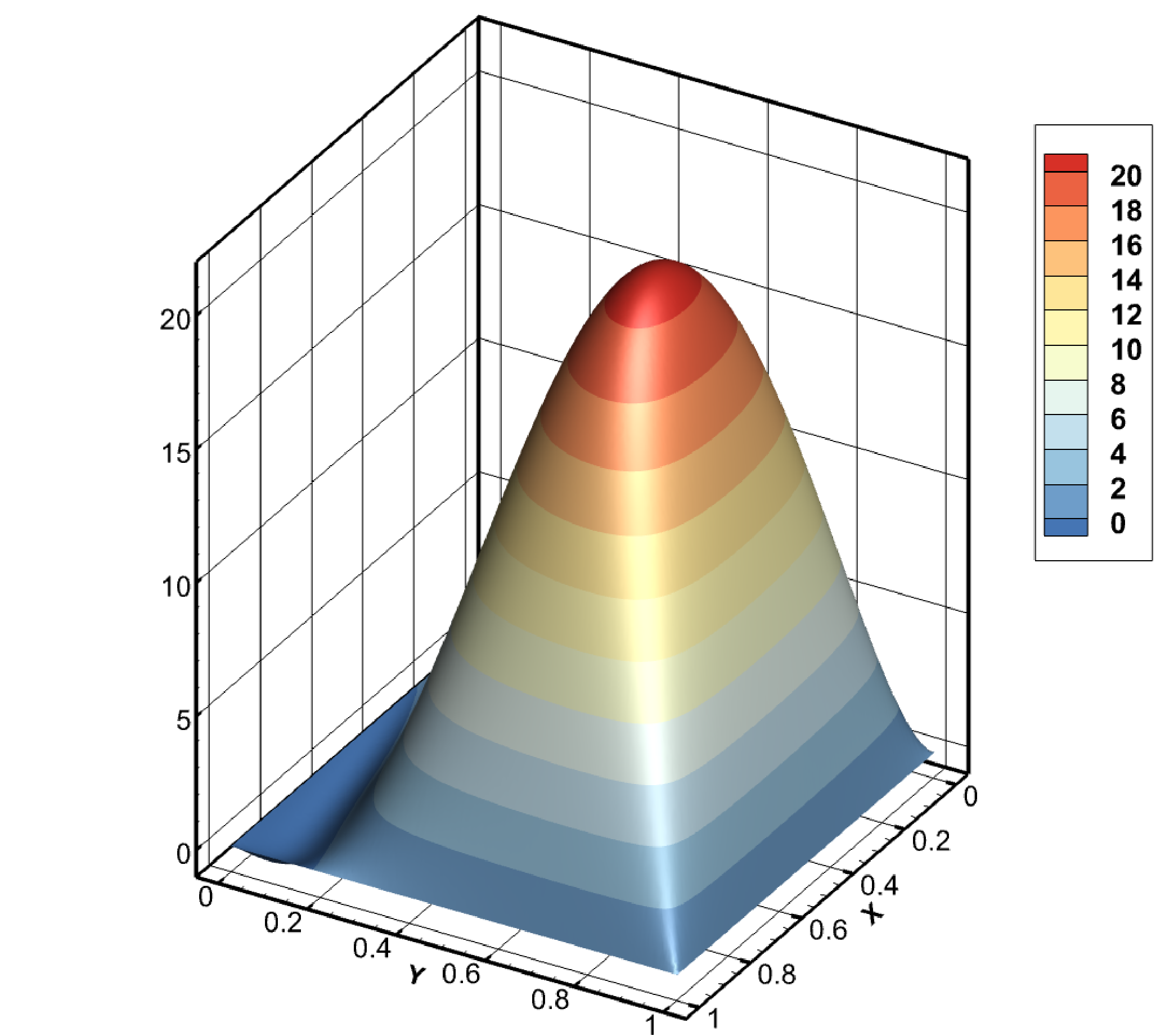}
			(a)
		\end{minipage}
		\hfill
		\begin{minipage}[c]{0.24\textwidth}
			\centering
			\includegraphics[width=\linewidth]{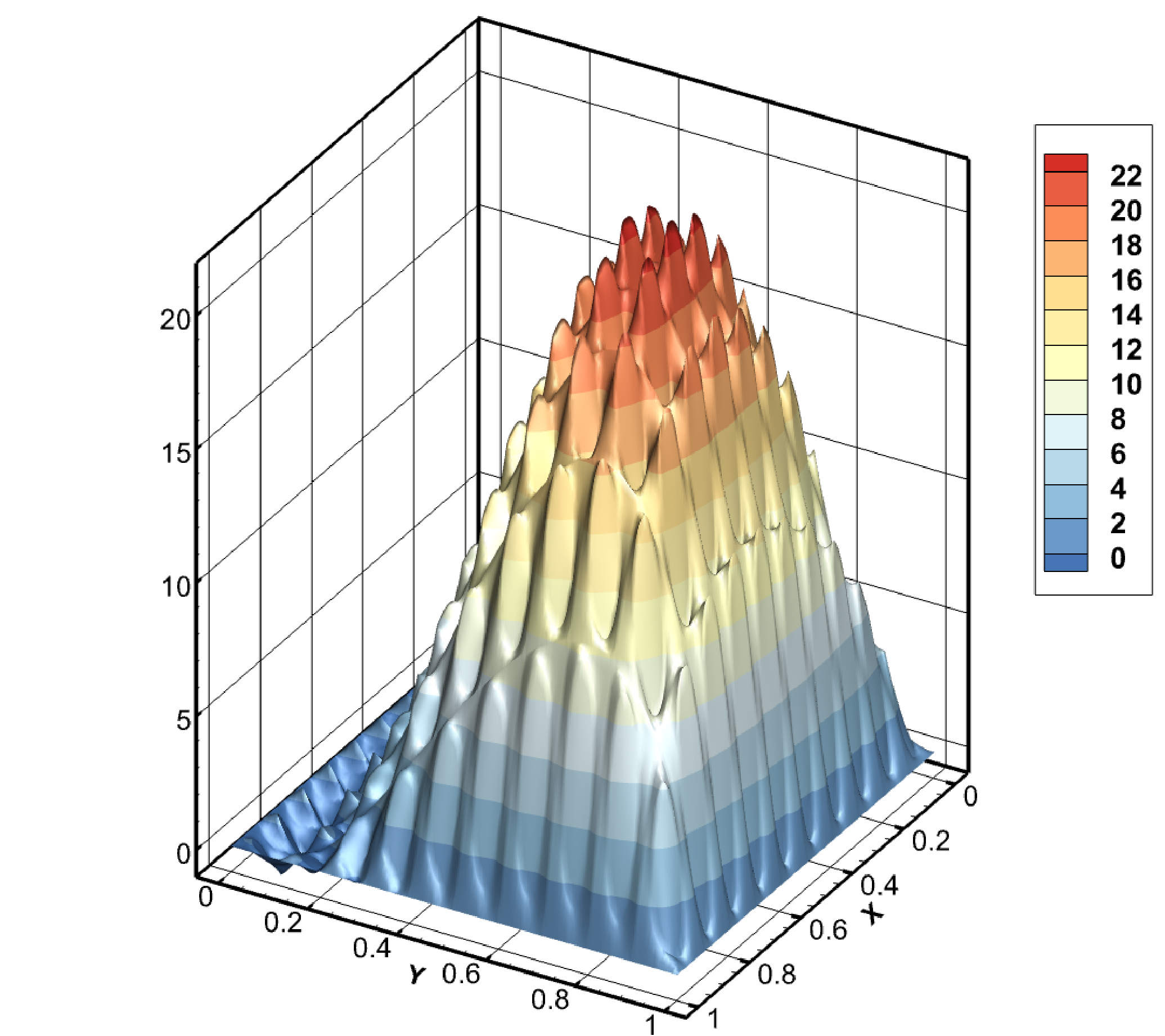}
			(b)
		\end{minipage}
		\hfill
		\begin{minipage}[c]{0.24\textwidth}
			\centering
			\includegraphics[width=\linewidth]{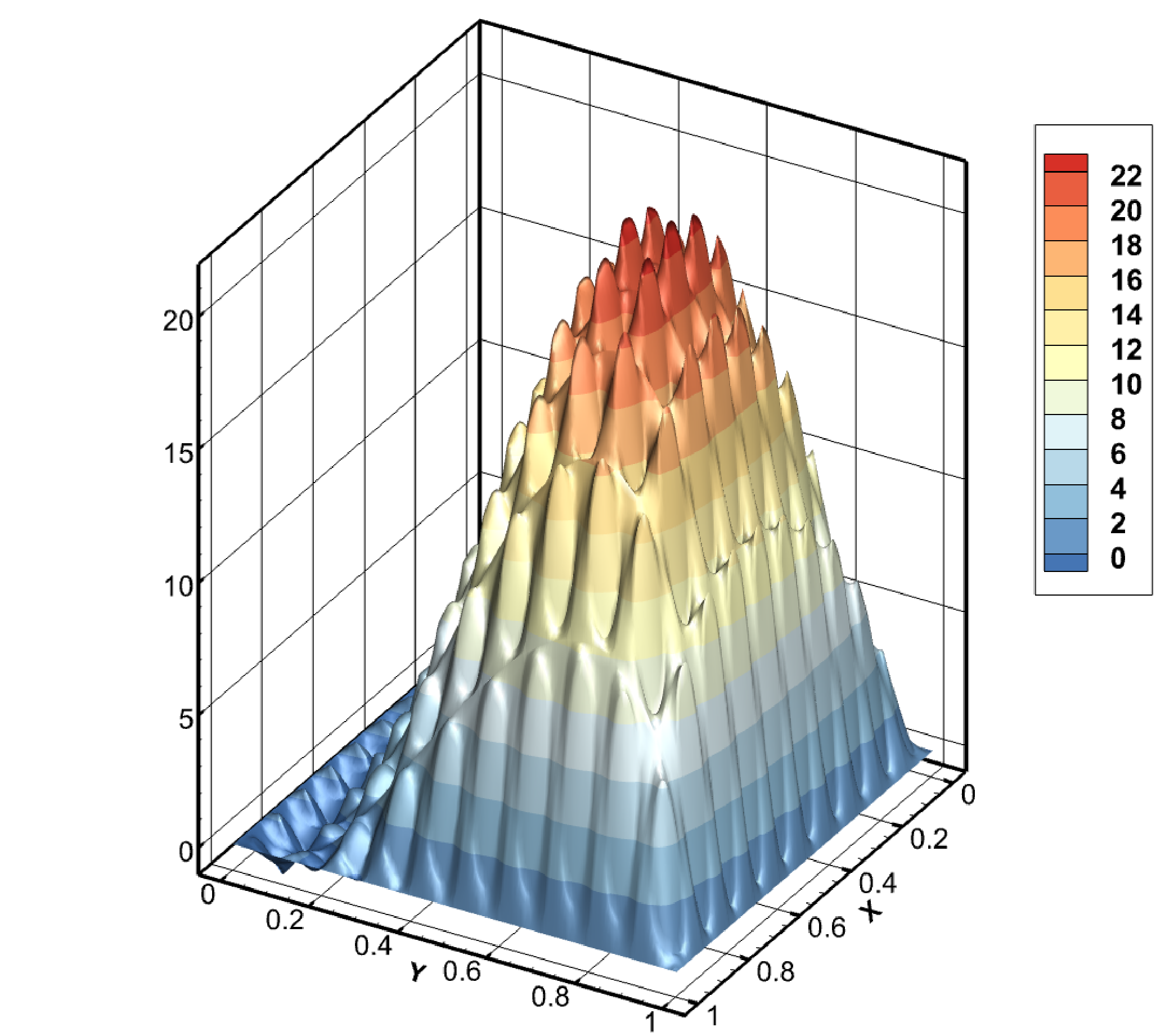}
			(c)
		\end{minipage}
		\hfill
		\begin{minipage}[c]{0.24\textwidth}
			\centering
			\includegraphics[width=\linewidth]{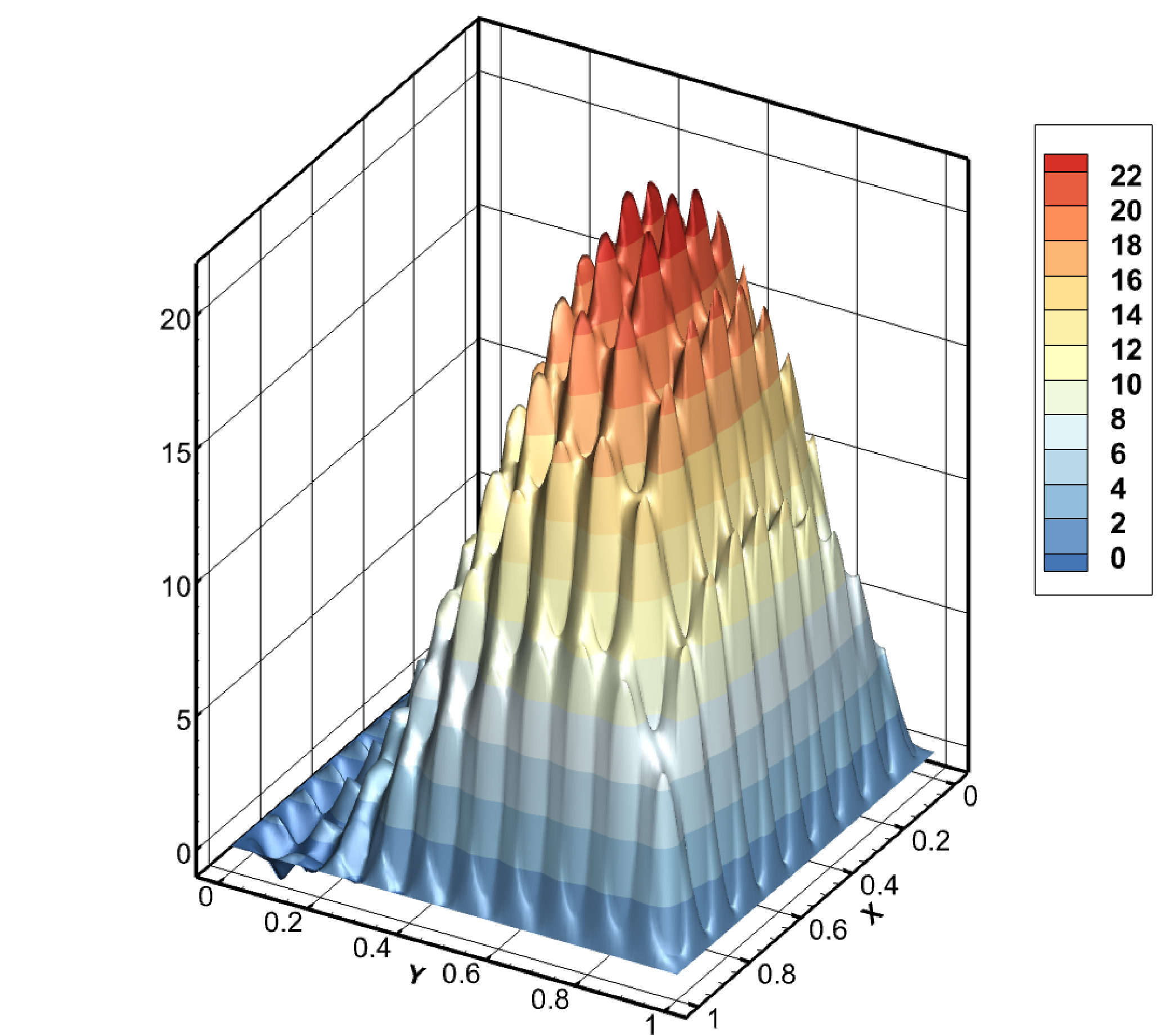}
			(d)
		\end{minipage}
		\caption{Second displacement field component under scale coupling: (a) $u_2^{(0)}$; (b) $u_2^{(1,\epsilon)}$; (c) $u_2^{(2,\epsilon)}$; (d) $u_2^{\epsilon}$.}\label{f9}
	\end{figure}

	After numerical calculations, Figs.\hspace{1mm}\ref{f2}-\ref{f9} display the simulative results for solutions $T^{(0)}$, $T^{(1,\epsilon)}$, $T^{(2,\epsilon)}$, $T^\epsilon$, and  $c^{(0)}$, $c^{(1,\epsilon)}$, $c^{(2,\epsilon)}$, $c^\epsilon$, and $u_1^{(0)}$, $u_1^{(1,\epsilon)}$, $u_1^{(2,\epsilon)}$, $u_1^\epsilon$, $u_2^{(0)}$, $u_2^{(1,\epsilon)}$, $u_2^{(2,\epsilon)}$, $u_2^\epsilon$, respectively. Among them, Figs.\hspace{1mm}\ref{f2}-\ref{f4} correspond to scale-separated material coefficients, while Figs.\hspace{1mm}\ref{f6}-\ref{f9} demonstrate the case without scale-separation. Additionally, the error values associated with the temperature increment, moisture, and displacement fields are provided in Tables\hspace{1mm}\ref{t4} and \ref{t5}, computed using the $L^2$ norm and $H^1$ semi-norm.
	\begin{table}[!htb]
		\centering
		\caption{The relative errors for scale-separated material parameters.}
		\label{t4}
		\begin{tabular}{cccccc}
			\hline
			\multicolumn{6}{c}{Temperature increment field} \\
			\hline
			$TerrorL^20$ & $TerrorL^21$ & $TerrorL^22$ & $TerrorH^10$ & $TerrorH^11$ & $TerrorH^12$ \\
			0.09888 & 0.09560 & 0.01302 & 0.86663 & 0.82266 & 0.08441 \\
			\hline
			\multicolumn{6}{c}{Moisture field} \\
			\hline
			$cerrorL^20$ & $cerrorL^21$ & $cerrorL^22$ & $cerrorH^10$ & $cerrorH^11$ & $cerrorH^12$ \\
			0.05672 & 0.04998 & 0.01137 & 0.69682 & 0.59505 & 0.08720 \\
			\hline
			\multicolumn{6}{c}{Displacement field} \\
			\hline
			$\bm{u}errorL^20$ & $\bm{u}errorL^21$ & $\bm{u}errorL^22$ & $\bm{u}errorH^10$ & $\bm{u}errorH^11$ & $\bm{u}errorH^12$ \\
			0.12365 & 0.04082 & 0.03929 & 0.77757 & 0.19190 & 0.19162 \\
			\hline
		\end{tabular}
	\end{table}
	\begin{table}[!htb]
		\centering
		\caption{The relative errors for scale-coupled material parameters.}
		\label{t5}
		\begin{tabular}{cccccc}
			\hline
			\multicolumn{6}{c}{Temperature increment field} \\
			\hline
			$TerrorL^20$ & $TerrorL^21$ & $TerrorL^22$ & $TerrorH^10$ & $TerrorH^11$ & $TerrorH^12$ \\
			0.09455 & 0.09115 & 0.01218 & 0.85758 & 0.81063 & 0.06621 \\
			\hline
			\multicolumn{6}{c}{Moisture field} \\
			\hline
			$cerrorL^20$ & $cerrorL^21$ & $cerrorL^22$ & $cerrorH^10$ & $cerrorH^11$ & $cerrorH^12$ \\
			0.05473 & 0.04778 & 0.01134 & 0.68460 & 0.57755 & 0.07891 \\
			\hline
			\multicolumn{6}{c}{Displacement field} \\
			\hline
			$\bm{u}errorL^20$ & $\bm{u}errorL^21$ & $\bm{u}errorL^22$ & $\bm{u}errorH^10$ & $\bm{u}errorH^11$ & $\bm{u}errorH^12$ \\
			0.07714 & 0.04235 & 0.04010 & 0.72242 & 0.21628 & 0.20540 \\
			\hline
		\end{tabular}
	\end{table}
	
	From Figs.\hspace{1mm}\ref{f2}-\ref{f9}, Table\hspace{1mm}\ref{t4}, and Table\hspace{1mm}\ref{t5}, it is evident that the HOMS solutions demonstrate significantly higher accuracy than the homogenized and LOMS solutions for the temperature increment and moisture fields, particularly in the $H^1$ semi-norm. For the displacement fields, the difference between the LOMS and HOMS solutions is not obvious. Indeed, the purpose of the HOMS method is to capture the micro-scale oscillating information arising from material heterogeneities. However, in this example, the displacement fields exhibit few highly oscillatory fluctuations. As a result, the accuracy of the LOMS solutions closely matches that of the HOMS solutions.
	
	\subsection{Example 2: 3D quasi-periodic composite structure}
	\label{sec:52}
\begin{figure}[!htb]
		\centering
		\begin{minipage}[c]{0.4\textwidth}
			\centering
			\includegraphics[width=50mm]{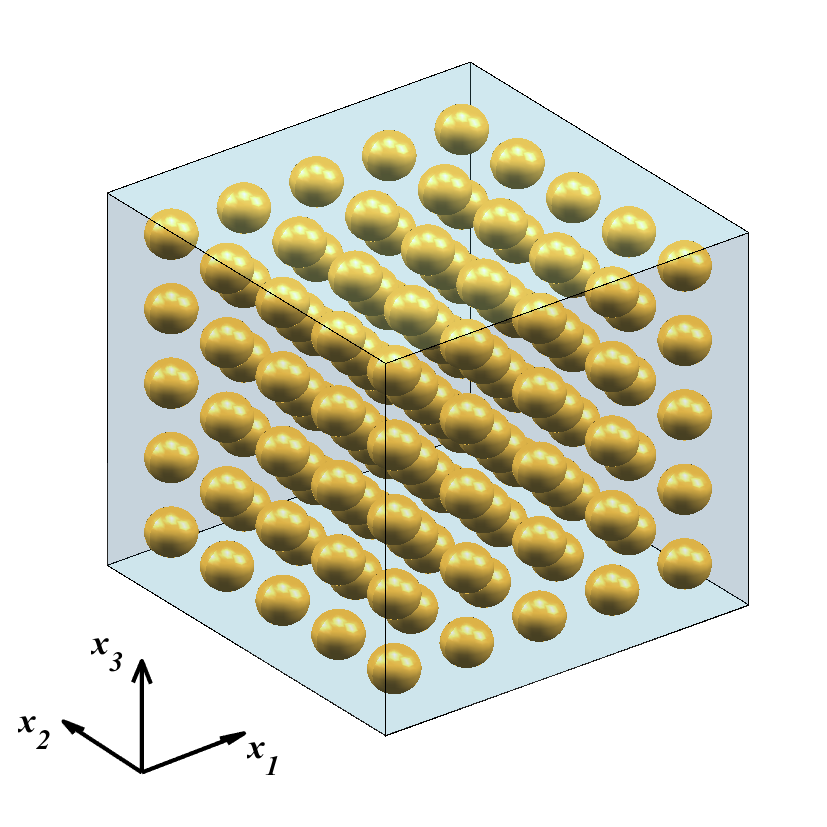}\\
			(a)
		\end{minipage}
		\begin{minipage}[c]{0.4\textwidth}
			\centering
			\includegraphics[width=50mm]{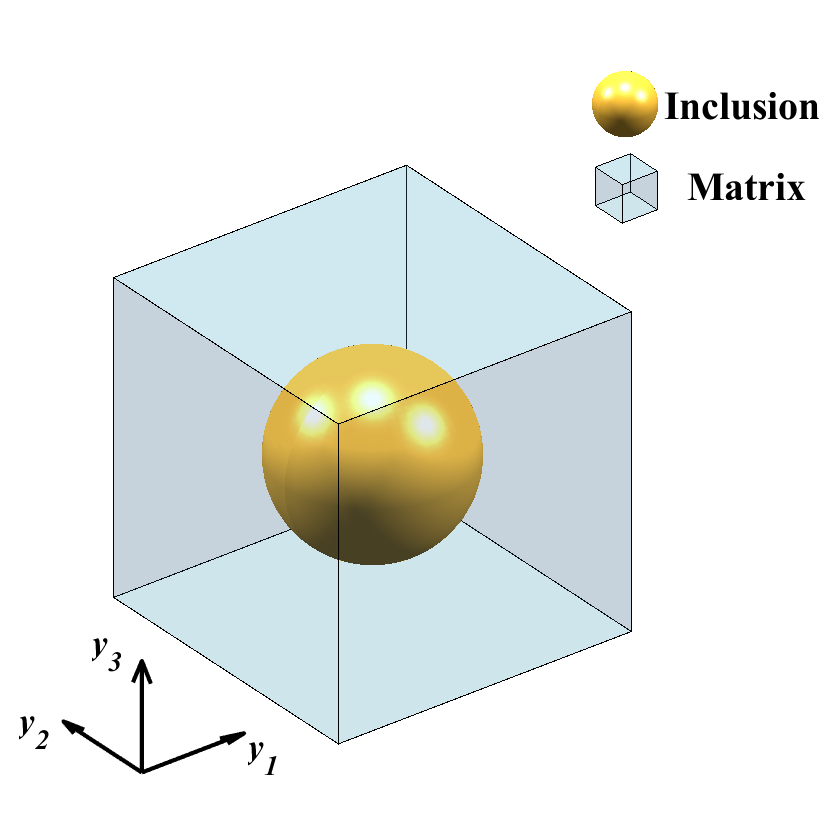}\\
			(b)
		\end{minipage}
		\caption{The schematic of composite structure: (a) composite structure $\Omega$; (b) PUC $Y$.}\label{f10:3D}
	\end{figure}

	In this example, a 3D quasi-periodic composite structure is investigated. The detailed macroscopic structure $\Omega$ and microscopic unit cell $Y$ are shown in Fig.\hspace{1mm}\ref{f10:3D}, where $\Omega= (x_1,x_2,x_3) = [0,1]^3 \mathrm{cm}^3$ and small periodic parameter $\epsilon=1/5$.
	
	This example also investigates both scale-separated and non-scale-separated cases. The scale-separated case employs the material parameters in Table\hspace{1mm}\ref{t1} with $\psi(\mathbf{x}) = 1 + x_3$, while the non-scale-separated case uses the material parameters in Table\hspace{1mm}\ref{t2} with $\psi(\mathbf{x}) = x_3$.
	
	Furthermore, the heat source, moisture source and body forces are defined by $h=500 \mathrm{J/(cm^{3}\cdot s)}$, $m=500 \mu\mathrm{ g/(cm^{3}\cdot s)}$ and $(f_1,f_2,f_3)=(0,0,-5000) \mathrm{N/cm^3}$, respectively. The boundary conditions on $\partial\Omega$ are prescribed as $\overline{T}=273.15 \mathrm{K}$, $\overline{c}=0 \mathrm{g/cm^3}$ and $\overline{\bm{u}}=0\mathrm{cm}$.
	
	The comparison of computational cost in Table\hspace{1mm}\ref{t6} reveals that, consistent with the findings in Example\hspace{1mm}1, the HOMS method substantially reduces computational cost compared to the precise FEM.
	\begin{table}[!htb]{\caption{Comparison of computational cost.}\label{t6}}
		\centering
		\begin{tabular}{cccc}
			\hline
			& Cell equations & Homogenized equations & Multi-scale equations\\
			\hline
			FEM nodes & 1212 & 4096 & 79021\\
			FEM elements & 6007 & 20250 & 480352 \\
			\hline
				& & HOMS method & precise FEM\\
				\hline
				\multicolumn{2}{c}{Computing time for case 1} & 65.244s & 155.019s \\
				\multicolumn{2}{c}{Computing time for case 2} & 2422.809s & 133.781s \\
				\hline
			\end{tabular}
		\end{table}
		
\begin{figure}[!htb]
			\centering
			\begin{minipage}[c]{0.24\textwidth}
				\centering
				\includegraphics[width=\linewidth]{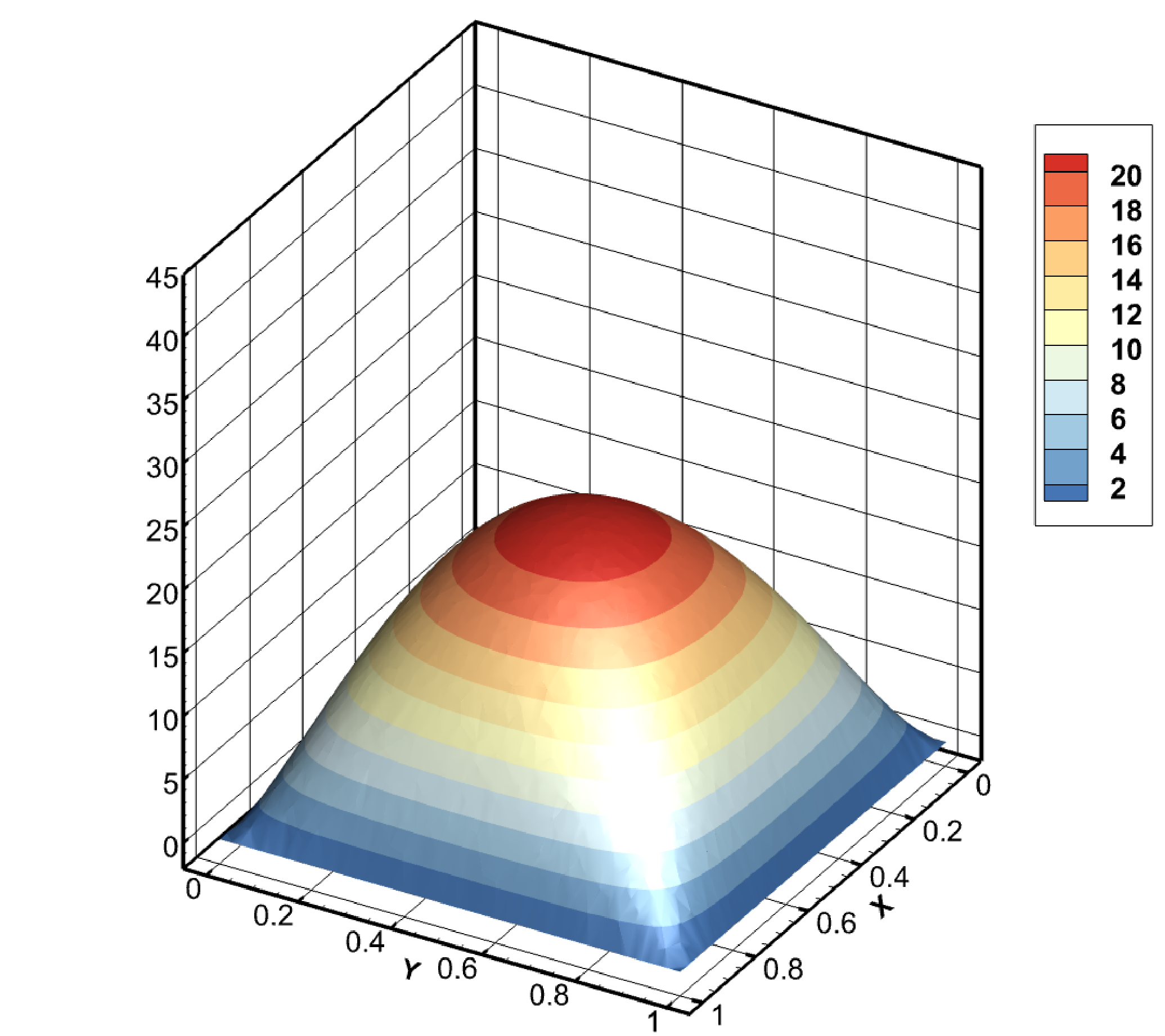}
				(a)
			\end{minipage}
			\hfill
			\begin{minipage}[c]{0.24\textwidth}
				\centering
				\includegraphics[width=\linewidth]{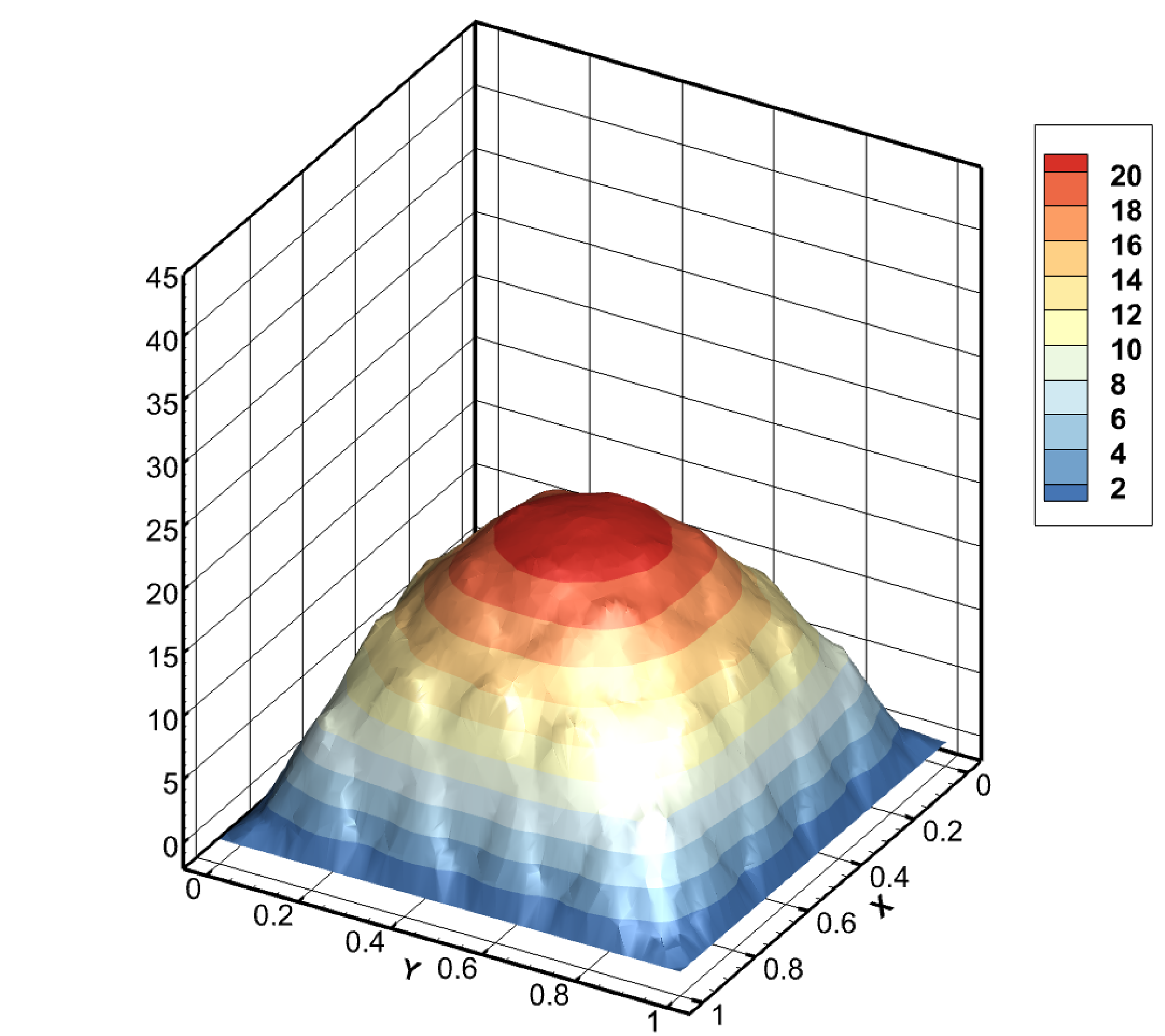}
				(b)
			\end{minipage}
			\hfill
			\begin{minipage}[c]{0.24\textwidth}
				\centering
				\includegraphics[width=\linewidth]{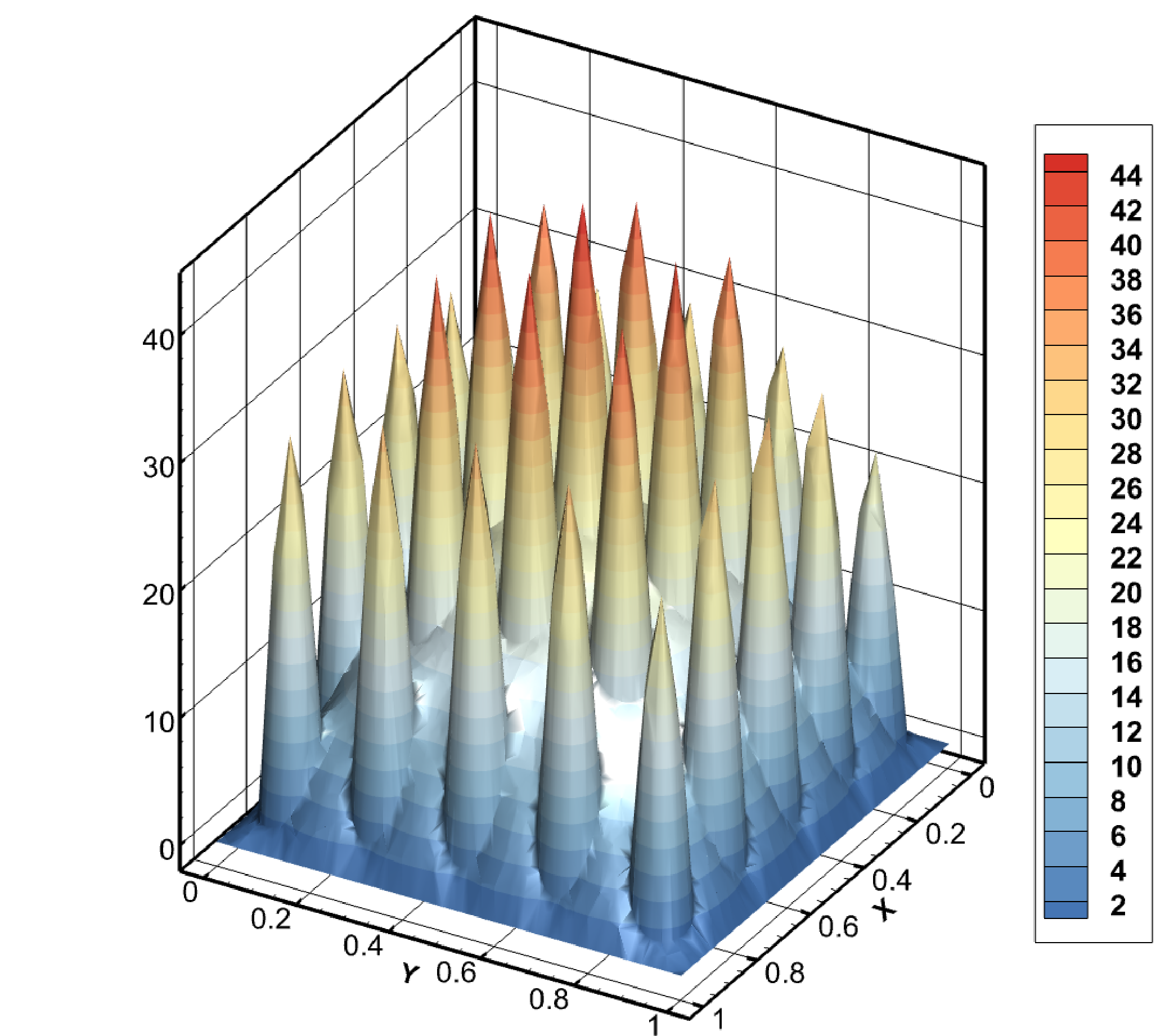}
				(c)
			\end{minipage}
			\hfill
			\begin{minipage}[c]{0.24\textwidth}
				\centering
				\includegraphics[width=\linewidth]{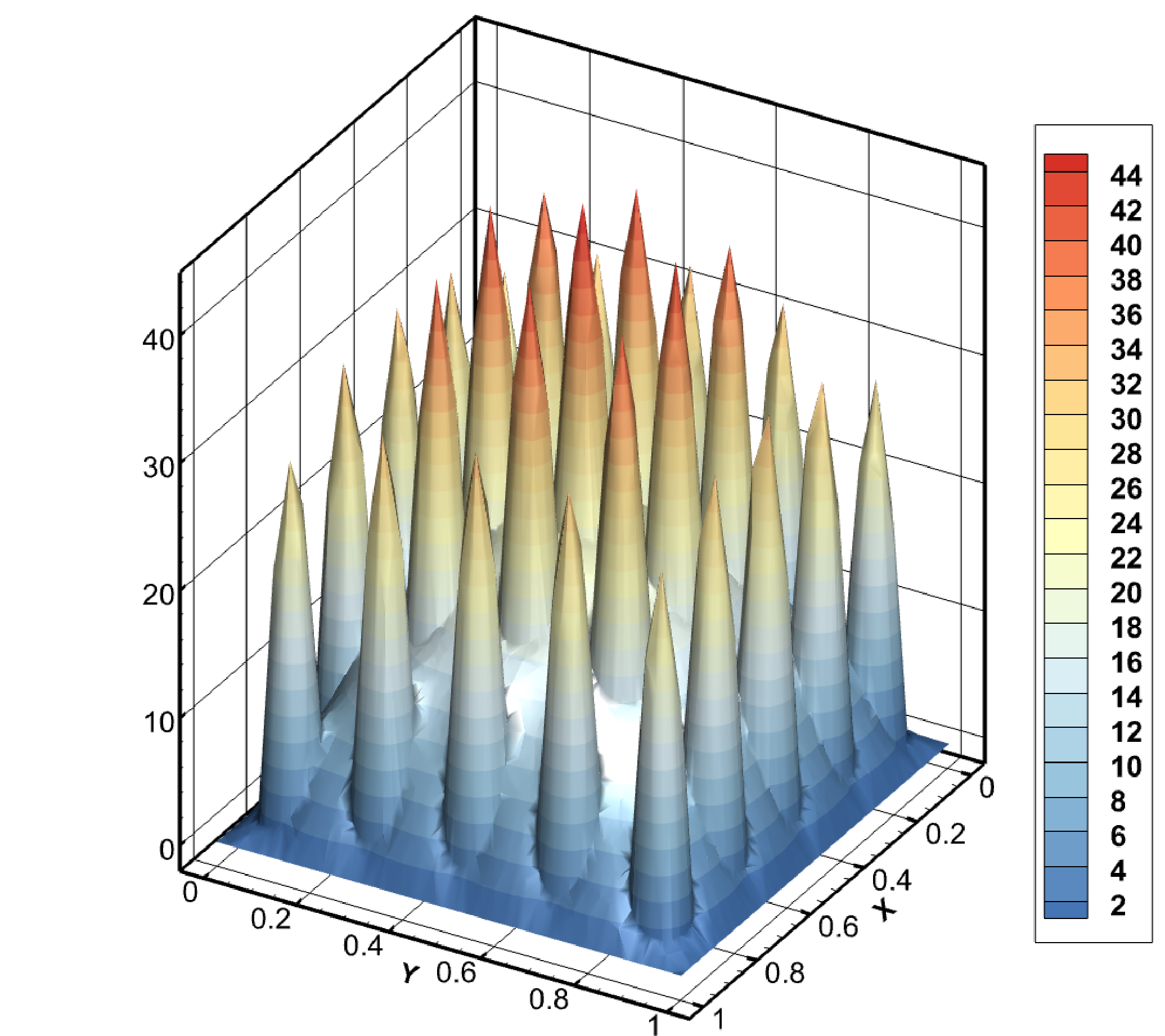}
				(d)
			\end{minipage}
			\caption{Temperature increment field in $x_3=0.3 \mathrm{cm}$ with scale separation: (a) $T^{(0)}$; (b) $T^{(1,\epsilon)}$; (c) $T^{(2,\epsilon)}$; (d) $T^{\epsilon}$.}\label{f11}
		\end{figure}

\begin{figure}[!htb]
			\centering
			\begin{minipage}[c]{0.24\textwidth}
				\centering
				\includegraphics[width=\linewidth]{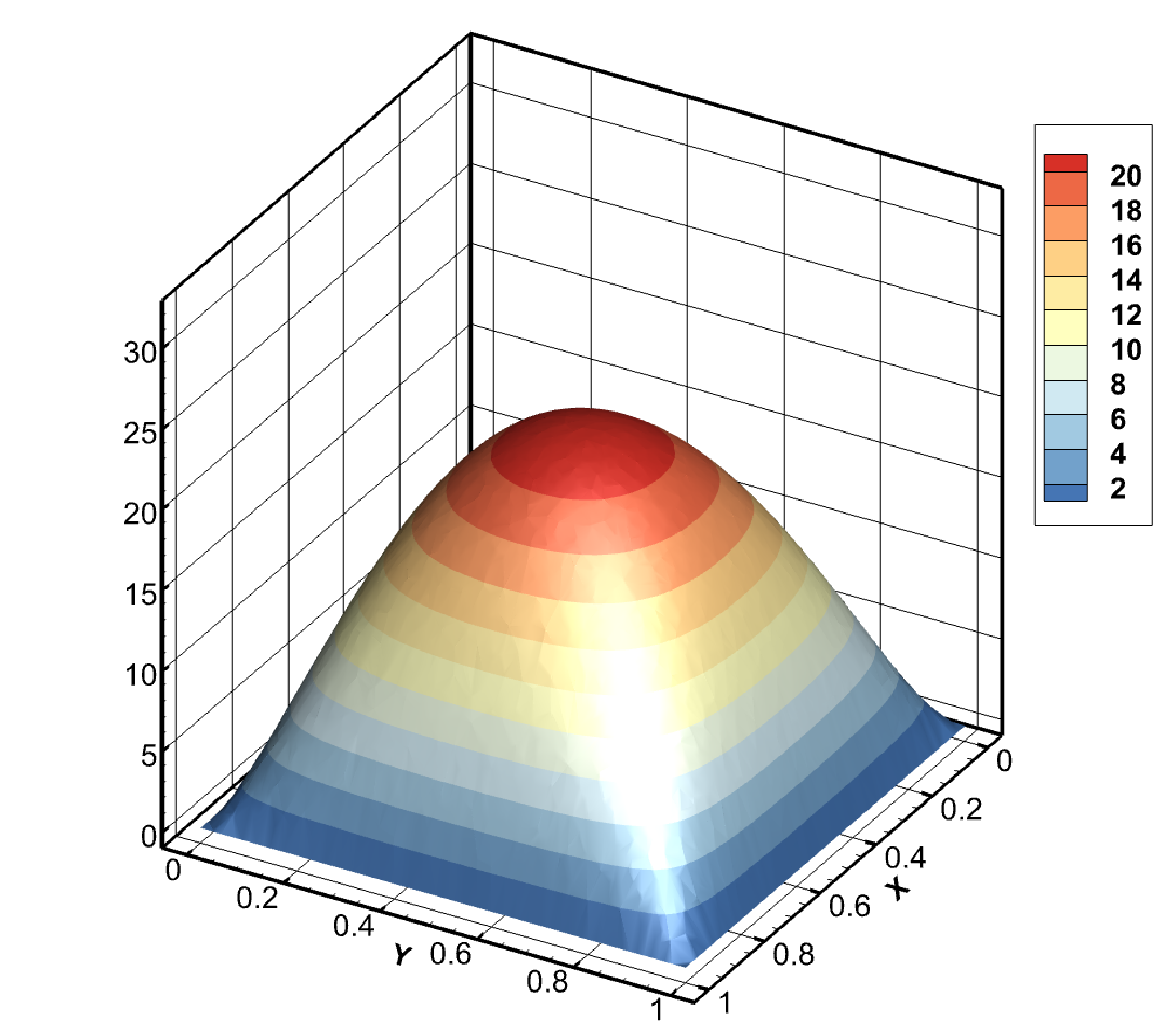}
				(a)
			\end{minipage}
			\hfill
			\begin{minipage}[c]{0.24\textwidth}
				\centering
				\includegraphics[width=\linewidth]{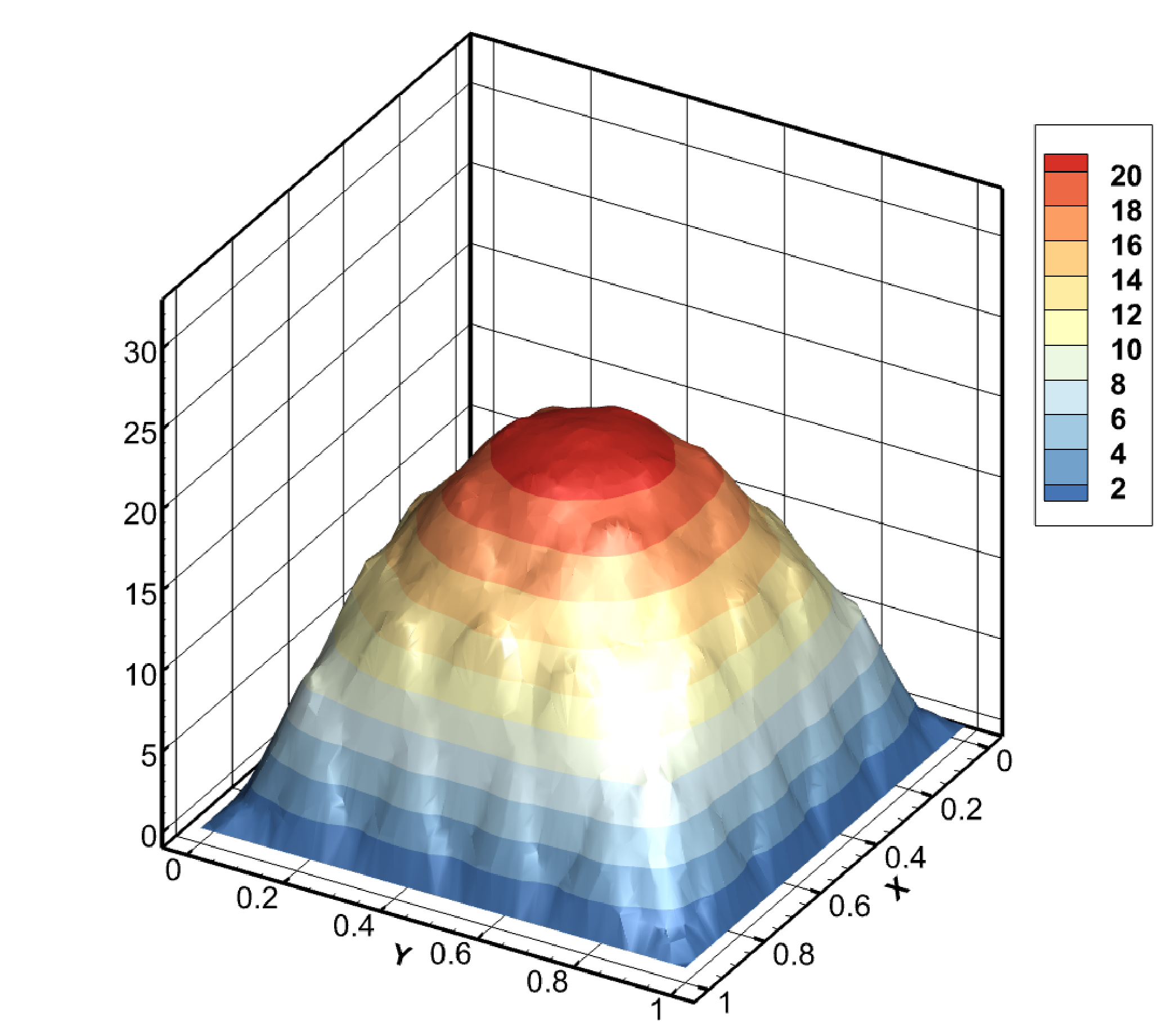}
				(b)
			\end{minipage}
			\hfill
			\begin{minipage}[c]{0.24\textwidth}
				\centering
				\includegraphics[width=\linewidth]{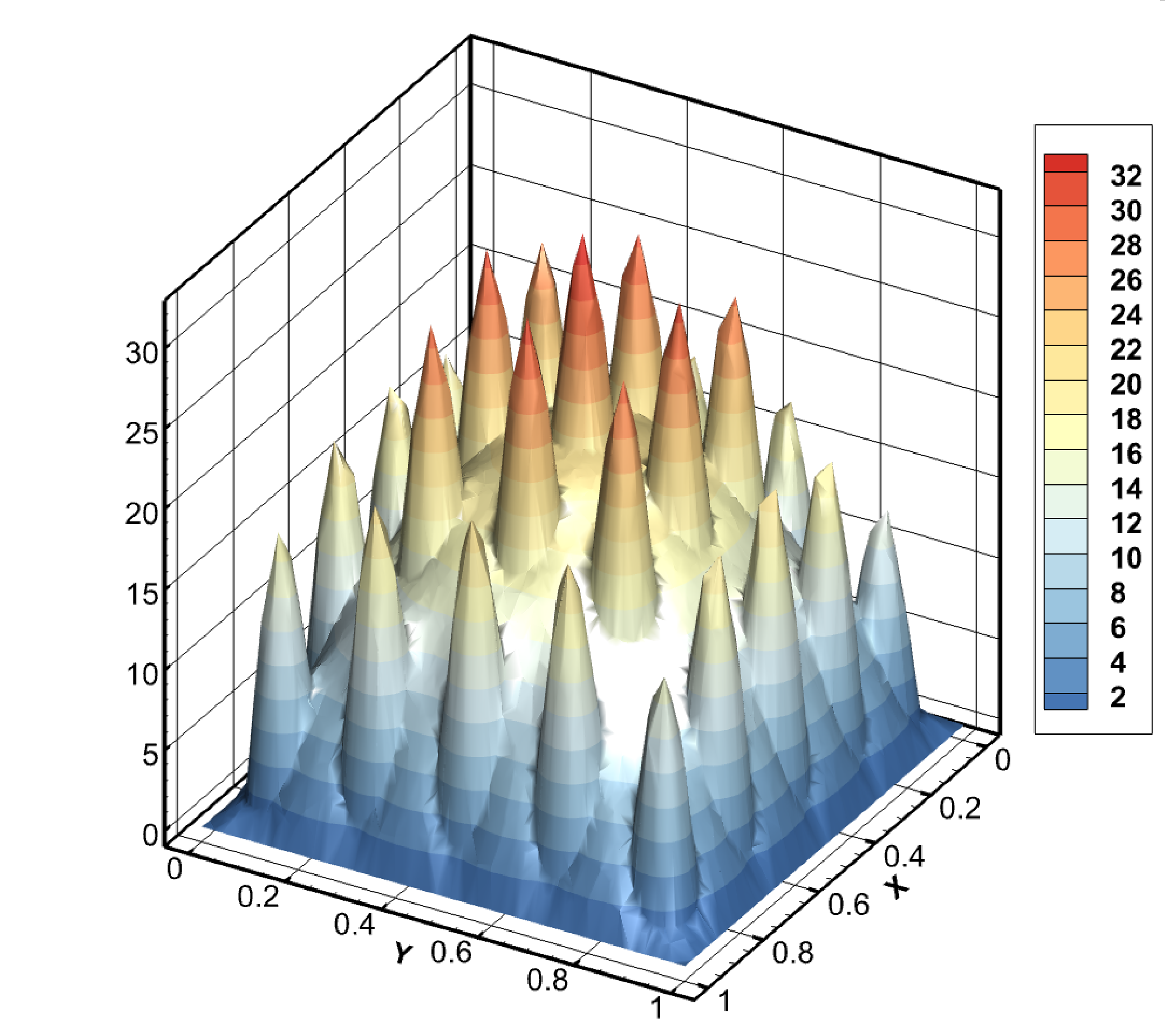}
				(c)
			\end{minipage}
			\hfill
			\begin{minipage}[c]{0.24\textwidth}
				\centering
				\includegraphics[width=\linewidth]{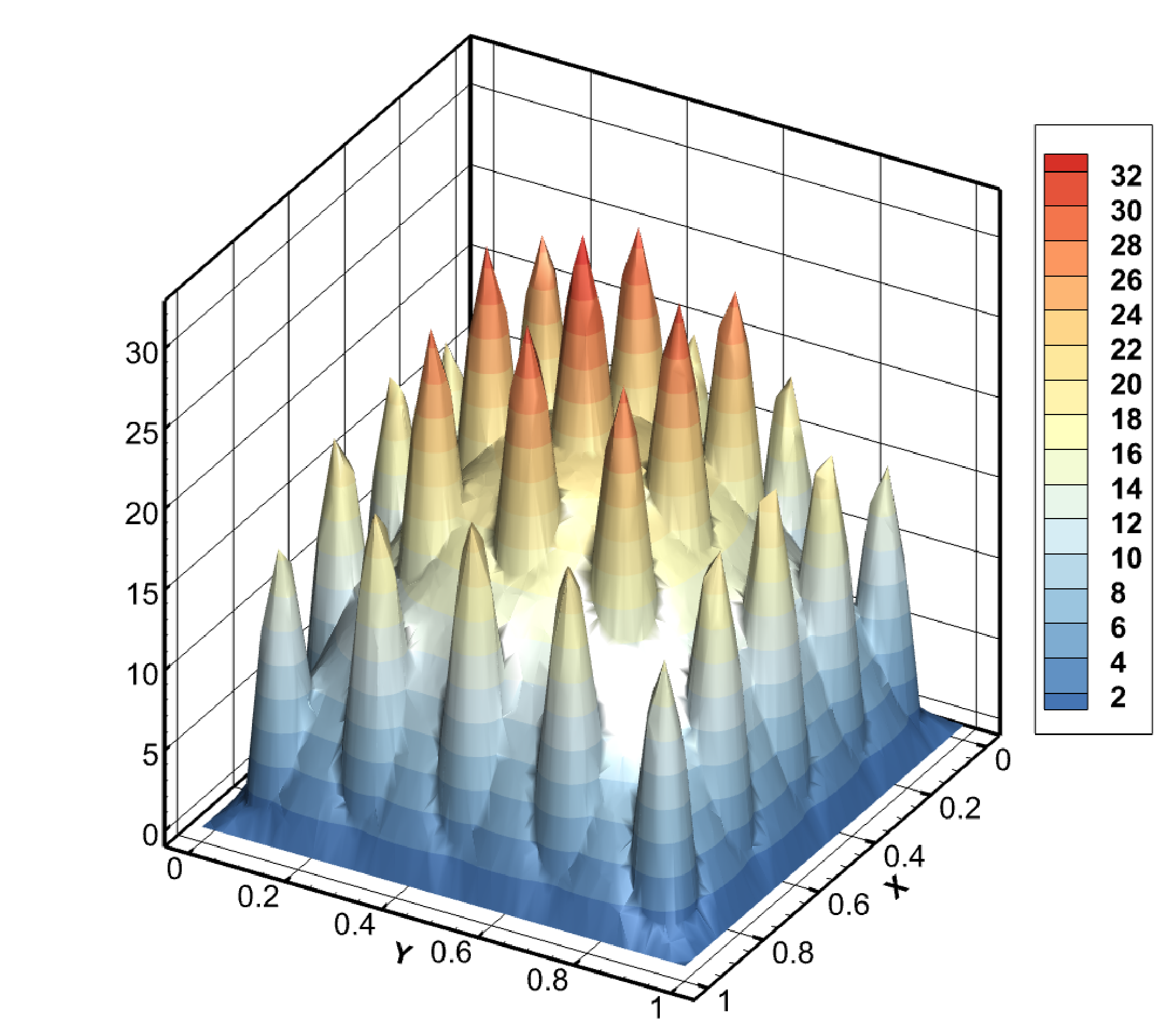}
				(d)
			\end{minipage}
			\caption{Moisture field in $x_3=0.3 \mathrm{cm}$ with scale separation: (a) $c^{(0)}$; (b) $c^{(1,\epsilon)}$; (c) $c^{(2,\epsilon)}$; (d) $c^{\epsilon}$.}\label{f12}
		\end{figure}

\begin{figure}[!htb]
			\centering
			\begin{minipage}[c]{0.24\textwidth}
				\centering
				\includegraphics[width=\linewidth]{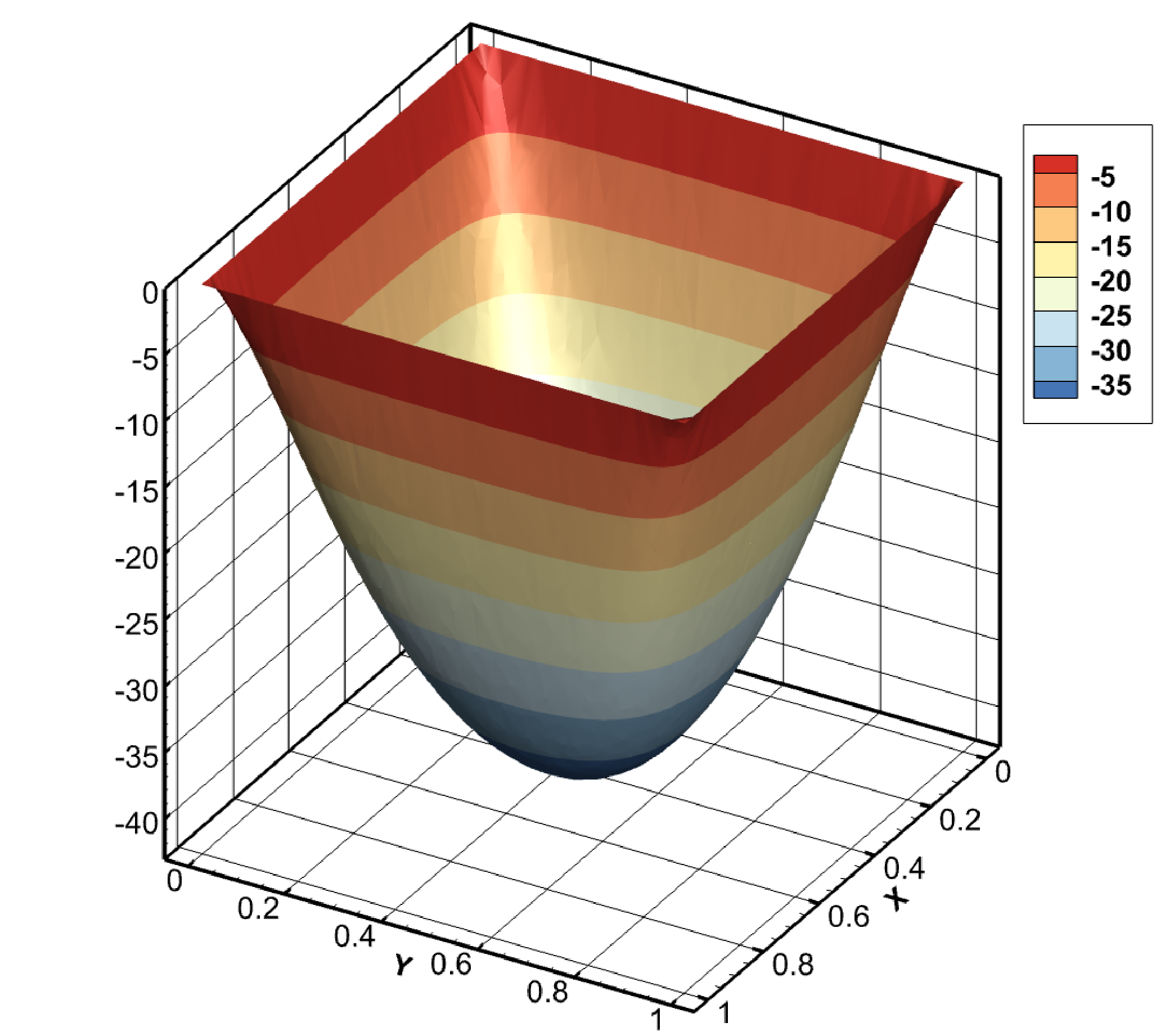}
				(a)
			\end{minipage}
			\hfill
			\begin{minipage}[c]{0.24\textwidth}
				\centering
				\includegraphics[width=\linewidth]{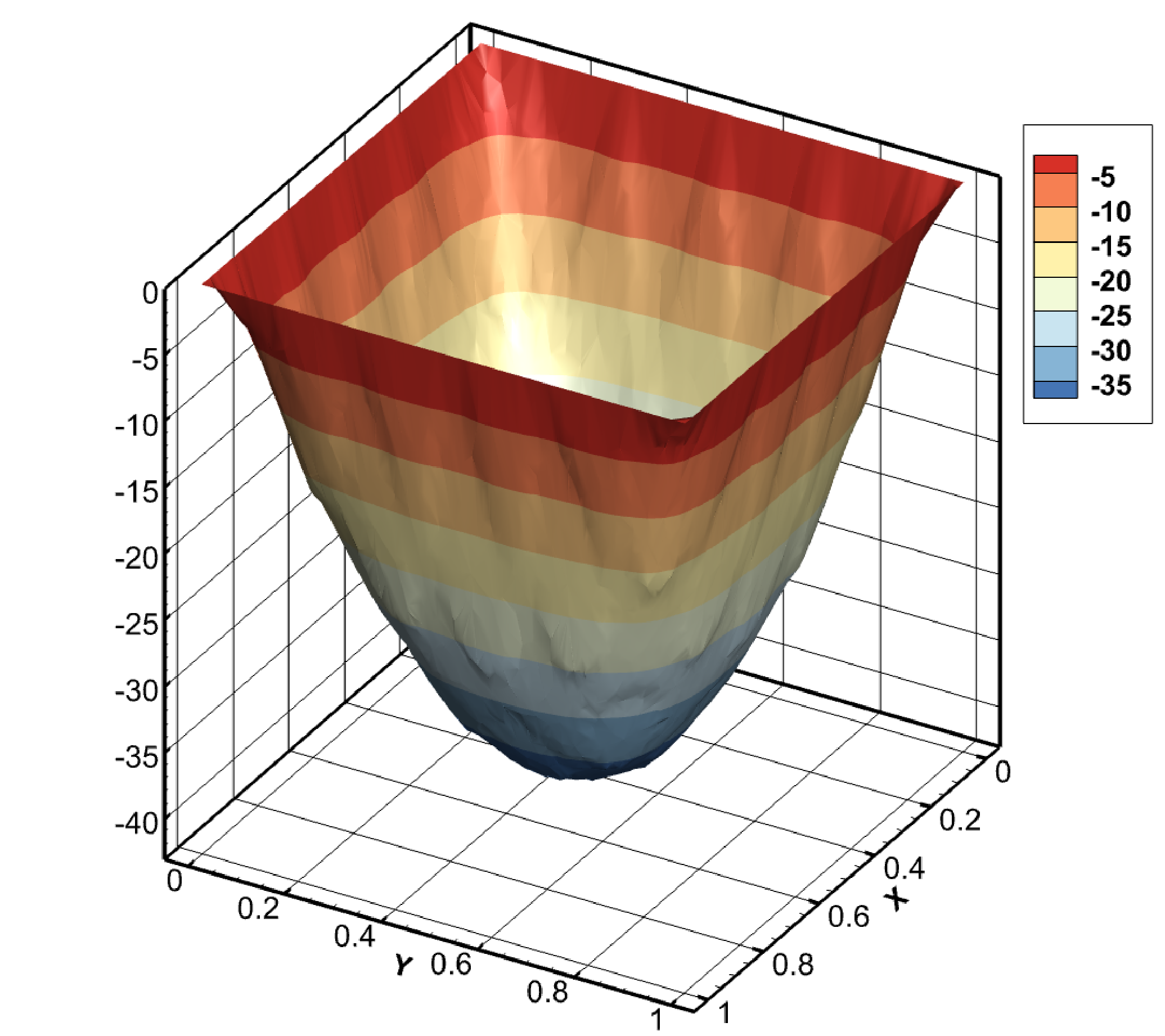}
				(b)
			\end{minipage}
			\hfill
			\begin{minipage}[c]{0.24\textwidth}
				\centering
				\includegraphics[width=\linewidth]{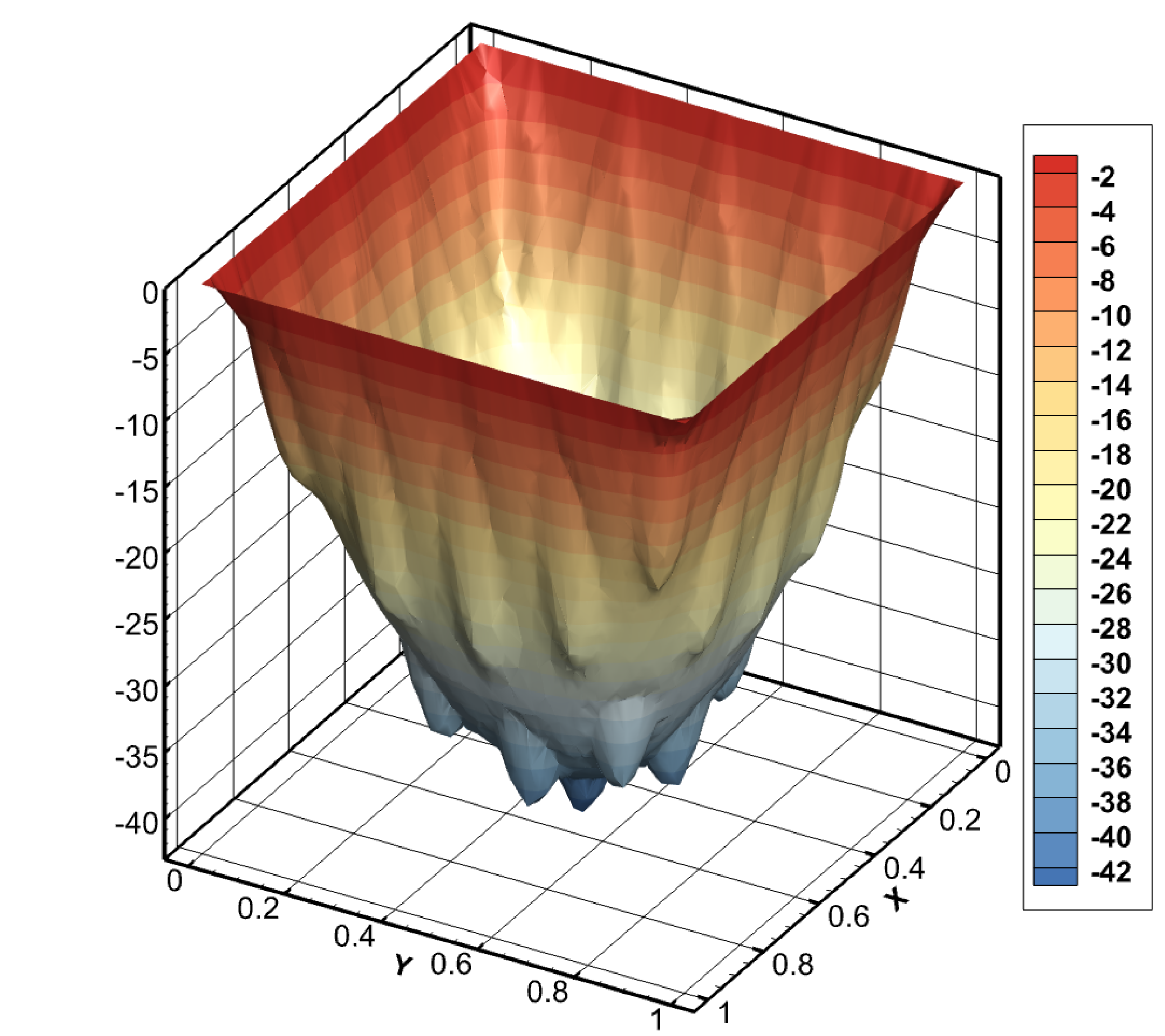}
				(c)
			\end{minipage}
			\hfill
			\begin{minipage}[c]{0.24\textwidth}
				\centering
				\includegraphics[width=\linewidth]{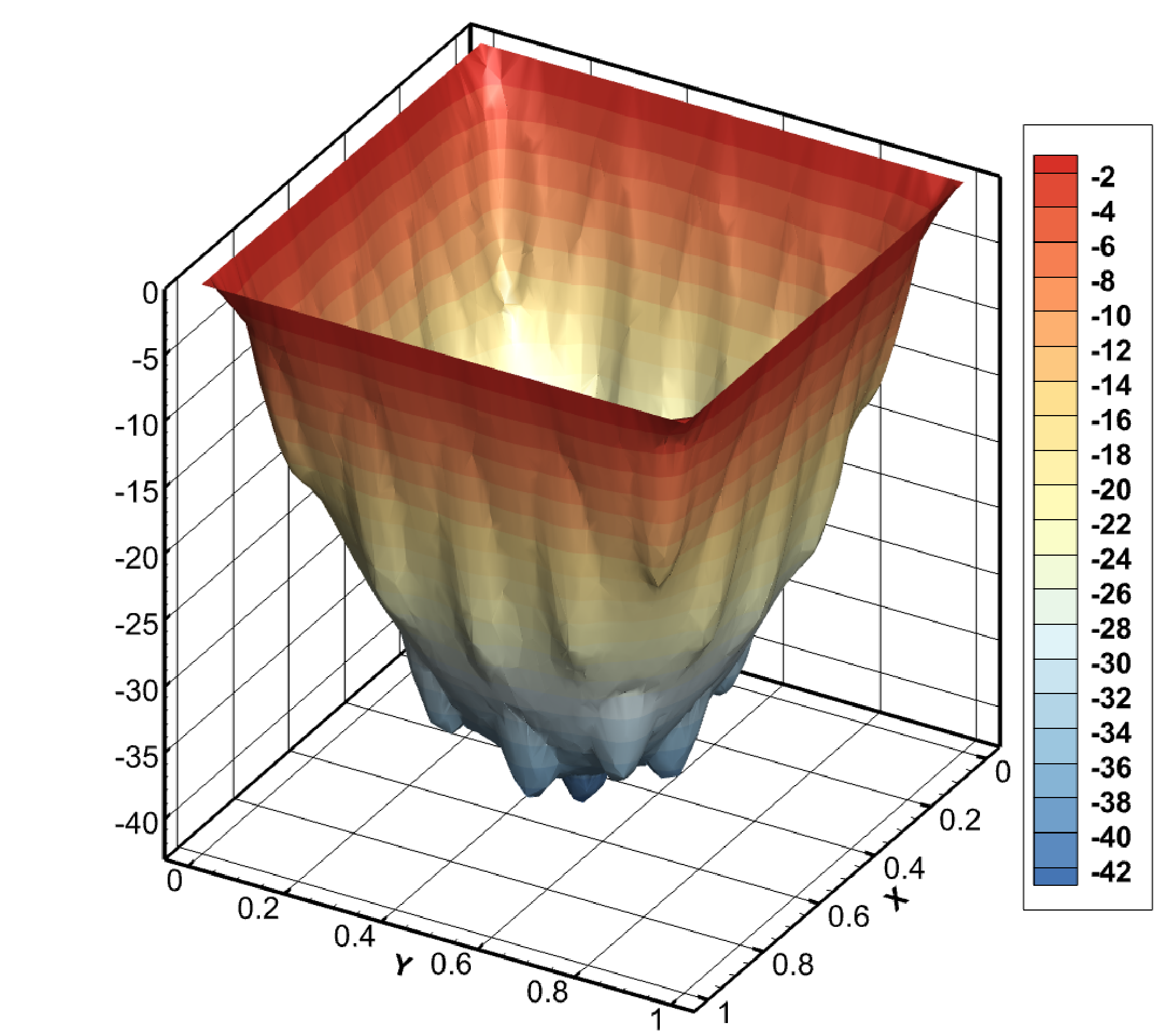}
				(d)
			\end{minipage}
			\caption{Third displacement field component in $x_3=0.3 \mathrm{cm}$ with scale separation: (a) $u_3^{(0)}$; (b) $u_3^{(1,\epsilon)}$; (c) $u_3^{(2,\epsilon)}$; (d) $u_3^{\epsilon}$.}\label{f15}
		\end{figure}

\begin{figure}[!htb]
			\centering
			\begin{minipage}[c]{0.24\textwidth}
				\centering
				\includegraphics[width=\linewidth]{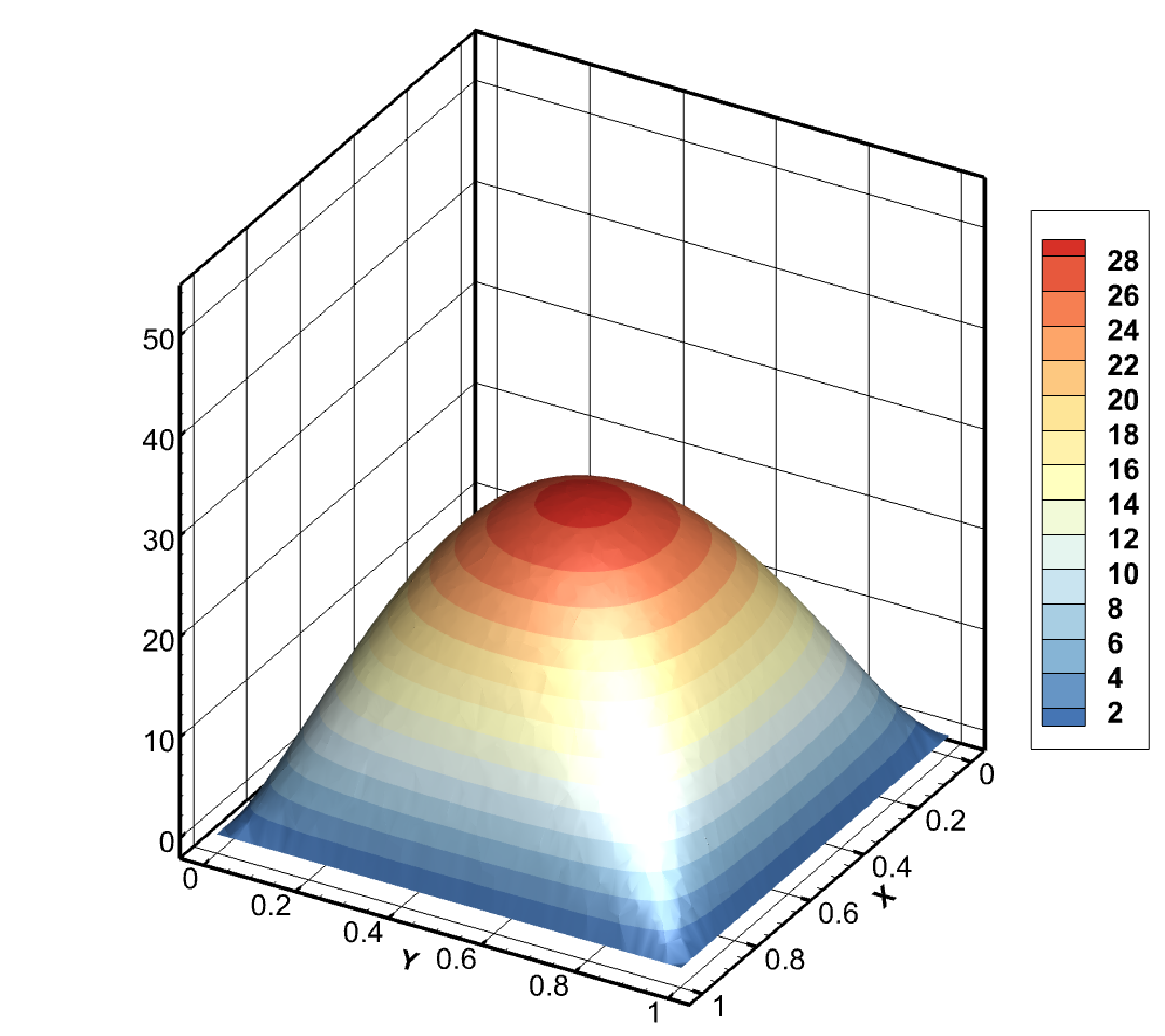}
				(a)
			\end{minipage}
			\hfill
			\begin{minipage}[c]{0.24\textwidth}
				\centering
				\includegraphics[width=\linewidth]{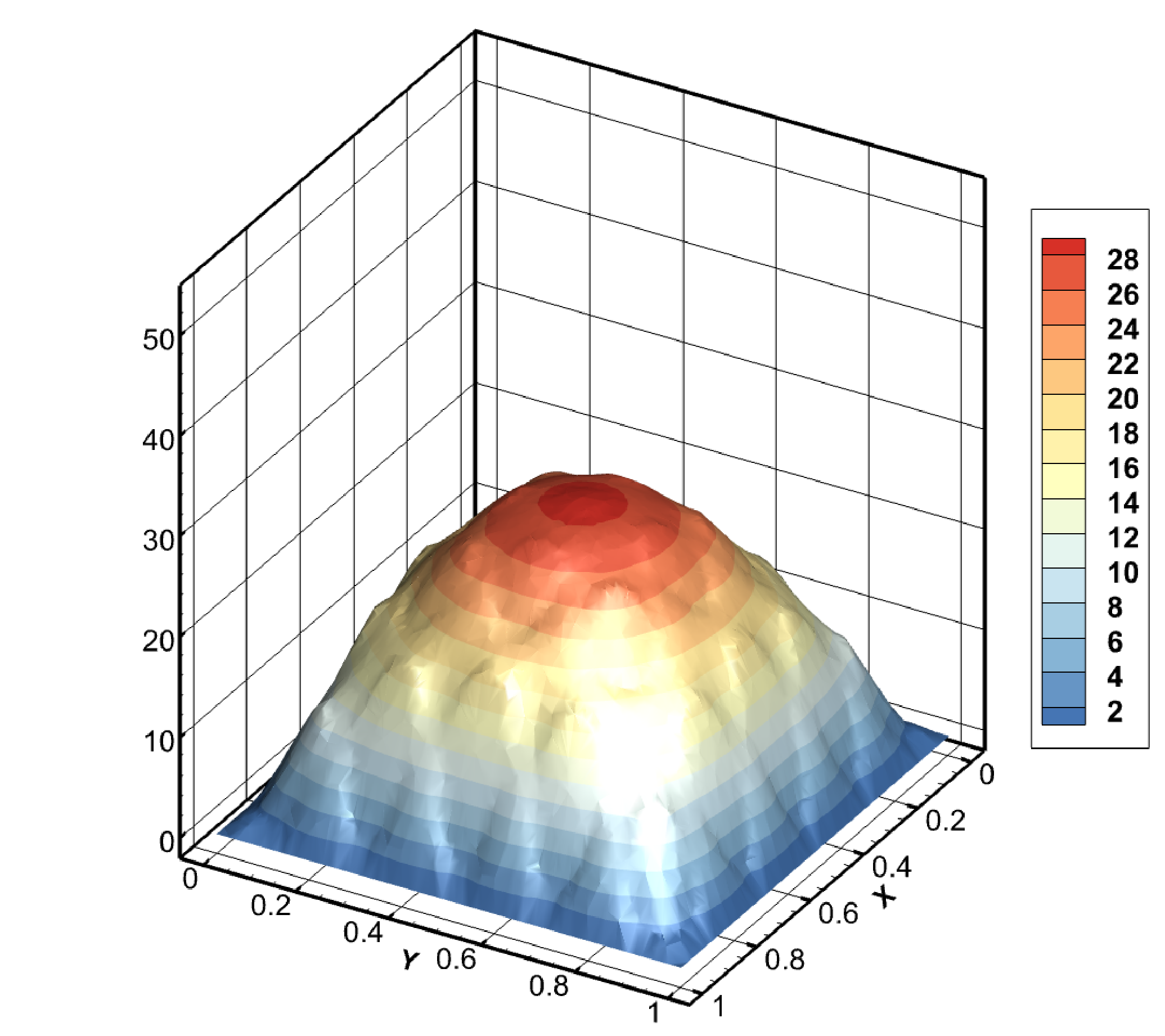}
				(b)
			\end{minipage}
			\hfill
			\begin{minipage}[c]{0.24\textwidth}
				\centering
				\includegraphics[width=\linewidth]{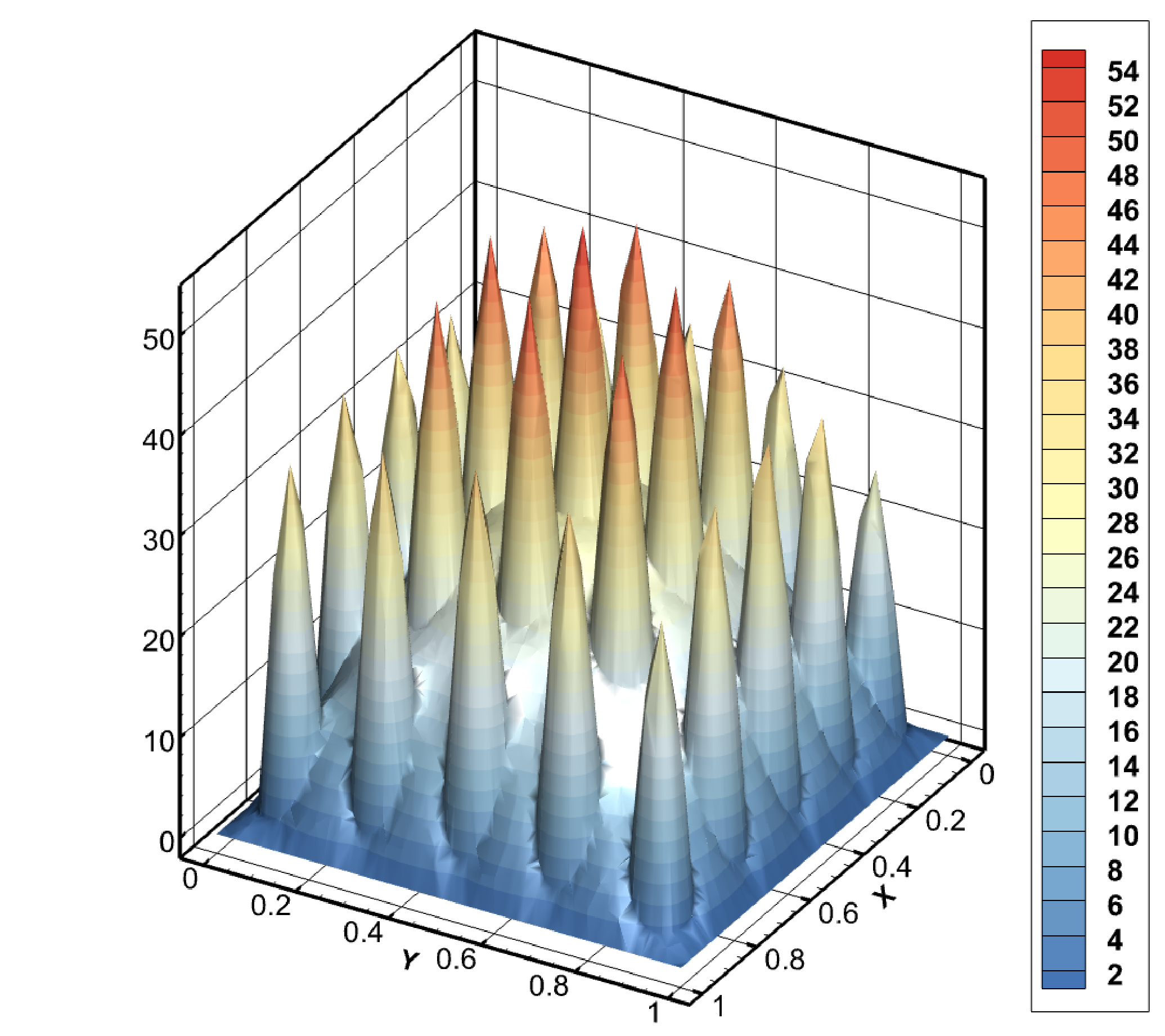}
				(c)
			\end{minipage}
			\hfill
			\begin{minipage}[c]{0.24\textwidth}
				\centering
				\includegraphics[width=\linewidth]{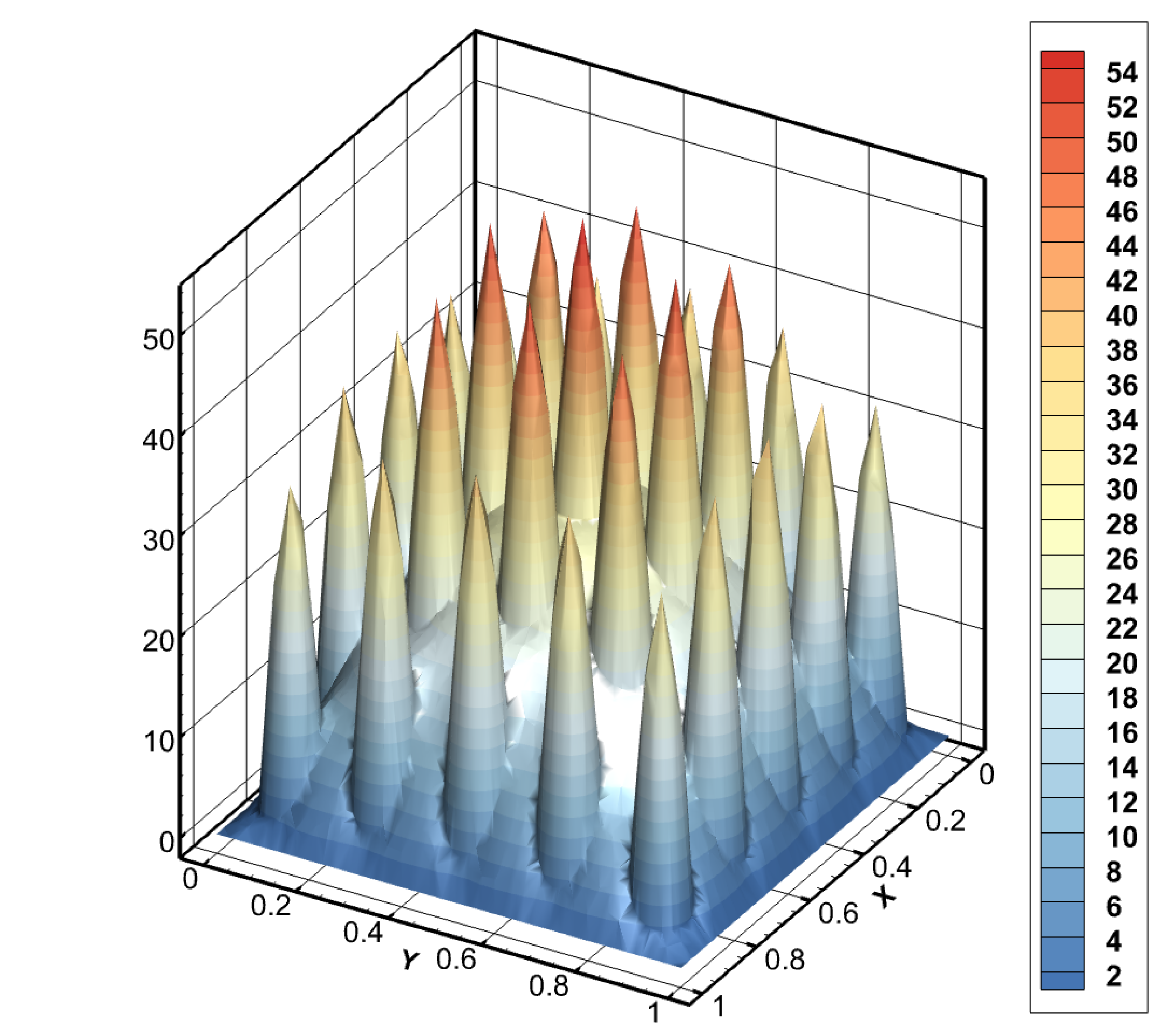}
				(d)
			\end{minipage}
			\caption{Temperature increment field in $x_3=0.3 \mathrm{cm}$ under scale coupling: (a) $T^{(0)}$; (b) $T^{(1,\epsilon)}$; (c) $T^{(2,\epsilon)}$; (d) $T^{\epsilon}$.}\label{f16}
		\end{figure}

\begin{figure}[!htb]
			\centering
			\begin{minipage}[c]{0.24\textwidth}
				\centering
				\includegraphics[width=\linewidth]{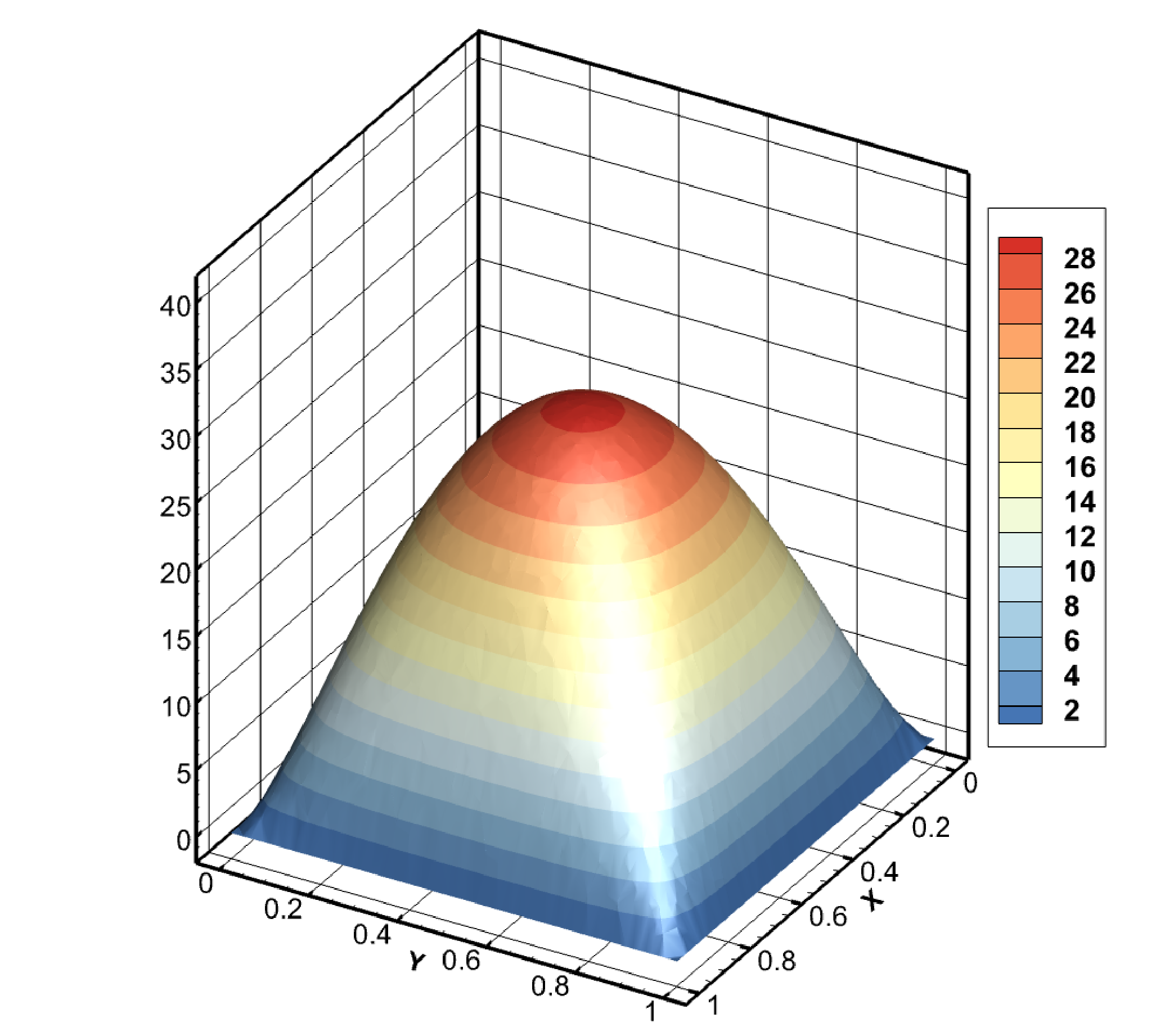}
				(a)
			\end{minipage}
			\hfill
			\begin{minipage}[c]{0.24\textwidth}
				\centering
				\includegraphics[width=\linewidth]{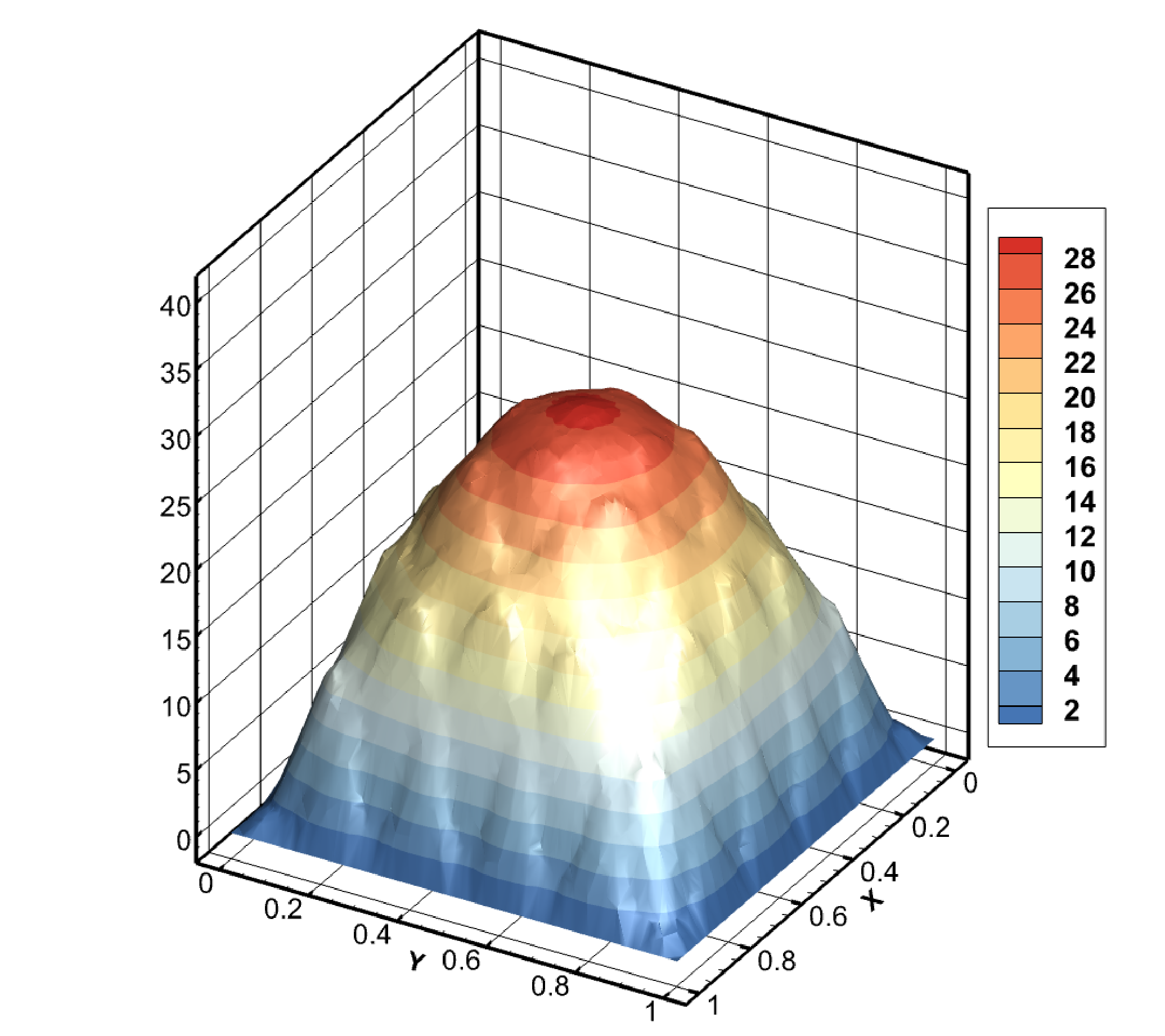}
				(b)
			\end{minipage}
			\hfill
			\begin{minipage}[c]{0.24\textwidth}
				\centering
				\includegraphics[width=\linewidth]{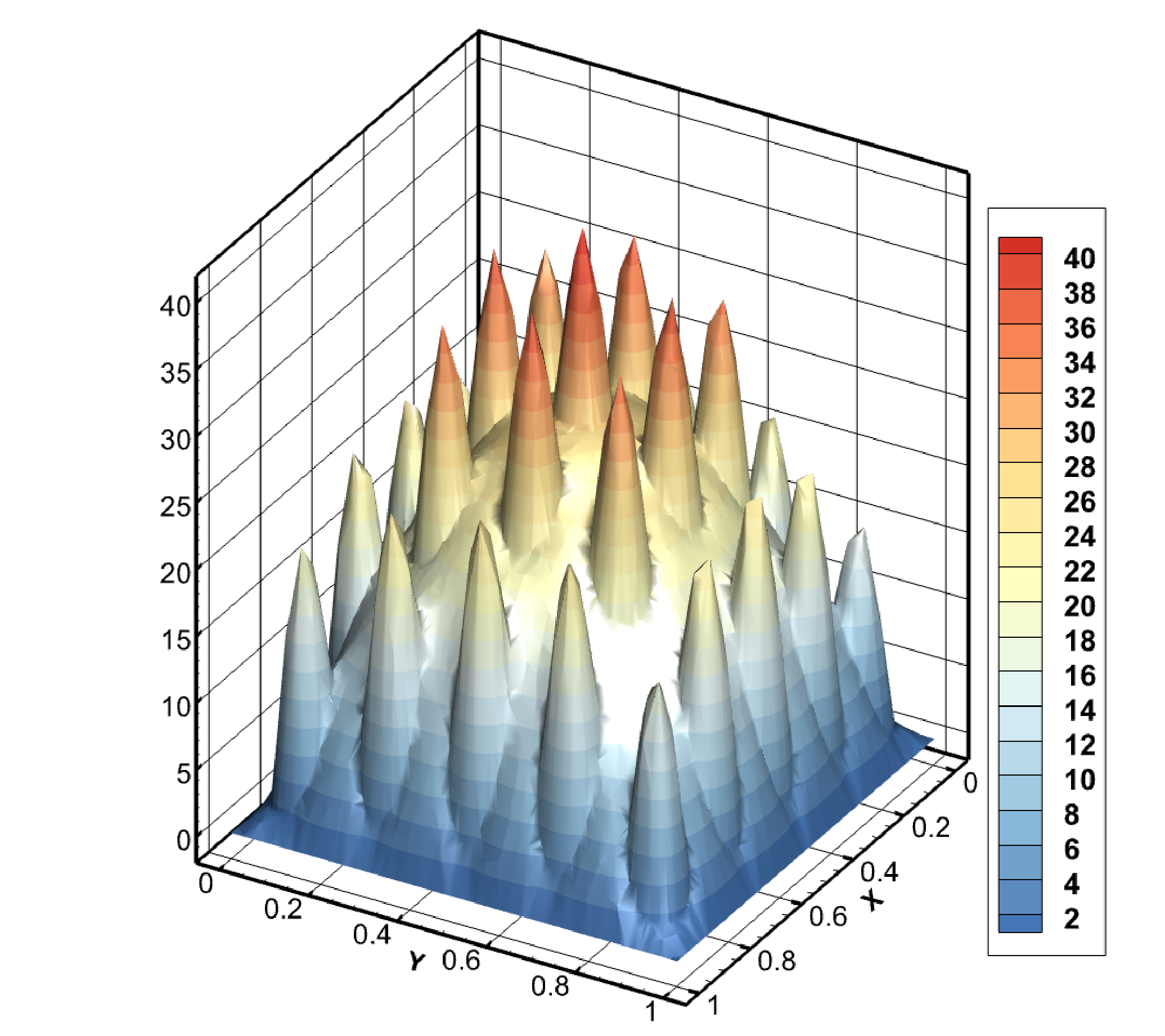}
				(c)
			\end{minipage}
			\hfill
			\begin{minipage}[c]{0.24\textwidth}
				\centering
				\includegraphics[width=\linewidth]{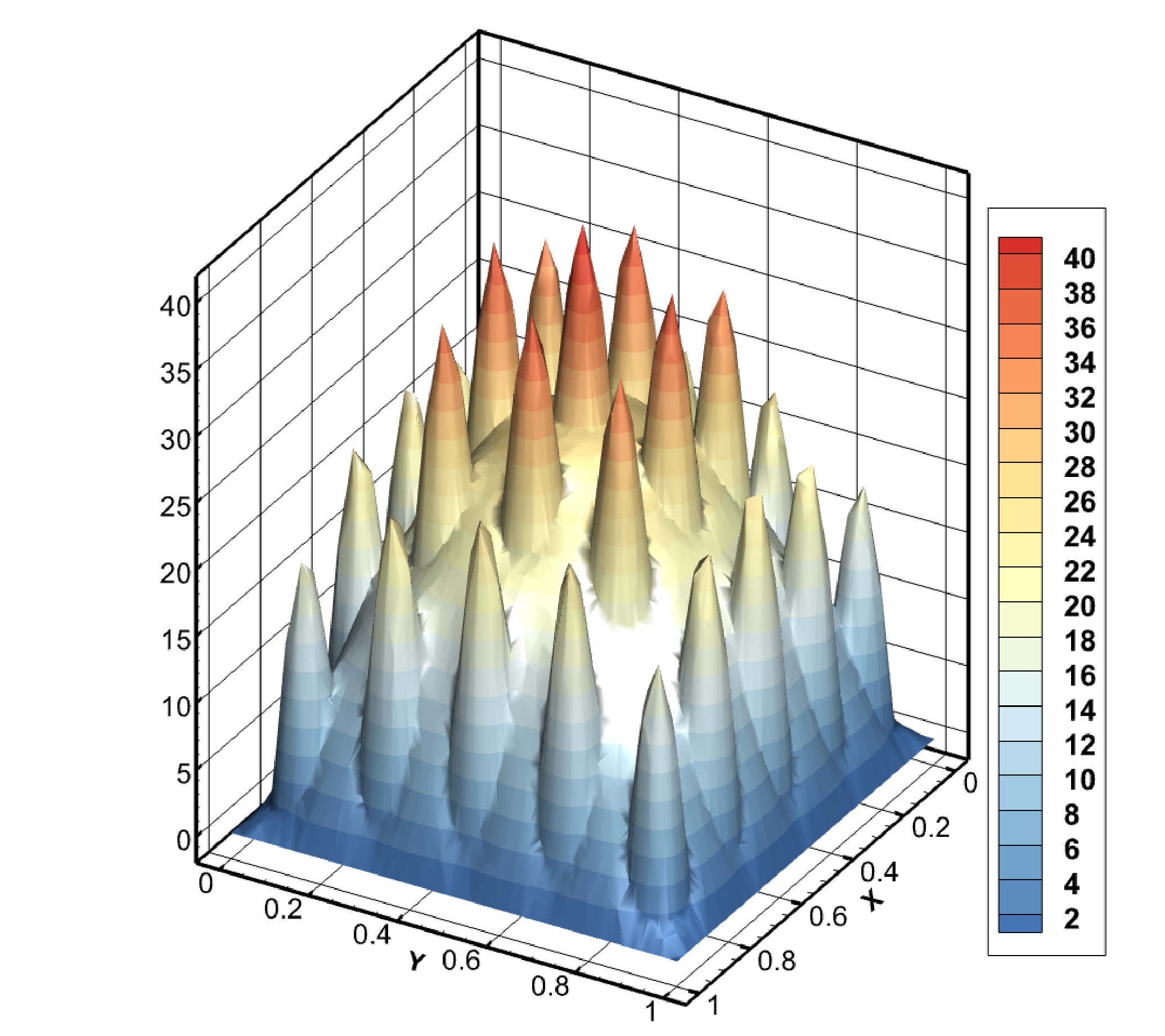}
				(d)
			\end{minipage}
			\caption{Moisture field in $x_3=0.3 \mathrm{cm}$ under scale coupling: (a) $c^{(0)}$; (b) $c^{(1,\epsilon)}$; (c) $c^{(2,\epsilon)}$; (d) $c^{\epsilon}$.}\label{f17}
		\end{figure}

\begin{figure}[!htb]
			\centering
			\begin{minipage}[c]{0.24\textwidth}
				\centering
				\includegraphics[width=\linewidth]{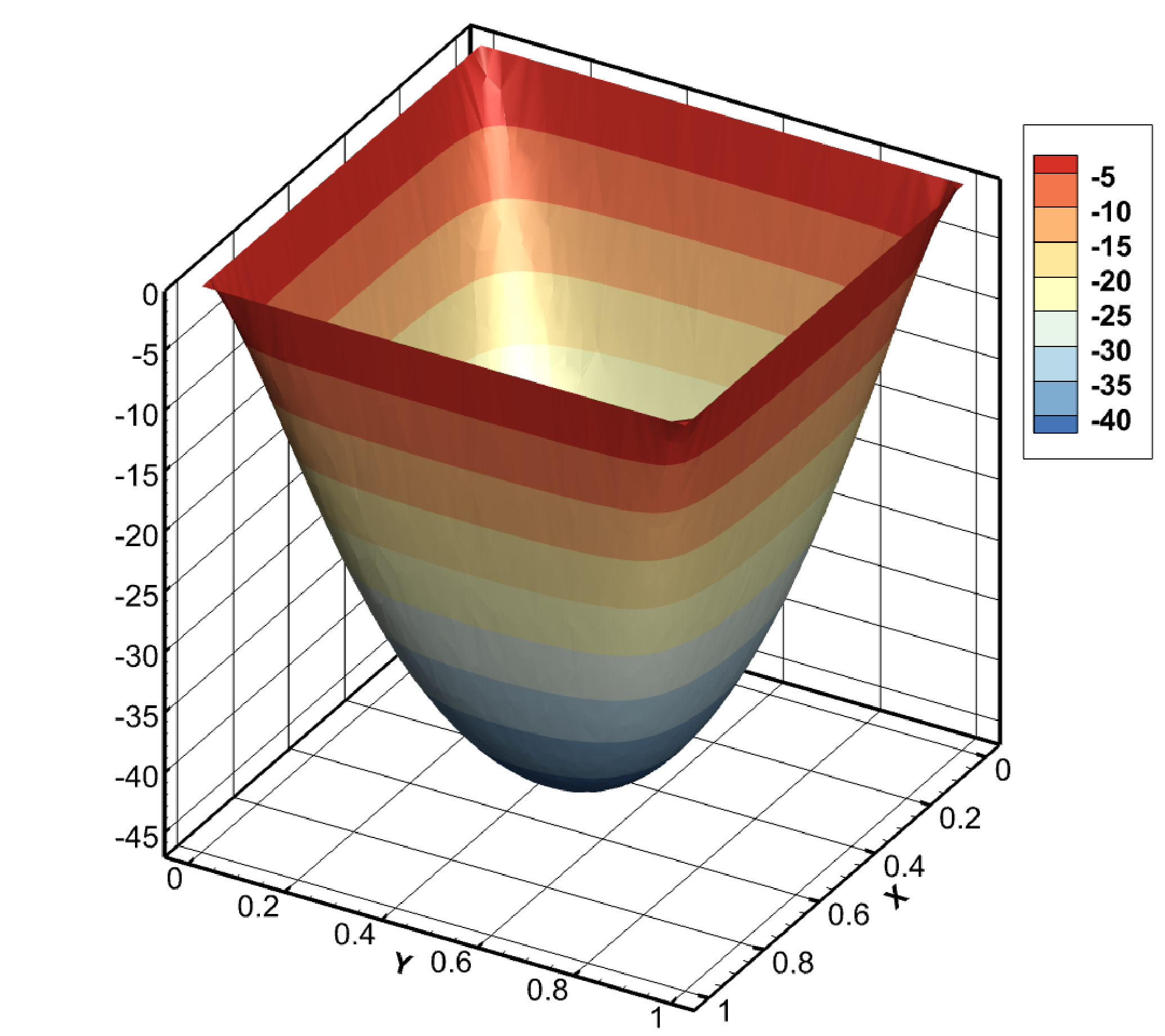}
				(a)
			\end{minipage}
			\hfill
			\begin{minipage}[c]{0.24\textwidth}
				\centering
				\includegraphics[width=\linewidth]{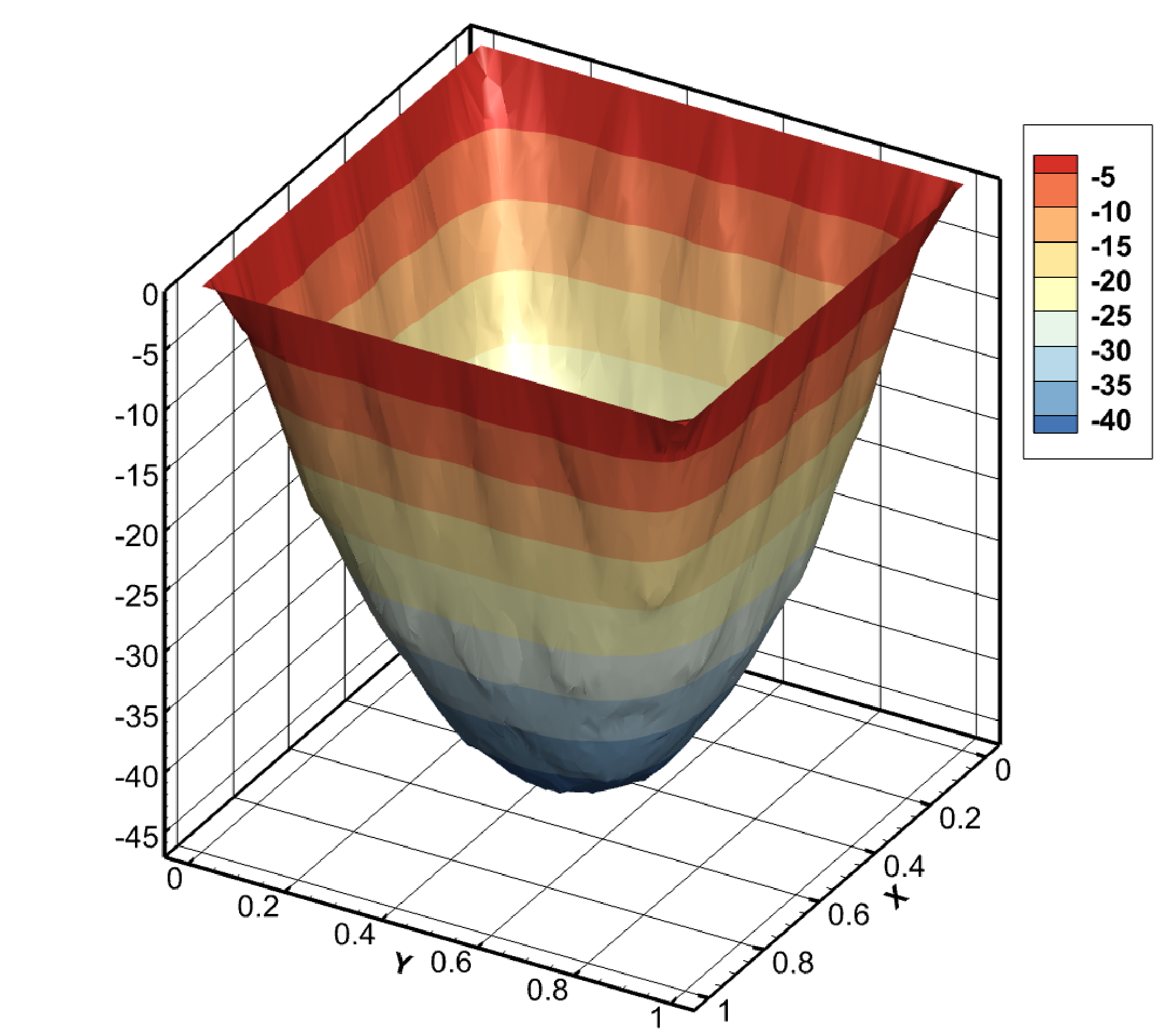}
				(b)
			\end{minipage}
			\hfill
			\begin{minipage}[c]{0.24\textwidth}
				\centering
				\includegraphics[width=\linewidth]{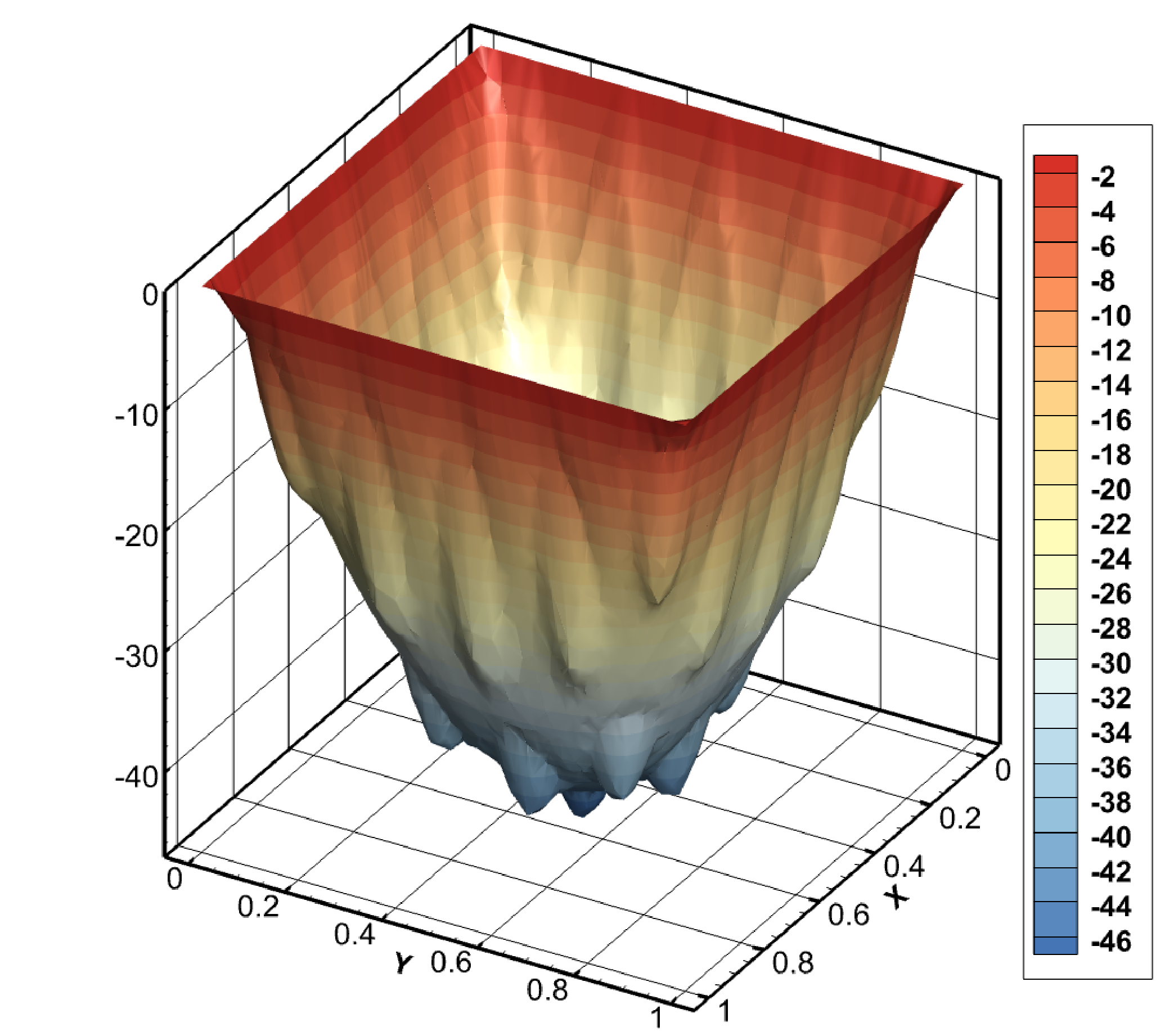}
				(c)
			\end{minipage}
			\hfill
			\begin{minipage}[c]{0.24\textwidth}
				\centering
				\includegraphics[width=\linewidth]{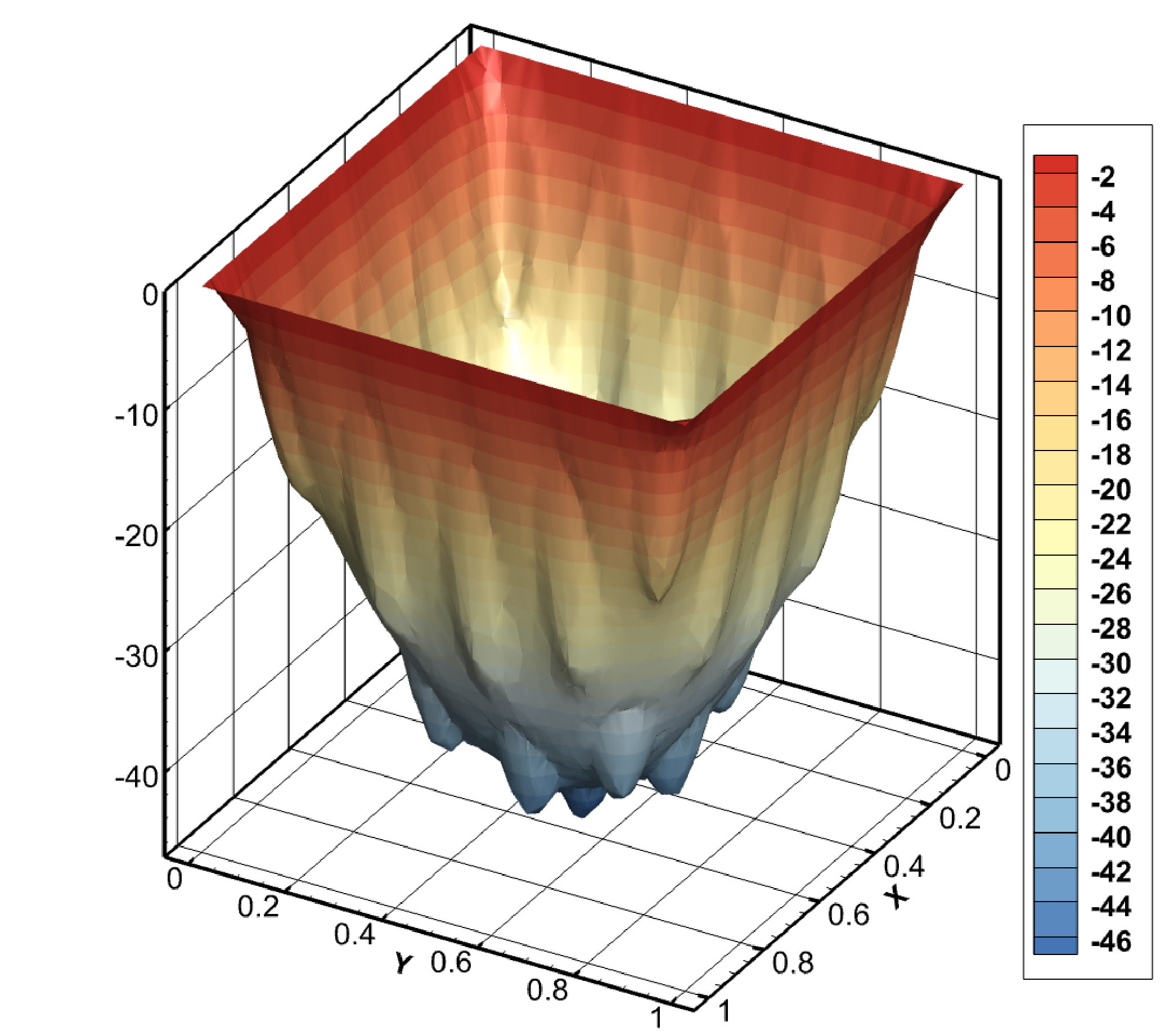}
				(d)
			\end{minipage}
			\caption{Third displacement field component in $x_3=0.3 \mathrm{cm}$ under scale coupling: (a) $u_3^{(0)}$; (b) $u_3^{(1,\epsilon)}$; (c) $u_3^{(2,\epsilon)}$; (d) $u_3^{\epsilon}$.}\label{f20}
		\end{figure}

		After numerical simulations, the results are presented in Figs.\hspace{1mm}\ref{f11}-\ref{f20}. Among them, Figs.\hspace{1mm}\ref{f11}-\ref{f15} illustrate the solutions for the case with scale-separated material coefficients, while Figs.\hspace{1mm}\ref{f16}-\ref{f20} correspond to the case without scale-separation. Furthermore, Tables\hspace{1mm}\ref{t7} and \ref{t8} present the numerical errors of different kinds of solutions.
		\begin{table}[!htb]
			\centering
			\caption{The relative errors for scale-separated material parameters.}
			\label{t7}
			\begin{tabular}{cccccc}
				\hline
				\multicolumn{6}{c}{Temperature increment field} \\
				\hline
				$TerrorL^20$ & $TerrorL^21$ & $TerrorL^22$ & $TerrorH^10$ & $TerrorH^11$ & $TerrorH^12$ \\
				0.27929 & 0.28013 & 0.03437 & 0.94076 & 0.93121 & 0.12766 \\
				\hline
				\multicolumn{6}{c}{Moisture field} \\
				\hline
				$cerrorL^20$ & $cerrorL^21$ & $cerrorL^22$ & $cerrorH^10$ & $cerrorH^11$ & $cerrorH^12$ \\
				0.14687 & 0.14648 & 0.02093 & 0.82318 & 0.80311 & 0.11502 \\
				\hline
				\multicolumn{6}{c}{Displacement field} \\
				\hline
				$\bm{u}errorL^20$ & $\bm{u}errorL^21$ & $\bm{u}errorL^22$ & $\bm{u}errorH^10$ & $\bm{u}errorH^11$ & $\bm{u}errorH^12$ \\
				0.05345 & 0.02731 & 0.02083 & 0.40714 & 0.22018 & 0.15868 \\
				\hline
			\end{tabular}
		\end{table}
		\begin{table}[!htb]
			\centering
			\caption{The relative errors for scale-coupled material parameters.}
			\label{t8}
			\begin{tabular}{cccccc}
				\hline
				\multicolumn{6}{c}{Temperature increment field} \\
				\hline
				$TerrorL^20$ & $TerrorL^21$ & $TerrorL^22$ & $TerrorH^10$ & $TerrorH^11$ & $TerrorH^12$ \\
				0.22774 & 0.22821 & 0.02952 & 0.91539 & 0.90360 & 0.12100 \\
				\hline
				\multicolumn{6}{c}{Moisture field} \\
				\hline
				$cerrorL^20$ & $cerrorL^21$ & $cerrorL^22$ & $cerrorH^10$ & $cerrorH^11$ & $cerrorH^12$ \\
				0.11843 & 0.11735 & 0.01810 & 0.76358 & 0.73929 & 0.11277 \\
				\hline
				\multicolumn{6}{c}{Displacement field} \\
				\hline
				$\bm{u}errorL^20$ & $\bm{u}errorL^21$ & $\bm{u}errorL^22$ & $\bm{u}errorH^10$ & $\bm{u}errorH^11$ & $\bm{u}errorH^12$ \\
				0.04919 & 0.02629 & 0.01713 & 0.38455 & 0.21311 & 0.12067 \\
				\hline
			\end{tabular}
		\end{table}
		
		As shown in Figs.\hspace{1mm}\ref{f11}-\ref{f20} and supported by the error results in Tables\hspace{1mm}\ref{t7} and \ref{t8}, only HOMS solutions closely match the precise FEM results, and successfully capture the micro-scale oscillations arising from the spatial heterogeneities within composite structures. By comparison, both the homogenized and LOMS solutions are insufficient to provide high-accuracy solutions for the multi-scale problem \eqref{eq:2.1}. Thus, the proposed HOMS approach demonstrates significantly higher computational precision relative to both the traditional homogenized method and the LOMS method.
		
		\subsection{Example 3: 3D quasi-periodic composite plate}
		\label{sec:53}
\begin{figure}[!htb]
			\centering
			\begin{minipage}[c]{0.45\textwidth}
				\centering
				\includegraphics[width=55mm]{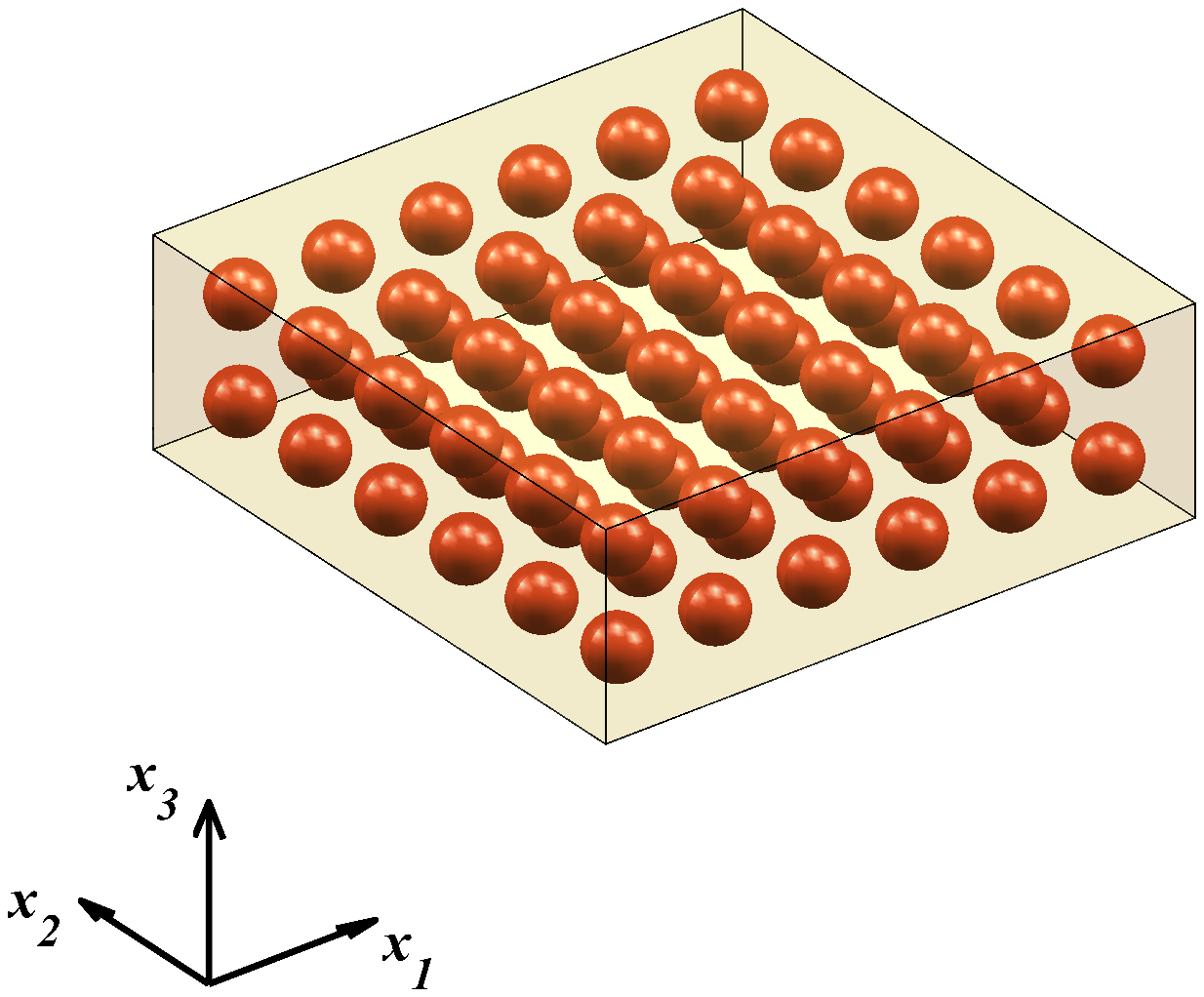}\\
				(a)
			\end{minipage}
			\begin{minipage}[c]{0.45\textwidth}
				\centering
				\includegraphics[width=55mm]{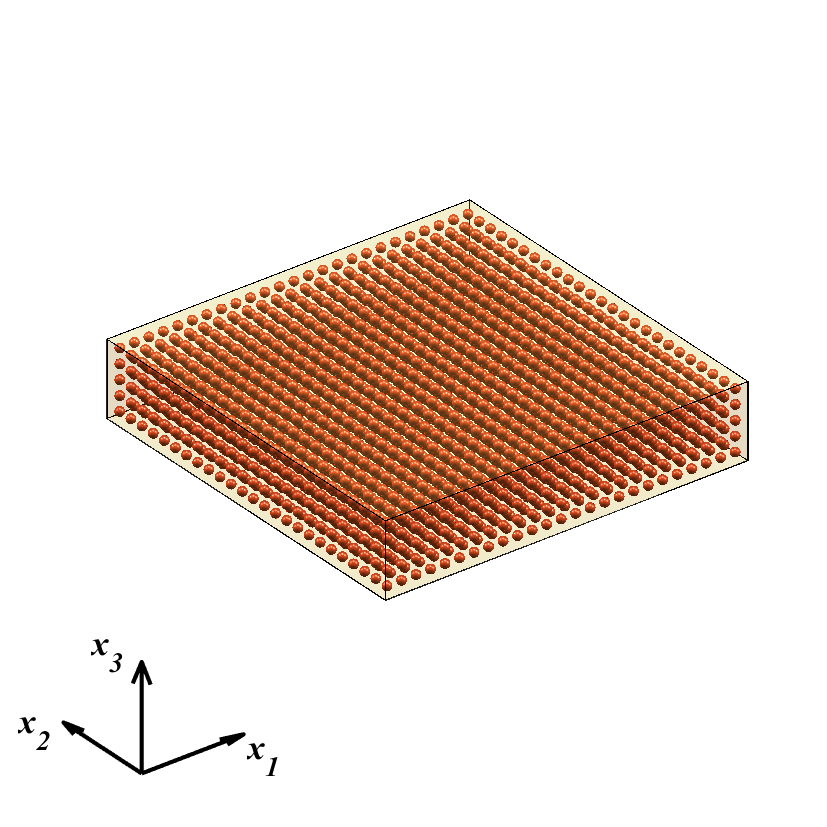}\\
				(b)
			\end{minipage}
			\begin{minipage}[c]{0.45\textwidth}
				\centering
				\includegraphics[width=55mm]{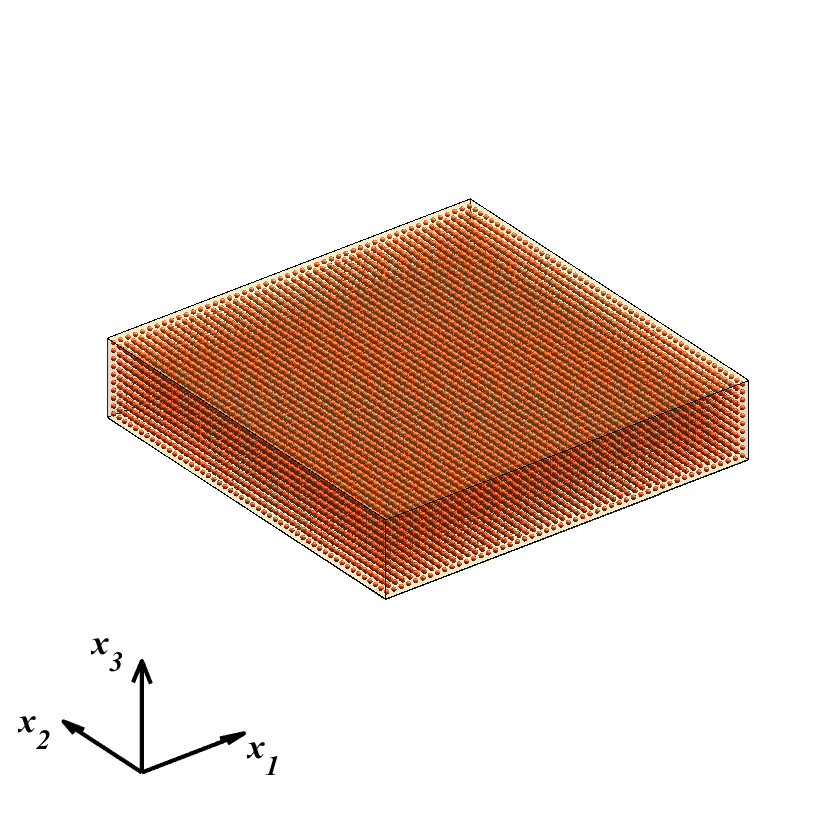}\\
				(c)
			\end{minipage}
			\begin{minipage}[c]{0.45\textwidth}
				\centering
				\includegraphics[width=55mm]{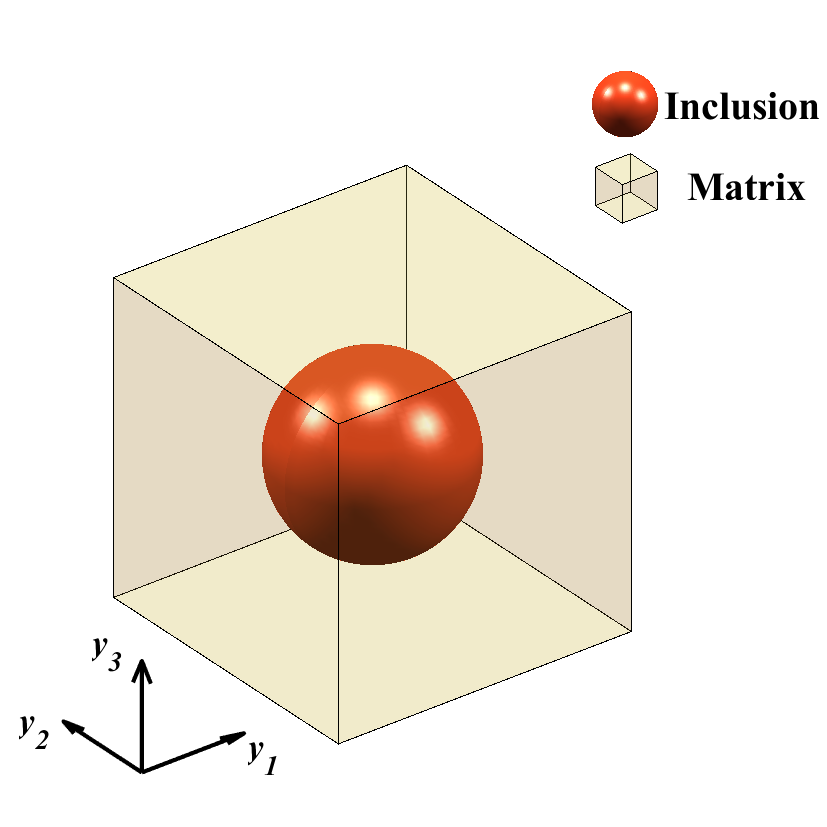}\\
				(d)
			\end{minipage}
			\caption{The schematic of composite structure: (a) composite structure $\Omega$ (Case 1); (b) composite structure $\Omega$ (Case 2); (c) composite structure $\Omega$ (Case 3);  (d) PUC $Y$.}\label{f21:3Dplate}
		\end{figure}

		In this example, we will validate the proposed HOMS numerical method for solving the multi-scale problem \eqref{eq:2.1} of a composite plate with quasi-periodic material properties. The structures of macroscopic domain $\Omega$ and PUC $Y$ are shown in Fig.\hspace{1mm}\ref{f21:3Dplate}. To examine the computational efficiency of the proposed method, three plate configurations are investigated:
		\begin{enumerate}
			\item[Case 1:] The macrostructure domain $\Omega$, shown in Fig.\hspace{1mm}\ref{f21:3Dplate}(a), is composed of $6\times6\times2$ unit cells, and occupies $\Omega=(x_1,x_2,x_3)=[0,1] \times [0,1] \times[0,1/3] \mathrm{cm}^3$. The small periodic parameter $\epsilon=1/6$.
			\item[Case 2:] The macrostructure domain $\Omega$, shown in Fig.\hspace{1mm}\ref{f21:3Dplate}(b), is composed of $25\times25\times5$ unit cells, and occupies $\Omega=(x_1,x_2,x_3)=[0,1] \times [0,1] \times[0,1/5] \mathrm{cm}^3$. The small periodic parameter $\epsilon=1/25$.
			\item[Case 3:] The macrostructure domain $\Omega$, shown in Fig.\hspace{1mm}\ref{f21:3Dplate}(c), is composed of $50\times50\times10$ unit cells, and occupies $\Omega=(x_1,x_2,x_3)=[0,1] \times [0,1] \times[0,1/5] \mathrm{cm}^3$. The small periodic parameter $\epsilon=1/50$.
		\end{enumerate}
		
		This example utilizes material parameters characterized by the quasi-periodic form $a(\mathbf{x},\mathbf{y}) = \tilde{a}(\mathbf{x}) + \hat{a}(\mathbf{y})$. The specific values are listed in Table\hspace{1mm}\ref{t2}, where $\psi(\mathbf{x}) = x_3$. The internal heat source, moisture source, body forces, and boundary conditions for this example are consistent with Example\hspace{1mm}2.
		
		Computational cost for the three cases are detailed in Tables\hspace{1mm}\ref{t9}-\ref{t11}. Since the HOMS method requires significantly fewer FEM nodes and elements than precise FEM, it significantly reduces computational resource requirement. As evidenced by Tables\hspace{1mm}\ref{t9} and \ref{t10}, the computational time for Case 2 using the HOMS method is 91 times that of Case 1, while the computational time for Case 2 using precise FEM is 435 times that of Case 1. This demonstrates that as the scale of composite structures increases, the computational time of the HOMS method grows at a much slower rate than precise FEM. Notably, Table\hspace{1mm}\ref{t11} shows that precise FEM failed to obtain high-accuracy solutions in Case 3 due to excessively fine mesh requirements. 
		
		In a summary, these results indicate that the proposed HOMS method exhibits progressively slower growth in computational time with increasing material scale. Moreover, it remains computationally feasible even when precise FEM approaches become infeasible. This confirms that the HOMS method reduces computational cost and enhances efficiency, which offer significant advantages for real-world engineering applications.
		\begin{table}[!htb]{\caption{Comparison of computational cost (Case 1: $6\times6\times2$ unit cells).}\label{t9}}
			\centering
			\begin{tabular}{cccc}
				\hline
				& Cell equations & Homogenized equations & Multi-scale equations\\
				\hline
				FEM nodes & 1212 & 1536 & 28059\\
				FEM elements & 6007 & 6750 & 165501 \\
				\hline
				& \multicolumn{2}{c}{HOMS method} & precise FEM\\
				\hline
				Computational time & \multicolumn{2}{c}{944.848s} & 26.604s \\
				\hline
			\end{tabular}
		\end{table}
		\begin{table}[!htb]{\caption{Comparison of computational cost (Case 2: $25\times25\times5$ unit cells).}\label{t10}}
			\centering
			\begin{tabular}{cccc}
				\hline
				& Cell equations & Homogenized equations & Multi-scale equations\\
				\hline
				FEM nodes & 1212 & 1024 & 2865732\\
				FEM elements & 6007 & 4050 & 18032867 \\
				\hline
				&\multicolumn{2}{c}{HOMS method} & precise FEM\\
				\hline
				Computational time & \multicolumn{2}{c}{86429.823s} & 11576.982s \\
				\hline
			\end{tabular}
		\end{table}
		\begin{table}[!htb]{\caption{Comparison of computational cost (Case 3: $50\times50\times10$ unit cells).}\label{t11}}
			\centering
			\begin{tabular}{cccc}
				\hline
				& Cell equations & Homogenized equations & Multi-scale equations \\
				\hline
				FEM nodes & 1212 & 1024 & $\approx 30300000 (estimated)$\\
				FEM elements & 6007 & 4050 & $\approx 150175000 (estimated)$ \\
				\hline
			\end{tabular}
		\end{table}
		
\begin{figure}[!htb]
			\centering
			\begin{minipage}[c]{0.24\textwidth}
				\centering
				\includegraphics[width=\linewidth]{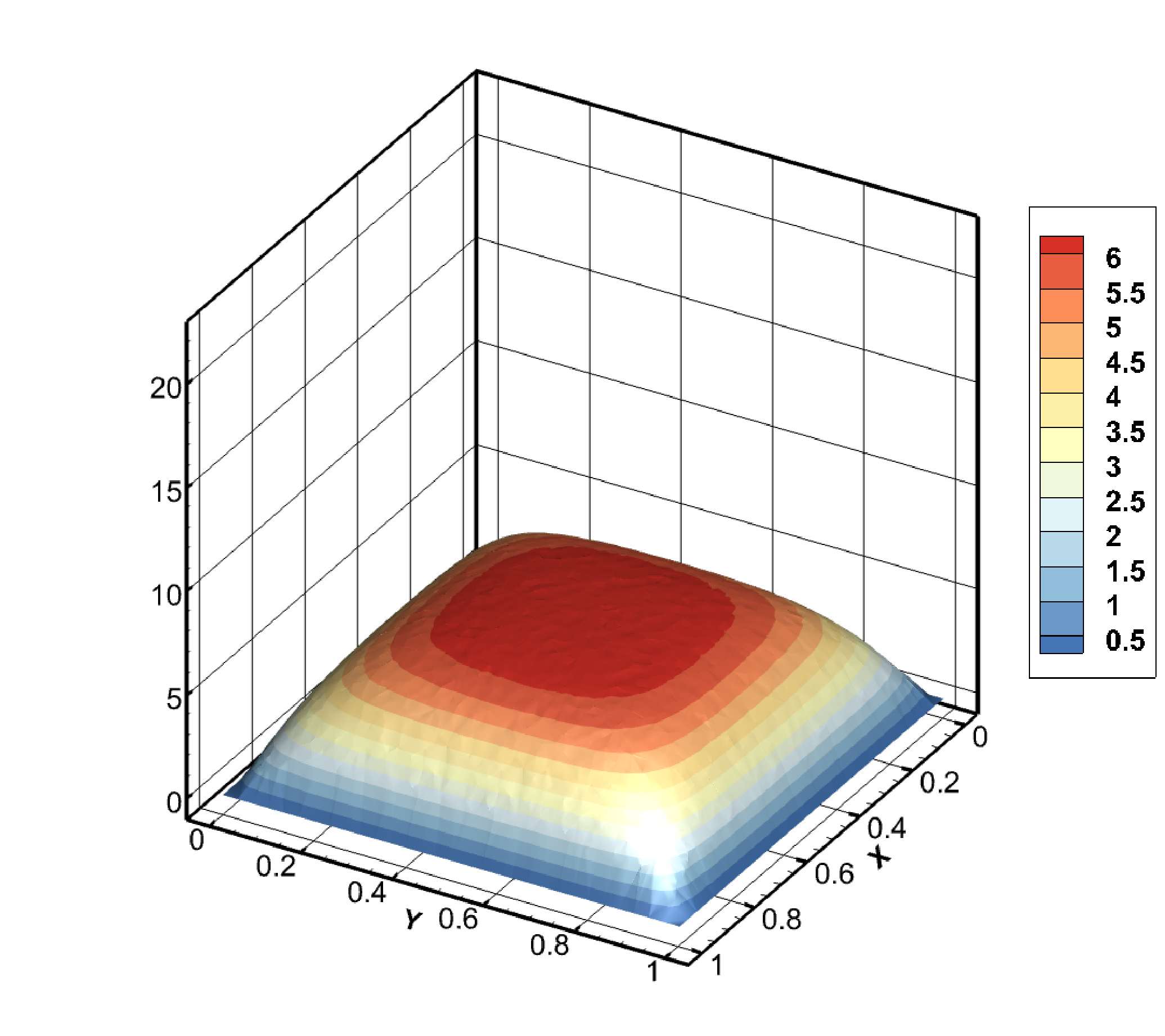}
				(a)
			\end{minipage}
			\hfill
			\begin{minipage}[c]{0.24\textwidth}
				\centering
				\includegraphics[width=\linewidth]{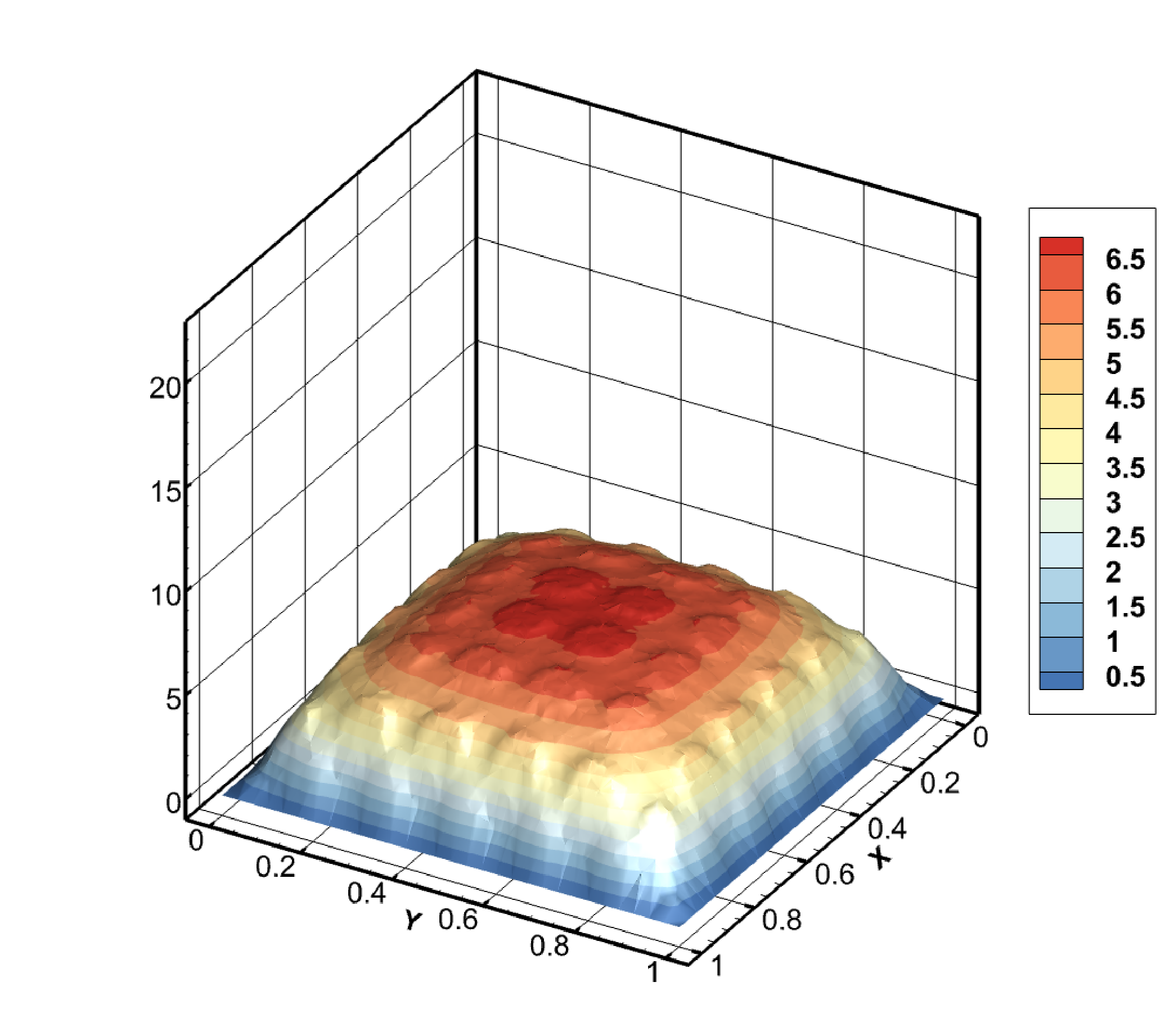}
				(b)
			\end{minipage}
			\hfill
			\begin{minipage}[c]{0.24\textwidth}
				\centering
				\includegraphics[width=\linewidth]{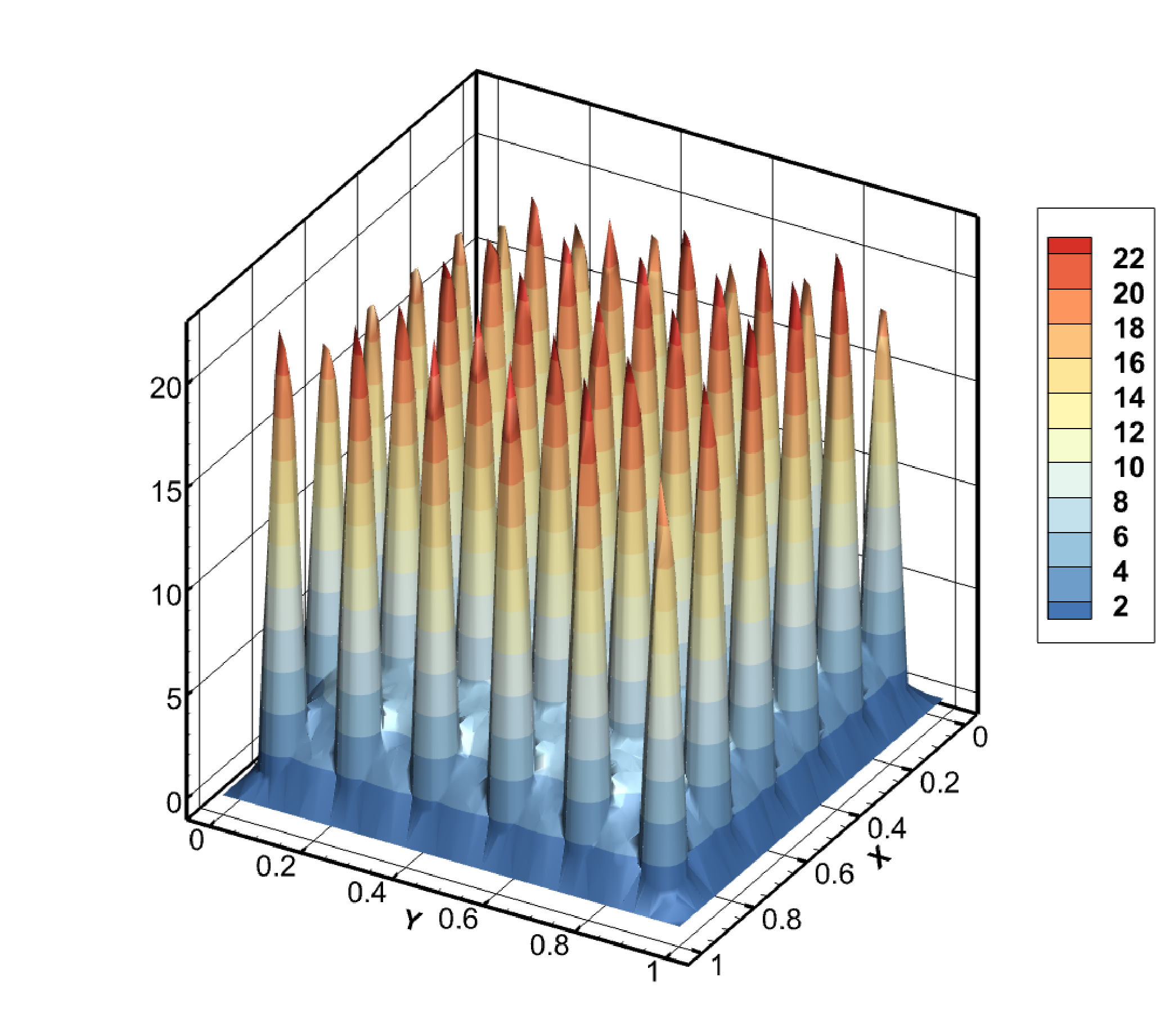}
				(c)
			\end{minipage}
			\hfill
			\begin{minipage}[c]{0.24\textwidth}
				\centering
				\includegraphics[width=\linewidth]{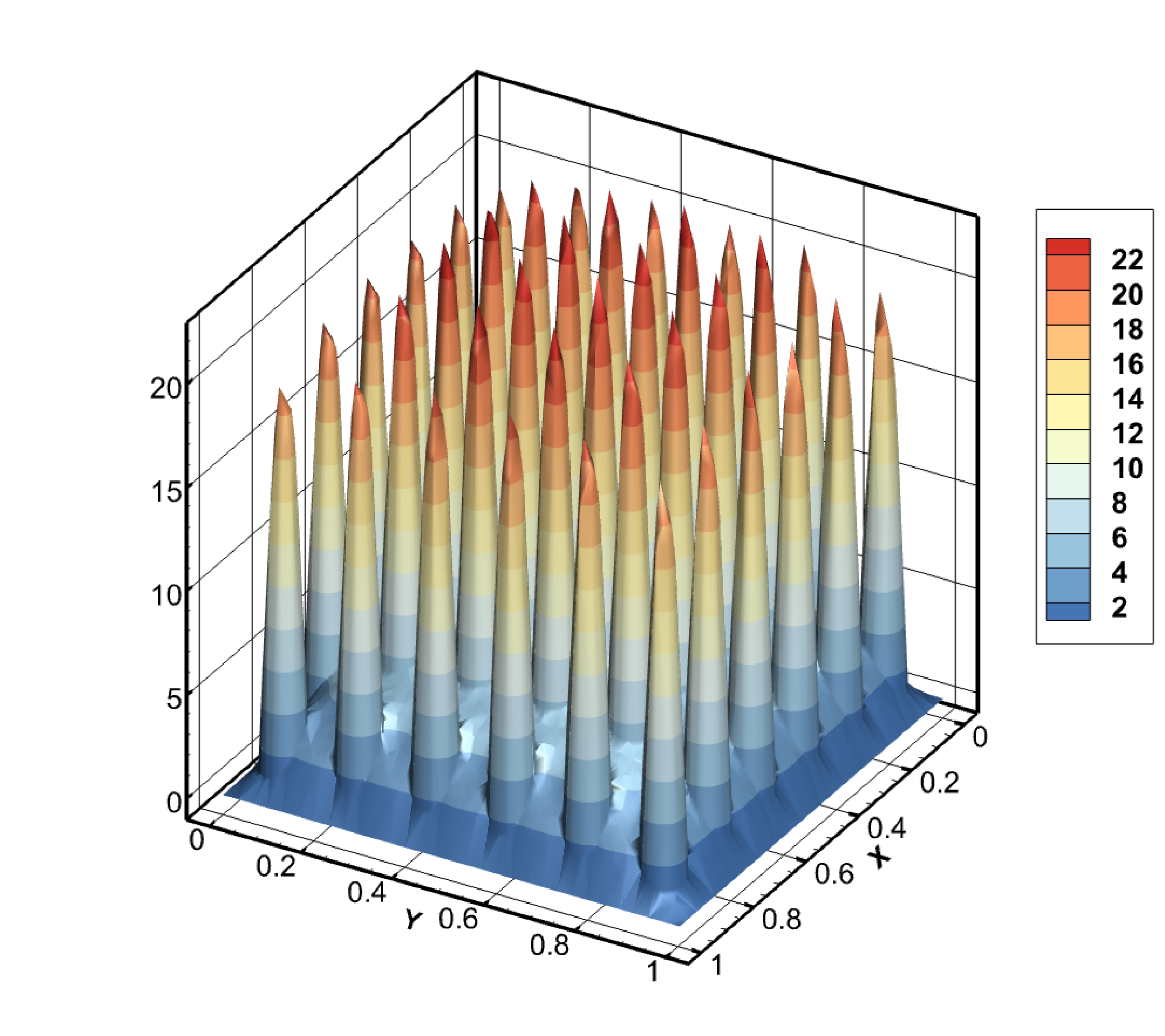}
				(d)
			\end{minipage}
			\caption{Case 1: Temperature increment field in $x_3=0.1 \mathrm{cm}$: (a) $T^{(0)}$; (b) $T^{(1,\epsilon)}$; (c) $T^{(2,\epsilon)}$; (d) $T^{\epsilon}$.}\label{f22}
		\end{figure}

\begin{figure}[!htb]
			\centering
			\begin{minipage}[c]{0.24\textwidth}
				\centering
				\includegraphics[width=\linewidth]{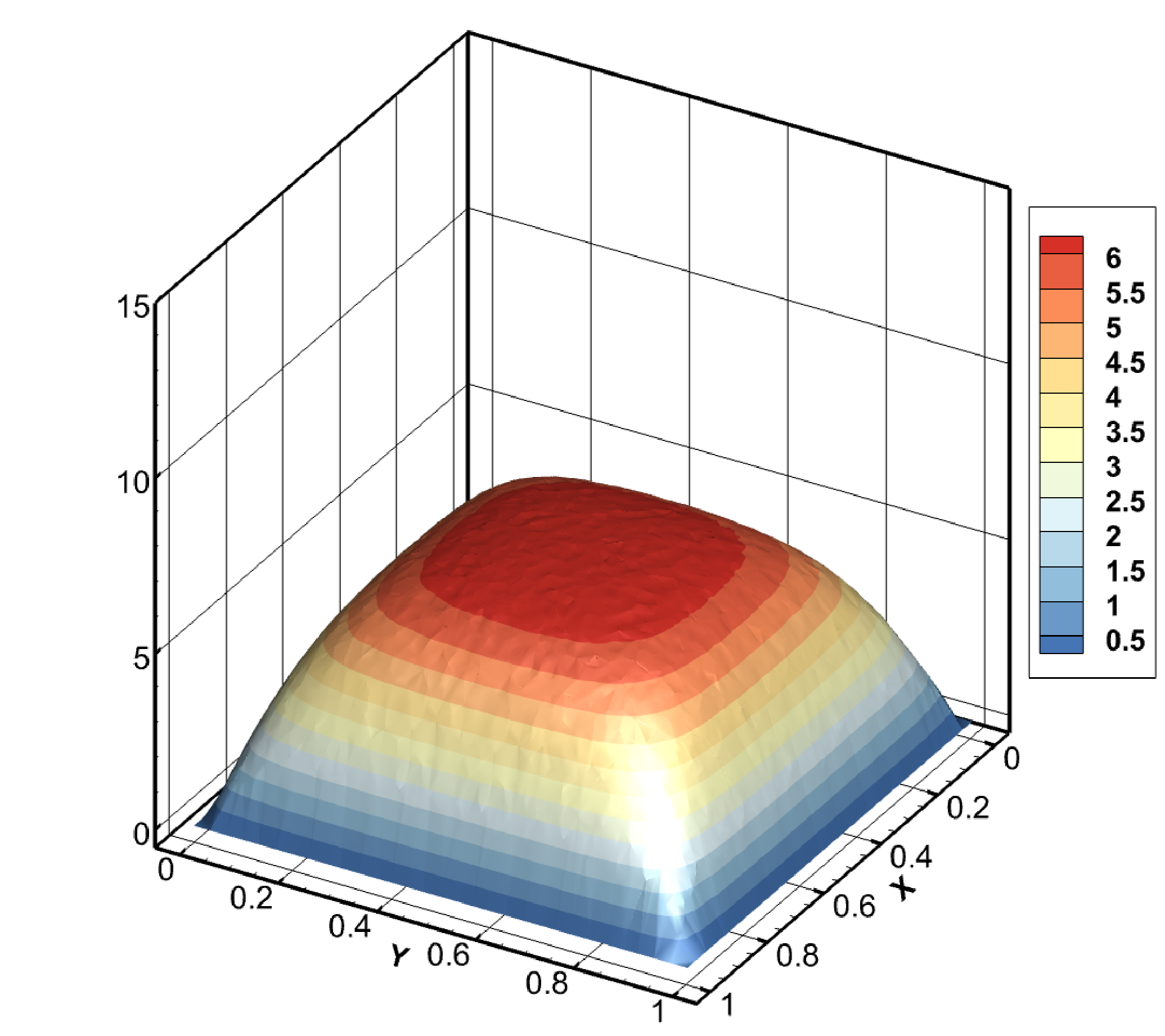}
				(a)
			\end{minipage}
			\hfill
			\begin{minipage}[c]{0.24\textwidth}
				\centering
				\includegraphics[width=\linewidth]{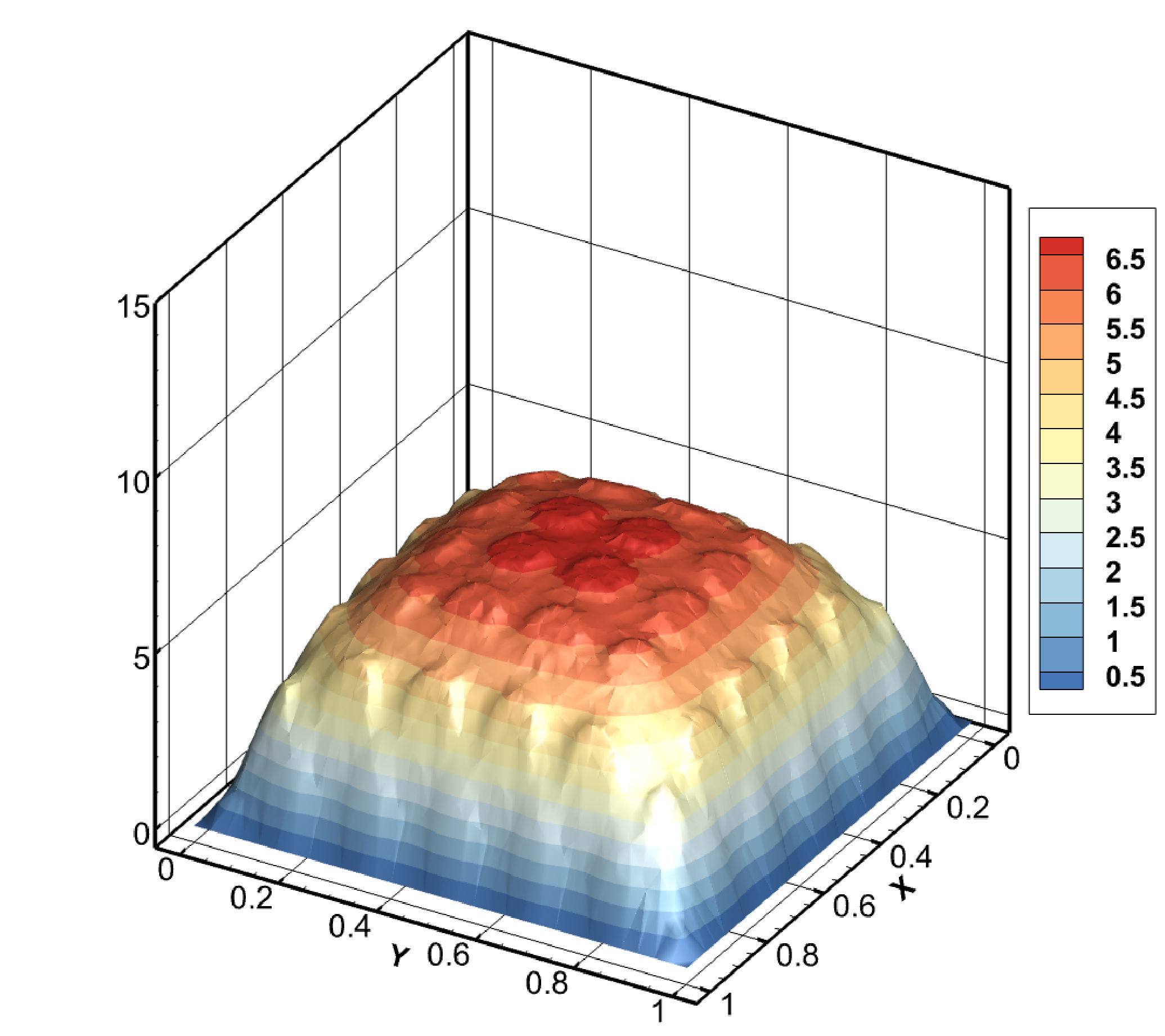}
				(b)
			\end{minipage}
			\hfill
			\begin{minipage}[c]{0.24\textwidth}
				\centering
				\includegraphics[width=\linewidth]{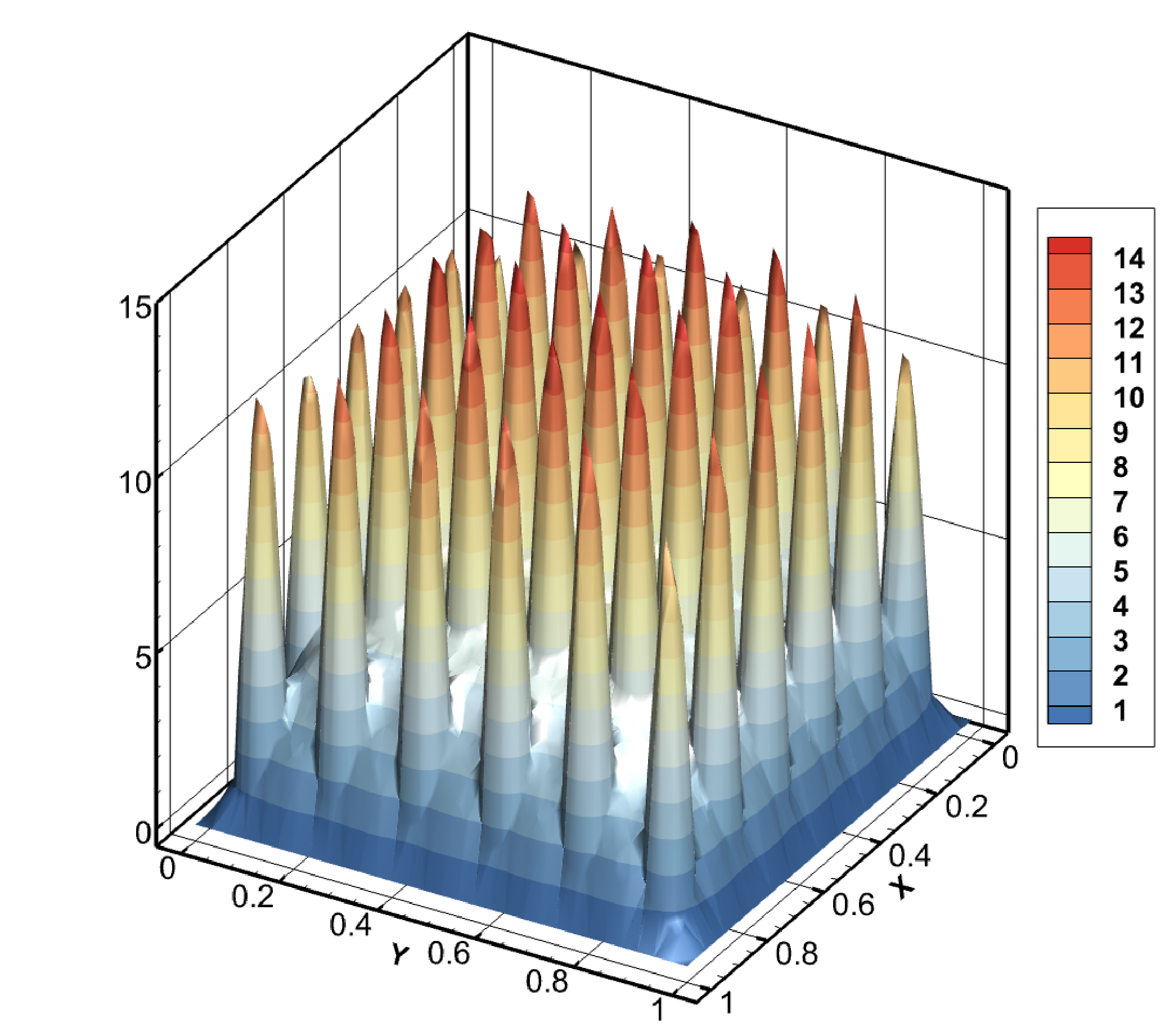}
				(c)
			\end{minipage}
			\hfill
			\begin{minipage}[c]{0.24\textwidth}
				\centering
				\includegraphics[width=\linewidth]{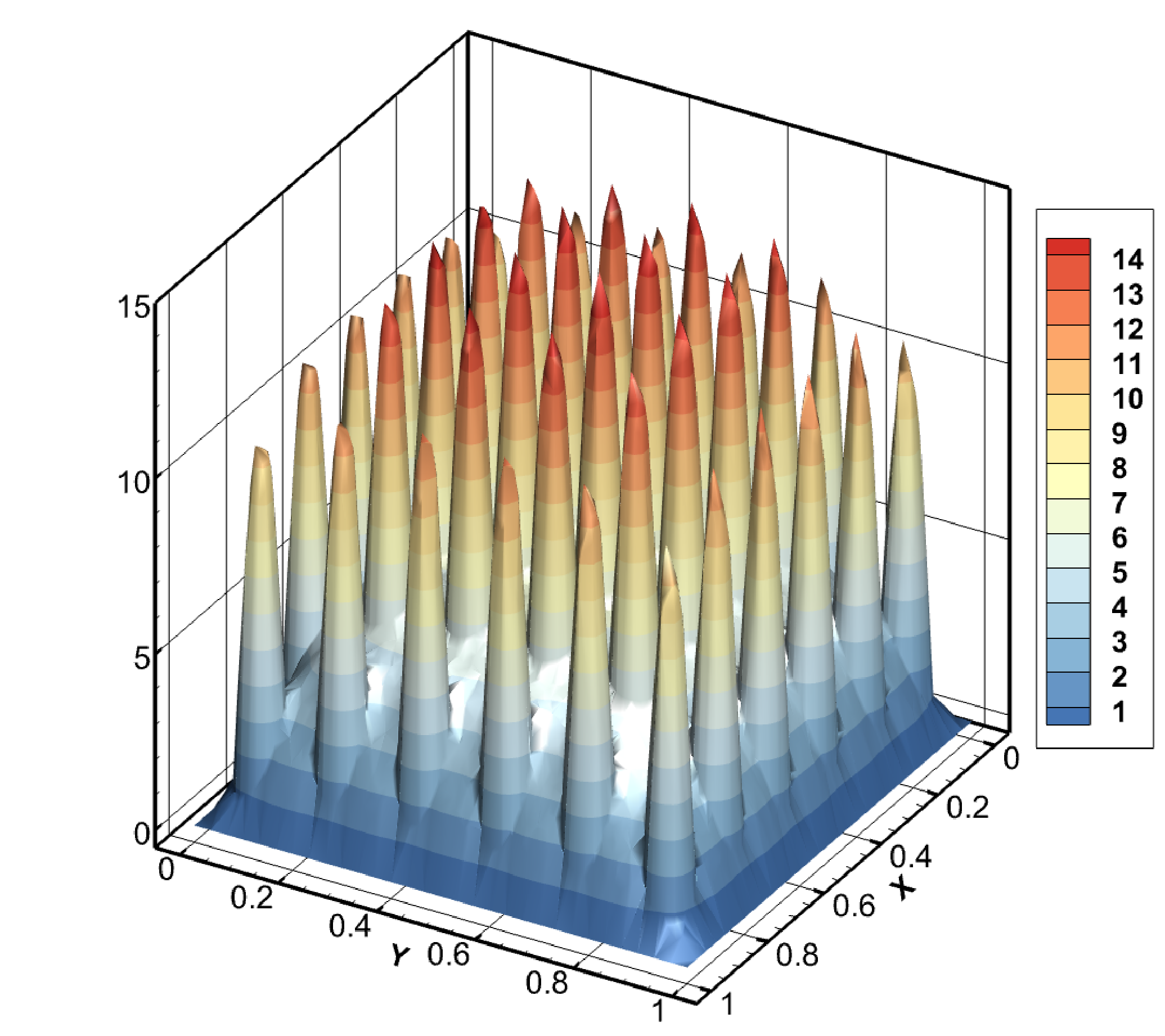}
				(d)
			\end{minipage}
			\caption{Case 1: Moisture field in $x_3=0.1 \mathrm{cm}$: (a) $c^{(0)}$; (b) $c^{(1,\epsilon)}$; (c) $c^{(2,\epsilon)}$; (d) $c^{\epsilon}$.}\label{f23}
		\end{figure}

\begin{figure}[!htb]
			\centering
			\begin{minipage}[c]{0.24\textwidth}
				\centering
				\includegraphics[width=\linewidth]{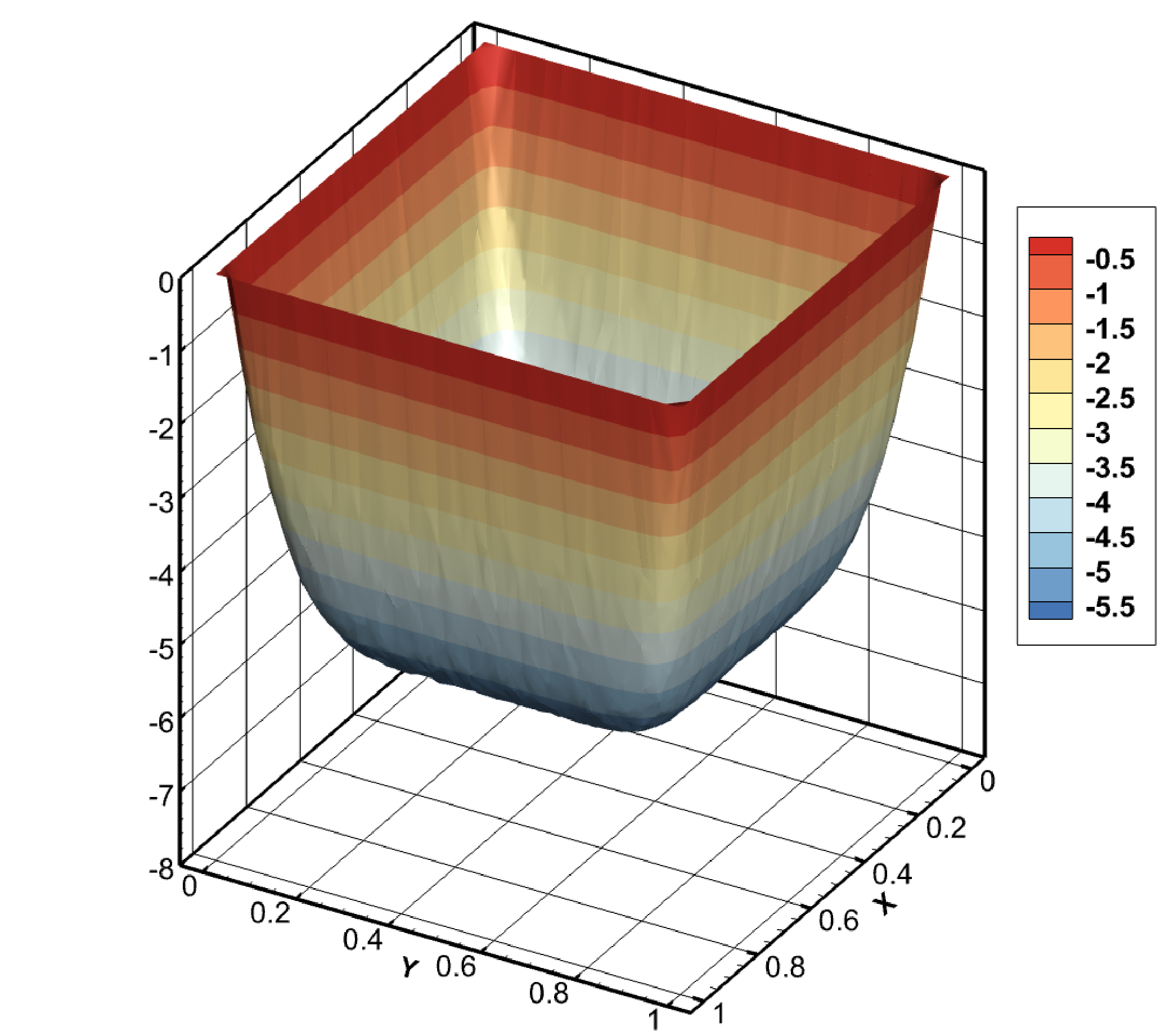}
				(a)
			\end{minipage}
			\hfill
			\begin{minipage}[c]{0.24\textwidth}
				\centering
				\includegraphics[width=\linewidth]{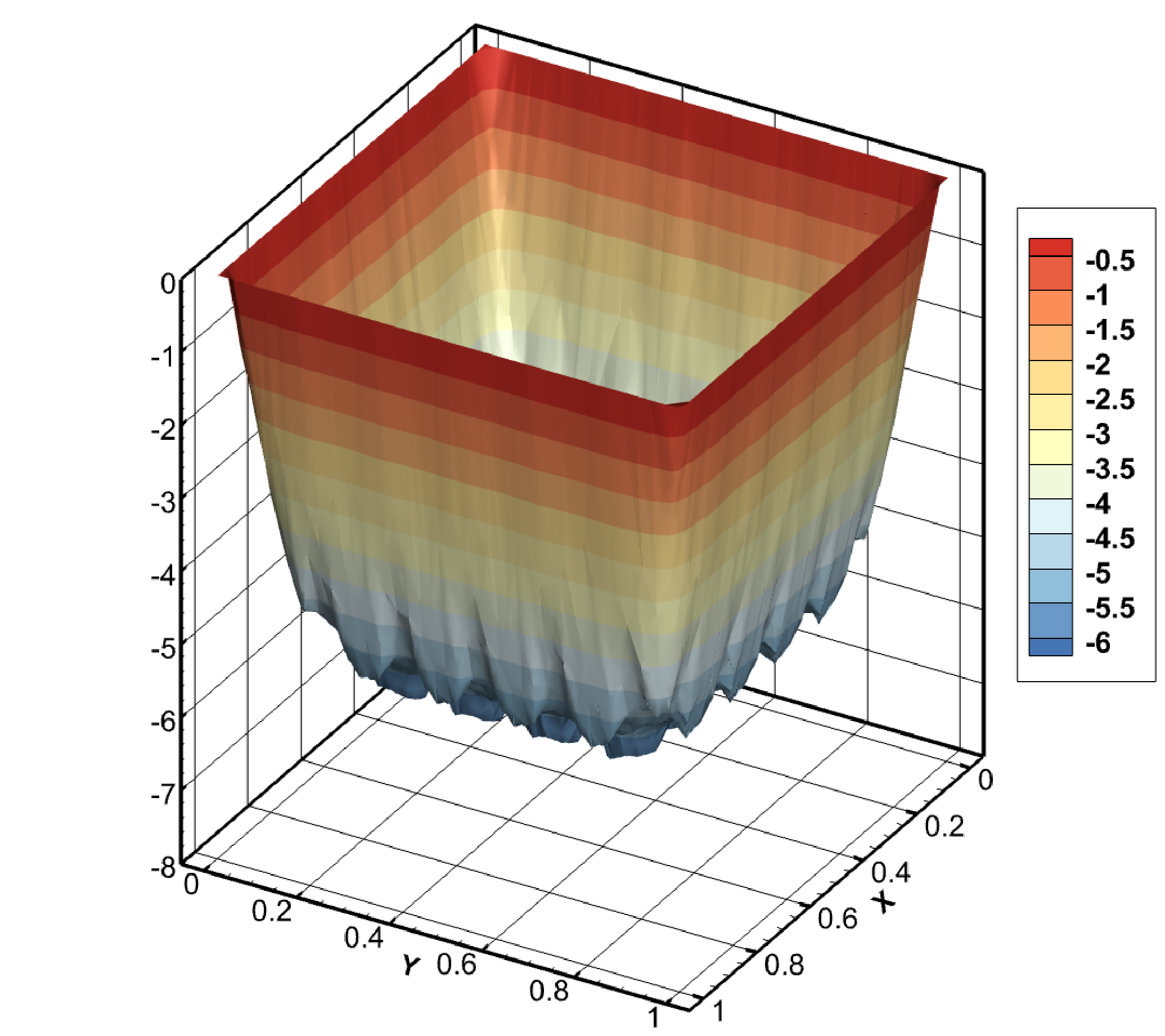}
				(b)
			\end{minipage}
			\hfill
			\begin{minipage}[c]{0.24\textwidth}
				\centering
				\includegraphics[width=\linewidth]{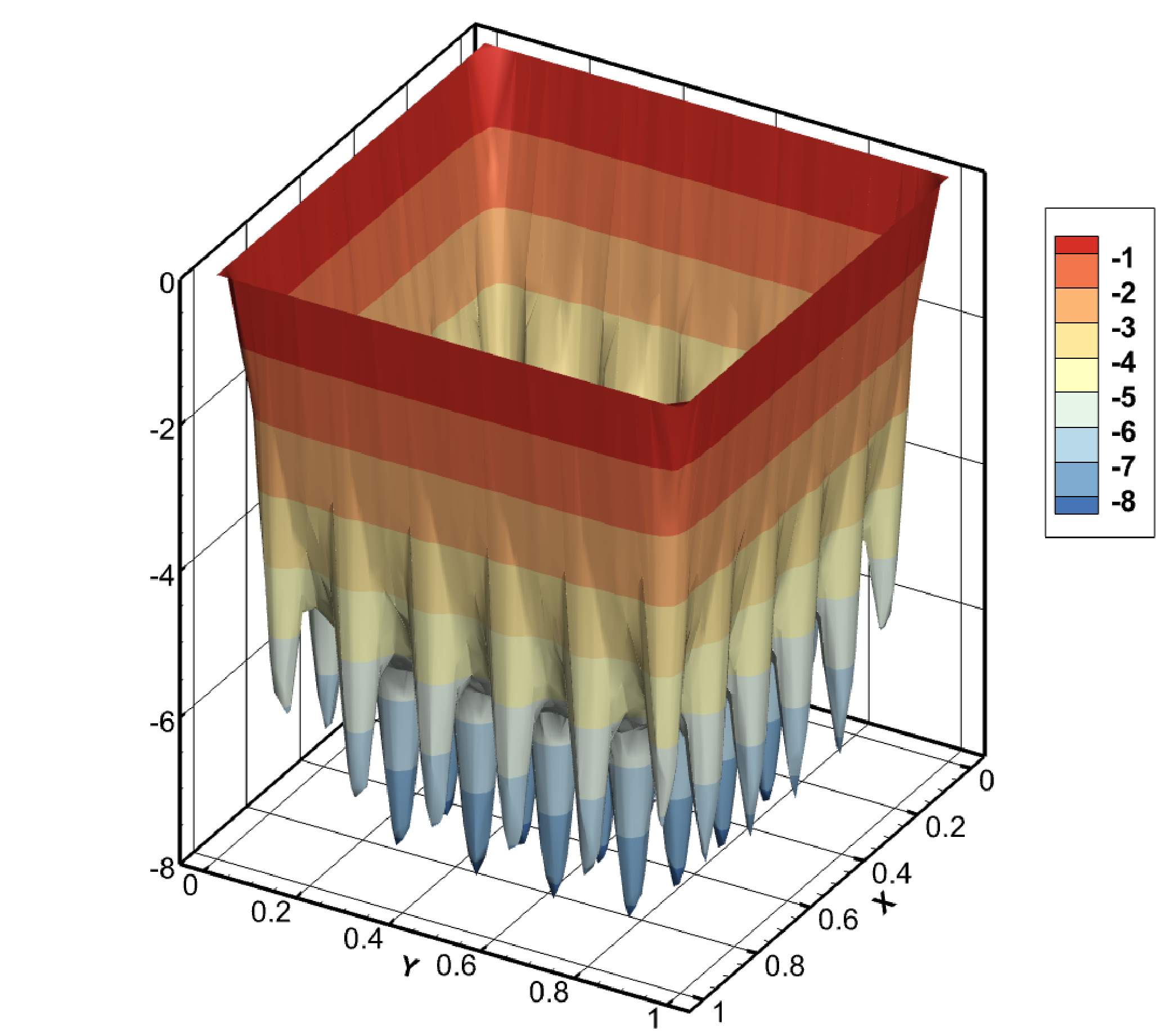}
				(c)
			\end{minipage}
			\hfill
			\begin{minipage}[c]{0.24\textwidth}
				\centering
				\includegraphics[width=\linewidth]{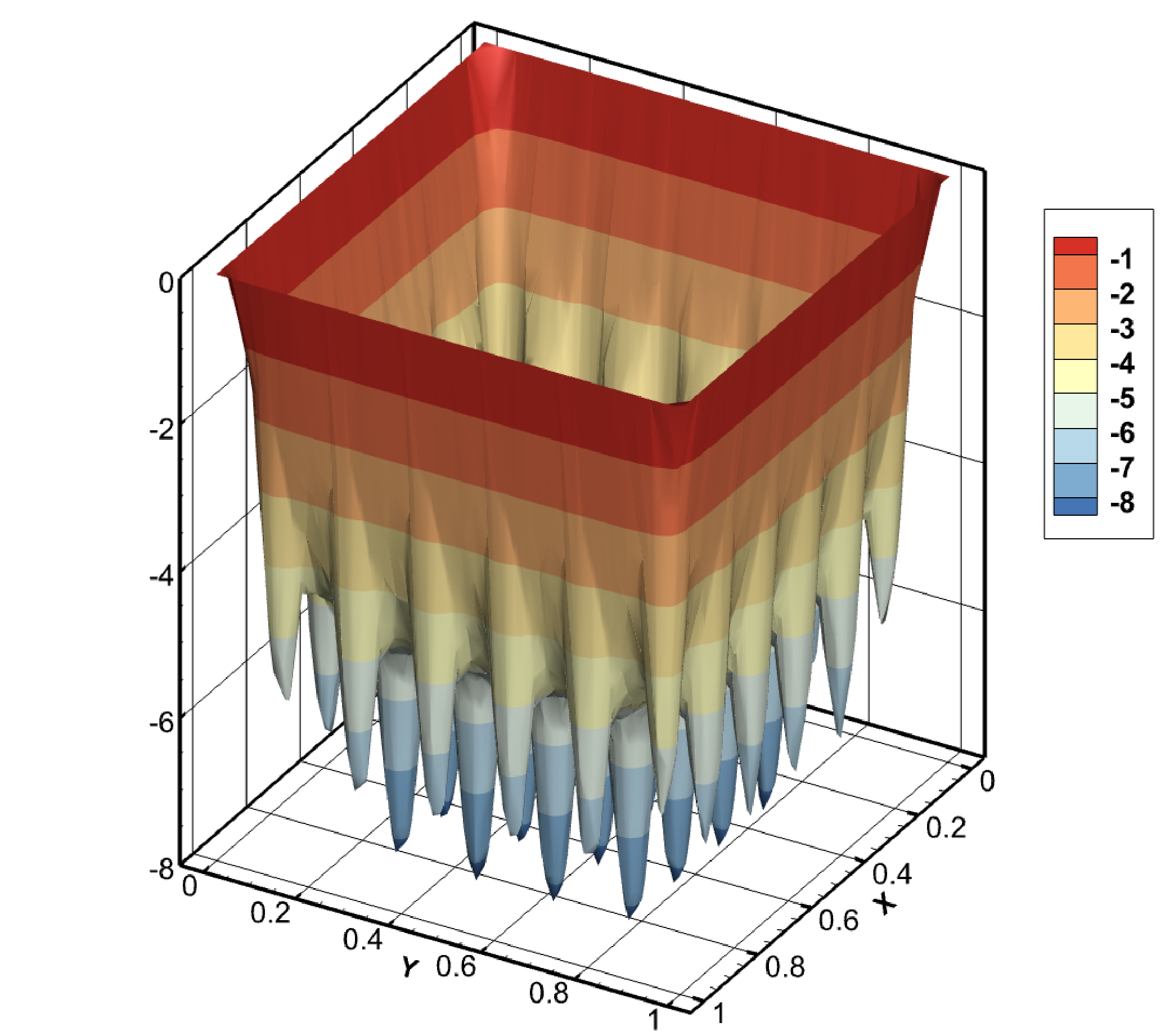}
				(d)
			\end{minipage}
			\caption{Case 1: Third displacement field component in $x_3\!\!=\!\!0.1 \mathrm{cm}$: (a) $u_3^{(0)}$; (b) $u_3^{(1,\epsilon)}$; (c) $u_3^{(2,\epsilon)}$; (d) $u_3^{\epsilon}$.}\label{f26}
		\end{figure}

\begin{figure}[!htb]
			\centering
			\begin{minipage}[c]{0.24\textwidth}
				\centering
				\includegraphics[width=\linewidth]{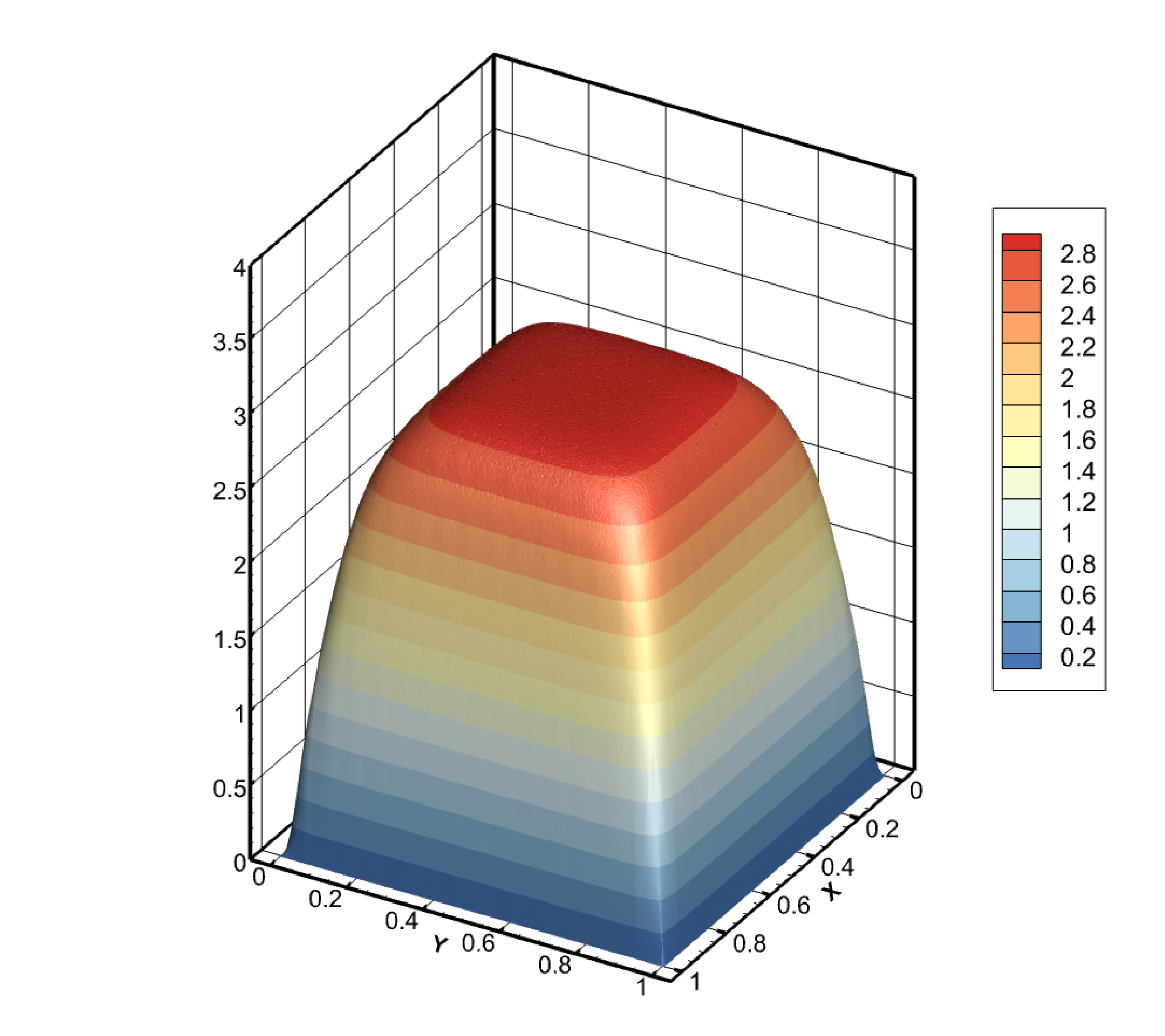}
				(a)
			\end{minipage}
			\hfill
			\begin{minipage}[c]{0.24\textwidth}
				\centering
				\includegraphics[width=\linewidth]{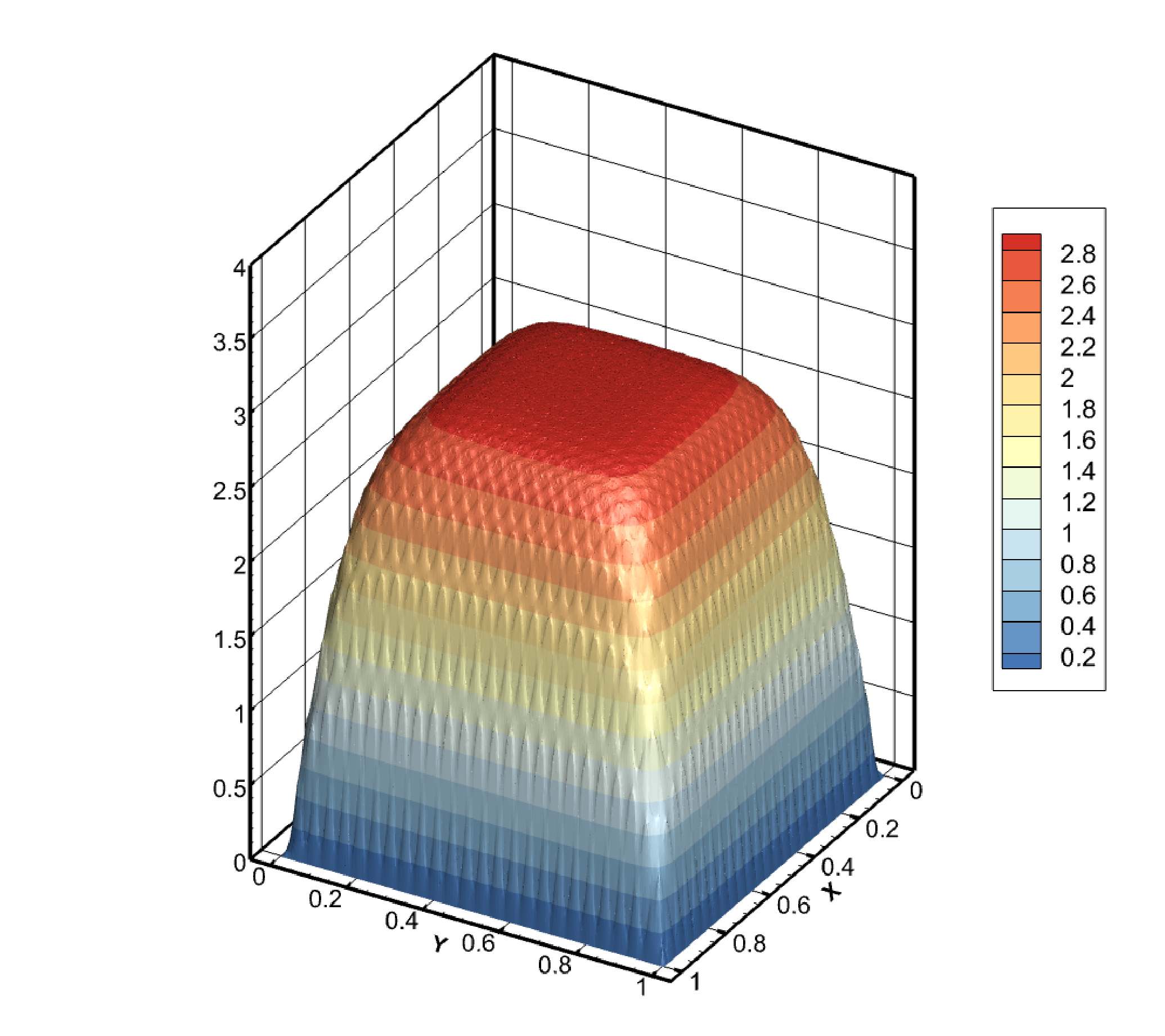}
				(b)
			\end{minipage}
			\hfill
			\begin{minipage}[c]{0.24\textwidth}
				\centering
				\includegraphics[width=\linewidth]{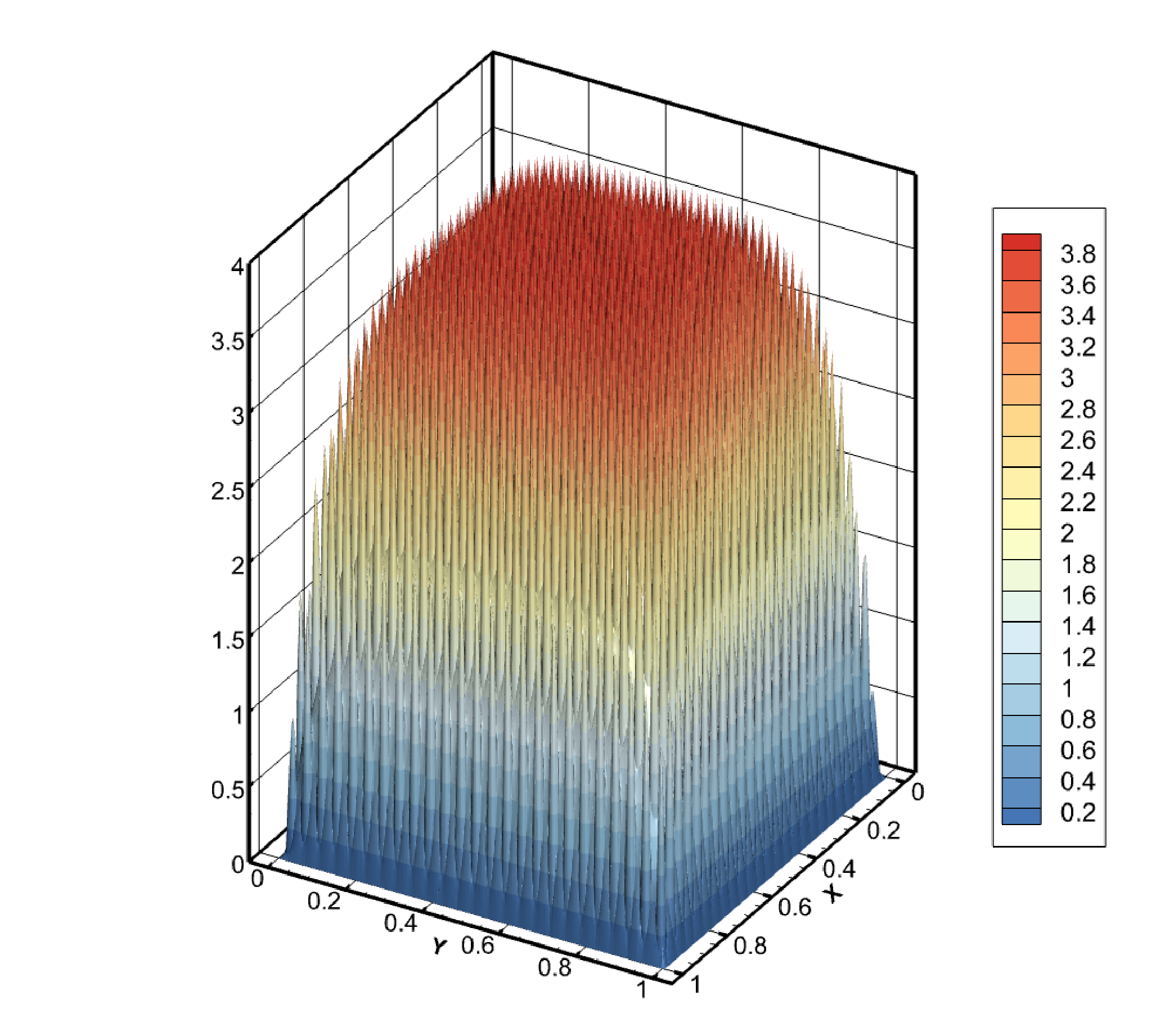}
				(c)
			\end{minipage}
			\hfill
			\begin{minipage}[c]{0.24\textwidth}
				\centering
				\includegraphics[width=\linewidth]{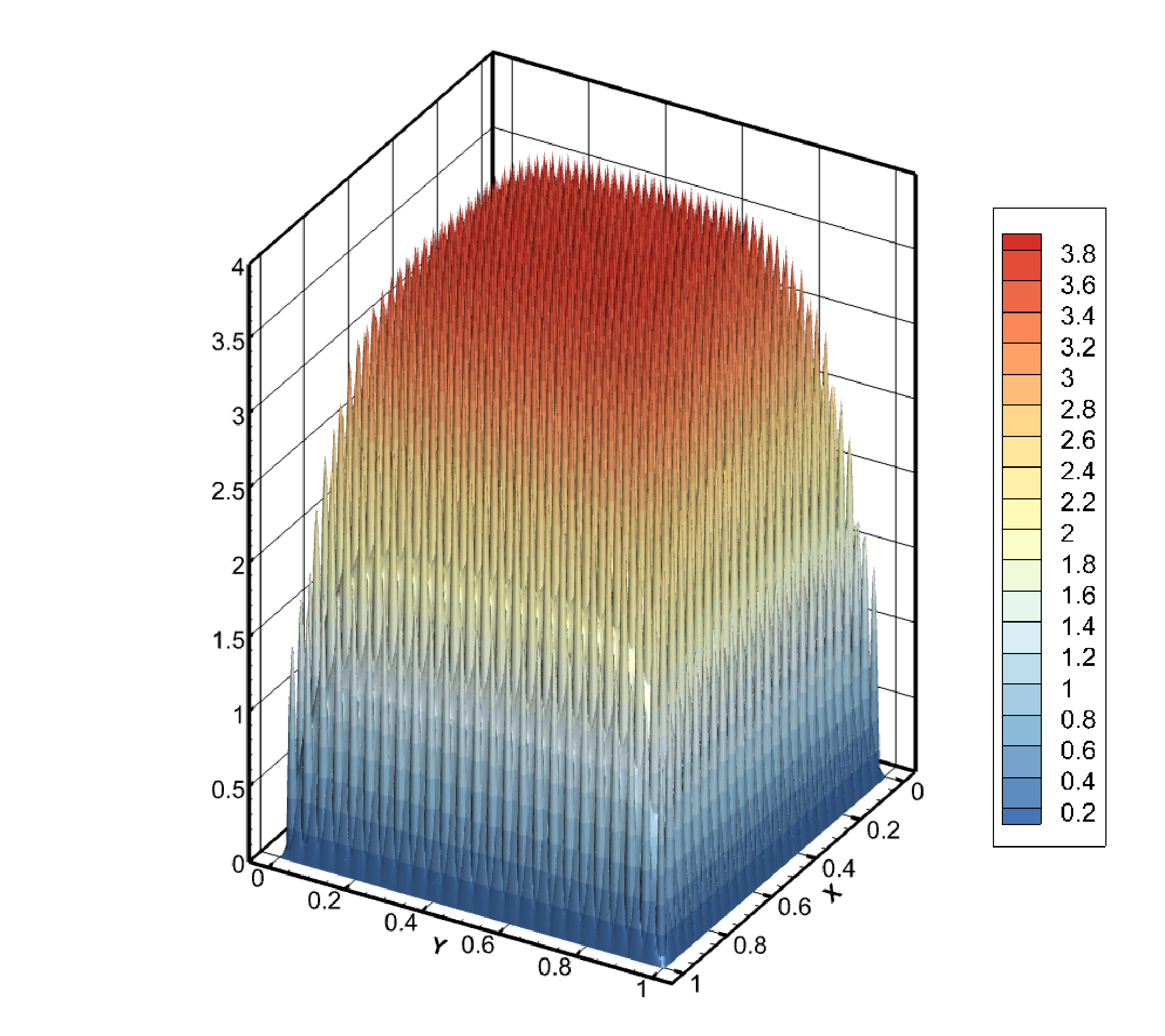}
				(d)
			\end{minipage}
			\caption{Case 2: Temperature increment field in $x_3=0.1 \mathrm{cm}$: (a) $T^{(0)}$; (b) $T^{(1,\epsilon)}$; (c) $T^{(2,\epsilon)}$; (d) $T^{\epsilon}$.}\label{f27}
		\end{figure}

\begin{figure}[!htb]
			\centering
			\begin{minipage}[c]{0.24\textwidth}
				\centering
				\includegraphics[width=\linewidth]{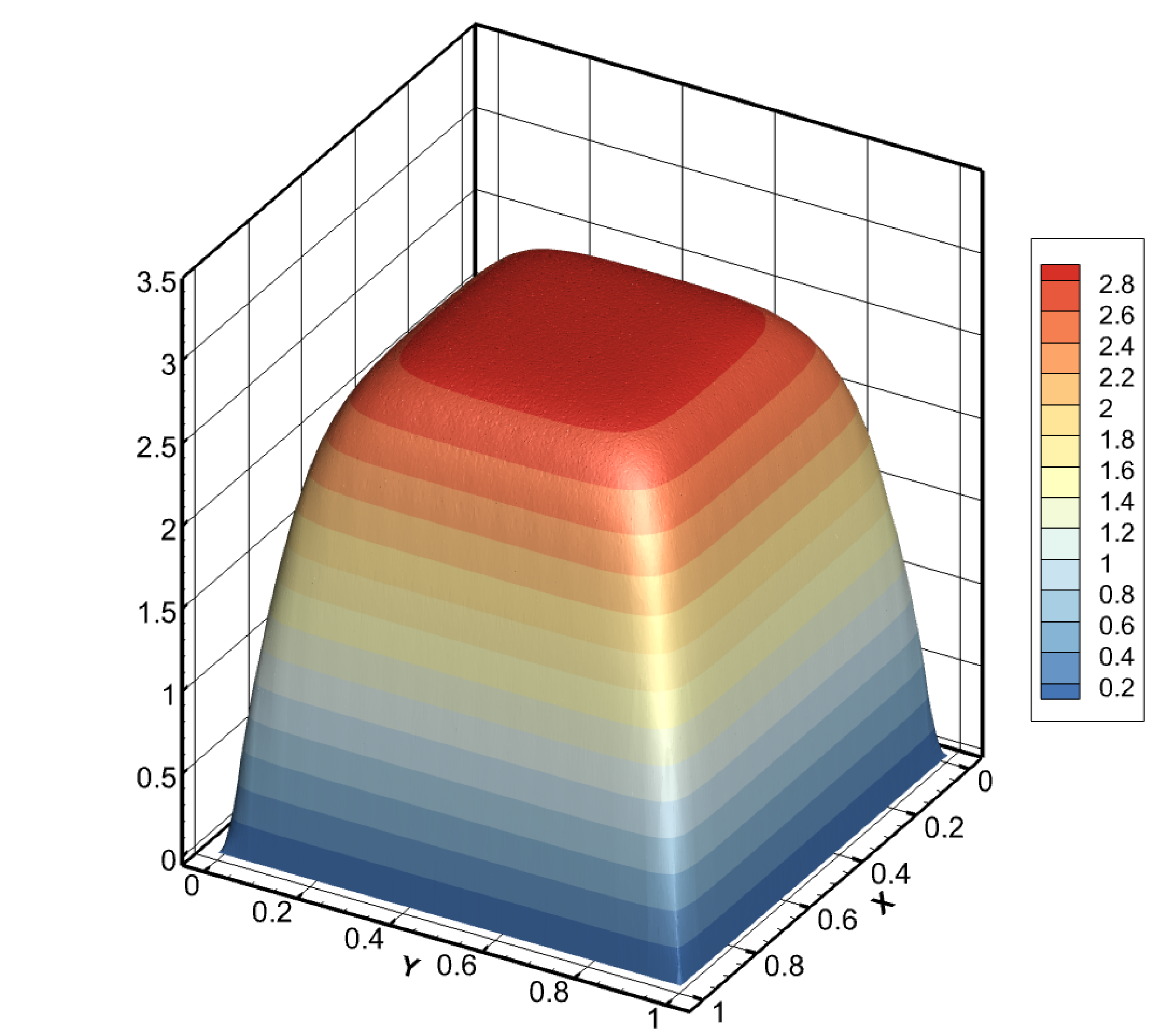}
				(a)
			\end{minipage}
			\hfill
			\begin{minipage}[c]{0.24\textwidth}
				\centering
				\includegraphics[width=\linewidth]{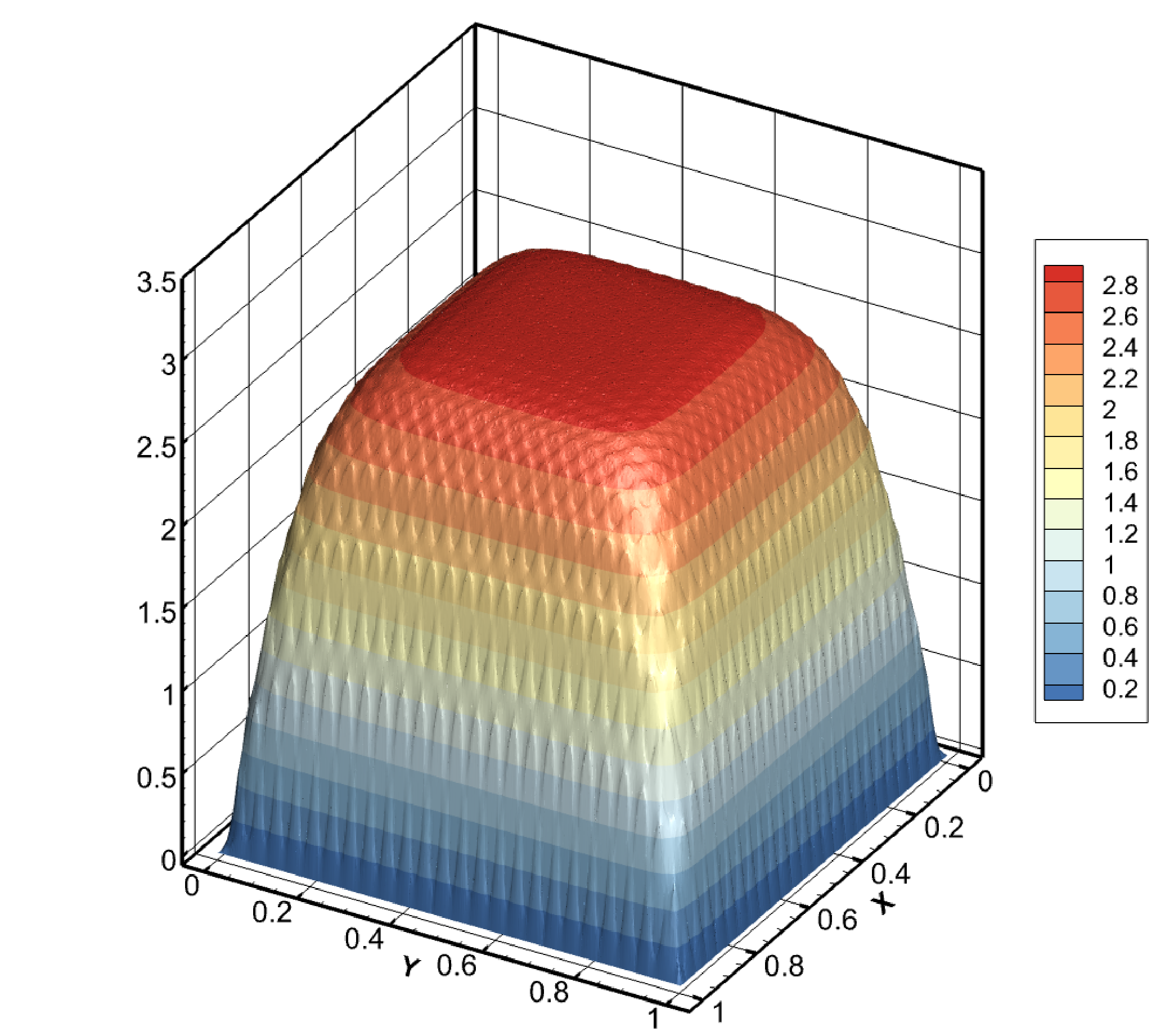}
				(b)
			\end{minipage}
			\hfill
			\begin{minipage}[c]{0.24\textwidth}
				\centering
				\includegraphics[width=\linewidth]{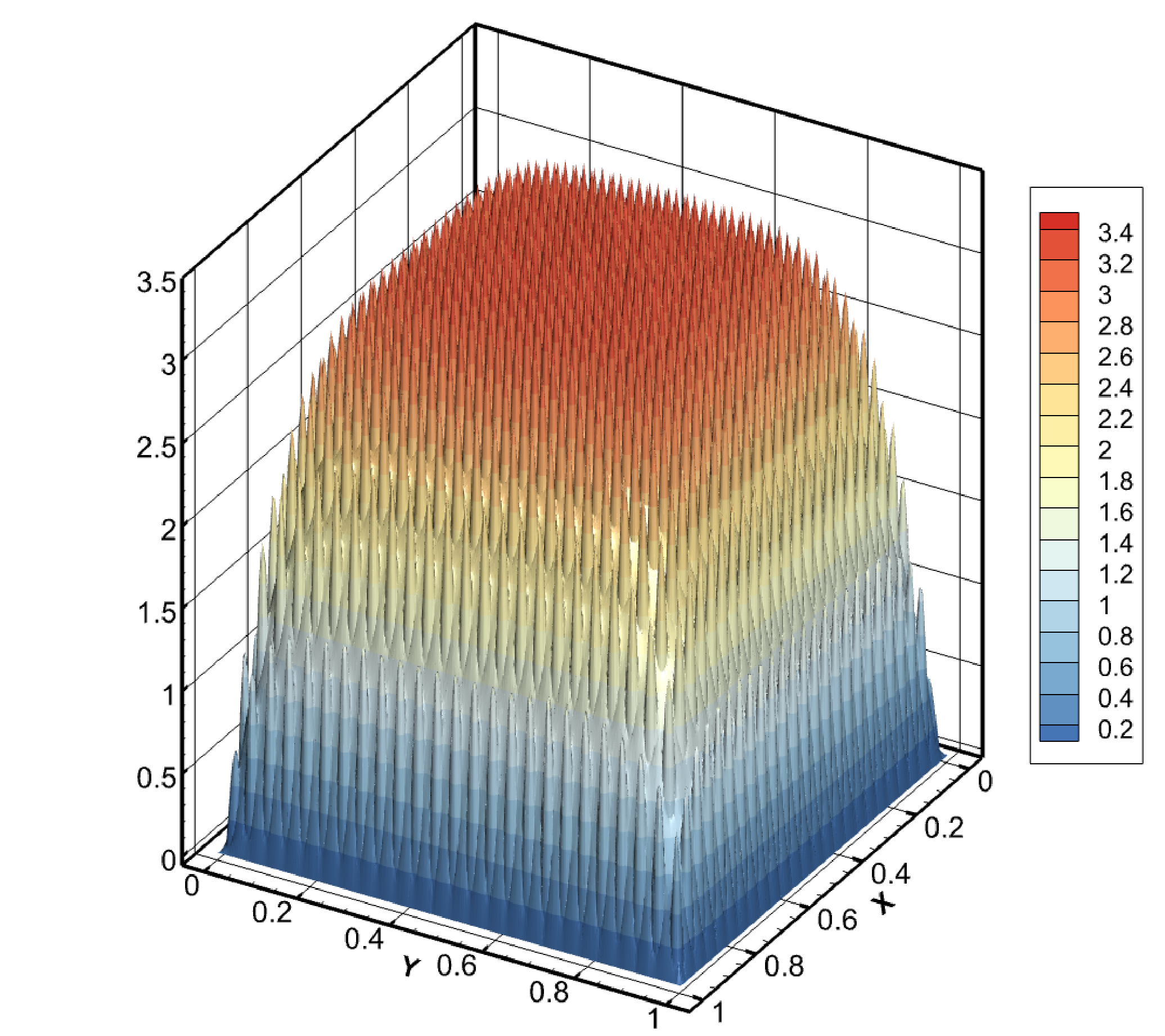}
				(c)
			\end{minipage}
			\hfill
			\begin{minipage}[c]{0.24\textwidth}
				\centering
				\includegraphics[width=\linewidth]{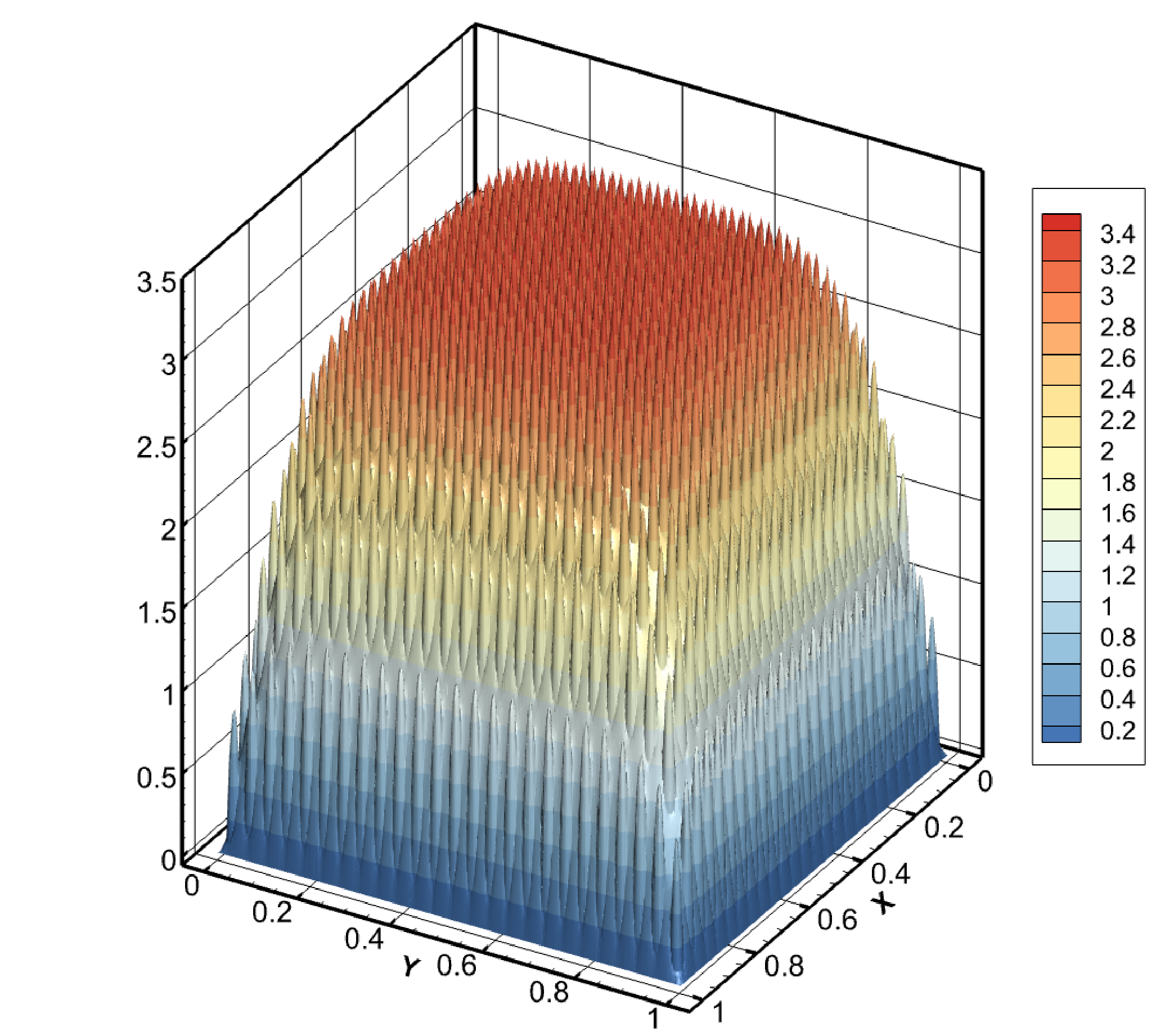}
				(d)
			\end{minipage}
			\caption{Case 2: Moisture field in $x_3=0.1 \mathrm{cm}$: (a) $c^{(0)}$; (b) $c^{(1,\epsilon)}$; (c) $c^{(2,\epsilon)}$; (d) $c^{\epsilon}$.}\label{f28}
		\end{figure}

\begin{figure}[!htb]
			\centering
			\begin{minipage}[c]{0.24\textwidth}
				\centering
				\includegraphics[width=\linewidth]{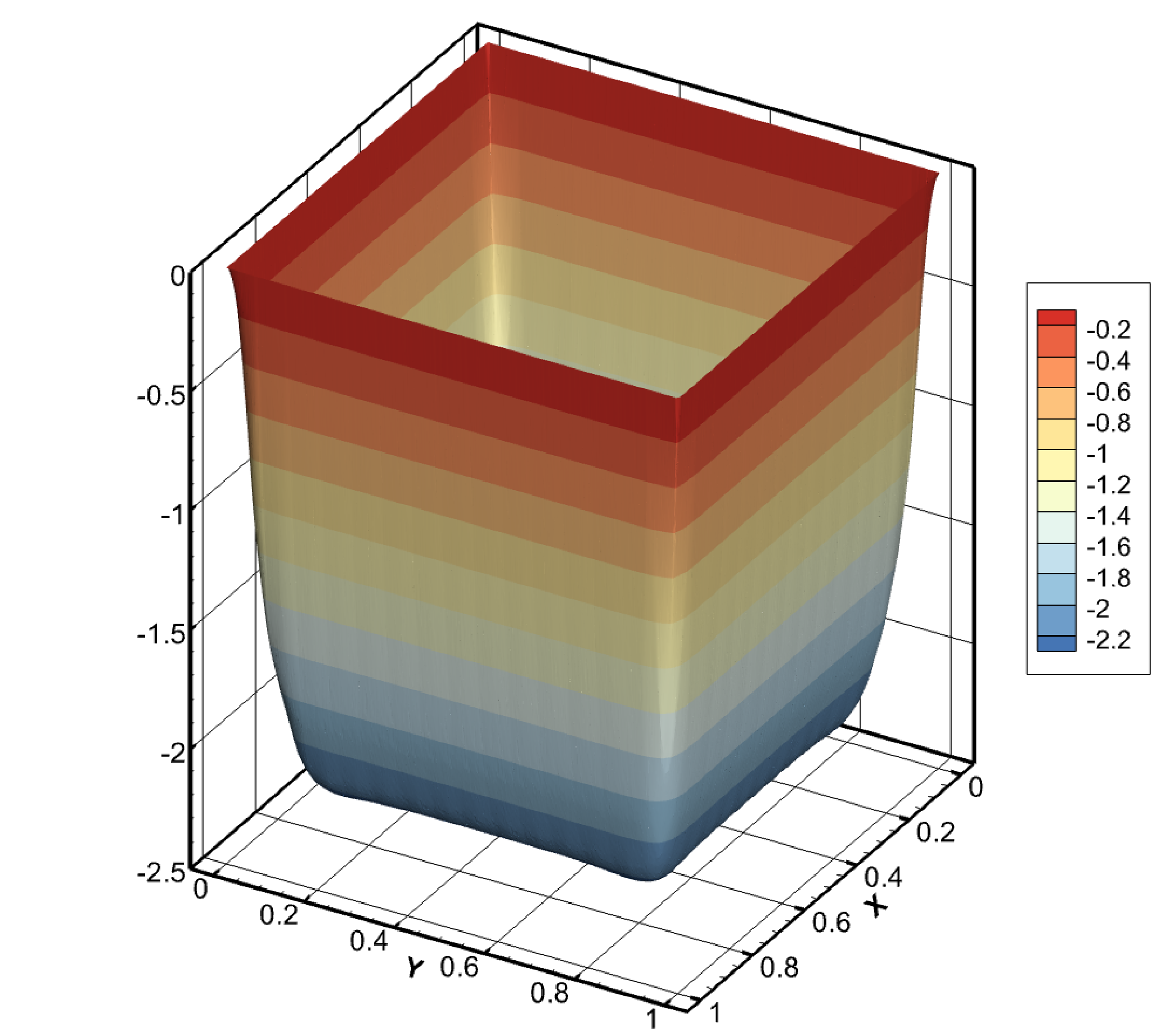}
				(a)
			\end{minipage}
			\hfill
			\begin{minipage}[c]{0.24\textwidth}
				\centering
				\includegraphics[width=\linewidth]{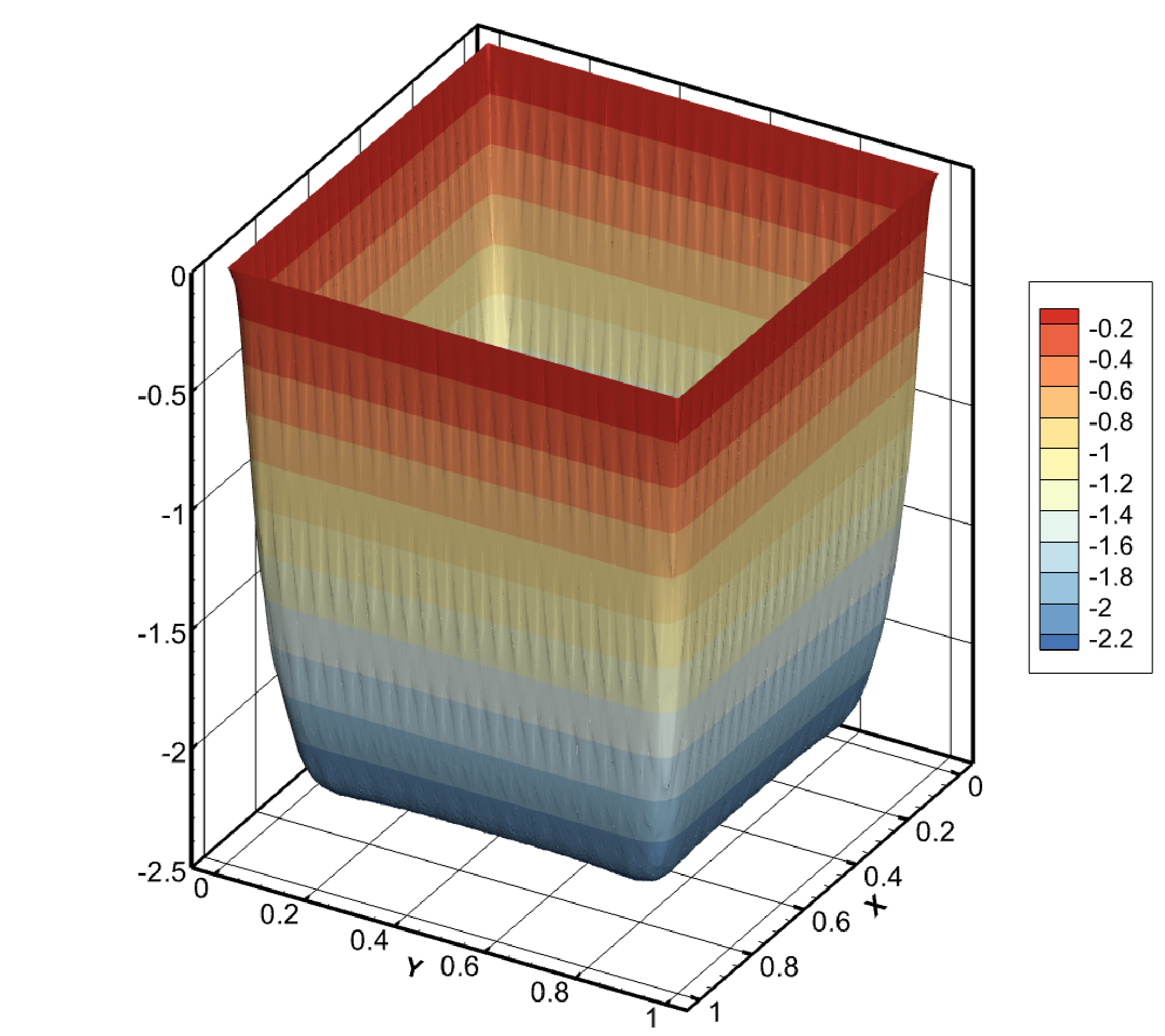}
				(b)
			\end{minipage}
			\hfill
			\begin{minipage}[c]{0.24\textwidth}
				\centering
				\includegraphics[width=\linewidth]{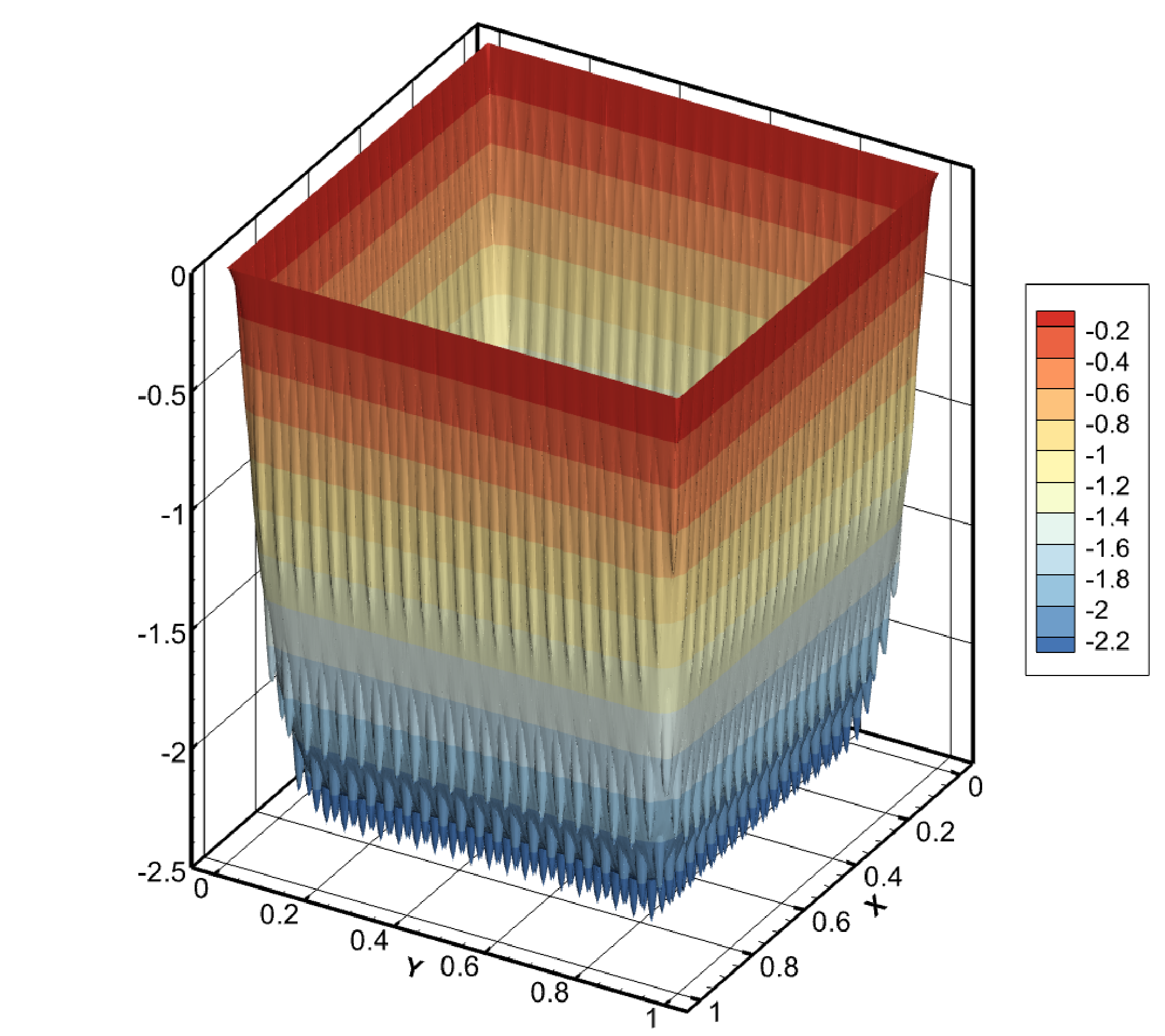}
				(c)
			\end{minipage}
			\hfill
			\begin{minipage}[c]{0.24\textwidth}
				\centering
				\includegraphics[width=\linewidth]{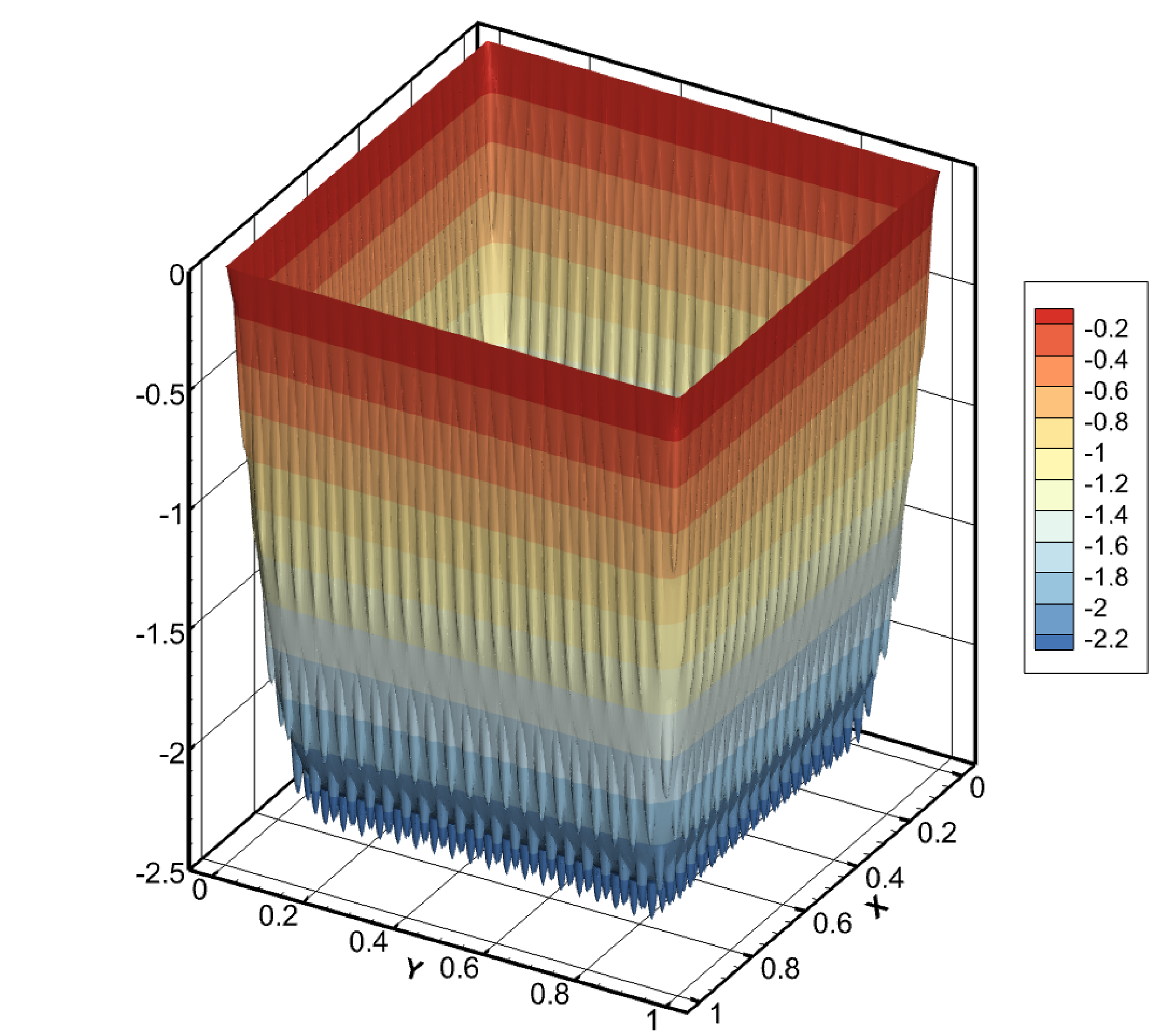}
				(d)
			\end{minipage}
			\caption{Case 2: Third displacement field component in $x_3\!\!=\!\!0.1 \mathrm{cm}$: (a) $u_3^{(0)}$; (b) $u_3^{(1,\epsilon)}$; (c) $u_3^{(2,\epsilon)}$; (d) $u_3^{\epsilon}$.}\label{f31}
		\end{figure}

\begin{figure}[!htb]
			\centering
			\begin{minipage}[c]{0.3\textwidth}
				\centering
				\includegraphics[width=\linewidth]{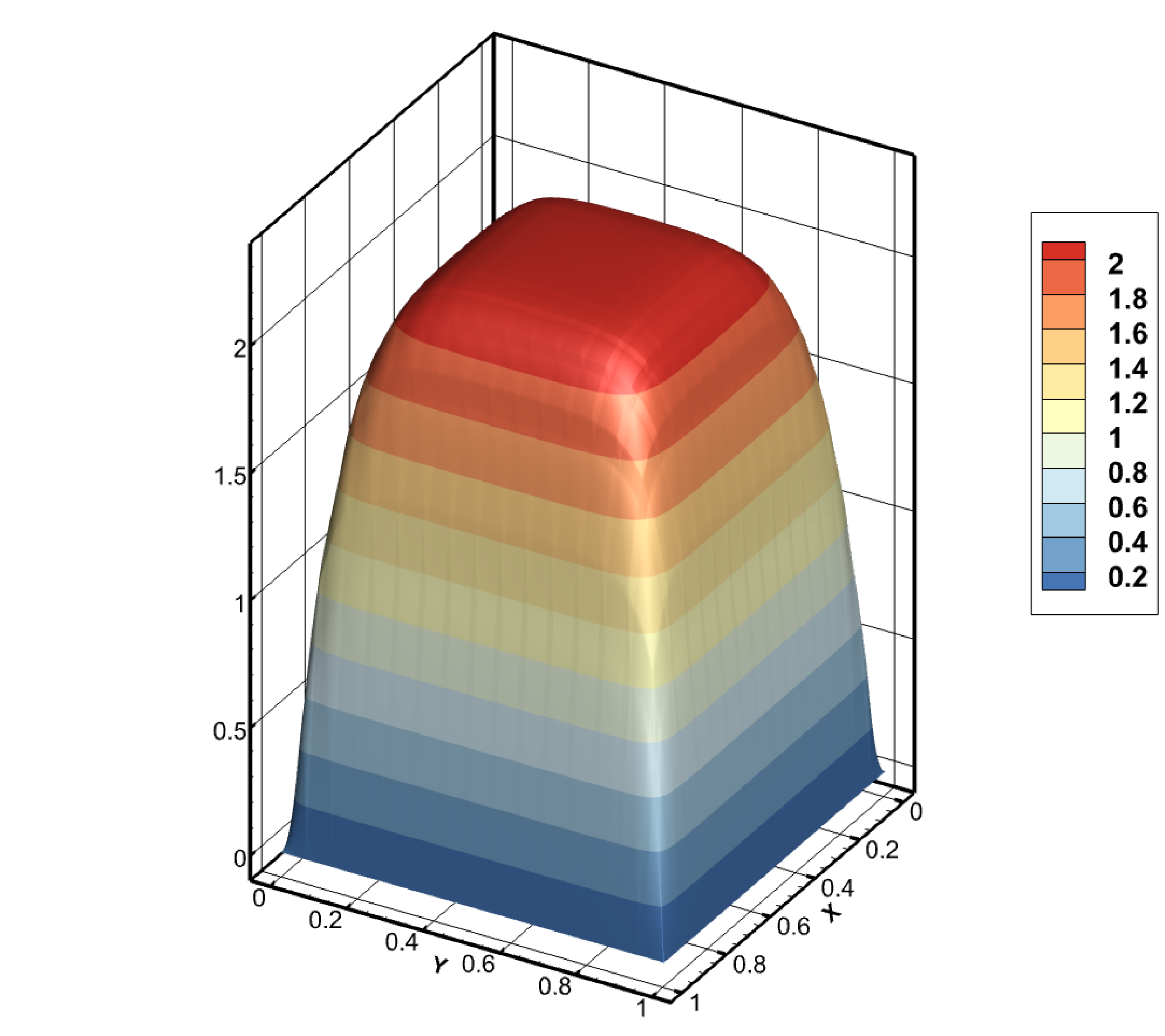}
				(a)
			\end{minipage}
			\hfill
			\begin{minipage}[c]{0.3\textwidth}
				\centering
				\includegraphics[width=\linewidth]{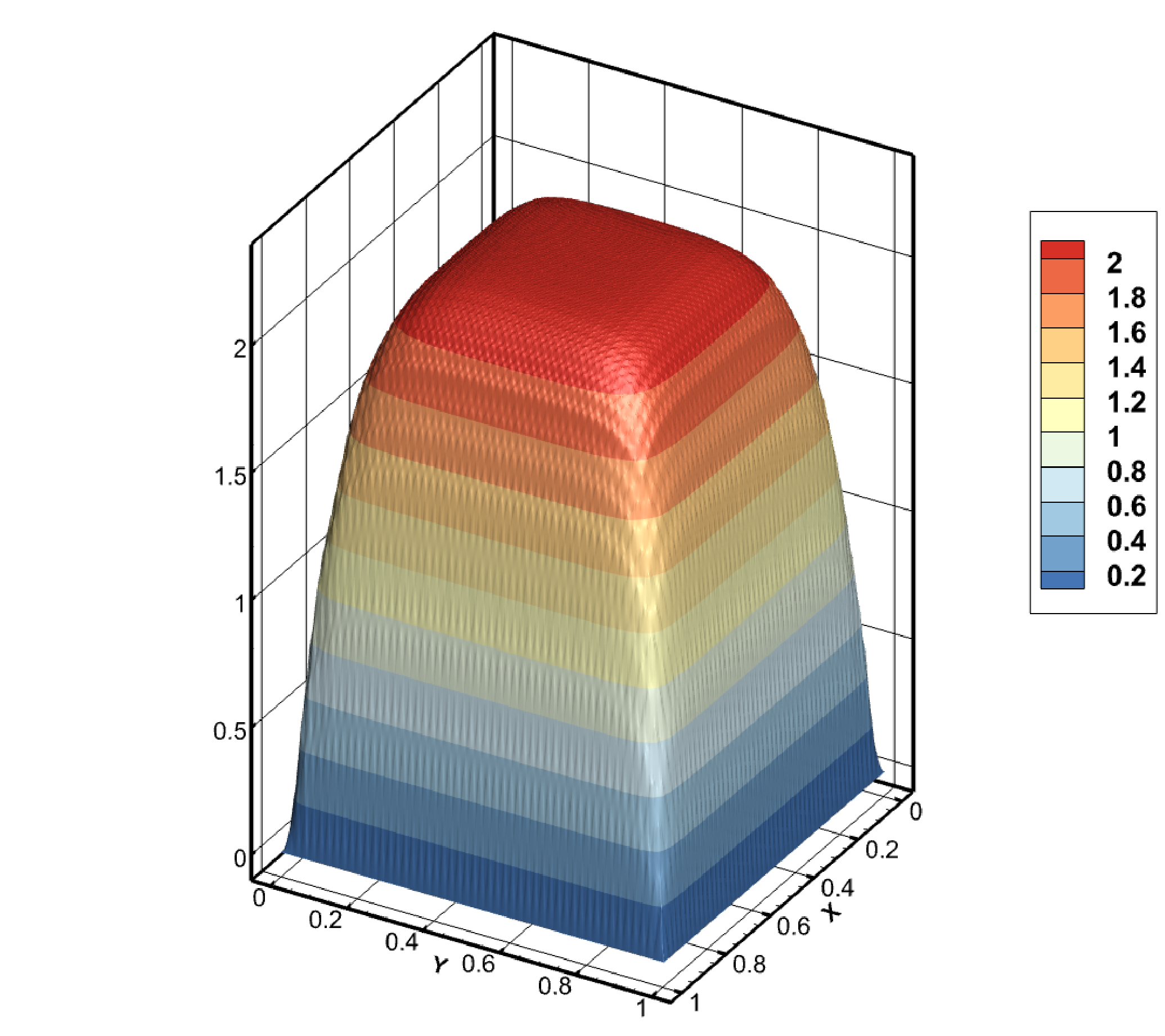}
				(b)
			\end{minipage}
			\hfill
			\begin{minipage}[c]{0.3\textwidth}
				\centering
				\includegraphics[width=\linewidth]{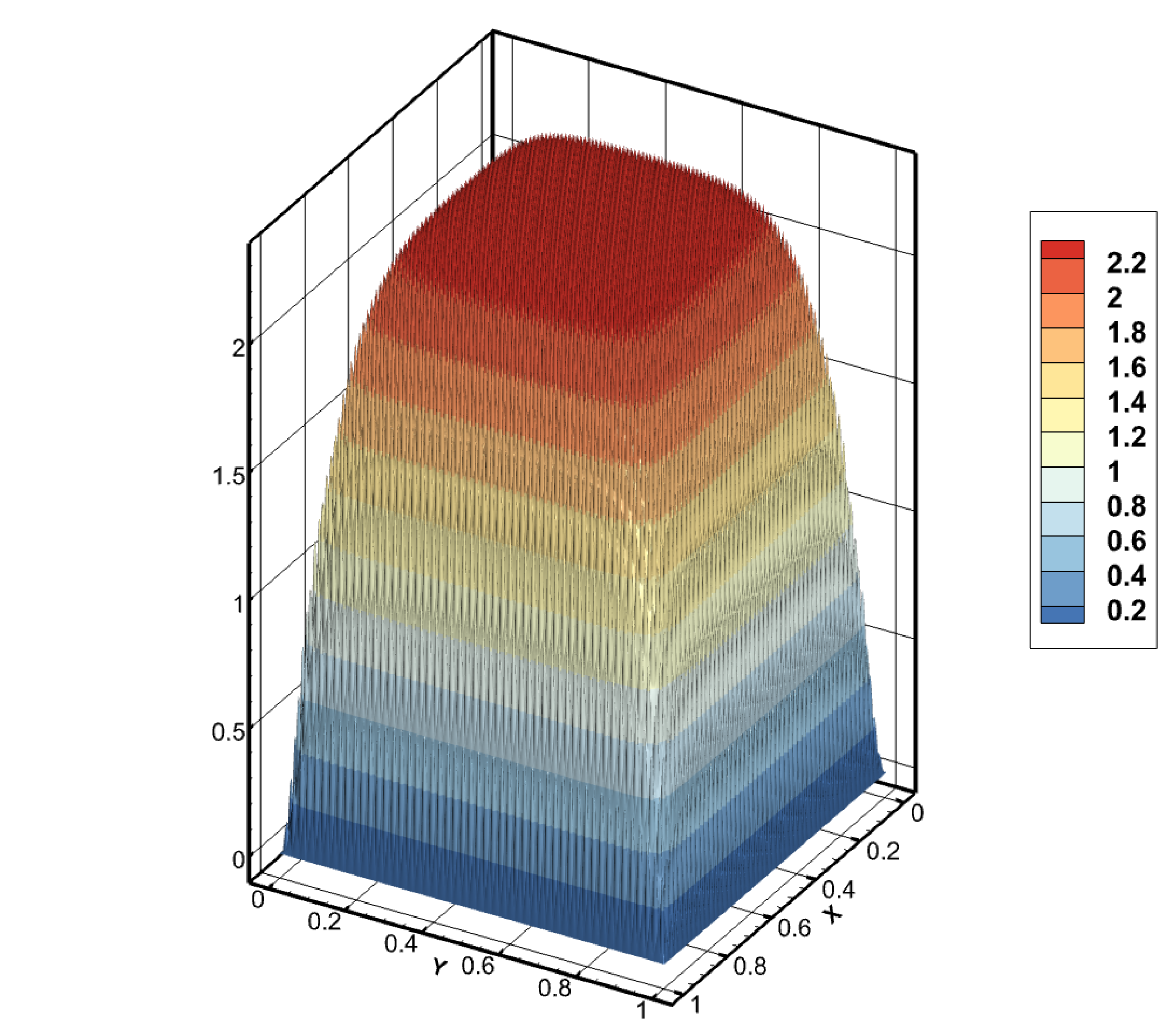}
				(c)
			\end{minipage}
			\caption{Case 3: Temperature increment field in $x_3=0.15 \mathrm{cm}$: (a) $T^{(0)}$; (b) $T^{(1,\epsilon)}$; (c) $T^{(2,\epsilon)}$.}\label{f32}
		\end{figure}

\begin{figure}[!htb]
			\centering
			\begin{minipage}[c]{0.3\textwidth}
				\centering
				\includegraphics[width=\linewidth]{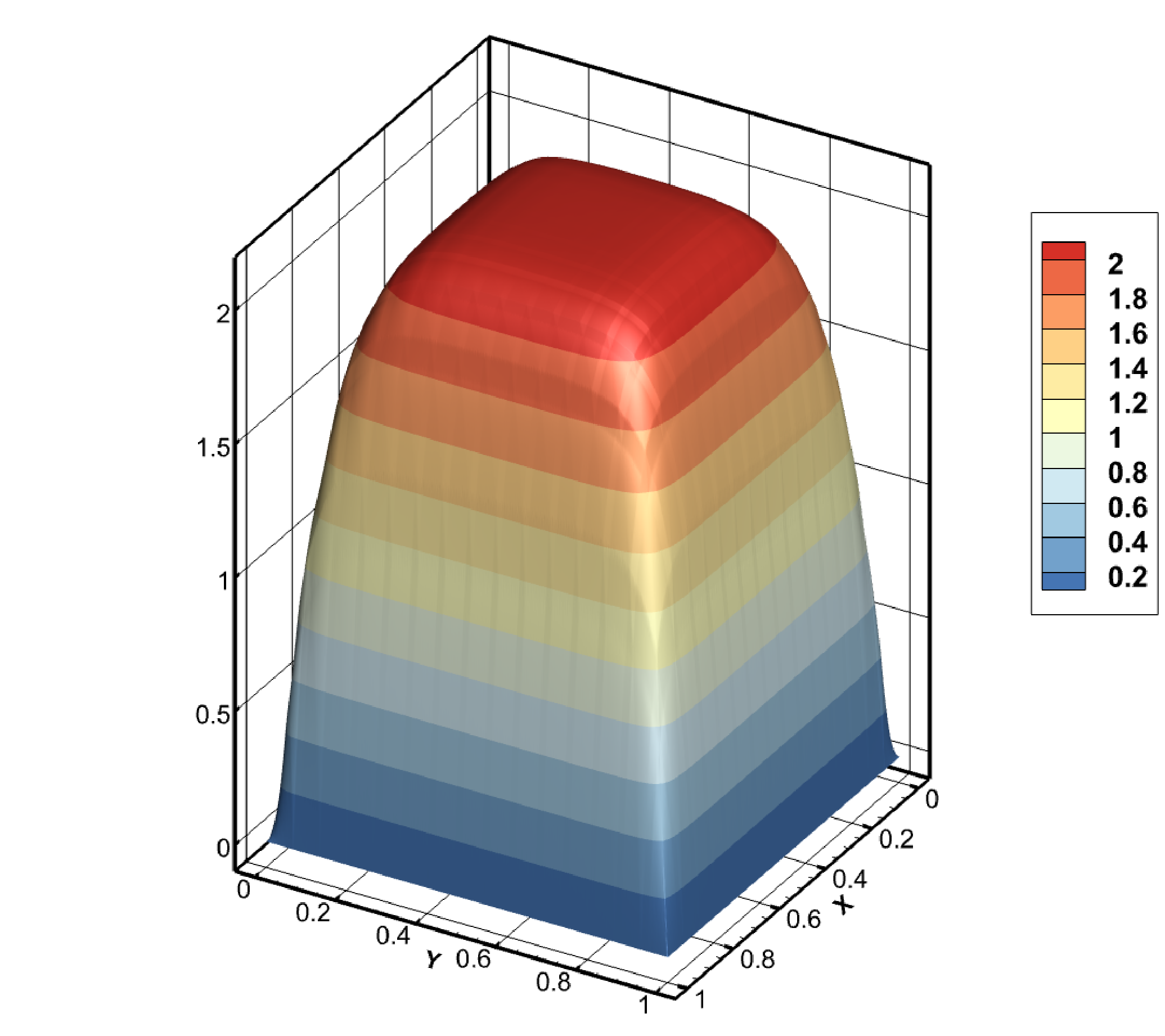}
				(a)
			\end{minipage}
			\hfill
			\begin{minipage}[c]{0.3\textwidth}
				\centering
				\includegraphics[width=\linewidth]{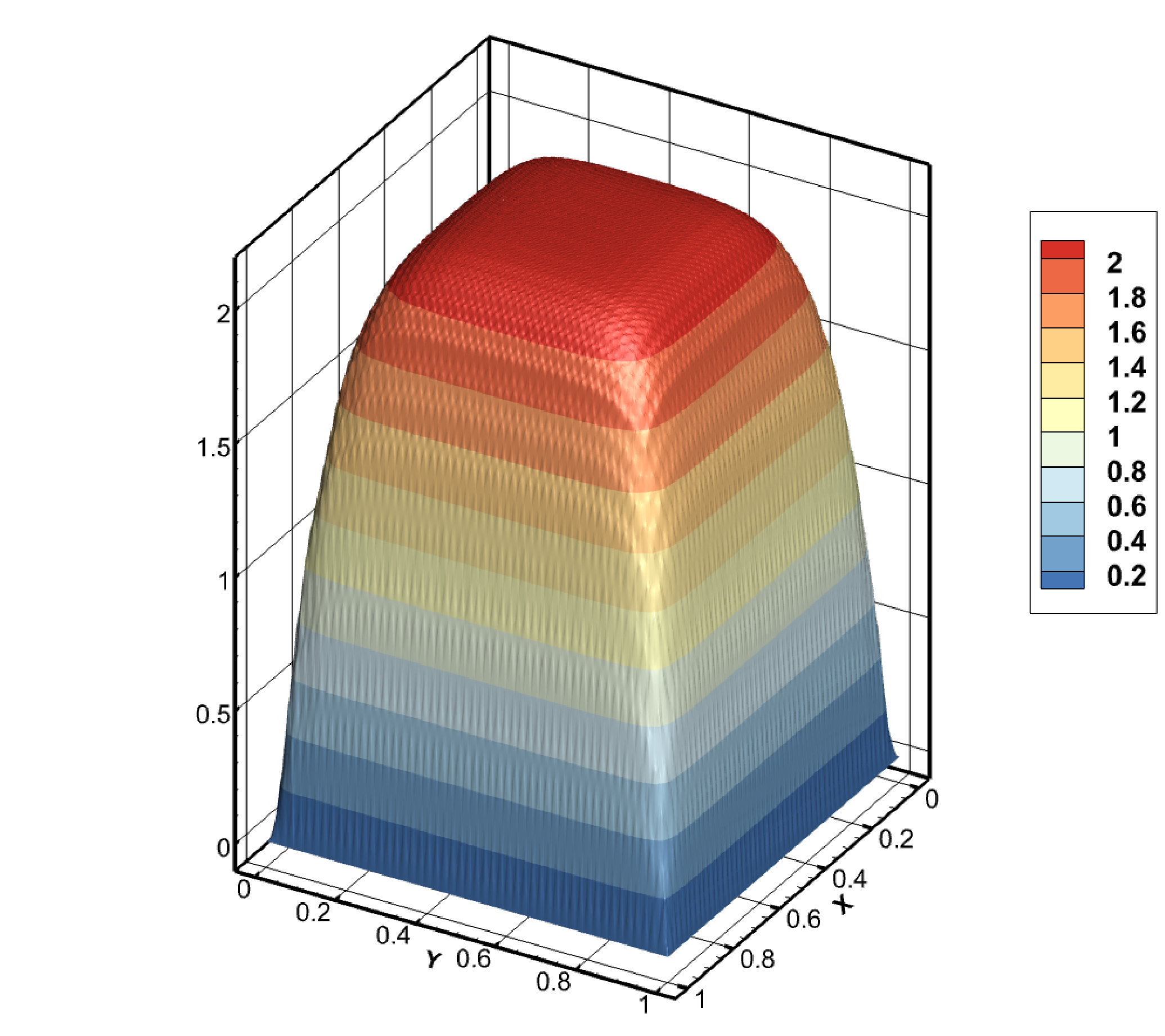}
				(b)
			\end{minipage}
			\hfill
			\begin{minipage}[c]{0.3\textwidth}
				\centering
				\includegraphics[width=\linewidth]{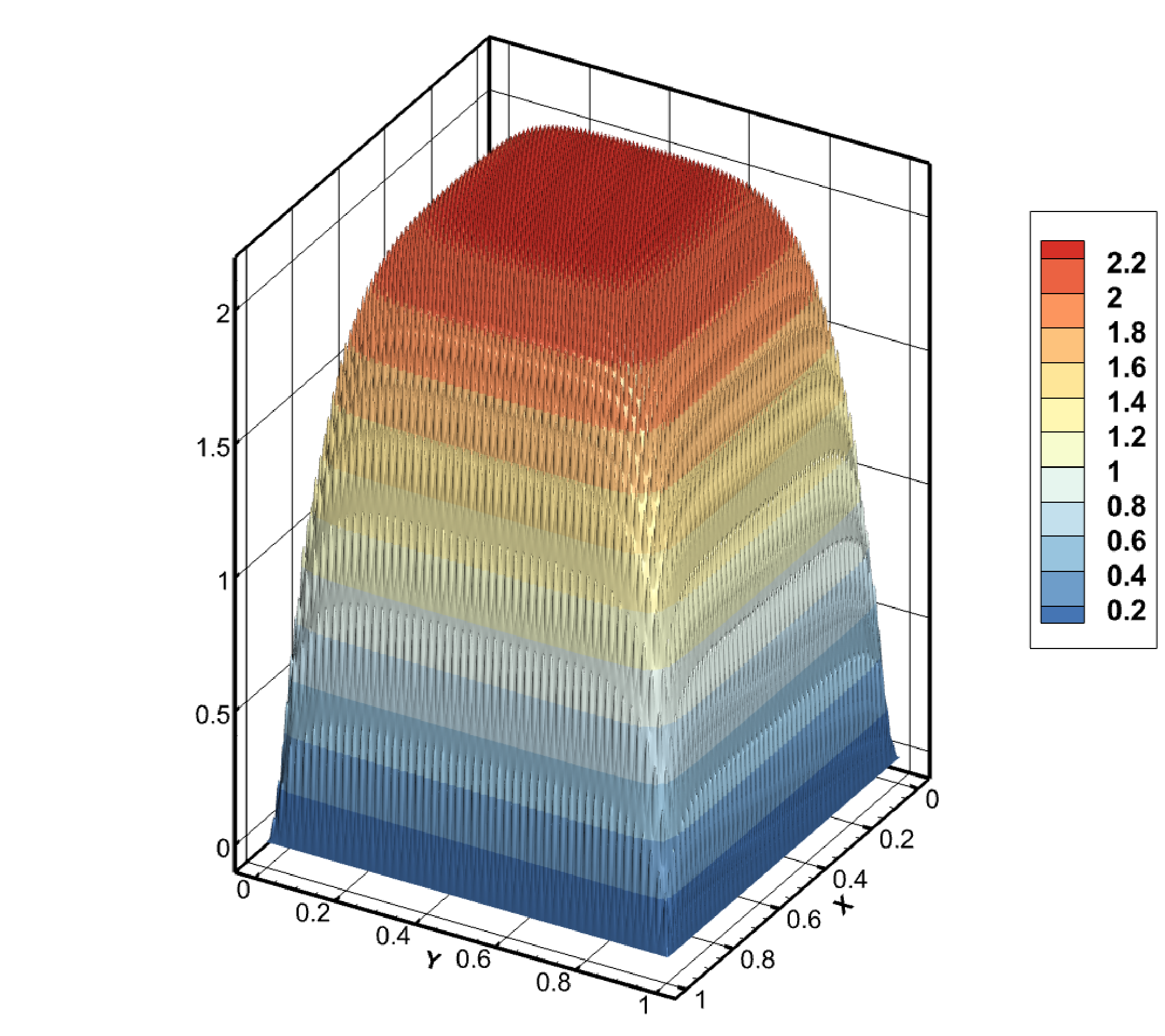}
				(c)
			\end{minipage}
			\caption{Case 3: Moisture field in $x_3=0.15 \mathrm{cm}$: (a) $c^{(0)}$; (b) $c^{(1,\epsilon)}$; (c) $c^{(2,\epsilon)}$.}\label{f33}
		\end{figure}

\begin{figure}[!htb]
			\centering
			\begin{minipage}[c]{0.3\textwidth}
				\centering
				\includegraphics[width=\linewidth]{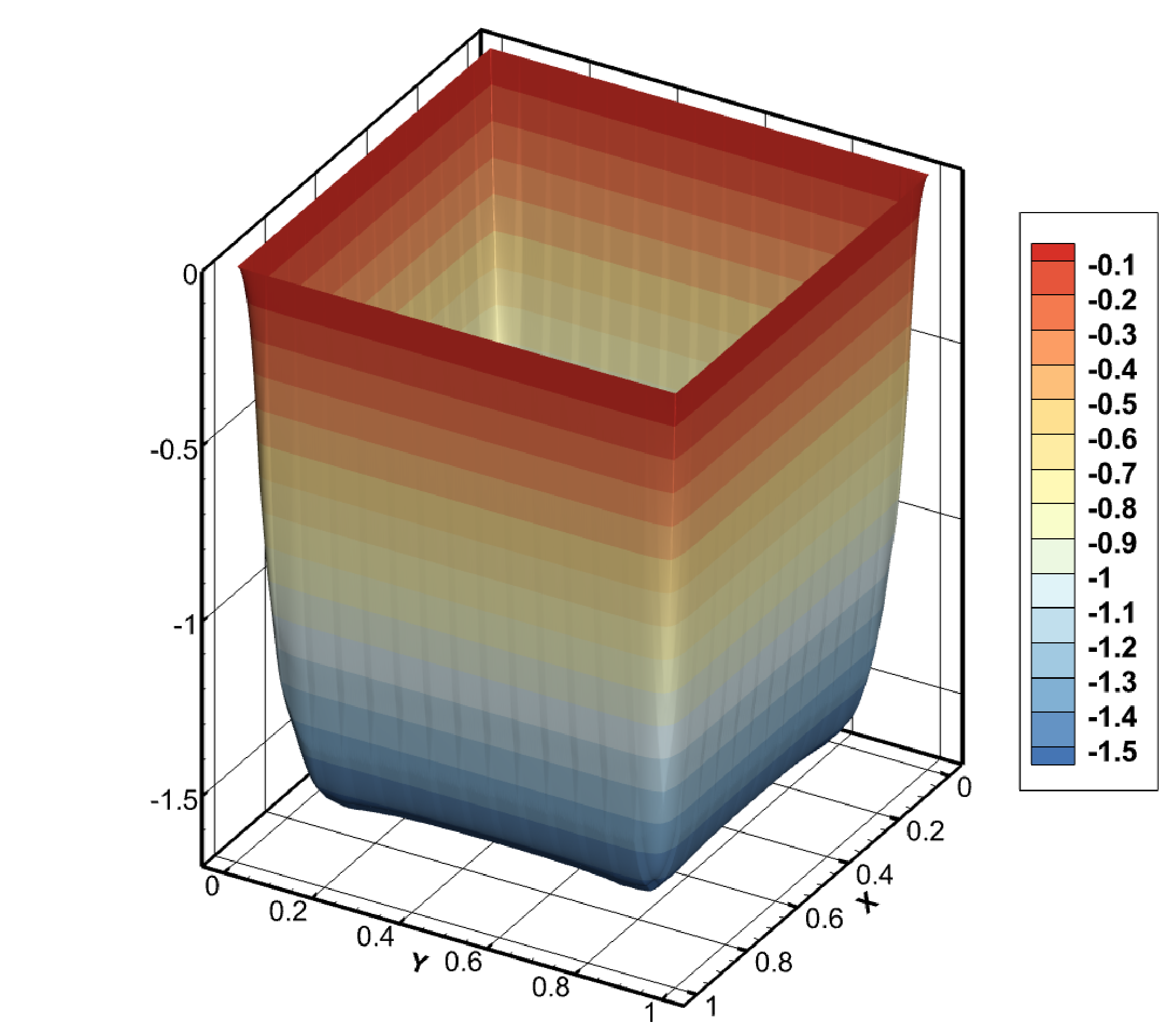}
				(a)
			\end{minipage}
			\hfill
			\begin{minipage}[c]{0.3\textwidth}
				\centering
				\includegraphics[width=\linewidth]{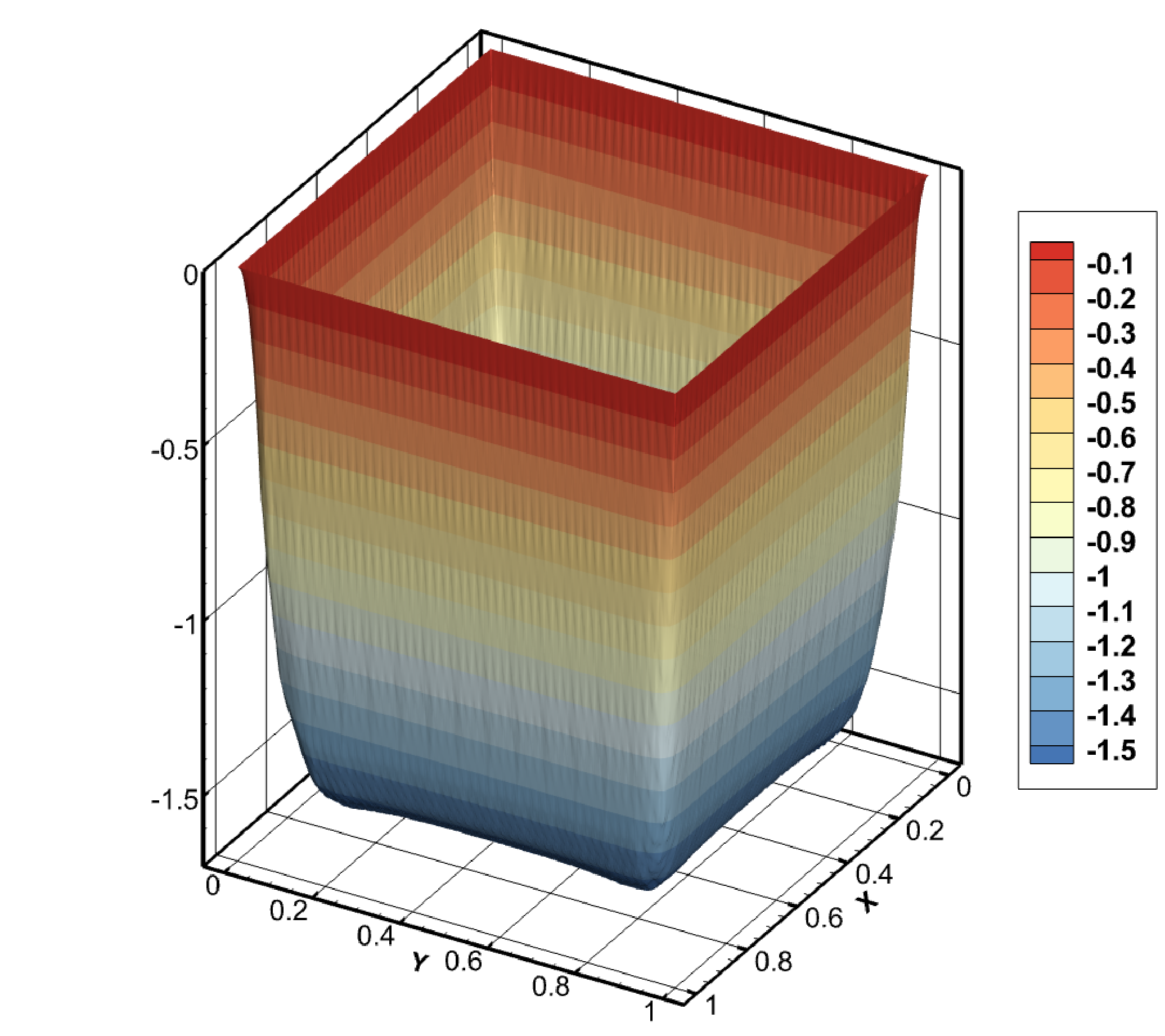}
				(b)
			\end{minipage}
			\hfill
			\begin{minipage}[c]{0.3\textwidth}
				\centering
				\includegraphics[width=\linewidth]{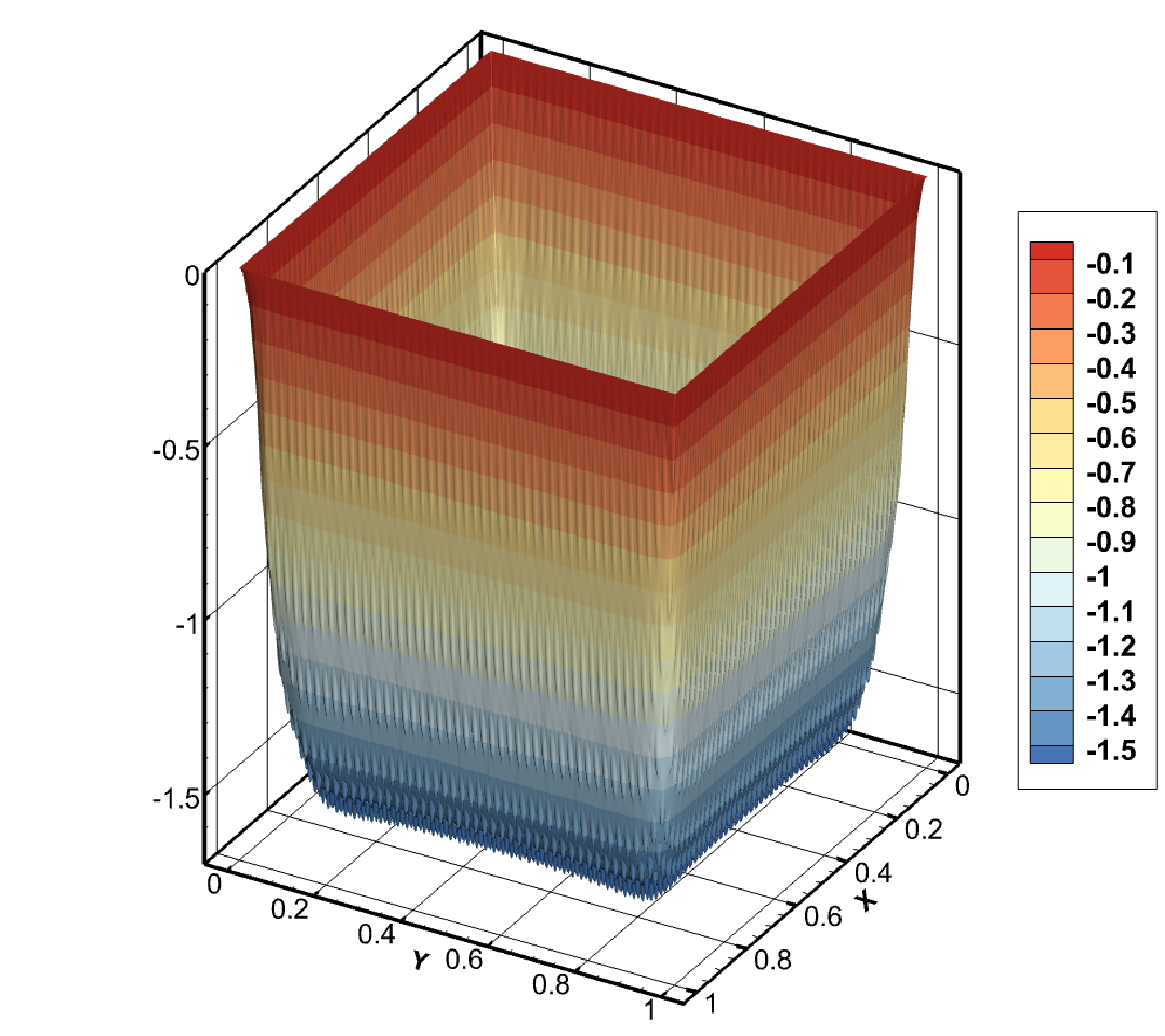}
				(c)
			\end{minipage}
			\caption{Case 3: Third displacement field component in $x_3=0.15 \mathrm{cm}$: (a) $u_3^{(0)}$; (b) $u_3^{(1,\epsilon)}$; (c) $u_3^{(2,\epsilon)}$.}\label{f36}
		\end{figure}

		After numerical computation, the computational results of Case 1 in Figs.\hspace{1mm}\ref{f22}-\ref{f26} with relative error data in Table\hspace{1mm}\ref{t12}, the results of Case 2 in Figs.\hspace{1mm}\ref{f27}-\ref{f31} with relative error data in Table\hspace{1mm}\ref{t13}, and the results of Case 3 in Figs.\hspace{1mm}\ref{f32}-\ref{f36}.
		\begin{table}[!htb]
			\centering
			\caption{The relative errors (Case 1: $6\!\times\!6\!\times\!2$ unit cells).}
			\label{t12}
			\begin{tabular}{cccccc}
				\hline
				\multicolumn{6}{c}{Temperature increment field} \\
				\hline
				$TerrorL^20$ & $TerrorL^21$ & $TerrorL^22$ & $TerrorH^10$ & $TerrorH^11$ & $TerrorH^12$ \\
				0.51371 & 0.51667 & 0.07944 & 0.96451 & 0.95587 & 0.16209 \\
				\hline
				\multicolumn{6}{c}{Moisture field} \\
				\hline
				$cerrorL^20$ & $cerrorL^21$ & $cerrorL^22$ & $cerrorH^10$ & $cerrorH^11$ & $cerrorH^12$ \\
				0.30177 & 0.30400 & 0.05219 & 0.89237 & 0.87514 & 0.15353 \\
				\hline
				\multicolumn{6}{c}{Displacement field} \\
				\hline
				$\bm{u}errorL^20$ & $\bm{u}errorL^21$ & $\bm{u}errorL^22$ & $\bm{u}errorH^10$ & $\bm{u}errorH^11$ & $\bm{u}errorH^12$ \\
				0.11890 & 0.10528 & 0.04306 & 0.57523 & 0.47770 & 0.16028 \\
				\hline
			\end{tabular}
		\end{table}
		\begin{table}[!htb]
			\centering
			\caption{The relative errors (Case 2: $25\!\times\!25\!\times\!5$ unit cells).}
			\label{t13}
			\begin{tabular}{cccccc}
				\hline
				\multicolumn{6}{c}{Temperature increment field} \\
				\hline
				$TerrorL^20$ & $TerrorL^21$ & $TerrorL^22$ & $TerrorH^10$ & $TerrorH^11$ & $TerrorH^12$ \\
				0.09642 & 0.09549 & 0.01907 & 0.84093 & 0.82644 & 0.16576 \\
				\hline
				\multicolumn{6}{c}{Moisture field} \\
				\hline
				$cerrorL^20$ & $cerrorL^21$ & $cerrorL^22$ & $cerrorH^10$ & $cerrorH^11$ & $cerrorH^12$ \\
				0.05091 & 0.04848 & 0.01214 & 0.62716 & 0.59557 & 0.13030 \\
				\hline
				\multicolumn{6}{c}{Displacement field} \\
				\hline
				$\bm{u}errorL^20$ & $\bm{u}errorL^21$ & $\bm{u}errorL^22$ & $\bm{u}errorH^10$ & $\bm{u}errorH^11$ & $\bm{u}errorH^12$ \\
				0.03806 & 0.02510 & 0.01904 & 0.43239 & 0.24303 & 0.14025 \\
				\hline
			\end{tabular}
		\end{table}
		
		From Figs.\hspace{1mm}\ref{f22}-\ref{f36}, Table\hspace{1mm}\ref{t12} and Table\hspace{1mm}\ref{t13}, it is evident that, consistent with previous results, only the HOMS solutions accurately capture the highly oscillating information at micro-scale, while LOMS solutions fail to achieve high-accuracy simulation for the temperature increment, moisture, and displacement fields. This demonstrates that the HOMS method remains effective in solving the H-T-M coupling problems \eqref{eq:2.1} for 3D quasi-periodic composite plates. Furthermore, this example confirms that the proposed approach not only substantially reduces computational cost but also delivers superior computational accuracy.
		
		\section{Conclusions}
		\label{sec:6}
		This study proposes a novel HOMS computational method for accurate and efficient simulation of complex material behaviors in quasi-periodic composite structures under H-T-M coupling effects, capable of accurately capturing multi-scale responses. This work makes three primary contributions: constructing the multi-scale asymptotic solutions to multi-field coupling governing equations incorporating higher-order correction terms using asymptotic homogenization; developing local and global error estimation for the multi-scale approximations; formulating an effective multi-scale numerical algorithm and deriving its corresponding convergent analysis. Numerical results demonstrate that the proposed HOMS computational method exhibits outstanding numerical accuracy and computational efficiency, regardless of whether the material parameters exhibit scale-separation characteristics or not. Crucially, only HOMS solutions accurately capture microscopic oscillatory information while meeting high-accuracy requirement of engineering simulation. Particularly for large-scale composite structures unattainable by classical finite element method, the proposed HOMS method provides an effective computational framework, which further validates the approach's theoretical soundness, algorithmic feasibility, and practical potential.
		
		To further advance the proposed computational framework, future research should address two critical aspects: First, the current work on H-T-M coupling problems in quasi-periodic composite structures is limited to the static linear case, necessitating extension of the present multi-scale computational framework to time-dependent nonlinear problems; second, given the growing engineering applications of composites with more than two-level spatial configurations, the presented two-scale computational methodology shall extend to three-scale spatial levels. Resolving these challenges will constitute crucial research directions.
		\section*{Acknowledgments}
		This research was supported by the National Natural Science Foundation of China (Nos.\hspace{1mm}12471387 and 12401523), Young Talent Fund of Association for Science and Technology in Shaanxi, China (No.\hspace{1mm}20220506), Xidian University Specially Funded Project for Interdisciplinary Exploration (No.\hspace{1mm}TZJH2024008), Fundamental Research Funds for the Central Universities (No.\hspace{1mm}QTZX25082), Innovation Capability Support Program of Shaanxi Province (No.\hspace{1mm}2024RS-CXTD-88), and also supported by the Center for high performance computing of Xidian University.
		{\section*{Appendix A. HOMS solutions for scale-separated parameters.}
			For quasi-periodic composite structures with scale-separation parameters, the multi-scale problem \eqref{eq:2.1} admits the following HOMS asymptotic solutions
			\begin{equation*}
				\begin{aligned}
					T^{\epsilon}(\mathbf{x}) & \approx T^{(0)}(\mathbf{x})
					+ \epsilon \mathcal{H}_{\alpha_1}(\mathbf{y}) \frac{\partial T^{(0)}(\mathbf{x})}{\partial x_{\alpha_1}}\\
					& + \epsilon^2 \Bigl( \mathcal{H}_{\alpha_1\alpha_2}(\mathbf{y}) \frac{\partial^2 T^{(0)}(\mathbf{x})}{\partial x_{\alpha_1} \partial x_{\alpha_2}}
					+ \frac{1}{\omega(\mathbf{x})} \frac{\partial \omega(\mathbf{x})}{\partial x_{\alpha_2}} \mathcal{\tilde{R}}_{\alpha_1\alpha_2}(\mathbf{y}) \frac{\partial T^{(0)}(\mathbf{x})}{\partial x_{\alpha_1}} \Bigr),
				\end{aligned}
			\end{equation*}		
			\begin{equation*}
				\begin{aligned}
					c^{\epsilon}(\mathbf{x}) & \approx c^{(0)}(\mathbf{x})
					+ \epsilon \mathcal{L}_{\alpha_1}(\mathbf{y}) \frac{\partial c^{(0)}(\mathbf{x})}{\partial x_{\alpha_1}}\\
					& + \epsilon^2 \Bigl( \mathcal{L}_{\alpha_1\alpha_2}(\mathbf{y}) \frac{\partial^2 c^{(0)}(\mathbf{x})}{\partial x_{\alpha_1} \partial x_{\alpha_2}}
					+ \frac{1}{\omega(\mathbf{x})} \frac{\partial \omega(\mathbf{x})}{\partial x_{\alpha_2}} \mathcal{\tilde{S}}_{\alpha_1\alpha_2}(\mathbf{y}) \frac{\partial c^{(0)}(\mathbf{x})}{\partial x_{\alpha_1}} \Bigr),
				\end{aligned}
			\end{equation*}		
			\begin{equation*}
				\begin{aligned}
					u_i^{\epsilon}(\mathbf{x}) &\approx u_i^{(0)}(\mathbf{x})
					+ \epsilon \Bigl( \mathcal{X}_{ih}^{\alpha_1}(\mathbf{y}) \frac{\partial u_h^{(0)}(\mathbf{x})}{\partial x_{\alpha_1}}
					- \omega(\mathbf{x}) \mathcal{\tilde{M}}_i(\mathbf{y}) T^{(0)}(\mathbf{x})
					- \omega(\mathbf{x}) \mathcal{\tilde{N}}_i(\mathbf{y}) c^{(0)}(\mathbf{x}) \Bigr) \\
					&+ \epsilon^2 \Bigl( \mathcal{P}_{ih}^{\alpha_1\alpha_2}(\mathbf{y}) \frac{\partial^2 u_h^{(0)}(\mathbf{x})}{\partial x_{\alpha_1} \partial x_{\alpha_2}}
					+ \frac{1}{\omega(\mathbf{x})} \frac{\partial \omega(\mathbf{x})}{\partial x_{\alpha_2}} \mathcal{\tilde{Q}}_{ih}^{\alpha_1\alpha_2}(\mathbf{y}) \frac{\partial u_h^{(0)}(\mathbf{x})}{\partial x_{\alpha_1}} + \frac{\partial \omega(\mathbf{x})}{\partial x_{\alpha_1}} \mathcal{\tilde{W}}_i^{\alpha_1}(\mathbf{y}) T^{(0)}(\mathbf{x})  \\
					&+ \omega(\mathbf{x}) \mathcal{\tilde{Z}}_i^{\alpha_1}(\mathbf{y}) \frac{\partial T^{(0)}(\mathbf{x})}{\partial x_{\alpha_1}} + \frac{\partial \omega(\mathbf{x})}{\partial x_{\alpha_1}} \mathcal{\tilde{F}}_i^{\alpha_1}(\mathbf{y}) c^{(0)}(\mathbf{x})
					+ \omega(\mathbf{x}) \mathcal{\tilde{G}}_i^{\alpha_1}(\mathbf{y}) \frac{\partial c^{(0)}(\mathbf{x})}{\partial x_{\alpha_1}} \Bigr),
				\end{aligned}
			\end{equation*}	
			where first-order cell functions satisfy the following new unit cell problems
			\begin{equation*}
				\begin{cases}
					\begin{aligned}
						& \frac{\partial}{\partial y_i}\Bigl(k_{ij}^*(\mathbf{y})\frac{\partial \mathcal{H}_{\alpha_1}}{\partial y_j}\Bigr)=-\frac{\partial k_{i\alpha_1}^*(\mathbf{y})}{\partial y_i},\quad \mathbf{y}\in Y, \\
						& \mathcal{H}_{\alpha_1}(\mathbf{y})=0, \quad \mathbf{y}\in\partial Y.
					\end{aligned}
				\end{cases}
			\end{equation*}		
			\begin{equation*}
				\begin{cases}
					\begin{aligned}
						& \frac{\partial}{\partial y_i}\Bigl(g_{ij}^*(\mathbf{y})\frac{\partial \mathcal{L}_{\alpha_1}}{\partial y_j}\Bigr)=-\frac{\partial g_{i\alpha_1}^*(\mathbf{y})}{\partial y_i}, \quad\mathbf{y}\in Y, \\
						& \mathcal{L}_{\alpha_1}(\mathbf{y})=0,\quad \mathbf{y}\in\partial Y.
					\end{aligned}
				\end{cases}
			\end{equation*}		
			\begin{equation*}
				\begin{cases}
					\begin{aligned}
						& \frac{\partial}{\partial y_j}\Bigl(D_{ijkl}^*(\mathbf{y})\frac{\partial \mathcal{X}_{kh}^{\alpha_1}}{\partial y_l}\Bigr)=-\frac{\partial D_{ijh\alpha_1}^*(\mathbf{y})}{\partial y_j},\quad \mathbf{y}\in Y, \\
						& \mathcal{X}_{kh}^{\alpha_1}(\mathbf{y})=0, \quad \mathbf{y}\in\partial Y.
					\end{aligned}
				\end{cases}
			\end{equation*}
			\begin{equation*}
				\begin{cases}
					\begin{aligned}
						& \frac{\partial}{\partial y_j}\Bigl(D_{ijkl}^*(\mathbf{y})\frac{\partial\mathcal{\tilde{M}}_k}{\partial y_l}\Bigr)=-\frac{\partial\bigl(D_{ijkl}^*(\mathbf{y})\alpha_{kl}^*(\mathbf{y})\bigr)}{\partial y_j},\quad \mathbf{y}\in Y, \\
						& \mathcal{\tilde{M}}_k(\mathbf{y})=0, \quad \mathbf{y}\in\partial Y.
					\end{aligned}
				\end{cases}
			\end{equation*}
			\begin{equation*}
				\begin{cases}
					\begin{aligned}
						& \frac{\partial}{\partial y_j}\Bigl(D_{ijkl}^*(\mathbf{y})\frac{\partial \mathcal{\tilde{N}}_k}{\partial y_l}\Bigr)=-\frac{\partial \bigl(D_{ijkl}^*(\mathbf{y})\beta_{kl}^*(\mathbf{y})\bigr)}{\partial y_j}, \quad \mathbf{y}\in Y, \\
						& \mathcal{\tilde{N}}_k((\mathbf{y})=0, \quad \mathbf{y}\in\partial Y.
					\end{aligned}
				\end{cases}
			\end{equation*}
			Moreover, the new macroscopic homogenized problem is defined as below
			\begin{equation*}
				\begin{cases}
					\begin{aligned}
						& -\frac{\partial}{\partial x_{i}}\Bigl(\omega(\mathbf{x})\hat{k}_{ij}^{*}\frac{\partial T^{(0)}}{\partial x_{j}}\Bigr)=h, \text{in }\Omega, \\
						& -\frac{\partial}{\partial x_{i}}\Bigl(\omega(\mathbf{x})\hat{g}_{ij}^{*}\frac{\partial c^{(0)}}{\partial x_{j}}\Bigr)=m, \text{in }\Omega, \\
						&\!-\frac{\partial}{\partial x_{j}}\!\Bigl(\!\omega(\mathbf{x})\hat{D}_{ijkl}^{*}\frac{\partial u_{k}^{(0)}}{\partial x_{l}}-\omega^{2}(\mathbf{x})\hat{A}_{ij}^{*} T^{(0)}-\omega^{2}(\mathbf{x})\hat{B}_{ij}^{*} c^{(0)}\!\Bigr)\!\!=\! f_{i}, \text{in }\Omega, \\
						& T^{(0)}(\mathbf{x}) = \overline{T}(\mathbf{x}),\;\text{on }\;\Gamma_{T},\\
						& \omega(\mathbf{x})\hat{k}_{ij}^{*}\bigl(\mathbf{x}\bigr) \frac{\partial T^{(0)}(\mathbf{x})}{\partial x_{j}}n_i=\overline{q}(\mathbf{x}),\;\text{on }\;\Gamma_{q},\\
						& c^{(0)}(\mathbf{x}) = \overline{c}(\mathbf{x}),\;\text{on }\;\Gamma_{c},\\
						& \omega(\mathbf{x})\hat{g}_{ij}^{*}\bigl(\mathbf{x}\bigr) \frac{\partial c^{(0)}(\mathbf{x})}{\partial x_{j}}n_i=\overline{d}(\mathbf{x}),\;\text{on }\;\Gamma_{d},\\
						& \bm{u}^{(0)}(\mathbf{x}) = \overline{\bm{u}}(\mathbf{x}),\;\text{on }\;\Gamma_{{u}},\\
						& \!\Bigl[\!\omega(\mathbf{x})\hat{D}_{ijkl}^{*}(\mathbf{x}) \frac{\partial u_{k}^{(0)}(\mathbf{x})}{\partial x_l}
						\!-\!\omega^{2}(\mathbf{x}) \hat{A}_{ij}^{*}(\mathbf{x}) T^{(0)}(\mathbf{x})
						\!-\!\omega^{2}(\mathbf{x}) \hat{B}_{ij}^{*}(\mathbf{x}) c^{(0)}(\mathbf{x})\!\Bigr]\! n_{j}\!=\!\overline{\sigma}_{i}(\mathbf{x}),\;\text{on }\;\Gamma_{\sigma},
					\end{aligned}
				\end{cases}
			\end{equation*}
			where the new homogenized material parameters are defined as
			\begin{equation*}
				\begin{aligned}
					& \hat{k}_{ij}(\mathbf{x})=\omega(\mathbf{x}) \frac{1}{|Y|} \int_{Y} \Bigl( k_{ij}^{*}(\mathbf{y})+k_{ik}^{*}(\mathbf{y})\frac{\partial \mathcal{H}^{j}(\mathbf{y})}{\partial y_{k}}\Bigr) dY=\omega(\mathbf{x})\hat{k}_{ij}^{*}, \\
					& \hat{g}_{ij}(\mathbf{x})=\omega(\mathbf{x})\frac{1}{|Y|} \int_{Y} \Bigl( g_{ij}^*(\mathbf{y})+g_{ik}^*(\mathbf{y})\frac{\partial \mathcal{L}^j(\mathbf{y})}{\partial y_k}\Bigr) dY=\omega(\mathbf{x})\hat{g}_{ij}^* ,\\
					& \hat{D}_{ijkl}(\mathbf{x})=\omega(\mathbf{x})\frac{1}{|Y|} \int_{Y} \Bigl( D_{ijkl}^{*}(\mathbf{y})+D_{ijmn}^{*}(\mathbf{y})\frac{\partial \mathcal{X}_{mk}^{l}(\mathbf{y})}{\partial y_{n}}\Bigr) dY=\omega(\mathbf{x})\hat{D}_{ijkl}^{*} ,\\
					& \hat{A}_{ij}(\mathbf{x})=\omega^2(\mathbf{x})\frac{1}{|Y|} \int_{Y}  D_{ijkl}^*(\mathbf{y}) \Bigl( \alpha_{kl}^*(\mathbf{y})+\frac{\partial \mathcal{\tilde{M}}_k(\mathbf{y})}{\partial y_l}\Bigr) dY=\omega^2(\mathbf{x})\hat{A}_{ij}^* ,\\
					& \hat{B}_{ij}(\mathbf{x})=\omega^2(\mathbf{x})\frac{1}{|Y|} \int_{Y}
					D_{ijkl}^{*}(\mathbf{y}) \Bigl(\beta_{kl}^*(\mathbf{y})+\frac{\partial \mathcal{\tilde{N}}_k(\mathbf{y})}{\partial y_{l}}\Bigr) dY=\omega^2(\mathbf{x})\hat{B}_{ij}^*.
				\end{aligned}
			\end{equation*}
			Moreover, the second-order cell functions satisfy the following new unit cell problems
			\begin{equation*}
				\begin{cases}
					\begin{aligned}
						& \frac{\partial}{\partial y_j}\Bigl(D_{ijkl}^*(\mathbf{y})\frac{\partial \mathcal{P}_{kh}^{\alpha_1\alpha_2}(\mathbf{y})}{\partial y_l}\Bigr)=\hat{D}_{i\alpha_1h\alpha_2}^*-D_{i\alpha_1h\alpha_2}^*(\mathbf{y})-D_{i\alpha_1kl}^*(\mathbf{y})\frac{\partial \mathcal{X}_{kh}^{\alpha_2}(\mathbf{y})}{\partial y_l} \\
						& -\frac{\partial}{\partial y_j}\bigl(D_{ijk\alpha_2}^*(\mathbf{y}) \mathcal{X}_{kh}^{\alpha_1}(\mathbf{y})\bigr),\quad \mathbf{y}\in Y ,\\
						& \mathcal{P}_{kh}^{\alpha_1\alpha_2}(\mathbf{y})=0,\quad \mathbf{y}\in\partial Y .
					\end{aligned}
				\end{cases}
			\end{equation*}
			\begin{equation*}			
				\begin{cases}
					\begin{aligned}
						& \frac{\partial}{\partial y_j}\Bigl(D_{ijkl}^*(\mathbf{y})\frac{\partial \mathcal{\tilde{Z}}_k^{\alpha_1}(\mathbf{y})}{\partial y_l}\Bigr)=D_{i\alpha_1kl}^*(\mathbf{y})\Bigl(\frac{\partial \mathcal{\tilde{M}}_k(\mathbf{y})}{\partial y_l}+\alpha_{kl}^*(\mathbf{y})\Bigr)-\hat{A}_{i\alpha_1}^*\\
						& +\frac{\partial}{\partial y_j}\bigl(D_{ijk\alpha_1}^*(\mathbf{y}) \mathcal{\tilde{M}}_k(\mathbf{y})+D_{ijkl}^*(\mathbf{y})\alpha_{kl}^*(\mathbf{y}) \mathcal{H}_{\alpha_1}(\mathbf{y})\bigr),\quad \mathbf{y}\in Y, \\
						&\mathcal{\tilde{Z}}_k^{\alpha_1}(\mathbf{y})=0,\quad \mathbf{y}\in\partial Y.				
					\end{aligned}
				\end{cases}
			\end{equation*}
			\begin{equation*}
				\begin{cases}
					\begin{aligned}
						& \frac{\partial}{\partial y_{j}}\Bigl(D_{ijkl}^{*}(\mathbf{y})\frac{\partial \mathcal{\tilde{G}}_{k}^{\alpha_{1}}(\mathbf{y})}{\partial y_{l}}\Bigr) = D_{i\alpha_{1}kl}^{*}(\mathbf{y})\Bigl(\frac{\partial \mathcal{\tilde{N}}_{k}(\mathbf{y})}{\partial y_{l}}+\beta_{kl}^{*}(\mathbf{y})\Bigr)-\hat{B}_{i\alpha_1}^*\\
						& +\frac{\partial}{\partial y_{j}}\bigl(D_{ijk\alpha_{1}}^{*}(\mathbf{y}) \mathcal{\tilde{N}}_{k}(\mathbf{y})+D_{ijkl}^{*}(\mathbf{y})\beta_{kl}^{*}(\mathbf{y}) \mathcal{L}_{\alpha_{1}}(\mathbf{y})\bigr),\quad \mathbf{y}\in Y, \\
						& \mathcal{\tilde{G}}_{k}^{\alpha_{1}}(\mathbf{y})=0,\quad \mathbf{y}\in\partial Y.
					\end{aligned}
				\end{cases}
			\end{equation*}
			\begin{equation*}
				\begin{cases}
					\begin{aligned}
						& \frac{\partial}{\partial y_j}\Bigl(D_{ijkl}^*(\mathbf{y})\frac{\partial \mathcal{\tilde{Q}}_{kh}^{\alpha_1\alpha_2}(\mathbf{y})}{\partial y_l}\Bigr)=\hat{D}_{i\alpha_2h\alpha_1}^*-D_{i\alpha_2h\alpha_1}^*(\mathbf{y})-D_{i\alpha_2kl}^*(\mathbf{y})\frac{\partial \mathcal{X}_{kh}^{\alpha_1}(\mathbf{y})}{\partial y_l},\; \mathbf{y}\in Y, \\
						& \mathcal{\tilde{Q}}_{kh}^{\alpha_1\alpha_2}(\mathbf{y})=0,\; \mathbf{y}\in\partial Y .
					\end{aligned}
				\end{cases}
			\end{equation*}
			\begin{equation*}
				\begin{cases}
					\begin{aligned}
						& \frac{\partial}{\partial y_j}\Bigl(D_{ijkl}^*(\mathbf{y})\frac{\partial \mathcal{\tilde{W}}_k^{\alpha_1}(\mathbf{y})}{\partial y_l}\Bigr)=2\Bigl(D_{i\alpha_1kl}^*(\mathbf{y})\frac{\partial \mathcal{\tilde{M}}_k(\mathbf{y})}{\partial y_l}+D_{i\alpha_1kl}^*(\mathbf{y})\alpha_{kl}^*(\mathbf{y})-\hat{A}_{i\alpha_1}^*\Bigr) \\
						& +\frac{\partial}{\partial y_j}\bigl(D_{ijk\alpha_1}^*(\mathbf{y}) \mathcal{\tilde{M}}_k(\mathbf{y})\bigr),\quad \mathbf{y}\in Y ,\\
						& \mathcal{\tilde{W}}_k^{\alpha_1}(\mathbf{y})=0,\quad \mathbf{y}\in\partial Y .
					\end{aligned}
				\end{cases}
			\end{equation*}
			\begin{equation*}
				\begin{cases}
					\begin{aligned}
						& \frac{\partial}{\partial y_j}\Bigl(D_{ijkl}^*(\mathbf{y})\frac{\partial \mathcal{\tilde{F}}_k^{\alpha_1}(\mathbf{y})}{\partial y_l}\Bigr)=2\Bigl(D_{i\alpha_1kl}^*(\mathbf{y})\frac{\partial \mathcal{\tilde{N}}_k(\mathbf{y})}{\partial y_l}+D_{i\alpha_1kl}^*(\mathbf{y})\beta_{kl}^*(\mathbf{y})-\hat{B}_{i\alpha_1}^*\Bigr) \\
						& +\frac{\partial}{\partial y_j}\bigl(D_{ijk\alpha_1}^*(\mathbf{y}) \mathcal{\tilde{N}}_k(\mathbf{y})\bigr),\quad \mathbf{y}\in Y ,\\
						& \mathcal{\tilde{F}}_k^{\alpha_1}(\mathbf{y})=0,\quad \mathbf{y}\in\partial Y .
					\end{aligned}
				\end{cases}
			\end{equation*}
			\begin{equation*}
				\begin{cases}
					\begin{aligned}
						& \frac{\partial}{\partial y_i}\Bigl(k_{ij}^*(\mathbf{y})\frac{\partial  \mathcal{H}_{\alpha_1\alpha_2}(\mathbf{y})}{\partial y_j}\Bigr)=\hat{k}_{\alpha_1\alpha_2}^*-k_{\alpha_1\alpha_2}^*(\mathbf{y})-k_{\alpha_1j}^*(\mathbf{y})\frac{\partial \mathcal{H}_{\alpha_2}(\mathbf{y})}{\partial y_j}-\frac{\partial}{\partial y_i}\bigl(k_{i\alpha_2}^*(\mathbf{y}) \mathcal{H}_{\alpha_1}(\mathbf{y})\bigr),\;\; \mathbf{y}\in Y, \\
						& \mathcal{H}_{\alpha_1\alpha_2}(\mathbf{y})=0,\quad \mathbf{y}\in\partial Y .
					\end{aligned}
				\end{cases}
			\end{equation*}
			\begin{equation*}
				\begin{cases}
					\begin{aligned}
						& \frac{\partial}{\partial y_i}\Bigl(g_{ij}^*(\mathbf{y})\frac{\partial \mathcal{L}_{\alpha_1\alpha_2}(\mathbf{y})}{\partial y_j}\Bigr)=\hat{g}_{\alpha_1\alpha_2}^*-g_{\alpha_1\alpha_2}^*(\mathbf{y})-g_{\alpha_1j}^*(\mathbf{y})\frac{\partial \mathcal{L}_{\alpha_2}(\mathbf{y})}{\partial y_j}-\frac{\partial}{\partial y_i}\bigl(g_{i\alpha_2}^*(\mathbf{y}) \mathcal{L}_{\alpha_1}(\mathbf{y})\bigr),\;\; \mathbf{y}\in Y, \\
						& \mathcal{L}_{\alpha_1\alpha_2}(\mathbf{y})=0,\quad \mathbf{y}\in\partial Y .
					\end{aligned}
				\end{cases}
			\end{equation*}
			\begin{equation*}
				\begin{cases}
					\begin{aligned}
						& \frac{\partial}{\partial y_i}\Bigl(k_{ij}^*(\mathbf{y})\frac{\partial \mathcal{\tilde{R}}_{\alpha_1\alpha_2}(\mathbf{y})}{\partial y_j}\Bigr)=\hat{k}_{\alpha_2\alpha_1}^*-k_{\alpha_2\alpha_1}^*(\mathbf{y})-k_{\alpha_2j}^*(\mathbf{y})\frac{\partial \mathcal{H}_{\alpha_1}(\mathbf{y})}{\partial y_j}, \quad \mathbf{y}\in Y, \\
						& \mathcal{\tilde{R}}_{\alpha_1\alpha_2}(\mathbf{y})=0, \quad \mathbf{y}\in\partial Y .
					\end{aligned}
				\end{cases}
			\end{equation*}
			\begin{equation*}
				\begin{cases}
					\begin{aligned}
						& \frac{\partial}{\partial y_i}\Bigl(g_{ij}^*(\mathbf{y})\frac{\partial \mathcal{\tilde{S}}_{\alpha_1\alpha_2}(\mathbf{y})}{\partial y_j}\Bigr)=\hat{g}_{\alpha_2\alpha_1}^*-g_{\alpha_2\alpha_1}^*(\mathbf{y})-g_{\alpha_2j}^*(\mathbf{y})\frac{\partial \mathcal{L}_{\alpha_1}(\mathbf{y})}{\partial y_j}, \quad \mathbf{y}\in Y ,\\
						& \mathcal{\tilde{S}}_{\alpha_1\alpha_2}(\mathbf{y})=0, \quad \mathbf{y} \in\partial Y .
					\end{aligned}
				\end{cases}
			\end{equation*}
		}
		
		{\section*{Appendix B. Supplemental expressions for Section 3.1.}
			\begin{equation*}
				\mathcal{A}_0(\mathbf{x},\mathbf{y}) = -\frac{\partial}{\partial x_i}\Bigl(\hat{k}_{ij}\frac{\partial T^{(0)}}{\partial x_j}\Bigr) + \frac{\partial}{\partial y_i}\Bigl[k_{ij} \frac{\partial}{\partial x_j}\Bigl(\mathcal{H}_{\alpha_1}\frac{\partial T^{(0)}}{\partial x_{\alpha_1}}\Bigr)\Bigr],
			\end{equation*}
			\begin{equation*}
				\mathcal{A}_1(\mathbf{x},\mathbf{y}) = \frac{\partial}{\partial x_i} \Bigl[ k_{ij} \frac{\partial}{\partial x_j} \Bigl( \mathcal{H}_{\alpha_1} \frac{\partial T^{(0)}}{\partial x_{\alpha_1}} \Bigr) \Bigr],
			\end{equation*}
			\begin{equation*}
				\mathcal{B}_0(\mathbf{x},\mathbf{y}) = -\frac{\partial}{\partial x_i} \Bigl( \hat{g}_{ij}\frac{\partial c^{(0)}}{\partial x_j} \Bigr) + \frac{\partial}{\partial y_i} \Bigl[ g_{ij} \frac{\partial}{\partial x_j} \Bigl( \mathcal{L}_{\alpha_1} \frac{\partial c^{(0)}}{\partial x_{\alpha_1}} \Bigr)\Bigr],
			\end{equation*}
			\begin{equation*}
				\mathcal{B}_1(\mathbf{x}, \mathbf{y}) = \frac{\partial}{\partial x_i}\Bigl[ g_{ij} \frac{\partial}{\partial x_j} \Bigl( \mathcal{L}_{\alpha_1} \frac{\partial c^{(0)}}{\partial x_{\alpha_1}} \Bigr)\Bigr],
			\end{equation*}
			\begin{equation*}
				\begin{small}
					\begin{aligned}
						\mathcal{C}_{0i} (\mathbf{x}, \mathbf{y}) & =  -\frac{\partial}{\partial x_j} \Bigl( \hat{D}_{ijkl} \frac{\partial u_k^{(0)}}{\partial x_l} - \hat{A}_{ij} T^{(0)} - \hat{B}_{ij}  c^{(0)} \Bigr) + \frac{\partial}{\partial y_j} \Bigl[ D_{ijkl} \frac{\partial}{\partial x_l} \Bigl( \mathcal{X}_{kh}^{\alpha_1} \frac{\partial u_h^{(0)}}{\partial x_{\alpha_1}} - \mathcal{M}_k T^{(0)} - \mathcal{N}_k c^{(0)} \Bigr)\Bigr] \\
						& - \frac{\partial}{\partial y_j} \Bigl( D_{ijkl} \alpha_{kl} \mathcal{H}_{\alpha_1} \Bigr) \frac{\partial T^{(0)}}{\partial x_{\alpha_1}}  + \frac{\partial}{\partial y_j} \Bigl( D_{ijkl} \beta_{kl} \mathcal{L}_{\alpha_1} \Bigr) \frac{\partial c^{(0)}}{\partial x_{\alpha_1}},
					\end{aligned}
				\end{small}
			\end{equation*}
			\begin{equation*}
				\begin{small}
					\begin{aligned}
						\mathcal{C}_{1i} (\mathbf{x}, \mathbf{y}) & =  \frac{\partial}{\partial x_j} \Bigl[ D_{ijkl} \frac{\partial}{\partial x_l} \Bigl(\mathcal{X}_{kh}^{\alpha_1} \frac{\partial u_h^{(0)}}{\partial x_{\alpha_1}} - \mathcal{M}_k  T^{(0)} - \mathcal{N}_k  c^{(0)} \Bigr)\Bigr] - \frac{\partial}{\partial x_j} \Bigl( D_{ijkl} \alpha_{kl} \mathcal{H}_{\alpha_1}  \frac{\partial T^{(0)}}{\partial x_{\alpha_1}} \Bigr) - \frac{\partial}{\partial x_j} \Bigl( D_{ijkl} \beta_{kl} \mathcal{L}_{\alpha_1}  \frac{\partial c^{(0)}}{\partial x_{\alpha_1}} \Bigr),
					\end{aligned}
				\end{small}
			\end{equation*}
			\begin{equation*}
				\begin{aligned}
					\mathcal{U}(\mathbf{x},\mathbf{y}) & =\frac{\partial}{\partial x_i}\Bigl[k_{ij}\frac{\partial}{\partial x_j}\Bigl(\mathcal{H}_{\alpha_1}\frac{\partial T^{(0)}}{\partial x_{\alpha_1}}\Bigr)\Bigr] +\frac{\partial}{\partial x_i}\Bigl[k_{ij}\Bigl(\frac{\partial \mathcal{H}_{\alpha_1\alpha_2}}{\partial y_j}\frac{\partial^2T^{(0)}}{\partial x_{\alpha_1}\partial x_{\alpha_2}}+\frac{\partial \mathcal{R}_{\alpha_1}}{\partial y_j}\frac{\partial T^{(0)}}{\partial x_{\alpha_1}}\Bigr)\Bigr] \\
					& +\frac{\partial}{\partial y_i}\Bigl[k_{ij}\frac{\partial}{\partial x_j}\Bigl(\mathcal{H}_{\alpha_1\alpha_2}\frac{\partial^2T^{(0)}}{\partial x_{\alpha_1}\partial x_{\alpha_2}}+\mathcal{R}_{\alpha_1}\frac{\partial T^{(0)}}{\partial x_{\alpha_1}}\Bigr)\Bigr] +\epsilon\frac{\partial}{\partial x_i}\Bigl[k_{ij}\frac{\partial}{\partial x_j}\Bigl(\mathcal{H}_{\alpha_1\alpha_2}\frac{\partial^2T^{(0)}}{\partial x_{\alpha_1}\partial x_{\alpha_2}}+\mathcal{R}_{\alpha_1}\frac{\partial T^{(0)}}{\partial x_{\alpha_1}}\Bigr)\Bigr],
				\end{aligned}
			\end{equation*}
			\begin{equation*}
				\begin{aligned}
					\mathcal{V}(\mathbf{x},\mathbf{y}) & =\frac{\partial}{\partial x_{i}}\Bigl[g_{ij}\frac{\partial}{\partial x_{j}}\Bigl(\mathcal{L}_{\alpha_{1}}\frac{\partial c^{(0)}}{\partial x_{\alpha_{1}}}\Bigr)\Bigr] +\frac{\partial}{\partial x_{i}}\Bigl[g_{ij}\Bigl(\frac{\partial \mathcal{L}_{\alpha_{1}\alpha_{2}}}{\partial y_{j}}\frac{\partial^{2}c^{(0)}}{\partial x_{\alpha_{1}}\partial x_{\alpha_{2}}}+\frac{\partial \mathcal{S}_{\alpha_{1}}}{\partial y_{j}}\frac{\partial c^{(0)}}{\partial x_{\alpha_{1}}}\Bigr)\Bigr] \\
					& +\frac{\partial}{\partial y_i}\Bigl[g_{ij}\frac{\partial}{\partial x_j}\Bigl(\mathcal{L}_{\alpha_1\alpha_2}\frac{\partial^2c^{(0)}}{\partial x_{\alpha_1}\partial x_{\alpha_2}}+\mathcal{S}_{\alpha_1}\frac{\partial c^{(0)}}{\partial x_{\alpha_1}}\Bigr)\Bigr] +\epsilon\frac{\partial}{\partial x_i}\Bigl[g_{ij}\frac{\partial}{\partial x_j}\Bigl(\mathcal{L}_{\alpha_1\alpha_2}\frac{\partial^2c^{(0)}}{\partial x_{\alpha_1}\partial x_{\alpha_2}}+\mathcal{S}_{\alpha_1}\frac{\partial c^{(0)}}{\partial x_{\alpha_1}}\Bigr)\Bigr],
				\end{aligned}
			\end{equation*}
			\begin{equation*}
				\begin{aligned}
					\mathcal{J}_{i}(\mathbf{x},\mathbf{y})
					&= \frac{\partial}{\partial x_j} \Bigl[ D_{ijkl} \Bigl( \frac{\partial}{\partial x_l} \bigl( \mathcal{X}_{kh}^{\alpha_1} \frac{\partial u_h^{(0)}}{\partial x_{\alpha_1}}
					- \mathcal{M}_k T^{(0)}
					- \mathcal{N}_k c^{(0)} \bigr) - \alpha_{kl} \mathcal{H}_{\alpha_1} \frac{\partial T^{(0)}}{\partial x_{\alpha_1}}
					- \beta_{kl} \mathcal{L}_{\alpha_1} \frac{\partial c^{(0)}}{\partial x_{\alpha_1}} \Bigr) \Bigr] \\
					& + \frac{\partial}{\partial x_j} \Bigl[ D_{ijkl} \bigl( \frac{\partial \mathcal{P}_{kh}^{\alpha_1\alpha_2}}{\partial y_l} \frac{\partial^2 u_h^{(0)}}{\partial x_{\alpha_1} \partial x_{\alpha_2}}
					+ \frac{\partial \mathcal{Q}_{kh}^{\alpha_1}}{\partial y_l} \frac{\partial u_h^{(0)}}{\partial x_{\alpha_1}} + \frac{\partial \mathcal{W}_k}{\partial y_l} T^{(0)}
					+ \frac{\partial \mathcal{Z}_{k}^{\alpha_1}}{\partial y_l} \frac{\partial T^{(0)}}{\partial x_{\alpha_1}} + \frac{\partial \mathcal{F}_k}{\partial y_l} c^{(0)} + \frac{\partial \mathcal{G}_k^{\alpha_1}}{\partial y_l} \frac{\partial c^{(0)}}{\partial x_{\alpha_1}} \bigr) \Bigr] \\
					& + \frac{\partial}{\partial y_j} \Bigl[D_{ijkl} \Bigl( \frac{\partial}{\partial x_l} \bigl( \mathcal{P}_{kh}^{\alpha_1\alpha_2} \frac{\partial^2 u_h^{(0)}}{\partial x_{\alpha_1} \partial x_{\alpha_2}} + \mathcal{Q}_{kh}^{\alpha_1} \frac{\partial u_h^{(0)}}{\partial x_{\alpha_1}} + \mathcal{W}_k T^{(0)}
					+ \mathcal{Z}_k^{\alpha_1} \frac{\partial T^{(0)}}{\partial x_{\alpha_1}}
					+ \mathcal{F}_k c^{(0)}
					+ \mathcal{G}_k^{\alpha_1} \frac{\partial c^{(0)}}{\partial x_{\alpha_1}} \bigr)
					\\
					& - \alpha_{kl} \bigl( \mathcal{H}_{\alpha_1\alpha_2} \frac{\partial^2 T^{(0)}}{\partial x_{\alpha_1} \partial x_{\alpha_2}} + \mathcal{R}_{\alpha_1} \frac{\partial T^{(0)}}{\partial x_{\alpha_1}} \bigr) - \beta_{kl} \bigl( \mathcal{L}_{\alpha_1\alpha_2} \frac{\partial^2 c^{(0)}}{\partial x_{\alpha_1} \partial x_{\alpha_2}}
					+ \mathcal{S}_{\alpha_1} \frac{\partial c^{(0)}}{\partial x_{\alpha_1}} \bigr) \Bigr) \Bigr] \\
					&+ \epsilon \frac{\partial}{\partial x_j} \Bigl[ D_{ijkl} \Bigl( \frac{\partial}{\partial x_l} \bigl( \mathcal{P}_{kh}^{\alpha_1\alpha_2} \frac{\partial^2 u_h^{(0)}}{\partial x_{\alpha_1} \partial x_{\alpha_2}}
					+ \mathcal{Q}_{kh}^{\alpha_1} \frac{\partial u_h^{(0)}}{\partial x_{\alpha_1}} + \mathcal{W}_k T^{(0)} + \mathcal{Z}_k^{\alpha_1} \frac{\partial T^{(0)}}{\partial x_{\alpha_1}}	+ \mathcal{F}_k c^{(0)}
					+ \mathcal{G}_k^{\alpha_1} \frac{\partial c^{(0)}}{\partial x_{\alpha_1}} \bigr) \\
					& - \alpha_{kl} \bigl( \mathcal{H}_{\alpha_1\alpha_2} \frac{\partial^2 T^{(0)}}{\partial x_{\alpha_1} \partial x_{\alpha_2}}
					+ \mathcal{R}_{\alpha_1} \frac{\partial T^{(0)}}{\partial x_{\alpha_1}} \bigr) - \beta_{kl} \bigl( \mathcal{L}_{\alpha_1\alpha_2} \frac{\partial^2 c^{(0)}}{\partial x_{\alpha_1} \partial x_{\alpha_2}}
					+ \mathcal{S}_{\alpha_1} \frac{\partial c^{(0)}}{\partial x_{\alpha_1}} \bigr) \Bigr) \Bigr],
				\end{aligned}
			\end{equation*}
			{\section*{Appendix C. Supplemental expressions for Section 3.2.}
				\begin{small}
				\begin{equation*}
					\begin{aligned}
						F_0^T & = \Bigl[\frac{\partial}{\partial x_i} \Bigl( k_{ij} \frac{\partial \mathcal{R}_{\alpha_1}}{\partial y_j} \Bigr) + \frac{\partial}{\partial y_i} \Bigl( k_{ij} \frac{\partial \mathcal{R}_{\alpha_1}}{\partial x_j} \Bigr) \Bigr] \frac{\partial T^{(0)}}{\partial x_{\alpha_1}} + \Bigl[ \frac{\partial}{\partial x_i} \Bigl( k_{i \alpha_2} \mathcal{H}_{\alpha_1} + k_{ij} \frac{\partial \mathcal{H}_{\alpha_1 \alpha_2}}{\partial y_j} \Bigr) + \frac{\partial}{\partial y_i} \Bigl( k_{i \alpha_2} \mathcal{R}_{\alpha_1} + k_{ij} \frac{\partial \mathcal{H}_{\alpha_1 \alpha_2}}{\partial x_j} \Bigr)\Bigr] \frac{\partial^2 T^{(0)}}{\partial x_{\alpha_1} \partial x_{\alpha_2}}  \\
						& + \epsilon \Bigl( \frac{\partial k_{ij}(\mathbf{x},\mathbf{y})}{\partial x_i} \frac{\partial \mathcal{H}_{\alpha_1 \alpha_2}}{\partial x_j} + \frac{\partial}{\partial x_i} \bigl( k_{i \alpha_2} \mathcal{R}_{\alpha_1} \bigr)\Bigr) \frac{\partial^2 T^{(0)}}{\partial x_{\alpha_1} \partial x_{\alpha_2}} + \Bigl( k_{\alpha_2 \alpha_3} \mathcal{H}_{\alpha_1} + k_{\alpha_3 j} \frac{\partial \mathcal{H}_{\alpha_1 \alpha_2}}{\partial y_j} + \frac{\partial}{\partial y_i} \bigl( k_{i \alpha_3} \mathcal{H}_{\alpha_1 \alpha_2} \bigr) \Bigr) \frac{\partial^3 T^{(0)}}{\partial x_{\alpha_1} \partial x_{\alpha_2} \partial x_{\alpha_3}}\\
						& + \epsilon \Bigl( \frac{\partial}{\partial x_i} \bigl( k_{i \alpha_3} \mathcal{H}_{\alpha_1 \alpha_2} \bigr) + k_{\alpha_2 \alpha_3} \mathcal{R}_{\alpha_1} \Bigr) \frac{\partial^3 T^{(0)}}{\partial x_{\alpha_1} \partial x_{\alpha_2} \partial x_{\alpha_3}} + \epsilon  k_{ij} \mathcal{H}_{\alpha_1 \alpha_2} \frac{\partial^4 T^{(0)}}{\partial x_{\alpha_1} \partial x_{\alpha_2} \partial x_j \partial x_i} ,
					\end{aligned}
				\end{equation*}
				\end{small}
				\begin{equation*}
					F_i^T = k_{ij}\frac{\partial \mathcal{H}_{\alpha_1} }{\partial x_j}  \frac{\partial T^{(0)}}{\partial x_{\alpha_1}} + \epsilon k_{ij} \frac{\partial \mathcal{R}_{\alpha_1}}{\partial x_j} \frac{\partial T^{(0)}}{\partial x_{\alpha_1}},
				\end{equation*}
				\begin{small}
				\begin{equation*}
					\begin{aligned}
						F_0^c & = \Bigl[\frac{\partial}{\partial x_i} \Bigl( g_{ij} \frac{\partial \mathcal{S}_{\alpha_1}}{\partial y_j} \Bigr) + \frac{\partial}{\partial y_i} \Bigl( g_{ij} \frac{\partial \mathcal{S}_{\alpha_1}}{\partial x_j} \Bigr) \Bigr] \frac{\partial c^{(0)}}{\partial x_{\alpha_1}} + \Bigl[ \frac{\partial}{\partial x_i} \Bigl( g_{i\alpha_2} \mathcal{L}_{\alpha_1}+ g_{ij}\frac{\partial \mathcal{L}_{\alpha_1 \alpha_2}}{\partial y_j} \Bigr) + \frac{\partial}{\partial y_i} \Bigl( g_{i\alpha_2} \mathcal{S}_{\alpha_1} + g_{ij} \frac{\partial \mathcal{L}_{\alpha_1 \alpha_2}}{\partial x_j} \Bigr) \Bigr] \frac{\partial^2 c^{(0)}}{\partial x_{\alpha_1} \partial x_{\alpha_2}} \\
						& + \epsilon \Bigl( \frac{\partial g_{ij}}{\partial x_i} \frac{\partial \mathcal{L}_{\alpha_1 \alpha_2}}{\partial x_j} + \frac{\partial}{\partial x_i} \bigl( g_{i\alpha_2} \mathcal{S}_{\alpha_1} \bigr) \Bigr) \frac{\partial^2 c^{(0)}}{\partial x_{\alpha_1} \partial x_{\alpha_2}} + \Bigl( g_{\alpha_2 \alpha_3} \mathcal{L}_{\alpha_1} + g_{\alpha_3 j} \frac{\partial \mathcal{L}_{\alpha_1 \alpha_2}}{\partial y_j} + \frac{\partial}{\partial y_i} \bigl( g_{ia_3} \mathcal{L}_{a_1 a_2} \bigr) \Bigr) \frac{\partial^3 c^{(0)}}{\partial x_{\alpha_1} \partial x_{a_2} \partial x_{a_3}} \\
						& + \epsilon \Bigl( \frac{\partial}{\partial x_i} \bigl( g_{i\alpha_3} \mathcal{L}_{\alpha_1 \alpha_2}\bigr) + g_{\alpha_2 \alpha_3} \mathcal{S}_{\alpha_1} \Bigr) \frac{\partial^3 c^{(0)}}{\partial x_{\alpha_1} \partial x_{\alpha_2} \partial x_{\alpha_3}} + \epsilon  g_{ij} \mathcal{L}_{\alpha_1 \alpha_2} \frac{\partial^4 c^{(0)}}{\partial x_{\alpha_1} \partial x_{\alpha_2} \partial x_j \partial x_i},
					\end{aligned}
				\end{equation*}
				\end{small}
				\begin{equation*}
					F_i^c = g_{ij} \frac{\partial \mathcal{L}_{\alpha_1}}{\partial x_j} \frac{\partial c^{(0)}}{\partial x_{\alpha_1}} + \epsilon g_{ij} \frac{\partial \mathcal{S}_{\alpha_1}}{\partial x_j} \frac{\partial c^{(0)}}{\partial x_{\alpha_1}},
				\end{equation*}
				\begin{small}
				\begin{equation*}
					\begin{aligned}
						F_{i0}^u & = \Bigl[ \frac{\partial}{\partial x_j} \bigl( D_{ijk\alpha_2} \mathcal{X}_{kh}^{\alpha_1} \bigr) + \frac{\partial}{\partial x_j} \bigl( D_{ijkl} \frac{\partial \mathcal{P}_{kh}^{\alpha_1 \alpha_2}}{\partial y_l} \bigr) + \frac{\partial}{\partial y_j} \bigl( D_{ijkl} \frac{\partial \mathcal{P}_{kh}^{\alpha_1 \alpha_2}}{\partial x_l} \bigr) + \frac{\partial}{\partial y_j} \bigl( D_{ijk\alpha_2} \mathcal{Q}_{kh}^{\alpha_1} \bigr)\Bigr] \frac{\partial^2 u_h^{(0)}}{\partial x_{\alpha_1} \partial x_{\alpha_2}} \\
						& + \frac{\partial}{\partial y_j} \Bigl( D_{ijkl} \frac{\partial \mathcal{Q}_{kh}^{\alpha_1}}{\partial x_l} \Bigr) \frac{\partial u_h^{(0)}}{\partial x_{\alpha_1}} + \epsilon \Bigl[ \frac{\partial}{\partial x_j} \bigl( D_{ijkl} \frac{\partial \mathcal{P}_{kh}^{\alpha_1 \alpha_2}}{\partial x_l} \bigr) + \frac{\partial}{\partial x_j} \bigl( D_{ijk\alpha_2} \mathcal{Q}_{kh}^{\alpha_1} \bigr) \Bigr] \frac{\partial^2 u_h^{(0)}}{\partial x_{\alpha_1} \partial x_{\alpha_2}} \\
						& + \Bigl[ D_{i\alpha_2 k \alpha_3} \mathcal{X}_{kh}^{\alpha_1} + D_{i \alpha_3 kl} \frac{\partial \mathcal{P}_{kh}^{\alpha_1 \alpha_2}}{\partial y_l} + \frac{\partial}{\partial y_j} \bigl( D_{ijk\alpha_3} \mathcal{P}_{kh}^{\alpha_1 \alpha_2} \bigr)\Bigr] \frac{\partial^3 u_h^{(0)}}{\partial x_{\alpha_1} \partial x_{\alpha_2} \partial x_{\alpha_3}} \\
						& + \epsilon \Bigl[ \frac{\partial}{\partial x_j} \bigl( D_{ijk\alpha_3} \mathcal{P}_{kh}^{\alpha_1 \alpha_2} \bigr) + D_{i \alpha_3 kl} \frac{\partial \mathcal{P}_{kh}^{\alpha_1 \alpha_2}}{\partial x_l} + D_{i \alpha_2 k \alpha_3} \mathcal{Q}_{kh}^{\alpha_1} \Bigr] \frac{\partial^3 u_h^{(0)}}{\partial x_{\alpha_1} \partial x_{\alpha_2} \partial x_{\alpha_3}} + \epsilon D_{i \alpha_3 k \alpha_4} \mathcal{P}_{kh}^{\alpha_1 \alpha_2} \frac{\partial^4 u_h^{(0)}}{\partial x_{\alpha_1} \partial x_{\alpha_2} \partial x_{\alpha_3} \partial x_{\alpha_4}},
					\end{aligned}
				\end{equation*}
				\end{small}
				\begin{equation*}
					F_{ij}^{u} = D_{ijkl} \Bigl( \frac{\partial \mathcal{X}_{kh}^{\alpha_1}}{\partial x_l} + \frac{\partial \mathcal{Q}_{kh}^{\alpha_1}}{\partial y_l} + \epsilon \frac{\partial \mathcal{Q}_{kh}^{\alpha_1}}{\partial x_l} \Bigr) \frac{\partial u_h^{(0)}}{\partial x_{\alpha_1}},
				\end{equation*}
				\begin{equation*}
					\begin{aligned}
						F_{i0}^{uT} & = \Bigl[ -\frac{\partial}{\partial x_j} \Bigl(D_{ijkl} \frac{\partial \mathcal{M}_k}{\partial x_l}\Bigr) + \frac{\partial}{\partial x_j} \Bigl(D_{ijkl}\frac{\partial \mathcal{W}_k}{\partial y_l}\Bigr) + \frac{\partial}{\partial y_j} \Bigl(D_{ijkl}\frac{\partial \mathcal{W}_k}{\partial x_l}\Bigr)\Bigr] T^{(0)} + \epsilon \frac{\partial}{\partial x_j} \Bigl(D_{ijkl}\frac{\partial \mathcal{W}_k}{\partial x_l}\Bigr) T^{(0)} \\
						& + \Bigl[ -D_{i\alpha_1 kl}\frac{\partial \mathcal{M}_k}{\partial x_l} + D_{i\alpha_1 kl}\frac{\partial \mathcal{W}_k}{\partial y_l} + \frac{\partial}{\partial y_j}\bigl(D_{ijk \alpha_1}\mathcal{W}_k\bigr) + \frac{\partial}{\partial y_j}\bigl(D_{ijkl}\frac{\partial \mathcal{Z}^{\alpha_1}_{k}}{\partial x_l}\bigr) - \frac{\partial}{\partial y_j}\bigl(D_{ijkl}\alpha_{kl}\mathcal{R}_{\alpha_1}\bigr)\Bigr] \frac{\partial T^{(0)}}{\partial x_{\alpha_l}}\\
						& + \epsilon D_{i\alpha_1 kl}\frac{\partial \mathcal{W}_k}{\partial x_l}\frac{\partial T^{(0)}}{\partial x_{\alpha_1}} + \frac{\partial}{\partial y_j}\bigl(D_{ijk\alpha_2}\mathcal{Z}^{\alpha_1}_{k} - D_{ijkl}\alpha_{kl}\mathcal{H}_{\alpha_1 \alpha_2}\bigr)\frac{\partial^2 T^{(0)}}{\partial x_{\alpha_1} \partial x_{\alpha_2}} \\
						& + \epsilon \frac{\partial}{\partial x_j}\bigl(D_{ijk\alpha_2}\mathcal{Z}^{\alpha_1}_{k} - D_{ijkl}\alpha_{kl}\mathcal{H}_{\alpha_1 \alpha_2}\bigr)\frac{\partial^2 T^{(0)}}{\partial x_{\alpha_1} \partial x_{\alpha_2}} + \epsilon \bigl(D_{i\alpha_2 k\alpha_3}\mathcal{Z}^{\alpha_1}_{k} - D_{i\alpha_3 kl}\alpha_{kl}\mathcal{H}_{\alpha_1 \alpha_2}\bigr)\frac{\partial^3 T^{(0)}}{\partial x_{\alpha_1} \partial x_{\alpha_2} \partial x_{\alpha_3}},
					\end{aligned}
				\end{equation*}
				\begin{equation*}
					\begin{aligned}
						F_{ij}^{uT} & \!=\!\!\Bigl[\!-D_{ijk\alpha_1} \mathcal{M}_k
						\!+\!D_{ijkl} \frac{\partial \mathcal{Z}_k^{\alpha_1}}{\partial y_l}\!-\!D_{ijkl} \alpha_{kl} \mathcal{H}_{\alpha_1}\!+\!\epsilon \Bigl(\!D_{ijk\alpha_1} \mathcal{W}_k\!+\! D_{ijkl} \bigl(\!\frac{\partial \mathcal{Z}_k^{\alpha_1}}{\partial x_l}
						\!-\!\alpha_{kl} \mathcal{R}_{\alpha_1}\!\bigr)\!\Bigr)\!\Bigr]\!\frac{\partial T^{(0)}}{\partial x_{\alpha_1}},
					\end{aligned}
				\end{equation*}
				\begin{equation*}
					\begin{aligned}
						F_{i0}^{uc} & = \Bigl[ -\frac{\partial}{\partial x_j} \Bigl(D_{ijkl}\frac{\partial \mathcal{N}_k}{\partial x_l}\Bigr) + \frac{\partial}{\partial x_j}\Bigl(D_{ijkl}\frac{\partial \mathcal{F}_k}{\partial y_l}\Bigr) + \frac{\partial}{\partial y_j}\Bigl(D_{ijkl}\frac{\partial \mathcal{F}_k}{\partial x_l}\Bigr)\Bigr] c^{(0)} + \epsilon \frac{\partial}{\partial x_j}\Bigl(D_{ijkl}\frac{\partial \mathcal{F}_k}{\partial x_l}\Bigr) c^{(0)} \\
						& + \Bigl[-D_{i\alpha_1 kl}\frac{\partial \mathcal{N}_k}{\partial x_l} + D_{i\alpha_1 kl}\frac{\partial \mathcal{F}_k}{\partial y_l} + \frac{\partial}{\partial y_j} \bigl(D_{ijk \alpha_1}\mathcal{F}_k\bigr) + \frac{\partial}{\partial y_j} \bigl(D_{ijkl}\frac{\partial \mathcal{G}^{\alpha_1}_{k}}{\partial x_l}\bigr) - \frac{\partial}{\partial y_j}\bigl(D_{ijkl}\alpha_{kl}\mathcal{S}_{\alpha_1}\bigr) \Bigr] \frac{\partial c^{(0)}}{\partial x_{\alpha_l}} \\
						& + \epsilon D_{i\alpha_1 kl}\frac{\partial \mathcal{F}_k}{\partial x_l}\frac{\partial c^{(0)}}{\partial x_{\alpha_1}} + \frac{\partial}{\partial y_j}\bigl(D_{ijk\alpha_2}\mathcal{G}^{\alpha_1}_{k} - D_{ijkl}\beta_{kl}\mathcal{L}_{\alpha_1 \alpha_2}\bigr)\frac{\partial^2 c^{(0)}}{\partial x_{\alpha_1} \partial x_{\alpha_2}} \\
						& + \epsilon \frac{\partial}{\partial x_j}\bigl(D_{ijk\alpha_2}\mathcal{G}^{\alpha_1}_{k} - D_{ijkl}\beta_{kl}\mathcal{L}_{\alpha_1 \alpha_2}\bigr)\frac{\partial^2 c^{(0)}}{\partial x_{\alpha_1} \partial x_{\alpha_2}} + \epsilon \bigl(D_{i\alpha_2 k\alpha_3}\mathcal{G}^{\alpha_1}_{k} - D_{i\alpha_3 kl}\beta_{kl}\mathcal{L}_{\alpha_1 \alpha_2}\bigr)\frac{\partial^3 c^{(0)}}{\partial x_{\alpha_1} \partial x_{\alpha_2} \partial x_{\alpha_3}},
					\end{aligned}
				\end{equation*}
				\begin{equation*}
					\begin{aligned}
						F_{ij}^{uc} & = \Bigl[ -D_{ijk\alpha_1} \mathcal{N}_k
						+ D_{ijkl} \frac{\partial \mathcal{G}_k^{\alpha_1}}{\partial y_l} - D_{ijkl} \beta_{kl} \mathcal{L}_{\alpha_1} + \epsilon \Bigl(D_{ijk\alpha_1} \mathcal{F}_k +  D_{ijkl} \bigl(\frac{\partial \mathcal{G}_k^{\alpha_1}}{\partial x_l} - \beta_{kl} \mathcal{S}_{\alpha_1}\bigr)\Bigr)\Bigr]\frac{\partial c^{(0)}}{\partial x_{\alpha_1}}.
					\end{aligned}
				\end{equation*}
				\bibliographystyle{SageV}
				\bibliography{paper}	
			\end{document}